\documentclass[11pt,reqno]{amsart}
\usepackage[a4paper,margin=1.1in]{geometry}

\usepackage{amssymb,mathtools,xcolor}
\usepackage{hyperref}
\hypersetup{colorlinks,linkcolor={red!50!black},citecolor={blue!50!black},urlcolor={blue!75!black}}
\hypersetup{pdftitle={Theta functions and intrinsic mirrors of Laurent phenomenon algebras},
  pdfauthor={JiaRui Fei},bookmarksnumbered=true,bookmarksdepth=3}
\numberwithin{equation}{section}
\allowdisplaybreaks

\newtheorem{theorem}{Theorem}[section]
\newtheorem{proposition}[theorem]{Proposition}
\newtheorem{lemma}[theorem]{Lemma}
\newtheorem{corollary}[theorem]{Corollary}
\newtheorem{boundarycorollary}{Corollary}[theorem]

\theoremstyle{definition}

\newtheorem{example}[theorem]{Example}
\theoremstyle{remark}

\DeclareMathOperator{\Newt}{Newt}
\DeclareMathOperator{\Supp}{Supp}
\DeclareMathOperator{\Spec}{Spec}
\DeclareMathOperator{\Frac}{Frac}
\DeclareMathOperator{\Pic}{Pic}
\DeclareMathOperator{\Cl}{Cl}
\DeclareMathOperator{\ord}{ord}
\DeclareMathOperator{\Tr}{Tr}

\begin{document}
	
\title[Theta functions for LP algebras]{Theta Functions and Intrinsic Mirrors of Laurent Phenomenon Algebras}
\author{JiaRui Fei}
\address{School of Mathematical Sciences, Shanghai Jiao Tong University, China}
\email{jiarui@sjtu.edu.cn}
\thanks{The author was supported in part by National Natural Science Foundation of China (No. 12571038)}

\subjclass[2020]{Primary 13F60, 14J33; Secondary 14M25, 14C20}
\date{}
\dedicatory{}
\keywords{Laurent phenomenon algebra, theta basis, Newton polytope,
intrinsic mirror, log Calabi--Yau variety, Gauss valuation, divisorial valuation,
valuative independence, reciprocity, atomicity}

\begin{abstract}
We construct bases of upper Laurent phenomenon algebras over a field
by minimizing Newton polytopes in a single seed. For any
translation-invariant total order on the exponents, each least exponent
$a$ that occurs determines a smallest Newton polytope $Q_a$.
Choose one element with Newton polytope $Q_a$ and least term $x^a$
for each such $a$. These elements form a basis. In every reachable
seed chart, every Gauss valuation takes a nonzero basis expansion to the
minimum of its values on the nonzero summands.

Suppose the algebra has Laurent coefficients and smooth spectrum over
an algebraically closed field of characteristic zero. After setting the
curve parameters to one, we compute theta products on the intrinsic
mirror from the exchange polynomials and give a finite presentation of
the mirror algebra. A seed has a polynomial corner if each
exchange polynomial can be multiplied by a unit of the coefficient ring
to give a polynomial with constant term one. At such a corner, order
the exponents by a weight with strictly positive, rationally independent
coordinates. For every lattice exponent $a$ and integral weight $p$,
the theta function indexed by $p$ has least $a$-weight
$\min_{b\in Q_a}\langle b,p\rangle$ in the corresponding mirror chart.
This reciprocity holds for cyclic exchanges and arbitrary pole orders.

For fixed polynomial-corner supports, with constant terms normalized
to one, we determine when the forward mirror is smooth and its
augmented mirror recovers the original algebra in the chosen torus
coordinates. The coefficient tuples for which both
conditions hold form a computable finite set. All their entries are positive integers.
\end{abstract}
\maketitle
\setcounter{tocdepth}{1}
{\linespread{0.89}\selectfont\tableofcontents}
\clearpage

\section*{Introduction}

Fomin and Zelevinsky proved that cluster variables are Laurent
polynomials in any cluster~\cite[Theorem~3.1]{FZ}.
Berenstein, Fomin and Zelevinsky defined the \emph{upper cluster algebra}
by intersecting all cluster Laurent rings~\cite[Definition~1.6]{BFZ}.
Lam and Pylyavskyy introduced \emph{Laurent phenomenon algebras}
(LP algebras), whose mutation rule allows irreducible exchanges with
more than two terms, and proved the Laurent phenomenon~\cite[Theorem~5.1]{LP}.
An LP seed specifies mutable variables and exchange polynomials;
mutation replaces one variable and updates the polynomials.
The \emph{upper LP algebra} $\mathcal U$, or full upper algebra,
is the intersection of the Laurent rings of all reachable seeds.

Exchange relations constrain Laurent support.
For example, if $xx'=1+u+v$ with $u,v$ invertible, the coefficient
of $x^{-m}$, $m>0$, must be divisible by $(1+u+v)^m$ for the
polynomial to be Laurent in $x',u,v$. For a fixed least exponent,
how small can the convex hull of the exponents be? Over a field we
prove that a smallest polytope exists, and normalized representatives
form a basis of $\mathcal U$.

For smooth affine log Calabi--Yau varieties with suitable seed charts
over an algebraically closed field of characteristic zero, the same
exchanges determine theta multiplication on the intrinsic mirror.
When does a mirror theta valuation recover the minimum of a linear
form on the smallest polytope? A seed has a \emph{polynomial corner}
if each exchange polynomial becomes a polynomial with nonzero
constant term after rescaling by a nonzero scalar and a Laurent monomial
in the invertible coefficient variables. For smooth full upper LP algebras
with Laurent coefficients, we prove equality at such a corner for strictly
positive ordering weights that separate the exponents. The exchanges may
be cyclic and the pole orders arbitrary.

\subsection*{Smallest Newton polytopes and a basis}
Fix a seed over a field $k$. We call its exchange polynomials the
\emph{ordinary exchanges}. Its coefficient variables are \emph{frozen},
hence fixed by mutation. The invertible ones are the \emph{designated
units}; the others have only nonnegative powers in each chart. With
\emph{Laurent coefficients}, all frozen variables are invertible.

The seed coordinates give the \emph{marking}. Write $M$ and
$L=\operatorname{Hom}(M,\mathbb Z)$ for the character and cocharacter
lattices of the torus where all coordinates are nonzero. The \emph{support}
$\Supp(f)$ of a nonzero Laurent polynomial consists of the exponents $a$
with nonzero coefficient at $x^a$; its \emph{Newton polytope} $\Newt(f)$
is their convex hull. Let $e(f)$ be the least exponent in a fixed
translation-invariant total order on $M$. A \emph{separating weight}
$\omega\in L_{\mathbb R}=L\otimes_{\mathbb Z}\mathbb R$, with distinct
pairings on distinct lattice points, gives such an order.

For every least exponent $a$ attained in the upper algebra, we can
choose $f_a$ with $e(f_a)=a$ such that
\[ Q_a:=\Newt(f_a)\subseteq\Newt(f)\qquad(f\in\mathcal U,\ e(f)=a). \]
The polytope $Q_a$ is unique, and $Q_{na}=nQ_a$ for every integer $n\ge1$.
Normalize each $f_a$ to have coefficient one at $x^a$. Any such
choices form a basis. In every reachable chart the \emph{Gauss valuation}
$\nu$ of a weight takes a Laurent polynomial to the least weight of
its terms. The basis has \emph{valuative independence} for all these
valuations, meaning that $\nu(\sum_a c_af_a)=\min_{c_a\ne0}\nu(f_a)$
for every nonzero sum (Theorem~\ref{thm:minimum-source-newton-polytope}).\footnote{For a related basis theorem on affine opens of log canonical
Calabi--Yau pairs, see Blum and Liu~\cite[Theorem~1.2]{BL}. Their
valuations have log discrepancy zero. For a prime divisor $E$ on a
normal birational model $\pi:X'\to X$, the \emph{log discrepancy} is
$1+\operatorname{ord}_E(K_{X'}-\pi^*(K_X+D))$, with compatible
canonical divisors and $K_X+D$ $\mathbb R$-Cartier. The pair
$(X,D)$ is log canonical if these numbers are all nonnegative.}
Each \emph{cluster monomial}, a product of nonnegative powers of the
mutable variables of a single seed, belongs to every such basis after
normalization at its least exponent, as does any regular Laurent monomial
in a reachable seed (Corollary~\ref{cor:reachable-monomials-minimum-basis}).

The LP starfish theorem reduces regularity to the \emph{first star},
the initial chart and its immediate neighbors~\cite[Theorem~1.3]{Branching}.
In each neighbor, regularity amounts to divisibility by powers of
one exchange polynomial, as in the opening example. The logarithmic
Cartier operator preserves these conditions in positive characteristic.
Specializing the finite linear systems that express them gives the
result in characteristic zero. Starting with any element of least
exponent $a$, finitely many divisibility tests determine $Q_a$.
The \emph{normal fan} of a polytope records its minimizing face in each
direction. With Laurent coefficients, if the degrees of two exchanges in
each other's mutable variables have product at least four at a reachable
seed, no finite fan refines the normal fans of all the $Q_a$ for any
fixed separating weight
(Theorem~\ref{thm:mutual-degrees-infinite-fan}).

A weight \emph{exposes} a vertex when its minimum is attained only there.
With Laurent coefficients, suppose a separating weight exposes a vertex
of $\Newt(g)$, for some regular $g$, that is negative in every mutable
coordinate. Then every $a\in M$ occurs as a least exponent and the upper
algebra is finitely generated. For each rational polytope $\Delta$, the
graded algebra of sections supported in $n\Delta$ has a finite
\emph{Khovanskii basis}: sections whose initial forms generate the
associated graded algebra. An \emph{initial form} records the degree
and least Laurent term (Theorem~\ref{thm:polyhedral-newton-sections}).

The first star also detects orders along frozen coordinate divisors.
Under the valuative independence hypotheses of
Theorem~\ref{thm:polynomial-star-descent}, these orders recover a basis
over the original coefficient ring from one with frozen variables inverted.
The full upper algebra is factorial over a unique factorization domain
(Proposition~\ref{prop:upper-factorial}). For the algebra presented by
ordinary exchanges, Theorem~\ref{thm:local-exact-exchange-strata}
computes the divisor class group, the group of Weil divisors modulo
principal divisors, under explicit domain and regularity hypotheses.

\subsection*{Theta multiplication from the exchanges}
Let $U$ be a smooth connected affine variety over an algebraically
closed field $k$ of characteristic zero. Assume its dense seed torus
is \emph{volume-compatible}, meaning that the wedge of the differentials
$d\log x_j=dx_j/x_j$ extends to a nowhere vanishing volume form on $U$.
Under this log Calabi--Yau condition the \emph{essential skeleton} of
$U$ identifies with the space $L_{\mathbb R}$ of Gauss weights in
the seed torus~\cite[\S2 and Lemma~2.2]{KY}. Keel and Yu construct an
algebra by counting rational curves. Its basis consists of \emph{theta
functions} indexed by the integral points $L$ of the skeleton.
Setting the numerical curve parameters to one gives the \emph{augmented
algebra} $B_{\mathbb Z}(U)$. The \emph{forward mirror} is
$Y=\Spec B_k(U)$, where $B_k(U)=k\otimes_{\mathbb Z}B_{\mathbb Z}(U)$.

Mandel recovers cluster theta multiplication from \emph{descendant
logarithmic Gromov--Witten invariants} associated with products of two or
three theta functions. These are virtual counts on mirror cluster log
pairs with boundary contacts and cotangent-line insertions
\cite[Theorems~1.1 and~1.5]{MandelLogGW}. They agree with
the Keel--Yu counts in the smooth affine cluster case
\cite[\S1.1 and Proposition~6.1]{MandelLogGW}.
Mandel and Ruddat prove the genus-zero correspondence between logarithmic
counts on smooth projective toric varieties and \emph{tropical curves},
balanced weighted piecewise linear graphs with integral slopes
\cite[Theorem~1.1]{MRDesc}.
Their multiplicity formulas compute cluster theta products
\cite[Example~1.3 and Theorem~5.1]{MRTropical}.

Suppose the mutable coordinates $x_i$ are regular on $U$, the designated
units are units, and the seed torus is the locus where all $x_i\ne0$.
Assume the first-star tori cover $U$ in codimension one and inverting the
other mutable variables gives the hypersurface $x_i y_i=F_i$ over the
Laurent ring of unchanged coordinates, with ordinary exchange $F_i$
(Section~\ref{sec:incidence}). We compute the \emph{trace} of a theta
product, its coefficient at the unit theta function, by polynomial
congruences from these possibly nonbinomial exchanges. Traces of two
and three theta functions determine each product through a finite
triangular system with diagonal entries one
(Theorems~\ref{thm:incidence} and~\ref{thm:corner-free-product}).
Finitely many traces present $B_{\mathbb Z}(U)$ and, over a field of
definition with computable field operations and equality, give equations
for the singular locus of $Y$ and test smoothness
(Theorem~\ref{thm:finite-forward-presentation} and
Corollary~\ref{cor:forward-jacobian}). Applications include minuscule
Grassmannians~\cite{Minuscule} and three-point conformal blocks~\cite{Fusion}.

For any smooth affine log Calabi--Yau variety $U$ over $k$ with a
dense split torus and trivial Picard group $\Pic(U)$, the group of line
bundles, the input and output labels determine each contributing curve's
intersections with the boundary divisors of a smooth projective
compactification, hence its numerical class. Augmented multiplication
therefore recovers multiplication over the monoid ring of effective
numerical curve classes, and
$B_k(U)$ is a finitely generated normal domain
(Corollaries~\ref{cor:full-integral-coefficient-gauge} and~\ref{cor:augmentation-normal}).

Under the first-star hypotheses, let $g\in k[U]$ have a Newton vertex
negative in every mutable coordinate. Its \emph{normal cone} $C$
consists of the weights attaining their minimum there. It determines
a Laurent chart on $Y$ with coordinates $z$ dual to the marking of $U$
(Proposition~\ref{prop:negative-vertex-chart}), in which the forward
theta basis is valuatively independent at every real Gauss weight
(Theorem~\ref{thm:forward-vi-main}). After augmentation, an integral
element is \emph{positive} if its local expansions at all generic points
of the skeleton have nonnegative coefficients. The theta functions are
precisely the nonzero positive elements that are not sums of two nonzero
positive elements
(Theorem~\ref{thm:intrinsic-atomicity}).

Proposition~\ref{prop:monic-integral-descent} gives an integral descent
criterion for theta bases of full upper LP algebras.

\subsection*{Reciprocity and change of seed}
Fock and Goncharov proposed that integral tropical points of one
cluster space index a canonical basis on the dual
space~\cite[Conjecture~4.1]{FG}. Gross, Hacking, Keel and Kontsevich
constructed cluster theta functions from scattering diagrams and
conjectured that the tropical evaluations obtained by interchanging
dual theta labels agree. They proved this when one label lies in the
cluster complex, the union of the cones with nonnegative mutable
coordinates in the tropical seed charts
\cite[Lemma~9.10(3) and Remark~9.11]{GHKK}. Cheung, Magee, Mandel and
Muller proved valuative independence and established this equality for
chiral dual cluster scattering diagrams whenever both valuations are finite. Their
Langlands dual formulation requires an integrality condition on the
exchange matrix~\cite[Theorem~1.1 and Corollary~4.20]{CMMM}.

Here the two quantities are the valuation of a forward theta function
and the minimum of a linear form on $Q_a$. Write $\Theta_p$ for
the Laurent expansion of the theta function of label $p\in L$ in
the chart determined by $C$, and $\nu^z_a$ for the Gauss valuation
of weight $a\in M$ there. Choose a separating weight $\omega$ in the interior of $C$ to
define $Q_a$. By Theorem~\ref{thm:reciprocal-polytope}, for each $a$
the forward valuations determine a lattice polytope $P_a$ satisfying
\[ P_a\subseteq Q_a,\qquad \nu^z_a(\Theta_p)=\min_{b\in P_a}\langle b,p\rangle\ge\min_{b\in Q_a}\langle b,p\rangle\qquad(p\in L). \]

Now let $U$ be the smooth spectrum of a full upper LP algebra with
Laurent coefficients, and suppose the marked seed has a polynomial corner.
The first-star hypotheses then hold,
and the nonnegative orthant gives a chart on $Y$. For every
separating weight $\omega$ strictly positive on all marked coordinate
characters, Theorem~\ref{thm:polynomial-corner-reciprocity} gives
$P_a=Q_a$ for every $a\in M$. In particular, any normalized
representative $f_a$ on $Q_a$ satisfies
\[ \nu^z_a(\Theta_p)=\nu^x_p(f_a)\qquad(a\in M,\ p\in L), \]
where $\nu^x_p$ is the Gauss valuation in the marked chart of $U$.
Both evaluation labels are unrestricted.

Return to $xx'=1+u+v$. At $a=(-m,0,0)$, $m>0$, the function
$(x')^m=x^{-m}(1+u+v)^m$ has smallest Newton polytope
$Q_a=a+m\operatorname{conv}(0,e_u,e_v)$, where $e_u,e_v$ are the
exponent vectors of $u,v$. At $p=(0,-n,1)$, $n>0$, the minimum
on this triangle is $-mn$. In the dual coordinates $X,U_1,U_2$,
Proposition~\ref{prop:axis} gives $\Theta_p=U_1^{-n}U_2(1+X)^n$,
whose least $a$-weight is again $-mn$. The exchange
$xx'=1+u+v+2uv$ gives an example with a singular forward mirror
(Section~\ref{sec:second-mirror}).

The proof of reciprocity fixes $p$ and lets $a$ vary. The minimum
of $p$ on $Q_a$ is then the minimum of $a$ on a lattice polytope
$T_p\subset L_{\mathbb R}$ (Proposition~\ref{prop:source-polar-polytope}).
The bound~\eqref{eq:forward-newton-bound} gives $\Newt(\Theta_p)\subseteq T_p$. To obtain
the opposite inclusion, we perturb $p$ and choose regular functions
whose least exponents form lattice bases both at the perturbed weight
and at $\omega$. They give an annular coordinate on a curve over a
function field. Riemann--Roch and the slope formula extend the annular
section across a disk (Proposition~\ref{prop:corner-disk}). The disk
contributes to the coefficient at the perturbed vertex. Passing to
the limit with an exposing weight fixed shows that every vertex of
$T_p$ occurs in $\Theta_p$.

For exchanges affine-linear in successive groups of variables, we also
construct representatives on $Q_a$ from principal minors under the
coprimality hypothesis of Theorem~\ref{thm:principal-minor-reciprocity}.
For mutually dependent nonlinear exchanges, Theorem~\ref{thm:nonlinear-simplex-reciprocity}
gives conditions under which $Q_a$ is a translated convex hull of
dilated exchange polytopes. Example~\ref{ex:cyclic-unequal-poles} shows
that this particular description can fail even when reciprocity holds.

Reciprocity also describes how the least exponent of a fixed section
changes under mutation. Fomin and Zelevinsky's separation formula expresses a cluster
variable in terms of its $F$-polynomial and its $g$-vector
\cite[Corollary~6.3]{FZIV}.\footnote{With \emph{principal coefficients},
the extended exchange matrix has an identity matrix below the initial
exchange matrix $(b_{ij})$. Write $x_j$ for the initial mutable
variables and $y_j$ for the coefficient variables. Setting all
$x_j=1$ in a cluster variable gives its $F$-polynomial. Its
$g$-vector is its multidegree for $\deg x_j=\epsilon_j$ and
$\deg y_j=-\sum_i b_{ij}\epsilon_i$, where the $\epsilon_i$ form the
standard basis~\cite[Definitions~3.1, 3.3 and \S6]{FZIV}.}
In the mirror chart of a seed $t$ with a polynomial corner,
$\Theta_p^t=z^pA_p^t$. Here $z^p$ is the character of exponent
$p$, and $A_p^t$ is a polynomial in the mutable mirror coordinates
with nonnegative integral coefficients and constant term one.
It plays the role of an $F$-polynomial. On
$U$, the least exponent $a$ of $f_a$ plays the role of a
$g$-vector. Suppose mutation at $i$ gives another polynomial corner
and the exchange Laurent polynomial $\widehat F_i=x_ix'_i$ is
itself a polynomial with nonzero constant term. Write $e_i$ for the
$i$-th coordinate cocharacter. The factor $H_i=A_{-e_i}^t$
gives the new mirror coordinate $z'_i=H_i/z_i$
(Proposition~\ref{prop:theta-factor-mutation}). With a strictly positive
separating weight in each seed, reciprocity gives the least exponent
$a'$ of the same $f_a$ in the new seed by
$a'_i=\nu^z_a(H_i)-a_i$ and $a'_j=a_j$ for $j\ne i$
(Corollary~\ref{cor:corner-least-exponent-change}).

\subsection*{Marked return and variation of the coefficients}
At a polynomial corner, reciprocity recovers the smallest Newton
polytopes from the forward mirror. For smooth $Y$, the \emph{reverse
theta functions} have Laurent expansions $\Phi_a$ in the marked
coordinates of $U$. Keel, White and Yu conjecture that the original
variety occurs as a fiber of the second mirror family, and separately
conjecture symmetry of the theta pairing
\cite[Conjectures~21.4 and~21.6]{KWY}. The family-level conjecture is
open in this generality. We study \emph{marked return}:
equality $B_k(Y)=k[U]$ inside the prescribed Laurent field.

The \emph{first mutation partners} are the variables $x'_i$ obtained
by mutating the initial seed once in each mutable direction. They already
determine the ordinary exchanges up to coefficient units in a fixed
marking (Proposition~\ref{prop:marked-exchange-rigidity}). They also
detect marked return. Let $a_i$ be the least exponent of $x'_i$
for the strictly positive separating weight $\omega$, and $c_i$
its coefficient there. Theorem~\ref{thm:first-reverse-partner-return}
gives
\[ B_k(Y)=k[U]\ \text{in the marked field}\quad\Longleftrightarrow\quad \Phi_{a_i}=c_i^{-1}x'_i\ \text{for every mutable }i. \]
When return holds, every normalized cluster monomial in a reachable
seed of $U$ equals the reverse theta function at its least exponent
(Proposition~\ref{prop:arithmetic-marked-return}).

Fix the exchange supports and normalize constant terms to one.
Only finitely many marked augmented forward algebras occur
(Corollary~\ref{cor:constructible-reciprocal-locus}). The tuples for which
the forward mirror is smooth and marked return holds form a computable
finite set of positive integral tuples. We may also allow diagonal torus
translations, which rescale the marked coordinates. Then these properties
hold exactly when the values of the \emph{coefficient characters}, the
invariant Laurent monomials in the exchange coefficients, agree with
their values on one tuple in the finite set
(Theorem~\ref{thm:finite-return-coefficients}).

In every mutable rank at least two, there are smooth varieties with
smooth forward mirror for which return fails after every diagonal
translation (Theorem~\ref{thm:smooth-both-sides-nonreturn}). The
forward theta basis itself can be insensitive to the coefficients.
In an irreducible complex family of smooth polynomial-corner sources
with fixed supports, suppose some coefficient character is nonconstant.
Outside any prescribed countable union of proper closed subsets of the
parameter space, we can then choose two fibers with identical marked
forward theta bases whose upper algebras cannot be identified by a
diagonal torus translation (Theorem~\ref{thm:generic-coefficient-obstruction}).

\section{Upper Laurent phenomenon algebras and their tori}\label{sec:upper}

\subsection{The first star of a seed}

Let $S$ be a unique factorization domain. We use \emph{LP seeds}
in the sense of Lam and Pylyavskyy \cite[\S2]{LP}. A seed $t$ of rank $r$
has algebraically independent mutable variables $x_1,\ldots,x_r$ and
irreducible \emph{exchange polynomials} $F_i\in S[x_1,\ldots,x_r]$.
Each $F_i$ is independent of $x_i$, and no mutable variable divides it.
The \emph{exchange Laurent polynomial} $\widehat F_i$ is $F_i$ times
a monomial with nonpositive exponents in the other mutable variables. Its exponents are
chosen so that, for every $j\ne i$, substituting $x_j=F_j/y$
in $\widehat F_i$ gives a Laurent polynomial not divisible by $F_j$,
where $y$ is an indeterminate. Mutation in direction $i$ gives
the seed $\mu_i t$ and replaces $x_i$ by $x_i'=\widehat F_i/x_i$. We call $F_i$
the ordinary exchange polynomial when we need to distinguish it from
$\widehat F_i$.

Let $L(t)=S[x_1^{\pm1},\ldots,x_r^{\pm1}]$. The \emph{first star}
of $t$ consists of its initial chart and its $r$ adjacent charts.
Their Laurent rings lie in the same rational function field. Write
\[ \mathcal U_\star(t)=L(t)\cap\bigcap_{i=1}^r L(\mu_i t). \]
The \emph{upper LP algebra} $\mathcal U$ is the intersection of the
Laurent rings of all reachable seeds. The LP starfish theorem identifies
$\mathcal U_\star(t)$ with $\mathcal U$ for every LP seed $t$
\cite[Theorem~1.3]{Branching}.
Our coefficient rings include Laurent polynomial rings
$k[u_1^{\pm1},\ldots,u_c^{\pm1}]$, abbreviated $k[u^{\pm1}]$,
and their polynomial extensions. A \emph{polynomial frozen variable}
is one of the polynomial generators of the coefficient ring, with
only nonnegative powers allowed in each chart.

\begin{proposition}\label{prop:upper-factorial}
The ring $\mathcal U$ is factorial, and every mutable seed variable
is prime in $\mathcal U$. For $C_i=S[x_j^{\pm1}:j\ne i]$ and
$s_i=\prod_{j\ne i}x_j$, one has
\[ \mathcal U[s_i^{-1}]=C_i[x_i,F_i/x_i] \simeq C_i[X,Y]/(XY-F_i). \]
The relative Laurent charts $\Spec L(t)$ and $\Spec L(\mu_i t)$,
viewed over $\Spec S$, cover every height-one
point of $\Spec\mathcal U$ over any coefficient UFD $S$. For Laurent
coefficients over a field these charts are tori over that field.
\end{proposition}
\begin{proof}
Factoriality and the description of the prime elements up to associates
follow from \cite[Proposition~2.2]{Minuscule}. We recall the factorization
argument of \cite{Pompili} for the intersection $\mathcal U_\star(t)$.
Write $L_0=L(t)$ and $L_i=L(\mu_i t)$.
In $L_i$, the element $x_i$ is associated to the prime $F_i$.
Indeed, inverting the other variables preserves irreducibility of $F_i$,
and adjoining the independent Laurent variable $x_i'$ preserves
primality. In $L_j$ for $j\ne i$, and in $L_0$, the element
$x_i$ is a unit. If $v_i$ denotes the order at $F_i$ in $L_i$,
then $v_i(x_j)=\delta_{ij}$. Inverting $F_i$ identifies $L_i[F_i^{-1}]$
with $L_0[F_i^{-1}]$. Thus a Laurent polynomial $a\in L_0$ belongs
to $\mathcal U$ precisely when $v_i(a)\ge0$ for each $i$.
In particular $(x_i)\mathcal U=\mathcal U\cap(F_i)L_i$ is prime.

For any prime factor $p$ in $L_0$, the element
$p^\circ=p\prod_i x_i^{-v_i(p)}$ belongs to $\mathcal U$ and has
all its $v_i$ equal to zero. It is prime there. Indeed, if it divides
$ab$, one of $a/p^\circ,b/p^\circ$ is Laurent in $L_0$ by
primality in $L_0$, and that quotient satisfies all the valuation
inequalities. Factoring $a\in\mathcal U\setminus\{0\}$ in $L_0$
and replacing each factor by $p^\circ$ leaves a unit of $S$
and the powers $x_i^{v_i(a)}$. This gives a factorization into
primes. The same argument works at any reachable seed.

The Laurent phenomenon gives $x_i'\in\mathcal U$.
After inverting $s_i$, the ring therefore contains $C_i[x_i,F_i/x_i]$.
It is contained in the intersection of $C_i[x_i^{\pm1}]$ and
$C_i[(F_i/x_i)^{\pm1}]$, which is exactly $C_i[x_i,F_i/x_i]$
by comparing powers of $x_i$. A height-one prime containing
$x_i$ equals the nonzero prime $(x_i)$. It contains neither
any $x_j$ with $j\ne i$ nor the first partner $x_i'$, since
these elements have $v_i$-value zero. Inverting them gives the $i$-th adjacent Laurent chart.
A prime containing no mutable variable lies in the initial chart.
For rank zero the upper algebra is the coefficient ring.
\end{proof}

\subsection{The algebra generated by the ordinary exchanges}

Let $S$ be a regular finitely generated factorial domain over a field
$k$, and let $x_1,\ldots,x_r$ be an LP seed over $S$ with
ordinary exchanges $F_i$. Consider the presentation
\[ A_0=S[x_1,\ldots,x_r,y_1,\ldots,y_r]/(x_i y_i-F_i:1\le i\le r). \]
For a cluster seed, the image of this presentation in the initial
Laurent ring is the lower bound of Berenstein, Fomin and Zelevinsky.
For a coprime acyclic seed they prove that the lower bound equals the
upper bound \cite[Theorem~1.18]{BFZ}. For an LP seed, the following
ideals describe the strata missed by the first star.
For a nonempty $I\subseteq\{1,\ldots,r\}$, write
\[ S_I=S[x_j^{\pm1}:j\notin I],\qquad J_I=(F_i|_{x_a=0\ (a\in I)}:i\in I)\subseteq S_I . \]
Here and below $x_I=0$ means $x_i=0$ for every $i\in I$.

\begin{theorem}\label{thm:local-exact-exchange-strata}
The map $A_0\to S[x_1^{\pm1},\ldots,x_r^{\pm1}]$ sending $y_i$
to $F_i/x_i$ is injective if and only if $A_0$ is a domain,
equivalently if and only if every $J_I$ is nonzero.
In that case $A_0$ is normal and Cohen--Macaulay.
For each $I$ with $|I|\ge2$ and each prime factor $q$ common to
the generators of $J_I$, taken up to a unit in $S_I$, let
$H_{I,q}$ be the closure in $\Spec A_0$ of
$\Spec((S_I/(q))[y_i:i\in I])$ in the stratum $x_I=0$ with the
other $x_j$ invertible.
These are the prime divisors missed by the initial and adjacent
Laurent charts. The divisor class group is
\[ \Cl(A_0)=\bigoplus_{I,q}\mathbb Z[H_{I,q}]. \]
If $A_0$ is a domain, the following conditions are equivalent:
\begin{enumerate}
\item $A_0=\mathcal U$.
\item $A_0$ is factorial.
\item Each $x_i$ is prime in $A_0$.
\item For every $I$ with $|I|\ge2$, the generators of $J_I$
      have unit gcd in $S_I$.
\end{enumerate}
Under these conditions
$\widehat F_i=F_i$ up to the coefficient-unit convention.
In particular, any $S$-subalgebra of $\mathcal U$ containing the
initial variables and their first mutation partners is $\mathcal U$.
\end{theorem}

\begin{proof}
Write $c=\dim S$. On the stratum with exactly the coordinates in
$I$ zero, eliminating $y_j$ for $j\notin I$ gives the ring
$(S_I/J_I)[y_i:i\in I]$. A nonempty such stratum has dimension
$c+r-\operatorname{ht}J_I$, by the dimension formula for affine
domains over a field. The initial Laurent stratum is integral of
dimension $c+r$, so $A_0$ has this dimension. The height theorem
in the regular ambient ring now gives height $r$ at every minimal
prime of the defining ideal. The ideal is locally a complete
intersection, being generated by $r$ elements of height $r$ in
a regular ring. Thus $A_0$ is Cohen--Macaulay with pure dimension
$c+r$ and no embedded primes. If $J_I=0$ for some nonempty $I$,
the closure of that stratum gives a component disjoint from the Laurent
open. If all $J_I\ne0$, every component meets this open. There is
then a single component, generically reduced and without embedded
primes, so $A_0$ is a domain. Localizing a domain at the nonzero
$x_i$ is injective, proving the first assertions.

Suppose henceforth that $A_0$ is a domain. The dimension formula
shows that each boundary divisor has generic point on a stratum of
the form $\Spec((S_I/(q))[y_i:i\in I])$, where $(q)$ is a
height-one prime containing $J_I$. Each such pair occurs once.
At this point eliminate $y_j$ for $j\notin I$ in the ambient
localization. The ambient regular local ring has parameters $q$
and $x_j$ for $j\in I$. In their cotangent space the coefficients
of the $x_j$ in the equations have matrix
\[ \left(y_i\delta_{ij} -\left.\frac{\partial F_i}{\partial x_j}\right|_{x_I=0} \right)_{i,j\in I}\pmod q . \]
Its determinant has coefficient one at $\prod_{i\in I}y_i$, and the
$y_i$ are algebraically independent over $\Frac(S_I/(q))$ at
this generic point. The equations therefore have independent linear
parts. The quotient local ring is regular of dimension one. The
Laurent open is regular, so Serre's criterion gives normality.

For $I=\{i\}$ the irreducible $F_i$ remains prime in $S_I$
because no mutable variable divides it. This gives one divisor
$E_i$. At its generic point $y_i$ is a unit and $F_i$ has
order one in the Laurent ring with coordinates $y_i,x_{j\ne i}$.
Hence $\ord_{E_i}(x_j)=\delta_{ij}$. The $i$-th adjacent chart has
this Laurent ring since $\widehat F_i/F_i$ is a monomial unit in
the unchanged coordinates. Every $H_{I,q}$ with $|I|\ge2$ is
missed by each chart, which requires all but at most one of the
initial coordinates to be units. Localizing $A_0$ at $\prod_i x_i$
gives a UFD whose units are precisely $s x^b$ with
$s\in S^\times$, $b\in\mathbb Z^r$. The divisor localization
sequence thus presents $\Cl(A_0)$ on the boundary divisors with
relations $\operatorname{div}(x_i)$. The identity matrix of the
orders along the $E_j$ eliminates these generators and gives the
claimed free group.

Consequently factoriality is equivalent to the absence of a common
prime factor on every stratum of size at least two. When there is
no such factor, $\operatorname{div}(x_i)=E_i$ and normality gives
$(x_i)=\{f\in A_0:\ord_{E_i}f\ge1\}$, a prime ideal. Conversely
a divisor $H_{I,q}$ would contain $(x_i)$ for $i\in I$ and
would be distinct from $E_i$, which precludes primality of $x_i$.
If the charts miss no divisor, their intersection equals $A_0$
by normality. This intersection is $\mathcal U_\star(t)=\mathcal U$.
Conversely $\mathcal U$ is factorial by
Proposition~\ref{prop:upper-factorial}.

Finally, a nontrivial normalizing power of $x_j$ in $\widehat F_i$
would require $F_j$ to divide the substitution
$F_i|_{x_j=F_j/t}$ in the LP normalization test. Reduction modulo
$F_j$ gives $F_j\mid F_i|_{x_j=0}$. This last polynomial is
nonzero, since $x_j\nmid F_i$, and is independent of $x_i$.
The divisibility holds in the polynomial ring by primality and
localization at mutable monomials, so $F_j$ is independent of
$x_i$ as well. It would then remain a common nonunit factor for
$J_{\{i,j\}}$, contrary to the gcd condition. Hence every $y_i$
is the corresponding first partner up to a coefficient unit. These
elements and the initial variables generate $A_0=\mathcal U$,
which proves the last assertion.
\end{proof}

\subsection{Recovering the exchanges in a marked torus}

Let $S$ be a unique factorization domain. Consider two LP seeds
with the same ordered mutable coordinates $x_1,\dots,x_r$ over
$S$, in the same rational function field. Write $F_i,G_i$ for
their respective ordinary exchanges, $x'_{i,F}=\widehat F_i/x_i$
for the first mutation partners of the first seed, and
$\mathcal U_F,\mathcal U_G$ for the full upper algebras.

\begin{proposition}\label{prop:marked-exchange-rigidity}
The following conditions are equivalent:
\begin{enumerate}
\item $x'_{i,F}\in\mathcal U_G$ for every $i$;
\item $\mathcal U_F\subseteq\mathcal U_G$;
\item $\mathcal U_F=\mathcal U_G$;
\item for each $i$ there is a unit $e_i\in S^\times$ with
      $F_i=e_iG_i$.
\end{enumerate}
\end{proposition}
\begin{proof}
Each $x'_{i,F}$ belongs to $\mathcal U_F$ by the Laurent
phenomenon. Conversely, assume the first condition and write
$C_i=S[x_j^{\pm1}:j\ne i]$. In the $i$-th adjacent Laurent
ring for $G$ the coefficient of $x'_{i,G}$ in $x'_{i,F}$ is
$\widehat F_i/\widehat G_i$. This coefficient belongs to $C_i$, so
$G_i\mid F_i$ there. Both ordinary exchanges remain
irreducible nonunits in $C_i$: each is irreducible before
localization and divisible by no mutable variable. Hence they
are associates in $C_i$. A unit there is a unit of $S$ times
a mutable Laurent monomial. Both exchanges have least degree zero
in each unchanged mutable variable, so this monomial has exponent
zero. This proves the fourth condition. Under that condition the
two adjacent Laurent rings agree at every index, and their initial
Laurent ring is the same. Thus their first-star intersections, and
hence their upper LP algebras, agree.
\end{proof}

Now let $S=k[u_1^{\pm1},\ldots,u_c^{\pm1}]$, where $k$ is
algebraically closed of characteristic zero. Let $M$ be the
character lattice of the full initial torus. A point
$t\in\operatorname{Hom}(M,k^\times)$ acts by translation,
sending the character $\chi^a$ to $t^a\chi^a$, where $t^a=t(a)$.
Let $D\subset M$ be generated by all
$a-b$ with $a,b\in\Supp(F_i)$ for some $i$. We take both
exponents from the same exchange, including its unit
coordinates.

\begin{corollary}\label{cor:diagonal-exchange-stabilizer}
The diagonal translations preserving $\mathcal U_F$ in its marked
Laurent field form $\operatorname{Hom}(M/D,k^\times)$.
If specified coefficient characters must be fixed, adjoin their
exponents to $D$ before forming the quotient.
\end{corollary}
\begin{proof}
Translation carries the first star to the star whose ordinary
exchanges are $t(F_i)$ up to scalar factors. Indeed the variables
and normalizing monomials change only by scalars in each of these
Laurent charts. By Proposition~\ref{prop:marked-exchange-rigidity}
the ring is preserved precisely when $t(F_i)=e_iF_i$ for
$e_i\in S^\times$. Translation preserves the finite support of
$F_i$, whereas a nonzero exponent of the monomial unit $e_i$ would
translate that support. Thus $e_i\in k^\times$, and comparison of
coefficients gives exactly $t^{a-b}=1$ for the indicated pairs.
\end{proof}

We use \emph{counting source} for a smooth connected affine log
Calabi--Yau variety with a volume-compatible dense split torus over
an algebraically closed field of characteristic zero.

Suppose the two marked seeds give counting sources. A torus translation
identifies a variety with its translate and preserves the
volume and the cocharacter labels. Transporting the compactification
and the evaluation point identifies the counting fibers. At curve-class
augmentation the two labeled forward algebras therefore agree.
Whenever their forward Laurent marking is defined by the same
specified theta labels, it agrees as well. A reverse construction
depending only on this marked forward algebra then gives the same
subalgebra in the fixed reverse Laurent field. If it recovers both
upper LP algebras in that field, the translation belongs to the
group of Corollary~\ref{cor:diagonal-exchange-stabilizer}.

We now take $S=k[u_1^{\pm1},\ldots,u_c^{\pm1}]$, where $k$ is
algebraically closed of characteristic zero. For an affine smooth
full upper algebra with marked torus $T$, the logarithmic volume form
$\Omega$ is the wedge of the logarithmic differentials of all the
coordinates on $T$, including the $u_j$. A mutation preserves
$\Omega$ up to sign, since $\widehat F_i$ is independent of $x_i$.
Thus $\Omega$ has no zeros or poles on any torus of the first star.
By Proposition~\ref{prop:upper-factorial}, it extends to a nowhere
vanishing regular top form on the smooth affine variety. The product $\prod_i x_i$ cuts out the complement of the initial torus.

\subsection{Gauss valuations in seed charts}

For a split torus $T=\Spec k[M]$, the group algebra $k[M]$
has basis the characters $\chi^a$, $a\in M$, with
$\chi^a\chi^b=\chi^{a+b}$. Write $L=\operatorname{Hom}(M,\mathbb Z)$.
For $f=\sum_a f_a\chi^a\ne0$ let
\[ \tau_f(w)=\min_{a:f_a\ne0}\langle a,w\rangle,\qquad w\in L_{\mathbb R}. \]
In torus coordinates $x=(x_1,\ldots,x_d)$, we write
$v_{x,w}(f)=\tau_f(w)$ for the \emph{Gauss valuation} of weight $w$.
We also call $\tau_f$ the \emph{tropicalization} of $f$ in this
chart. The \emph{initial form} at $w$ is the sum of the terms of
least $w$-weight. The initial forms of two nonzero Laurent polynomials have nonzero product,
so their values add under multiplication. The valuation extends to the
rational function field $K=k(x_1,\ldots,x_d)$ by subtraction on quotients.
We call $f$ \emph{coercive} if $-\tau_f$ is strictly positive away
from the origin, equivalently if $0$ lies in the interior of $\Newt(f)$.
The following form of the boundedness argument in \cite[Lemma~2.9
and Corollary~15.11]{KY} gives an explicit polytope.

\begin{proposition}\label{prop:coercive}
Suppose $U$ is smooth affine with dense torus $T$, and that
the logarithmic torus volume form extends to a nowhere vanishing
volume form on $U$. Adjoin $1$ to any finite set of algebra
generators of $k[U]$. The convex hull of their Laurent supports
contains $0$ in its interior in $M_{\mathbb R}$. Over an infinite
ground field a general linear combination $f$ of the generators and $1$ has
$\rho=-\tau_f$ strictly positive off the origin, with bounded
polyhedral sublevel sets.
\end{proposition}
\begin{proof}
Otherwise rational separation gives a nonzero primitive $v\in L$
whose pairing with every exponent is nonnegative. The corresponding valuation ring contains $k[U]$.
On a smooth toric model containing the ray of $v$, the logarithmic
form $\Omega$ has order $-1$ along its divisor. Yet the rational
map from this model to the affine variety $U$ is defined near the
generic point of the divisor, since all the generators belong to its
valuation ring. The pullback of the regular form on $U$ is regular
there, a contradiction. Avoiding the finitely many linear equations for
cancellation of an exponent gives $f$ with the asserted Newton
polytope. Its negative support minimum dominates a positive multiple
of a norm, by compactness of the unit sphere. The assertion about
sublevel sets follows.
\end{proof}

\begin{lemma}\label{lem:boundary-gauss-change}
Suppose two systems of torus coordinates differ by
$x_j'=G/x_j$, where $G$ is a nonzero rational function in the
coordinates other than $x_j$. Set $w'_j=v_{x,w}(G)-w_j$ and
$w'_i=w_i$ for $i\ne j$. Then $v_{x,w}=v_{x',w'}$ on $K$.
\end{lemma}

\begin{proof}
Write $E=k(x_i:i\ne j)$ and $v=v_{x,w}$. Grouping a Laurent
polynomial by its $x_j$-exponents gives
\[ v\!\left(\sum_r b_r x_j^r\right) =\min_{b_r\ne0}\bigl(v(b_r)+r w_j\bigr)\qquad (b_r\in E). \]
For rational $b_r$, this follows by clearing a common denominator in
$E$. Substituting $(x'_j)^r=G^r x_j^{-r}$ gives the same formula
with $x'_j,w'_j$ in place of $x_j,w_j$. Taking quotients proves
the assertion.
\end{proof}

For an LP mutation, take $G=\widehat F_j$, including its monomial
denominator. Iterating the lemma transports a Gauss valuation along any
finite sequence of mutations. Under a monomial change of coordinates,
the new weights are the values of the old weight on the new exponents.

\section{Newton polytopes of regular functions}\label{sec:source-polytopes}

\subsection{The smallest Newton polytope at a given least exponent}
\label{subsec:smallest-newton-polytope}

We allow both invertible and polynomial frozen variables. Over a
field $k$, take the coefficient ring
$S=k[u_1^{\pm1},\ldots,u_c^{\pm1},s_1,\ldots,s_l]$ and set
$C_i=S[x_j^{\pm1}:j\ne i]$. For each $i$, choose
a nonzero $E_i\in C_i$ and write $y_i=E_i/x_i$. We work with the
first-star intersection
\[ A=S[x_1^{\pm1},\ldots,x_r^{\pm1}]\cap\bigcap_{i=1}^r C_i[y_i^{\pm1}]. \]
For an LP seed one may take $E_i=\widehat F_i$ or $F_i$, since
they differ by a unit of $C_i$. Then $A=\mathcal U_\star(t)=\mathcal U$.
The arguments in this section apply to the stated intersection even
when an $E_i$ is reducible or a unit. All the rings lie in the common
rational function field. The initial coordinates embed $A$ in $k[M]$,
where $M$ is the character lattice of the torus with coordinates
$x_i,u_j,s_j$. The $s_j$ have nonnegative exponents in every chart.
Fix a translation-invariant total order on $M$, so $a<b$
implies $a+c<b+c$ for every $c\in M$. For a nonzero Laurent polynomial $f$, write
$e(f)$ for its \emph{least exponent}. The set of exponents attained
in $A$ is $\Gamma=\{e(f):0\ne f\in A\}$. A weight
$\omega\in M^*_{\mathbb R}$ is \emph{separating} if its pairings
with distinct lattice points are distinct. Such a weight orders $M$
by increasing pairing. We write $\nu_p=v_{(x,u,s),p}$ for the Gauss
valuation in the initial chart and use $\nu^x_p$ when more than one
chart is in play. By a \emph{section supported in} a polytope
$\Delta$ we mean an element of $A$ with support in $\Delta$.
A \emph{monic} section has least coefficient one.
Sums of polytopes below are Minkowski sums, so $P+Q=\{p+q:p\in P,\ q\in Q\}$.

\begin{theorem}\label{thm:minimum-source-newton-polytope}
For each $a\in\Gamma$, some $f\in A$ with $e(f)=a$ has Newton
polytope $Q_a$ contained in the Newton polytope of every other such
function. These smallest polytopes satisfy
\[ Q_{na}=nQ_a\quad(n\ge1),\qquad Q_{a+b}\subseteq Q_a+Q_b . \]
Choose any $f_a\in A$ supported in $Q_a$ with coefficient one
at $x^a$. The $f_a$, for $a\in\Gamma$, form a $k$-basis.
For every real Gauss weight $p$ and every nonzero finite sum,
\[ \nu_p\Bigl(\sum_a c_a f_a\Bigr)=\min_{c_a\ne0}\nu_p(f_a). \]
For an LP seed the same basis satisfies the formula in every
reachable seed chart.
\end{theorem}
\begin{proof}
For a finite set $E\subset M$, we first characterize the sections
supported in $E$ with a prescribed least term. If $j>0$, the
coefficient $h$ of $x_i^{-j}$ must belong to $E_i^j C_i$ for
the section to be Laurent in the $i$-th adjacent chart. This follows
by substituting $x_i=E_i/y_i$. Conversely, divisibility of
every negative slice suffices. Write $h=E_i^j h'$.
Take the coefficients of $h'$ as auxiliary unknowns. The
coordinate minima and maxima of the quotient are bounded by those
of the slice minus the corresponding extrema of the divisor,
since Newton polytopes add under multiplication. Polynomial
frozen exponents must also be nonnegative in each chart.
Equating coefficients in each slice and prescribing the least term
gives a finite linear system over $k$.

First suppose $k$ is perfect of characteristic $p>0$.
For $q=p^N$, we express the $N$-fold \emph{Cartier operator} on functions
using the torus logarithmic volume form:
\[ C_q\Bigl(\sum_m c_mx^m\Bigr)=\sum_{m\in qM}c_m^{1/q}x^{m/q}. \]
It preserves $A$. Indeed the slice of degree $-qj$ is
$E_i^{qj}h'$. Its contribution in degree $-j$ is
$E_i^j C_q(h')$, where the operator on the right uses
the unchanged coordinates. Nonnegative frozen exponents remain
nonnegative. Over an imperfect field we may make this construction
over its algebraic closure and descend solvability of the finite
linear system to $k$. Benito, Muller, Rajchgot and Smith give the
corresponding splitting for upper cluster algebras
\cite[Theorem~3.7]{BMRS}.

Fix $a\in\Gamma$ and start with any section of least exponent $a$.
Among the lattice subpolytopes of its Newton polytope that support
such a section, choose one minimal under inclusion and call it $Q$.
Let $f$ be a section on $Q$ of least exponent $a$.
Let $h$ be any other section at $a$, and write
$R=\Newt(h)$. Normalize both least coefficients to one. If
$Q\not\subseteq R$, there is an integral covector $v$ such
that $\min_Q v<\min_R v$. For large $q$, every lattice point of
$(1-1/q)Q+(1/q)R$ belongs to $Q$, by compactness and discreteness.
The strict inequality of minima excludes the face $Q^v$ where
$v$ attains its minimum on $Q$.
Yet $C_q(f^{q-1}h)$ is supported in this smaller set and has least
term $x^a$. This contradicts the choice of
$Q$ in positive characteristic.

In characteristic zero use the fixed finite system for the allowed
support $(Q\cap M)\setminus Q^v$, the prescribed coefficient at
$a$, and zero coefficients below it. If it is inconsistent,
include a linear combination of its equations giving $0=1$
in a finitely generated $\mathbb Z$-algebra of definition.
Adjoin the nonzero scalar coefficients of the exchange polynomials
and their inverses, together with the coefficients of the first-star
expressions for $f$ and $h$. By dominance and constructibility
over $\Spec\mathbb Z$, this subalgebra admits specializations to
algebraically closed fields of arbitrarily large characteristic $p$.
The support of each exchange survives specialization, so the same
bounds on the auxiliary coefficients define the specialized system.
The supports of $f,h$ may shrink, but their least coefficients stay
one and their chart identities persist. In particular,
$\Supp(f^{p-1}h)\subseteq (p-1)Q+R$ after specialization.
The Cartier construction with $q=p$ therefore gives a solution of
the specialized system, contrary to $0=1$. This proves $Q\subseteq R$ over
every field and proves uniqueness of $Q_a$.

Products give $Q_{a+b}\subseteq Q_a+Q_b$. To prove the reverse
inclusion for dilation, choose a monic representative $u$ at
$na$ with polytope $S=Q_{na}$. Write $q=nt+r$,
$0\le r<n$. The section $C_q(u^t f^r)$ has least exponent
$a$ and support in $(tS+rQ_a)/q$. These bodies tend to
$S/n$, so for large $q$ all their lattice points lie in
$(S/n)\cap M$. The same finite-system argument in characteristic
zero, and descent in positive characteristic, give a section at
$a$ supported there. Hence $nQ_a\subseteq S$.

If $a$ is the least exponent in the support of a Laurent polynomial,
it is also the least lattice point of its Newton polytope. Indeed, for
any lattice point $b$ in the polytope, clearing denominators in a
convex combination of support points expresses a positive multiple of
$b-a$ as a sum of differences that are nonnegative in the chosen
order. Hence $b\ge a$. In particular the least term of $f_a$
is $x^a$. For $0\ne h\in A$ subtract the appropriate multiple of
$f_{e(h)}$. The support remains inside $\Newt(h)$ and the
least exponent increases, so successive subtraction terminates.
The least exponents of the $f_a$ are distinct, proving
independence. In the resulting expansion every summand's polytope
lies in the Newton polytope of the sum. The opposite convex-hull
inclusion follows by addition. Taking support minima proves the
valuation formula. Under a seed mutation a Gauss valuation is
transported to a Gauss valuation, by
Lemma~\ref{lem:boundary-gauss-change}. Iterating proves the last claim.
\end{proof}

\begin{corollary}\label{cor:reachable-monomials-minimum-basis}
Suppose $A$ is a full upper LP algebra. Let $h\in A$ be a Laurent
monomial in the coordinates of a reachable seed, and let $a=e(h)$
in the initial chart. In every basis of
Theorem~\ref{thm:minimum-source-newton-polytope}, the element $f_a$
is $h$ divided by its coefficient at $x^a$.
\end{corollary}
\begin{proof}
Expand $h$ in the chosen basis. In the reachable chart the minimum
formula holds for every Gauss weight, so the Newton polytope of $h$
is the convex hull of the Newton polytopes of the nonzero summands.
It is a single point. Thus every summand is a scalar multiple of the
same monomial, and linear independence leaves just one. Its least
exponent and coefficient in the initial chart give the assertion.
\end{proof}

For a full upper LP algebra, the smallest polytopes also determine the
least exponents in another seed. Let $T$ carry Gauss weights from the
initial seed $t$ to a reachable seed $t'$, as in
Lemma~\ref{lem:boundary-gauss-change}. For a section $f_a$ on $Q_a$,
the minimum on $\Newt_{t'}(f_a)$ at $w'$ is
$\min_{b\in Q_a}\langle b,T^{-1}(w')\rangle$.
At a separating weight $w'$, the slope of this function is the new
least exponent $a'$. In particular it depends only on $Q_a$.
The new polytope is again the smallest one for its least exponent.
Indeed, complete $f_a$ to a basis as in
Theorem~\ref{thm:minimum-source-newton-polytope}. Its members have
distinct least exponents in $t'$, since cancelling two equal least
terms would contradict valuative independence at $w'$. Every
element of least exponent $a'$ therefore has $f_a$ in its basis
expansion. Valuative independence at all weights in $t'$ gives the
claimed inclusion of Newton polytopes.

We can recover $Q_a$ from any section at $a$. Starting with its
Newton polytope, remove a vertex other than $a$ and take the convex
hull of the remaining lattice points. Accept the smaller hull if its
normalized first-star system has a solution. If removal of a vertex
fails, every subsequent hull supporting a section at $a$ must still
contain that vertex. Each lattice point therefore needs at most one
test. When no further vertex can be removed, every proper lattice
subpolytope omits an indispensable vertex. By the theorem, the final
hull is $Q_a$.

For a polytope $\Delta$, let $V(\Delta)$ be the space of elements
of $A$ supported in $\Delta$. The representatives with
$Q_b\subseteq\Delta$ form a basis of $V(\Delta)$.
Their distinct least exponents give a triangular invertible
coefficient matrix on these indices. In particular there is one
normalized representative on $Q_a$ whose coefficient at every
other $b$ with $Q_b\subseteq Q_a$ is zero.

\subsection{Polytopes dual to the smallest Newton polytopes}

With Laurent coefficients, suppose that some $g\in A$ has a
vertex $v$ negative in every mutable coordinate. For any $a\in
M$, choose $l\ge0$ with $a_i-lv_i\ge0$ at all mutable indices.
After normalizing its leading coefficient,
$R_a=x^{a-lv}g^l$ has least exponent $a$ throughout the
interior of the normal cone $C$ of $v$. Thus $\Gamma=M$
at every separating weight in that interior. The inclusion
$Q_a\subseteq\Newt(R_a)$ makes $a$ its unique minimizing
vertex there. Comparing the smallest polytopes at any two such
weights now gives the same $Q_a$.

Write $L=M^*$. For each
$p\in L$, the minimum of $p$ on $Q_a$ varies with $a$ as the
support function of a lattice polytope. We use the minimum convention.
The next lemma is the one-sided form of Thurston's integer seminorm
theorem \cite{Thurston}. For the support-function argument see \cite[Lemma~0.1 and the proof
of Theorem~1]{deLaSalle}. The proof below uses the differentiability
argument of Vershynin recorded in \cite[Lemma~4]{KapovichInt}.

\begin{lemma}\label{lem:integral-support}
An integer-valued superadditive function on a lattice, homogeneous
under positive integral dilation, is the minimum support function
of a unique lattice polytope.
\end{lemma}
\begin{proof}
Write $h$ for the function and extend it to rational vectors by
homogeneity. Superadditivity bounds $h(x+y)-h(x)$ between $h(y)$
and $-h(-y)$.
Decomposition into signed coordinate vectors gives a Lipschitz bound,
so the extension to real vectors is homogeneous and concave.
At a differentiability point $v$, take lattice points $u_N$
within bounded distance of $Nv$. For a basis vector $e_i$,
differentiability at $v$ gives
$h(u_N+e_i)-h(u_N)\to\partial_i h(v)$. Indeed, write
$u_N=Nv+r_N$ with $r_N$ bounded. Homogeneity and differentiability
give $h(u_N)=Nh(v)+Dh(v)(r_N)+o(1)$, and replacing $r_N$ by
$r_N+e_i$ gives the corresponding expansion for $h(u_N+e_i)$.
The differences are integers. The gradients are consequently integral and bounded,
so only finitely many occur. Each gives a linear upper bound with
equality at its differentiability point. Density of those points
and continuity show that their minimum equals $h$ everywhere.
Their convex hull is the required polytope.
\end{proof}

\begin{proposition}\label{prop:source-polar-polytope}
For every $p\in L$ there is a unique lattice polytope $T_p$
in $L_{\mathbb R}$ such that
\begin{equation}\label{eq:dual-polytope-minima}
\min_{e\in T_p}\langle a,e\rangle=\min_{b\in Q_a}\langle b,p\rangle \qquad(a\in M).
\end{equation}
It is computable from finitely many first-star linear systems.
One has $T_{np}=nT_p$ for $n\ge1$ and
$T_{p+p'}\subseteq T_p+T_{p'}$.
\end{proposition}
\begin{proof}
Write $d_p(a)=\min_{b\in Q_a}\langle b,p\rangle$.
Theorem~\ref{thm:minimum-source-newton-polytope} makes $d_p$
superadditive and homogeneous in $a$. Apply
Lemma~\ref{lem:integral-support}. The dilation and sum statements
now follow from homogeneity and superadditivity in $p$.
For a given $a$ the explicit function $R_a$ above supplies a
finite set of possible exponents for computing $Q_a$. To obtain $T_p$,
use the coordinate box with lower and upper bounds $d_p(e_i)$
and $-d_p(-e_i)$. In each full-dimensional chamber cut out by
the differences of its lattice points, an interior lattice weight
has distinct pairings with those points. Its value identifies the
minimizing vertex. Every vertex occurs, and the procedure also
covers a box consisting of one point.
\end{proof}

\subsection{Section algebras and finite generation}

The condition $Q_a\subseteq\Delta$ determines exactly which least
exponents occur among sections supported in $\Delta$. We use the
dual polytopes to express this condition by finitely many inequalities.
The hypothesis on the negative vertex below has content even in rank
three. Speyer constructs a first-star intersection over a field that
is not finitely generated \cite{Speyer}. The theorem therefore
excludes such a vertex for that example.

\begin{theorem}\label{thm:polyhedral-newton-sections}
Let $A$ be the first-star intersection above with Laurent coefficient
ring $k[u^{\pm1}]$. Suppose some $g\in A$ has a Newton vertex
negative in every mutable coordinate. Fix a separating weight $\omega$
in the interior of its normal cone. Let $\Delta$ be a nonempty
rational polytope given by inequalities $\langle b,p_j\rangle\ge c_j$,
$1\le j\le l$, with integral normals and offsets, where each equation
for the affine hull is written as two opposite inequalities.
Write $\operatorname{Vert}(P)$ for the vertex set of a polytope $P$.
Then the set
\[ D_\omega(\Delta) =\{a\in M_{\mathbb R}: \langle a,v\rangle\ge c_j\quad (v\in\operatorname{Vert}(T_{p_j}),\ 1\le j\le l)\} \]
is a rational polytope, possibly empty. For $n>0$, the least exponents of the nonzero
elements of $V(n\Delta)$ are precisely $nD_\omega(\Delta)\cap M$.
The graded algebra $\bigoplus_{n\ge0}V(n\Delta)t^n$ is finitely
generated. It has a finite \emph{Khovanskii basis}, and its associated
graded algebra is a normal affine semigroup algebra. In particular,
$A$ is finitely generated over $k$.
\end{theorem}
\begin{proof}
For $n\ge1$, an element of $V(n\Delta)$ with least exponent $a$
exists if and only if $Q_a\subseteq n\Delta$. This is equivalent to
\[ \langle a,v\rangle\ge nc_j\qquad (v\in\operatorname{Vert}(T_{p_j}),\ 1\le j\le l). \]
At $n=1$ these inequalities define a rational polyhedron.
Every rational solution lies in $\Delta$: clearing denominators
reduces this to the integral case, where $a\in Q_a$. Hence this
set is bounded, by density of its rational points. The \emph{value
semigroup} consists of the pairs $(n,e(f))$ with $0\ne f\in V(n\Delta)$.
The rational polytope just obtained is the degree-one slice of the
closed cone generated by these pairs.
At degree zero the homogeneous inequalities also force $a=0$:
for an integral solution they imply $Q_a$ lies in the recession
cone of $\Delta$, the set of directions $v$ with
$\Delta+\mathbb R_{\ge0}v\subseteq\Delta$. Since $\Delta$ is
bounded, this cone is $\{0\}$. Rational polyhedrality
gives the same conclusion for real solutions.

The value semigroup is thus all lattice points of the indicated
rational polyhedral cone in degrees $n\ge0$. Choose its finite
\emph{Hilbert basis}, its minimal set of additive generators, by
Gordan's lemma. For each generator choose a
homogeneous section with that least exponent and normalize its
least coefficient to one. The associated graded algebra for degree
and least exponent is the semigroup algebra, since every nonzero
initial space is one-dimensional. Saturation of the semigroup gives
normality. Given a nonzero homogeneous section, choose a product
of the lifts of the same degree and least exponent. Subtract the
scalar multiple that cancels the least term. Each subtraction strictly
increases the least exponent in the finite set $n\Delta\cap M$,
so the lifts generate the section algebra. Finally, choose a
polytope $\Delta$ containing $0$ in its interior.
Every element of $A$ then occurs in some $V(n\Delta)$, and
evaluation at $t=1$ maps the section algebra onto $A$.
\end{proof}

\begin{corollary}\label{cor:polynomial-polyhedral-newton-sections}
Let $A$ be the first-star intersection over
$k[u_1^{\pm1},\ldots,u_c^{\pm1},s_1,\ldots,s_l]$. If its
localization $A^\circ=A[(s_1\cdots s_l)^{-1}]$ admits a regular
function with a Newton vertex negative in all mutable coordinates,
then its original polynomial-frozen section algebras and $A$
itself are finitely generated. The section algebras satisfy the
polyhedral value-semigroup description of
Theorem~\ref{thm:polyhedral-newton-sections} after imposing the
following inequalities on an exponent $b\in\Delta$. Here $b_{s_j}$
and $b_i$ are its $s_j$- and $x_i$-coordinates, and
$\operatorname{ord}_{s_j}(E_i)$ is the least exponent of $s_j$
in $E_i$.
\[ b_{s_j}\ge0,\qquad b_{s_j}+\operatorname{ord}_{s_j}(E_i)b_i\ge0 \qquad(i\text{ mutable},\ 1\le j\le l). \]
Here the polytopes $T_p$ are computed in the localized first-star
intersection, whether or not an original exchange becomes a unit.
\end{corollary}
\begin{proof}
Localization commutes with this finite intersection. Indeed, for
an element of every localized chart ring, a common power of
$s_1\cdots s_l$ clears the frozen denominators in all the charts.
Let $H$ be the cone cut out by the stated inequalities.
For $f\in A^\circ$, polynomiality in the initial chart
requires precisely the first group of inequalities on its support.
In mutable slice $-h<0$ in direction $i$, write its coefficient
as $E_i^h B$, where $B$ is Laurent in the unchanged
coordinates, including the localized frozen variables. For each
$j$ its frozen order is the minimum frozen exponent of the
original slice minus $h\operatorname{ord}_{s_j}(E_i)$.
Thus $B$ is polynomial in all frozen variables precisely when
the second inequalities hold on that slice. In nonnegative mutable
slices the adjacent coefficient is obtained by multiplication by
a nonnegative power of $E_i$, which introduces no frozen
poles. The finite first star now identifies the sections in $A$
supported in $n\Delta$ with those in $A^\circ$ supported in
$n(\Delta\cap H)$. The theorem applies to this rational polytope.
If the intersection is empty, the section algebra is $k$ in
degree zero. Taking an exhausting polytope again gives finite
generation of $A$.
\end{proof}

For a projective log canonical Calabi--Yau pair under the hypotheses
of \cite[Theorem~1.2]{BL}, Blum and Liu construct a basis adapted to
all valuations of log discrepancy zero. They also show that the
resulting tropical functions are independent of the choice of such a
basis \cite[Lemma~3.20]{BL}. In the LP setting the Cartier argument
gives the smallest Newton polytopes directly from the exchange
polynomials, including in the presence of polynomial frozen variables.

\section{Intrinsic multiplication from the exchange polynomials}
\label{sec:incidence}

\subsection{First-star charts for the curve counts}\label{subsec:counting-hypotheses}

To compute traces of theta products, we will express the condition that
a curve extends to the variety $U$ in terms of the exchange polynomials.
Let $U$ be a smooth connected affine log Calabi--Yau variety over
an algebraically closed field $k$ of characteristic zero. Fix an LP
seed with initial torus $T$, mutable coordinates $x_1,\ldots,x_r$
and designated units $u_1,\ldots,u_c$. For a regular function $h$
on $U$, write $U_h=\{h\ne0\}$ for its principal open subset.
Throughout this section we assume the following.
\begin{enumerate}
\item The $x_i$ are regular on $U$, the $u_j$ are units on $U$,
and $T=U_{x_1\cdots x_r}$. The logarithmic volume form on $T$
extends to a nowhere vanishing volume form on $U$.
\item The union $U_\star$ of $T$ and its $r$ adjacent tori
has complement of codimension at least two in $U$.
\item For each $i$, inverting the other mutable coordinates gives
\[ U_{\prod_{j\ne i}x_j}=\Spec C_i[X,Y]/(XY-F_i),\qquad X=x_i,\quad Y=F_i/x_i. \]
Here $F_i$ is the ordinary exchange polynomial and
$C_i=k[u_1^{\pm1},\ldots,u_c^{\pm1},x_j^{\pm1}\ (j\ne i)]$.
\end{enumerate}
The associated exchange Laurent polynomial differs from $F_i$ by a
unit in $C_i$, so $X\ne0$ and $Y\ne0$ give the initial and
adjacent seed tori in this open set. Normality and the codimension
assumption give $k[U]=\mathcal U_\star(t)$. The cocharacter lattice
$L=\mathbb Z^{r+c}$ identifies with the integral points of the
essential skeleton \cite[Lemma~2.2]{KY}. For $P\in L$, write
$P(v)$ for the order of a character $v=\chi^b$ at $P$, namely
$\langle b,P\rangle$.

\subsubsection*{Smooth upper LP sources}
Suppose $U$ is the smooth spectrum of a full upper LP algebra over
$k[u^{\pm1}]$ and its marked seed has a polynomial corner. Then the
preceding hypotheses hold. Indeed, at a polynomial corner we may
multiply each ordinary exchange $F_i$ by a scalar and a Laurent
monomial in the designated units to obtain a polynomial $G_i$
with $G_i(0)=1$. Each $G_i/x_i$ is a first
mutation partner times a coefficient unit and a monomial in the
other mutables with nonnegative exponents, so it is regular. For
a strictly positive weight, the product $\prod_i(G_i/x_i)$ has
a unique minimizing exponent, with every mutable entry equal to $-1$.
Theorem~\ref{thm:polyhedral-newton-sections} therefore gives finite
generation. Proposition~\ref{prop:upper-factorial} gives the
one-exchange charts and their coverage in codimension one. As in
Section~\ref{sec:upper}, the logarithmic torus form extends without
zeros or poles. This gives the log Calabi--Yau condition by the
criterion preceding \cite[Definition~2.1]{KY}.

\subsection{Equations for the trace}

In the intrinsic algebra the monomial $q^\beta$ records the numerical
curve class $\beta$. Setting all $q^\beta=1$ gives $B_{\mathbb Z}(U)$,
with basis $(\theta_m)_{m\in L}$ and unit $\theta_0$. Its trace
$\Tr$ takes the coefficient of this unit.
For $a\ge2$, write $\mathsf T(P_1,\ldots,P_a)=\Tr(\theta_{P_1}\cdots\theta_{P_a})$.
Each $P_\ell$ prescribes a \emph{contact order}, the vector of
orders of the torus coordinates at the $\ell$-th marked point.
By the Frobenius structure theorem, this integer counts maps from a
generic marked projective line with these contact orders and with
a further marked point mapping to a general point of $U$
\cite[Theorem~1.2 and Proposition~3.13]{KY}.\footnote{In the smooth affine
cluster case, the corresponding descendant insertion is
$\psi_{\mathrm{ev}}^{a-2}$ for these $a+1$ marks, where
$\psi_{\mathrm{ev}}$ is the first Chern class of the cotangent line at
the evaluation mark. Adding Mandel's free mark raises the exponent to
$a-1$ without changing the invariant, by the fundamental class axiom
\cite[Definition~1.2 and Remark~1.4]{MandelLogGW}. With $a+1$ marks and
nonzero input labels, the corresponding tropical curves have an evaluation vertex of
valency $a+1$, including its contracted marked end
\cite[Example~1.3]{MRTropical}.}

Take the evaluation mark to be $\infty$, and write $p_1,\ldots,p_a$
for the contact marks in $\mathbb A^1$. We may set $p_1=0$ and
$p_2=1$. The remaining marks and the coordinates
$(b_1,\ldots,b_r,b^u_1,\ldots,b^u_c)$ of the evaluation point are
algebraically independent. All degrees below are taken over their
rational function field. Define the mutable \emph{defects} and the total orders of the units by
\[ \delta_i=-\sum_{\ell=1}^{a}P_\ell(x_i),\qquad \epsilon_j=\sum_{\ell=1}^{a}P_\ell(u_j). \]
If any $\delta_i<0$ or any $\epsilon_j\ne0$, then the trace vanishes.
Indeed, away from the marked contacts a regular coordinate can only
have zeros, and a unit has neither zeros nor poles. The divisor of
each pulled-back coordinate on $\mathbb P^1$ has degree zero.

Assume henceforth that $\delta_i\ge0$ and $\epsilon_j=0$.
Use $t$ as the affine coordinate on the marked line.
Take monic polynomials $A_i$ of degree $\delta_i$, including
$A_i=1$ when $\delta_i=0$. Write
\[ X_i(t)=b_iA_i(t)\prod_\ell(t-p_\ell)^{P_\ell(x_i)},\qquad U_j(t)=b^u_j\prod_\ell(t-p_\ell)^{P_\ell(u_j)}. \]
For each $i$, substitute these functions in $F_i$ and clear
the marked poles to form $N_i(t)$. More precisely, multiply the
substitution by $\prod_\ell(t-p_\ell)^{e_{\ell i}}$, with
$e_{\ell i}=\max(0,-\min_{\mu\in\Supp(F_i)}\langle\mu,P_\ell\rangle)$.
The exponents of the mutable variables in $F_i$ are nonnegative,
so $N_i$ is a polynomial. The clearing factors are units at the
unmarked roots of the $A_i$.

\begin{theorem}\label{thm:incidence}
The trace $\mathsf T(P_1,\ldots,P_a)$ is the degree of the reduced
solution scheme of
\begin{equation}\label{eq:trace-congruences}
N_i\equiv0\pmod{A_i}\qquad(1\le i\le r)
\end{equation}
on the open locus defined by the following conditions.
\begin{enumerate}
\item The $A_i$ are squarefree and pairwise coprime, and none
 vanishes at a marked point.
\item At each root of $A_i$, the regular value of
      $F_i(X,U)/X_i$ is nonzero.
\item If $P_\ell=d e_i$ with $d>0$ and $e_i$ a mutable
 coordinate vector, the unchanged coordinates at $p_\ell$
 lie outside the hypersurface $F_i=0$.
\end{enumerate}
The reduced scheme is finite \'{e}tale over the generic parameter
field. Its equations involve $\sum_i\delta_i$ unknown coefficients
and the same number of scalar remainder equations.
\end{theorem}

\begin{proof}
Choose a smooth complete toric fan for $T$ containing the positive
mutable rays and the rays of the nonzero $P_\ell$. Remove the
toric strata of codimension at least two. In the open divisor orbit
of the $i$-th positive mutable ray, blow up the locus $F_i=0$.
These centers are disjoint. Each nonempty center is smooth, since a
singular point of $F_i=0$ would give a singular point of $XY=F_i$
at $X=Y=0$. In the blowup chart with coordinate $Y=F_i/X$ we
recover the specified one-exchange model. The union of the initial
and adjacent tori is therefore $U_\star$, with its original
coordinates.

Enlarge $U_\star$ by the open divisor orbits of the contact rays.
For a positive mutable contact ray take the part of the strict
transform lying over $F_i\ne0$. At every other contact ray take the
whole open orbit. Equivalently, adjoin the corresponding one-ray
toric opens, with the exchange center deleted on a positive mutable
ray. The resulting open subvariety of the blowup is smooth. Denote
it by $W$ and its boundary by $D=W\setminus U_\star$. The
components of $D$ are smooth and disjoint. The torus volume
trivializes $K_W(D)$. Along a generic point of an exceptional
divisor the blowup discrepancy is one and its center lies on a
single toric pole of order one, so that the pulled-back volume has
order zero there.

The congruences and the first two conditions extend the
rational map $(X(t),U(t))$ to $U_\star$ away from its contacts.
At a root of $A_i$ it lies in the adjacent $i$-th torus, with
coordinate $F_i/X_i$. At a nonzero contact, write $P_\ell=dv$
with $v$ primitive. Characters perpendicular to $v$ have nonzero
finite values at $p_\ell$. A character pairing to one with $v$
has order $d$. These orders extend the map across
$p_\ell$ into the open orbit belonging to $v$ with the required
multiplicity. At a positive mutable ray, condition (3) places the image in the chosen open orbit. A zero contact is an interior marked
point. Conversely, a map into this pair meeting the interior
coordinate divisors simply, separately and away from the marks
has the indicated unique expressions for its coordinates. This
description holds in families. The marked factors specify the pole
bounds on the fixed projective line, the evaluation specifies the
leading coefficients, and division by a monic $A_i$ is defined
over the parameter ring. Thus the equations describe an open part
of the fixed-domain space of maps with prescribed contacts.

Take a smooth projective simple normal crossings compactification
of $U$, extracting the divisors of the finitely many contact
valuations by modifications of the boundary as in
\cite[\S1, definition of the Frobenius multilinear form]{KY}.
Remove the other boundary divisors and the intersections of the contact divisors. Denote the resulting open by $Y'$.
The pairs $W$ and $Y'$ have a common open subset whose complements
have codimension at least two. On the interiors this follows from
the hypothesis on $U_\star$. At the generic point of each remaining
boundary divisor both local rings are the discrete valuation ring (DVR)
of the same marked monomial valuation. Indeed, the generators of an affine neighborhood on either model
belong to the common valuation ring, so the birational maps in both
directions are regular near the generic points. Their restrictions
identify the boundary divisors there.

We show that every generic solution lies in this common open.
Fix the marked domain and allow the evaluation point $b$ to vary.
The equations above define a family of finite type. Consider a
component of its reduction which dominates $T$. At a general
smooth point of this component, generic smoothness gives surjectivity
of the evaluation differential. For the corresponding map $g$,
fixed-contact infinitesimal deformations give sections of
$V=g^*T_W(-\log D)$, where $T_W(-\log D)$ is the sheaf of vector
fields tangent to $D$. In a coordinate normal to a contact divisor
this says that a variation of $(t-p_\ell)^d$ times a unit is still
divisible by $(t-p_\ell)^d$. The evaluation of sections of $V$
is surjective at $\infty$, and $\deg V=0$ since $K_W(D)$
is trivial. Splitting $V$ on $\mathbb P^1$ first shows that
each summand has nonnegative degree and then that $V$ is trivial.

Local deformations with fixed contacts on the marked line lift across a square-zero thickening in smooth coordinates,
using a unit multiple of $(t-p_\ell)^d$ for the normal parameter.
Differences of lifts on overlaps are sections of $V$ tensored with
the square-zero ideal. Thus the obstruction lies in $H^1(V)=0$,
and the tangent space is $H^0(V)$. Since the equations describe
this open mapping functor, it is smooth of dimension $\dim W$
near $g$. Thus every component dominating $T$ has dimension
$\dim T$, and the complement of this smooth locus in that component
has nondense evaluation image. Evaluation at any fixed interior domain point is
\'{e}tale there. With the interior point varying, evaluation is
smooth of relative dimension one. The incidence of such maps with
a fixed subset of codimension two in $W$ therefore has dimension
at most $\dim W-1$. At a boundary mark the evaluation differential
is the surjection from $H^0(V)$ to the fiber of $V$, followed by
the surjection onto the tangent space of its boundary divisor.
Consequently the maps sending this mark into a proper closed subset
of that divisor have dimension at most $\dim W-1$. The same
estimate holds at an interior mark by interior evaluation.
The excluded subsets do not dominate $T$. There are only finitely many dominant components, so their generic fibers avoid all these subsets.

In the other direction, the maps counted on the proper
compactification avoid the complement of the common open and meet
the reduced union of the interior coordinate divisors transversely
away from their marks \cite[Lemmas~3.6, 3.9 and~3.12]{KY}. On
$U_\star$ the zero divisor of $x_i$ is reduced, since it is cut
out by the smooth hypersurface $F_i=0$ in the adjacent torus.
Distinct initial coordinates cannot vanish together there.
We may therefore apply these avoidance and transversality statements
to the closures of the interior divisors and their pairwise
intersections. Every counted map belongs to the incidence scheme.
Only finitely many classes occur \cite[Lemma~3.8]{KY}. The two
generic reduced fibers consequently agree. Their degree, by
\cite[Theorem~1.2 and Proposition~3.13]{KY}, is the required trace.
\end{proof}

The degree in Theorem~\ref{thm:incidence} is computable by
ordinary polynomial arithmetic over an effective characteristic-zero
field of definition. Use discriminants, evaluations at the marks and
resultants for the three open conditions. For condition (2) take the
resultant of $A_i$ and the quotient of $N_i$ by $A_i$ for each
$i$ with $\delta_i>0$. Saturate the remainder ideal by the
product of the open factors. The resulting algebra over the generic
parameter field is finite-dimensional. The rank of its trace
pairing is its reduced degree. Indeed, after passing to an algebraic
closure, multiplication by $h$ on a local Artin factor of length
$l$ has trace $l h(\xi)$, where $\xi$ is its closed point.
In characteristic zero the pairing therefore has rank one on each local factor.

\subsection{Recovering the structure constants from traces}

Choose a coercive regular function $f$ as in Proposition~\ref{prop:coercive}.
Its Newton polytope contains $0$ in its interior. If $\tau_f$ is its
minimum support function, the sublevel sets of $\rho_f=-\tau_f$ are compact
rational polytopes.

\begin{theorem}\label{thm:corner-free-product}
Let $f$ be coercive, and write $\mathbb N=\mathbb Z_{\ge0}$.
For $p,q\in L$, form the finite set
\[ S=\{s\in L:\ \rho_f(s)\le\rho_f(p)+\rho_f(q),\ s-p-q\in\mathbb N^r\times\{0\}^{c}\}. \]
Order $S$ by any linear extension of coordinatewise order in the
mutable coordinates. Then
\begin{equation}\label{eq:theta-trace-system}
\theta_p\theta_q=\sum_{t\in S} n_t\theta_t,\qquad \sum_{t\in S}\mathsf T(t,-s)n_t=\mathsf T(p,q,-s)\quad (s\in S).
\end{equation}
The coefficient matrix is lower triangular with diagonal entries one.
Moreover $n_{p+q}=1$.
Consequently Theorem~\ref{thm:incidence} computes all the
integers $n_t$ from the exchange polynomials and the Laurent
support of $f$, for every pair of integral labels.
\end{theorem}

\begin{proof}
Let $g$ be a nonzero regular function. Whenever $\theta_s$ occurs in
$\theta_p\theta_q$, convexity for the principal divisor
$-\operatorname{div}(g)$ gives the \emph{regular-function inequality}
of Keel and Yu~\cite[Theorem~15.7(2)]{KY},
\begin{equation}\label{eq:regular-function-inequality}
\tau_g(s)\ge\tau_g(p)+\tau_g(q).
\end{equation}
Apply this first to $f$, then to the
mutable $x_i$, and to the units $u_j$ and $u_j^{-1}$.
These inequalities bound the support. Taking the trace after
multiplication by $\theta_{-s}$ gives the linear system. A nonzero entry $\mathsf T(t,-s)$ requires $t\le s$ in the
mutable coordinates, by the defect condition. If the labels of a trace
sum to zero, then all $A_i=1$. At a contact $P_\ell=d e_i$ with $d>0$ and $1\le i\le r$,
the unchanged coordinates are nonzero scalar multiples of the
algebraically independent evaluation parameters over the field of the
marks, so $F_i\ne0$. Thus every open condition of
Theorem~\ref{thm:incidence} holds and the trace is one, including the
constant map when all contacts vanish. This gives the diagonal entries.
Taking $s=p+q$
in the system leaves only $n_{p+q}$ on the left, and the three
labels on the right sum to zero. Hence $n_{p+q}=1$.
\end{proof}

When the seed is given over an effective field, one can also find
$f$ without a list of generators of $k[U]$. Intersect the space
of Laurent polynomials supported in $[-b,b]^{r+c}$ with each
adjacent Laurent ring. To do this on the $i$-th chart, expand in
powers of $x_i$. For $d>0$ the coefficient of $x_i^{-d}$
must be divisible in $C_i$ by $F_i^d$, since $x_i=F_i/Y$.
These are linear conditions on the coefficients in the box. They
are computed by reduction modulo $(F_i^d)$ in the unchanged
Laurent variables. Normality of $U$ and codimension-two coverage
identify this intersection with the regular functions supported in
the box. Increase $b$ until the supports of a basis of the
intersection have $0$ in the interior of their joint convex hull.
Proposition~\ref{prop:coercive} gives termination. A rational linear
combination chosen off the finitely many coefficient cancellation
hyperplanes then gives $f$. This constructs the multiplication from the ordinary exchanges over any effective field of definition.

\begin{lemma}\label{lem:coercive-theta-generators}
Choose a coercive regular
function $f$ and set $\rho=-\tau_f$. The theta functions
indexed by the Hilbert bases of the maximal cones of linearity
of $\rho$ generate the algebra over $\mathbb Z$.
\end{lemma}
\begin{proof}
Each cone is pointed because $\rho$ is positive away from zero.
Write a label $m$ as a sum of Hilbert generators in one cone.
The product of their theta functions has coefficient one at $m$.
Every other label has the form $s=m+d$, where $d$ is a nonzero
nonnegative mutable shift with zero unit coordinates.
Moreover $\rho(s)\le\rho(m)$, by the regular-function
inequality and linearity on the cone. Order the labels by increasing
$\rho$, and within each level by decreasing sum of mutable
coordinates, breaking remaining ties arbitrarily. The levels are
finite and $\rho$ is nonnegative and integral. Every other label in the product precedes $m$.
Subtracting the corresponding terms proves generation by induction.
\end{proof}

\subsection{A presentation and the singular locus of the forward mirror}

Work in the augmented algebra of Theorem~\ref{thm:corner-free-product}.
Choose a coercive regular function $f$, write $\rho=-\tau_f$, and
let $\Sigma$ be its fan of linearity in the label lattice $L$.
Set $s(p)=\sum_{i=1}^r p_i$, summing only mutable coordinates.
Choose an integer $N$ such that
$w(p)=N\rho(p)-s(p)$ is strictly positive for $p\ne0$.
It suffices to check the primitive rays of $\Sigma$.
Let $H$ be the union of the Hilbert bases of the maximal cone monoids.
Give $X_h$ weight $w(h)$ in the polynomial ring
$\mathbb Z[X_h:h\in H]$.

Write $z^h$ for the basis character of the group ring $\mathbb Z[L]$
with exponent $h$. Form the following finite set of homogeneous relations. For each
maximal cone take integral binomial generators of the relations
among its characters $z^h$, $h\in H$ in that cone. For each
inclusion-minimal subset of $H$ contained in no common cone, take the product
of its variables. Let $D$ be the maximum of zero, the variable
weights, and the weights of these relations, and set
$S_D=\{p\in L:w(p)\le D\}$.

\begin{theorem}\label{thm:finite-forward-presentation}
The traces of pairs and triples with labels in $S_D\cup(-S_D)$
determine a finite presentation of the augmented integral algebra
on the generators $\theta_h$, $h\in H$. Each homogeneous
relation above lifts to a relation among these generators by subtracting
terms of strictly smaller weight. The lifts generate the defining ideal.
The associated graded algebra is the toric face ring of $\Sigma$.
It has basis $(e_p)_{p\in L}$ and product $e_pe_q=e_{p+q}$
when $p,q$ belong to a common cone of $\Sigma$, and zero otherwise.
\end{theorem}
\begin{proof}
Every label $b$ occurring in $\theta_p\theta_q$ is $p+q+d$ with
$d$ a nonnegative mutable shift and zero unit part. The
regular-function inequality for $f$ gives $w(b)\le w(p)+w(q)-s(d)$.
Thus every nonzero shift strictly lowers the weight. Theorem~\ref{thm:corner-free-product} gives coefficient one at $p+q$. Its weight is
$w(p)+w(q)$ exactly when the minima of the support of $f$
at $p$ and $q$ have a common minimizing exponent, equivalently
when the two labels share a cone. Filtering by the spans of labels
of weight at most $j$ gives the asserted graded product.

The relations given above present the toric face ring over $\mathbb Z$.
Indeed, every monomial whose indices have no common cone vanishes. If two
surviving monomials have the same sum of labels, all their indices
lie in the smallest face containing this sum. Their difference
therefore belongs to a single cone's toric ideal. Grouping by the
sum of labels now gives all relations.

For each $p\in S_D$, decompose $p$ in the Hilbert generators of
a cone containing it. The corresponding product lifts $e_p$ and
has total weight $w(p)$. Subtract its terms of lower weight by
induction to obtain a polynomial of weight at most $w(p)$
representing $\theta_p$. The same argument in every degree proves
generation. Evaluate each relation of degree $j\le D$ in the
counting algebra and use these polynomials to subtract its value,
which has weight strictly less than $j$. The lifted relations
generate the kernel. For a relation, its top homogeneous part
belongs to the toric face ideal. Express this part in the homogeneous
generators at its degree and subtract the corresponding multiples of
the lifts. The top degree decreases, so the process terminates.

In each evaluation just described, multiply one generator at a time.
The weights of the two factors have sum at most $D$, and each product term lies
in $S_D$. In the trace system~\eqref{eq:theta-trace-system} it suffices to take
$b\in S_D$ with $b-p-q$ a nonnegative mutable shift and zero
unit part. This set contains the support of the product. The traces of pairs
give a unitriangular matrix on it. Its entries use only the
labels in the statement. In particular every incidence scheme used
has at most $3\max_{p\in S_D}\sum_{i=1}^r|p_i|$ unknown monic
coefficients.
\end{proof}

\begin{corollary}\label{cor:forward-jacobian}
Assume the hypotheses of Section~\ref{subsec:counting-hypotheses}
and suppose the seed is given over an effective field of definition
of characteristic zero. The ordinary
exchanges determine equations for the singular locus of the augmented
forward mirror and an algorithm to decide whether it is smooth.
\end{corollary}
\begin{proof}
The search following Theorem~\ref{thm:corner-free-product} gives a
coercive function. Theorem~\ref{thm:incidence} computes the traces
needed for the presentation of Theorem~\ref{thm:finite-forward-presentation}.
Over the ground field the algebra is a domain by the least-label
product rule. Its dimension is $n=\operatorname{rank}L$. Indeed, the rank of the
filtration through weight $j$ counts the lattice points with
$w(p)\le j$, which grows on the order of $j^n$. The weighted
presentation realizes this filtration by polynomial degree.
By the Jacobian criterion, adjoining the minors of size $|H|-n$
of the Jacobian of the lifted relations gives the ideal defining the
singular locus. We interpret a minor of size zero as one. The mirror
is geometrically smooth precisely when the resulting ideal contains one.
\end{proof}

One can also bound the integer coefficients needed in these
products before solving the incidence systems. By positivity and
$\mathsf T(b,-b)=1$, the coefficient at $b$ in
$\theta_p\theta_q$ is at most $\mathsf T(p,q,-b)$.
For each trace in Theorem~\ref{thm:finite-forward-presentation} impose its open
conditions by adjoining one variable inverse to the product of
the open factors. The resulting polynomial system has as many
equations as unknowns. Its reduced generic degree is bounded by
the product of the total degrees of these equations by the affine B\'ezout bound, with zero degrees replaced by one. Vanishing traces
with negative defect or unbalanced unit orders require no system.

\subsection{Changing the seed in the incidence equations}

Mutation induces a change of coordinates on the generic solution
schemes of the incidence equations~\eqref{eq:trace-congruences}. Write $E=\widehat F_i$ and
$x_i'=E/x_i$. By Lemma~\ref{lem:boundary-gauss-change}, the Gauss
valuation of weight $P$ has weight $P'$ in the adjacent chart,
where $P'(x_i')=\tau_E(P)-P(x_i)$ and the other coordinates agree.

\begin{proposition}\label{prop:incidence-mutation}
Suppose the two seeds give first-star charts on the same smooth
variety $U$ as in Theorem~\ref{thm:incidence}. Identify their generic
parameter fields by $b_i'=E(b)/b_i$, with all other evaluation
coordinates and all marks fixed. For contact orders $P_1,\ldots,P_a$
and their transforms $P_1',\ldots,P_a'$, the reduced generic
incidence schemes are isomorphic. On the locus of solutions,
the unchanged monic polynomials are the same. The polynomial for the
mutated variable is
\[ A_i'(t)= \frac{E(X(t),U(t))} {b_i' X_i(t)\prod_{\ell}(t-p_\ell)^{P_\ell'(x_i')}} . \]
\end{proposition}
\begin{proof}
Choose a single boundary refinement representing the finitely many
marked divisorial valuations. For each star the proof of
Theorem~\ref{thm:incidence} identifies the geometric generic fiber
with the maps in this compactification with those contacts and
evaluation. Intersect the two generic conditions on the domain and
evaluation. In particular the contact points avoid proper closed
subsets of their marked orbits. Thus the orders of $x_i'$ at the
marks are its divisorial orders $P_\ell'(x_i')$. If the residue
of its initial form vanishes on a proper subset of such an orbit,
that subset is among those avoided. At a zero contact the mark is
an interior point and we use ordinary evaluation there.

On the common maps $E(X,U)/X_i$ is the pullback of the regular
function $x_i'$. After removing the prescribed marked factors it
is a polynomial in $t$, and its leading coefficient is $b_i'$
by evaluation at infinity. All its remaining zeros are simple,
away from the marks and the zeros of the unchanged polynomials, by
the second star's generic transversality. The inverse exchange, with the same choice of coefficient units,
recovers the original functions. These inverse bijections of geometric
points commute with the Galois action. The two reduced generic schemes
are finite \'{e}tale, so the bijections identify them over the common
parameter field. In particular, a negative defect in either seed makes
both fibers empty.
\end{proof}

The normalizing monomial in $E$ can cancel zeros of the substituted
numerator at zeros of unchanged coordinates. The formula records these
cancellations before the remaining zeros are assigned to the mutant.

\subsection{Products with a split torus}

\begin{lemma}\label{lem:intrinsic-torus-product}
Let $V$ be a smooth connected affine log Calabi--Yau variety over an
algebraically closed field of characteristic zero with a volume-compatible
dense split torus. Let $S$ be a split torus with cocharacter lattice $N$.
For the product volume form, there is an isomorphism
\[ B_{\mathbb Z}(V\times S)\simeq B_{\mathbb Z}(V)\otimes_{\mathbb Z}\mathbb Z[N],\qquad \theta_{(p,a)}\longmapsto\theta_p\otimes z^a, \]
where $p$ is an integral point of the essential skeleton of $V$ and
$z^a$ is the group algebra element indexed by $a\in N$.
\end{lemma}
\begin{proof}
It suffices to adjoin one copy of $\mathbb G_m$ at a time. Choose a
smooth projective compactification $Y$ of $V$ with simple normal
crossings boundary and use $Y\times\mathbb P^1$ for the product.
On the marked projective line in the Frobenius count, let $a_i$ be
the order of the last coordinate at the mark $x_i$. A rational
function with divisor $\sum_i a_i[x_i]$ exists precisely when
$\sum_i a_i=0$. Its value at the evaluation mark determines it
uniquely, and its degree is $\sum_i\max(a_i,0)$. A zero contact
in either factor is an ordinary interior mark in that factor.

Projection therefore identifies the generic counting fibre for the
product with the fibre for $V$ whenever the orders balance and the
last degree has this value. In all other cases the product fibre is empty.
These counts are taken in the product compactification, with curve
classes pushed forward from any boundary refinement used to represent
the contacts \cite[Proposition~17.3(1)]{KY}. An effective numerical
class on $Y\times\mathbb P^1$ is specified by its class on $Y$ and
its degree on $\mathbb P^1$, and every such pair is effective.
Thus the Frobenius forms agree, class by class, with those of the
tensor product for the two factors. They determine multiplication
over the effective coefficient ring \cite[Theorem~1.2]{KY}.
Setting the curve monomials to one gives the asserted formula, since
the augmented algebra of $\mathbb G_m$ is its cocharacter group
algebra \cite[Lemma~9.4]{KY}.
\end{proof}

\section{Laurent charts for the forward theta functions}
\label{sec:polynomial-corner}

\subsection{The chart at a polynomial corner}

We keep the hypotheses of Theorem~\ref{thm:incidence}. Write
$v=(x_1,\ldots,x_r,u_1,\ldots,u_c)$ for the coordinates on the
initial torus and $L=\mathbb Z^{r+c}$ for its cocharacter lattice.
Set $n=r+c$. We say that the marked seed has a \emph{polynomial
corner} if for each mutable $i$ there is a Laurent monomial
$\rho_i$ in the designated units, times a nonzero scalar, such that
$G_i=\rho_iF_i$ belongs to $k[v_1,\ldots,v_n]$ and $G_i(0)\ne0$.

\begin{theorem}\label{thm:corner}
Suppose the marked seed has a polynomial corner.
Set $z_j=\theta_{e_j}$ for the standard basis of $L$.
The elements $z_j$ are algebraically independent, and localization gives
\[ B_{\mathbb Z}(U)[(z_1\cdots z_n)^{-1}]=\mathbb Z[z_1^{\pm1},\ldots,z_n^{\pm1}]. \]
In this Laurent ring the image $\Theta_m$ of $\theta_m$ is
\[ \Theta_m=z^m A_m(z_1,\ldots,z_r),\qquad A_m=\sum_{0\le d_i\le D_i(m)} c_{m,d}\,z^{(d,0)},\qquad D_i(m)=-\tau_{G_i}(m). \]
Here $d\in\mathbb N^r$, and $(d,0)$ appends $c$ zero unit
coordinates to $d$. The coefficients $c_{m,d}$ are nonnegative integers and $c_{m,0}=1$.
In particular $\Theta_m=z^m$ on the nonnegative orthant.
\end{theorem}

\begin{proof}
Theorem~\ref{thm:corner-free-product} gives coefficient one at
$\theta_{p+q}$ in $\theta_p\theta_q$. All other labels have the
form $p+q+(d,0)$ with $d\in\mathbb N^r\setminus\{0\}$.

To bound these corrections from above, we use the regular function $t_i=G_i/x_i$.
Indeed $F_i/\widehat F_i$ is a monomial with nonnegative
exponents in the other mutable coordinates, up to a coefficient
unit, and $\widehat F_i/x_i$ is the actual LP mutation partner.
The Laurent phenomenon makes this partner regular on the union of
the seed tori in Theorem~\ref{thm:incidence}. It extends to $U$
by normality and the codimension hypothesis.
For every $a$ in the nonnegative orthant we have
$\tau_{t_i}(a)=-a_i$, since $G_i$ has a nonzero constant
term. If $p,q$ lie in that orthant, so does every label occurring in their product. Applying
convexity for $t_i$ gives $-p_i-q_i-d_i\ge-p_i-q_i$ for all
$i$. The resulting multiplication rule is
$\theta_p\theta_q=\theta_{p+q}$. In particular the $z_j$
generate a polynomial subalgebra, since its distinct monomials have
distinct basis labels.

Now fix $m$, and take $b=N(1,\ldots,1)$ with $b+m$ in the
positive orthant. Every possible label in the product
$a=b+m+(d,0)$ of $\theta_b\theta_m$ is again in that orthant.
Applying convexity to $t_i$ gives
$-b_i-m_i-d_i\ge-b_i+\tau_{G_i}(m)-m_i$, or $d_i\le D_i(m)$.
The translated product is therefore a finite sum of the stated
monomials $z^{b+m+(d,0)}$, with nonnegative integral coefficients
and with coefficient one when $d=0$. It follows that the localization is generated by the $z_j^{\pm1}$.
Multiplication by $z_1\cdots z_n$ is injective. Indeed, order the
labels first by the sum of their mutable coordinates and then
lexicographically. Translation preserves this order, and any
nonzero shift $(d,0)$ strictly increases it. In a finite nonzero
theta sum the least label thus contributes the unique least term
after multiplication by $\theta_{(1,\ldots,1)}$, with unchanged
coefficient. It follows that $B_{\mathbb Z}(U)$ embeds in the localization.
Dividing the translated product by $z^b$ gives the formula.
\end{proof}

We call $A_m=z^{-m}\Theta_m$ the \emph{normalized polynomial factor}
of the theta function in this chart. For a chosen chart at a marked
seed $t$, we write $\Theta_m^t$ and $A_m^t$ when the seed varies.
The factor $A_m^t$ plays the role of an $F$-polynomial.
In a cluster algebra, Fomin and Zelevinsky obtain the $F$-polynomial
by setting the initial cluster variables to one in the algebra with
principal coefficients~\cite[Definition~3.3]{FZIV}. Here the variables
of $A_m^t$ are the mutable coordinates on the forward mirror.

\begin{corollary}\label{cor:trace-reconstruction}
For $m\in L$, let $\mathcal D_m$ be the integral box
$0\le d_i\le D_i(m)$. Take $b=N(1,\ldots,1)$ with
$N=1+\max(0,-m_1,\ldots,-m_n)$, and write $s_d=b+m+(d,0)$.
The coefficients of $A_m$ are the unique solution of
\[ \sum_{e\in\mathcal D_m}\mathsf T(s_e,-s_d)c_{m,e} =\mathsf T(b,m,-s_d) \qquad(d\in\mathcal D_m). \]
In any linear extension of componentwise order this is a unit lower
triangular integral system. If $s=\prod_i(1+D_i(m))$, its evaluation
uses at most $s$ traces of triples and $s(s-1)/2$ off-diagonal
traces of pairs, each with at most $\sum_iD_i(m)$ polynomial
coefficients to determine in Theorem~\ref{thm:incidence}.
\end{corollary}

\begin{proof}
By Theorem~\ref{thm:corner}, the translated product is
$\theta_b\theta_m=\sum_e c_{m,e}\theta_{s_e}$.
Taking the trace after multiplication by $\theta_{-s_d}$ gives
the equations. The mutable defect of $(s_e,-s_d)$ is $d-e$,
and that of $(b,m,-s_d)$ is $d$. A negative defect gives
zero, while a zero-defect pair has degree one as above. The number
of coefficients in each incidence system is the sum of its defects.
\end{proof}

Over an effective field the incidence equations therefore determine
any prescribed Laurent expansion. To compute a product, multiply the
two expansions. Positivity and the leading coefficient one ensure that the Laurent
support contains every label in the theta expansion, so comparison of coefficients
recovers the theta expansion.

\subsection{Coefficients on coordinate axes and faces}

\begin{lemma}\label{lem:single-axis-trace}
Fix a mutable index $i$ under the hypotheses of
Theorem~\ref{thm:incidence}. Let $G_i$ be any Laurent monomial
multiple of $F_i$ in the unchanged coordinates, up to a nonzero scalar.
Suppose the contacts $P_1,\ldots,P_a$ have total defect $k e_i$
with $k\ge0$, and total order zero on the units. Then
\[ \mathsf T(P_1,\ldots,P_a)=\binom{h}{k},\qquad h=-\sum_{\alpha=1}^a\tau_{G_i}(P_\alpha). \]
Here $h\ge0$, and the binomial coefficient is zero for $k>h$.
\end{lemma}

\begin{proof}
Let $T'$ be the torus on all coordinates except $x_i$.
For marks $p_\alpha$ and evaluation $b'\in T'$ at infinity,
these coordinates are prescribed by the incidence equations. Away
from the marks they have the form $b'g(t)$, where every
coordinate of $g(t)$ is a product of powers of $t-p_\alpha$.
The total order of each unchanged coordinate is zero, so $g(\infty)=1$.
Since $G_i$ is independent of $x_i$, the rational function
$H(t)=G_i(b'g(t))$ is now fixed. Its order at $p_\alpha$
is $\tau_{G_i}(P_\alpha)$ for generic $b'$. Indeed, the terms
of minimum order have distinct characters of $b'$ and nonzero
coefficients. Also $H(\infty)=G_i(b')\ne0$. It follows by degree
on $\mathbb P^1$ that $H$ has $h$ zeros off the marks,
counted with multiplicity.

These zeros are simple for generic $b'$. The hypersurface
$G_i=0$ in $T'$ is smooth, because a singular point would give
a singular point with $X=Y=0$ on the smooth one-exchange open
$XY=F_i$. Over the complement of the marks the map
$(b',t)\mapsto b'g(t)$ is smooth. The inverse image of the
hypersurface is therefore smooth of dimension $\dim T'$.
Its nondominant components miss a generic $b'$, and its dominant
components project generically \'{e}tale to $T'$ in characteristic
zero. Its geometric generic fiber consists of distinct simple zeros.

Only the monic polynomial for $x_i$ is undetermined. Its roots
must be $k$ among the $h$ zeros, by divisibility at an unmarked
root. Conversely every such subset gives a monic polynomial, and
the adjacent quotient at each selected root is nonzero by simplicity.
It remains to check the marked incidence conditions. At a mark on
the positive $i$-axis the unchanged residues are generic and
$G_i\ne0$. At a mark on another positive mutable axis $j$, write
$F_j$ as a polynomial in $x_i$. Its coefficients are independent
of $x_j$, and all coordinates other than $x_i,x_j$ have order
zero there. Their leading residues are independent evaluation
parameters times nonzero constants. At least one coefficient is
therefore nonzero at the mark. The value of $x_i$ at the
mark is the independent evaluation parameter $b_i$ times a nonzero
scalar over the field of the selected roots. Hence $F_j\ne0$
for generic $b_i$, simultaneously for all the finitely many
choices. Theorem~\ref{thm:incidence} now counts each subset once.
\end{proof}

\begin{proposition}\label{prop:axis}
For every $m\in L$ and every mutable $i$,
\[ A_m(0,\ldots,0,z_i,0,\ldots,0)=(1+z_i)^{D_i(m)}. \]
In particular $D_i(m)$ is the exact degree of $A_m$ in $z_i$.
For the normalized factor $A_{-e_j}=z_j\Theta_{-e_j}$, the exponent on
the $i$th axis is $\deg_{v_j}G_i$, for mutable and designated
unit coordinates alike.
\end{proposition}

\begin{proof}
In Corollary~\ref{cor:trace-reconstruction}, let $s_0=b+m$ and take $d=a e_i$,
where $0\le a\le D_i(m)$. The only terms on the left have $e=l e_i$, $0\le l\le a$.
The polynomial $G_i$ is independent of $x_i$, so changing a
contact by a multiple of $e_i$ leaves its tropical value unchanged.
Moreover $\tau_{G_i}(b)=\tau_{G_i}(s_0)=0$.
Writing $w=-\tau_{G_i}(-s_0)$, Lemma~\ref{lem:single-axis-trace}
gives precisely
\[ \sum_{l=0}^a \binom{w}{a-l} c_{m,l e_i} =\binom{w+D_i(m)}a. \]
The diagonal entries equal one. Comparing coefficients in
$(1+y)^w(1+y)^{D_i(m)}=(1+y)^{w+D_i(m)}$ shows that
$c_{m,l e_i}=\binom{D_i(m)}l$. The upper bound of
Theorem~\ref{thm:corner} rules out higher powers, and the coefficient
at the bound is one. Finally
$-\tau_{G_i}(-e_j)=\deg_{v_j}G_i$.
\end{proof}

\begin{corollary}\label{cor:independent-face}
Let $I\subseteq\{1,\ldots,r\}$, and suppose that $G_i$ is
independent of $x_j$ for all $i,j\in I$. Group the indices
$i\in I$ by their normalized exchanges $H=G_i/G_i(0)$.
For each distinct $H$, set $I_H=\{i\in I:G_i/G_i(0)=H\}$ and
$L_H=1+\sum_{i\in I_H}z_i$. Then
\[ \left.A_m\right|_{z_j=0\ (j\notin I,\ 1\le j\le r)}=\prod_H L_H^{-\tau_H(m)}. \]
\end{corollary}

\begin{proof}
Each $H$ is irreducible. Indeed $F_i$ is irreducible in the
Laurent coefficient ring with polynomial mutable variables. Inverting
the other mutable coordinates preserves irreducibility, and $H(0)=1$
excludes coordinate factors on passing back to the polynomial ring.

Take contacts with nonnegative defect $\delta$ supported on $I$
and with balanced unit orders. The unchanged coordinates give a fixed
map $b'g(t)$. As in Lemma~\ref{lem:single-axis-trace}, the function
$H(b'g(t))$ has $h_H=-\sum_\alpha\tau_H(P_\alpha)$ simple zeros
away from the marks. For distinct $H$ these zero sets are disjoint
for general $b'$. Indeed, the two hypersurfaces intersect in
codimension at least two on the torus of unchanged coordinates.
Their inverse image under the smooth map $(b',t)\mapsto b'g(t)$
has dimension less than the dimension of that torus.

Each zero of $H(b'g(t))$ may be left unused or assigned to one
$A_i(t)$ with $i\in I_H$. Pairwise coprimality forbids assigning
it twice. Every allocation of the prescribed degrees satisfies the
remaining conditions of Theorem~\ref{thm:incidence}. At a positive
axis in $I$ the exchange depends only on the generic unchanged
residues. At a positive axis $j\notin I$, expand $F_j$ in the
variables indexed by $I$. Some coefficient is nonzero at the
unchanged residues. The evaluation parameters $b_i$, $i\in I$,
are algebraically independent over the field of the selected roots,
so the exchange is nonzero at this mark as well. Thus the trace is
$[z^\delta]\prod_H L_H^{h_H}$, where $[z^\delta]$ extracts the
coefficient of $z^\delta$.

In Corollary~\ref{cor:trace-reconstruction}, take $d$ supported on
$I$. The defect condition restricts the summation to $e\le d$.
Write $s_0=b+m$ and $K=\prod_H L_H^{-\tau_H(-s_0)}$.
Each $H$ is independent of the active coordinates, so the pair
traces give the coefficients of $K$. The triple traces give those
of $K\prod_H L_H^{-\tau_H(m)}$, since
$\tau_H(b)=\tau_H(s_0)=0$. The constant term of $K$ is one,
so the triangular system determines the coefficients of the
restriction. For $i\in I_H$ its proposed degree in $z_i$ is
$-\tau_H(m)=D_i(m)$, as required.
\end{proof}

\subsection{A Laurent chart from a negative Newton vertex}

We now give a criterion that also applies to seeds without a polynomial
corner. Suppose a regular function $g\in k[U]$ has a
vertex $a$ of its Newton polytope with $a_i<0$ for each mutable
coordinate. Let $C\subset L_{\mathbb R}$ be the full-dimensional
normal cone on which $a$ attains the minimum, and set
$\delta_g(m)=\langle a,m\rangle-\tau_g(m)$. In the following
statement $z^p$ denotes the character with exponent $p$ of the
torus dual to the marked source torus. We call $C$ a \emph{straight cone}
for the resulting chart since $\theta_p=z^p$ for $p\in C\cap L$.

\begin{proposition}\label{prop:negative-vertex-chart}
Under the hypotheses of Theorem~\ref{thm:incidence}, a choice of $g,a$
as above determines a \emph{marked Laurent realization}, an embedding of
$B_{\mathbb Z}(U)$ in the Laurent ring $\mathbb Z[L]$ of the dual torus. It sends
$\theta_p$ to $z^p$ for $p\in C\cap L$, and sends $\theta_m$ to
a Laurent polynomial with coefficient one at $z^m$, nonnegative
integral coefficients, and support in
\[ m+\bigl\{(d,0):d\in\mathbb N^r,\ \textstyle\sum_i(-a_i)d_i\le\delta_g(m)\bigr\}. \]
For any $b\in L\cap\operatorname{int}C$, the realization identifies
$B_{\mathbb Z}(U)[\theta_b^{-1}]$ with $\mathbb Z[L]$.
\end{proposition}
\begin{proof}
By Theorem~\ref{thm:corner-free-product}, every product has
coefficient one at the sum of its labels, and all corrections are
nonnegative mutable shifts. Apply the regular-function
inequality to $g$ in $\theta_p\theta_m$, for $p\in C$. For each
label $s=p+m+(d,0)$ in that product, it gives
\[ \langle a,p\rangle+\tau_g(m)\le\tau_g(s) \le\langle a,p+m\rangle+\sum_i a_i d_i . \]
This bounds $d$ as asserted. For $m\in C$ the bound is zero.
Thus the straight products on $C$ identify a copy of
$\mathbb Z[C\cap L]$ inside $B_{\mathbb Z}(U)$. Multiplication by $\theta_b$ is
injective by the coefficient-one least term rule, for the order
used in the proof of Theorem~\ref{thm:corner}. Inverting $\theta_b$
in the semigroup ring gives $\mathbb Z[L]$, since for any $m$ one has
$m+Nb\in C$ for sufficiently large $N$.
The finite bound on the shifts is independent of $p$. Choosing
$p=Nb$ so that every label in the product lies in $C$, we obtain
$\theta_b^N\theta_m$ as a positive integer sum of straight
characters. Division by $z^{Nb}$ gives the claimed Laurent
polynomial and proves the localization identity.
\end{proof}

The support bound in Proposition~\ref{prop:negative-vertex-chart} gives
a polynomial factor $A_m=z^{-m}\Theta_m$ with nonnegative integral
coefficients and constant term one in this chart as well. The same face
formula applies. Take an independent set $I$ as in
Corollary~\ref{cor:independent-face}. Choose a separating weight in
$\operatorname{int}C$. Divide each $F_i$ for $i\in I$ by its
least term at this weight. For each distinct quotient $H$, let
$I_H$ consist of the indices that give $H$. Then
\[ A_m\big|_{z_j=0\ (j\notin I,\ 1\le j\le r)}=\prod_H L_H^{-\tau_H(m)},\qquad L_H=1+\sum_{i\in I_H}z_i. \]
In particular, when $I=\{i\}$, the restriction is
$(1+z_i)^{-\tau_H(m)}$.

Indeed, the unchanged coordinates have total contact order zero.
Dividing an exchange by its least term therefore preserves its
total contact order. Distinct normalized exchanges are nonassociate
irreducibles on the unchanged torus, so their zeros can be assigned
as in Corollary~\ref{cor:independent-face}. On an open subcone of
$C$ about the chosen weight, every $H$ has least term one.
For fixed $d$ supported on $I$, choose a lattice point $b$
sufficiently deep in that subcone that $b+m$, all labels in
$\theta_b\theta_m$, and all $b+m+e$ with $0\le e\le d$ lie
in it. This is possible by the uniform shift bound of
Proposition~\ref{prop:negative-vertex-chart}. Since each $H$ has
tropical value zero at $b$ and $b+m$, the same triangular trace
equations give the coefficient at $d$.

\subsection{Seed criteria for a negative Newton vertex}

Let $A$ be the full upper algebra of a marked seed over
$k[u_1^{\pm1},\ldots,u_c^{\pm1}]$, with mutable coordinates
$x_1,\ldots,x_r$ and positive mutable rank.
Write $F_i$ for the ordinary exchange at $i$. The functions
$g_i=F_i/x_i$ are regular. Each is the normalized first mutant
times a monomial with nonnegative mutable exponents and a coefficient
unit. We characterize the seeds for which a product of these functions has
the vertex required in
Proposition~\ref{prop:negative-vertex-chart}.

\begin{proposition}\label{prop:seed-negative-vertex}
Starting with $S=\varnothing$, adjoin an index $i$ whenever
$F_i\big|_{x_j=0\ (j\notin S)}\ne0$ as a polynomial in the remaining mutable variables over
$k[u_1^{\pm1},\ldots,u_c^{\pm1}]$. Repeat until $S$ stops
growing. The following are equivalent.
\begin{enumerate}
\item The final $S$ contains all mutable indices.
\item Some product $\prod_i g_i^{n_i}$, $n_i\in\mathbb N$, has
a Newton vertex strictly negative in all mutable coordinates.
\item There is no nonempty set $J$ of mutable indices such that
$F_i\in(x_j:j\in J)$ for every $i\in J$.
\item The ordinary exchange presentation
\[ A_0=k[u_1^{\pm1},\ldots,u_c^{\pm1}][x_1,\ldots,x_r,y_1,\ldots,y_r]/ (x_i y_i-F_i:1\le i\le r) \]
is a domain.
\end{enumerate}
When these conditions hold the mutable projection of the vertex
can be chosen to be $(-1,\ldots,-1)$. In particular the full upper
algebra is finitely generated when $A_0$ is a domain, by
Theorem~\ref{thm:polyhedral-newton-sections}.
Under the hypotheses of Theorem~\ref{thm:incidence}, integrality of
$A_0$ therefore gives the marked Laurent chart of
Proposition~\ref{prop:negative-vertex-chart}.
\end{proposition}
\begin{proof}
If the procedure stops early, take $J$ to be the complement of
its final set. Conversely, no member of a set $J$ as in the third
condition can ever be adjoined. This proves the equivalence of the
first and third conditions. For each nonempty $J$, first quotient the polynomial ring
$k[u^{\pm1}][x_1,\ldots,x_r]$ by $(x_j:j\in J)$. The quotient
is a polynomial ring over $k[u^{\pm1}]$, so it embeds in the Laurent
ring obtained by inverting the remaining mutable coordinates. The third condition says exactly that
every exchange-stratum ideal $J_I$ in
Theorem~\ref{thm:local-exact-exchange-strata} is nonzero. That theorem
identifies this with integrality of $A_0$.

Order the indices by their time of addition, breaking ties arbitrarily.
Under the first condition $F_i$ has a term involving only earlier
mutable variables. Give the $x_i$ successively increasing positive
weights, with each new weight greater than the weight of any exchange
term involving only earlier variables. Initially give the units weight
zero. For each $F_i$ some earlier-variable term then has smaller
weight than any term involving a variable at or after $i$. A small
generic perturbation in all coordinates preserves these strict
inequalities and chooses a unique minimizing term for each $F_i$.
Let $M_{ji}$ be the exponent of $x_j$ in that term. The matrix
$M$ is strictly triangular in the chosen order. With $I_r$ the
identity matrix and $\mathbf1$ the column of ones, the vector
$n=(I_r+M+\cdots+M^{r-1})\mathbf1$ has positive integral entries. The product of the $g_i^{n_i}$ has
unique minimizing exponent with mutable projection
$(M-I_r)n=-\mathbf1$, as required.

For the converse expose the vertex of a product by a generic
weight in the interior of its normal cone, and define $M$ from the
selected numerator terms as above. Strict negativity says
$Mn<n$ componentwise, so every $n_i>0$. An edge
$i\to j$ with $M_{ji}>0$ implies $n_i<n_j$, so there are
no directed cycles. Order the vertices so that the outgoing
neighbors of each index precede it. The selected term of $F_i$
then uses only earlier mutable variables, giving the first condition.
\end{proof}

The normalized first mutations can give further vertices. The following test applies to any finite collection of regular
functions, including coordinates from more distant seeds.

\begin{proposition}\label{prop:seed-negative-partners}
Fix nonzero regular Laurent functions $h_1,\ldots,h_t$. Take a
full-dimensional chamber of the common refinement of the normal fans
of their Newton polytopes, and let $D$ be the $r\times t$ matrix
of the mutable projections of their minimizing exponents in the
chamber. Exactly one of the following holds.
\begin{enumerate}
\item For some $n_j\in\mathbb Z_{\ge0}$, the product
$\prod_{j=1}^t h_j^{n_j}$ has an exposed vertex with mutable
projection $-\mathbf 1$ throughout the chamber after multiplication
by a monomial with nonnegative exponents in the initial mutable variables.
\item There is a nonzero $\lambda\in\mathbb Z_{\ge0}^r$ with
$D^t\lambda\ge0$.
\end{enumerate}
Thus a finite calculation on the normalized first partners tests
every product of those partners. At a lattice-separating weight $w$
in the chamber the second condition rules out a $w$-minimizing exponent
negative in every mutable coordinate for \emph{every} regular function
of $A$ if the minimizing exponents
of the $h_j$, the initial mutable variables and the unit characters
and their inverses generate the value semigroup of $A$ at $w$.
\end{proposition}
\begin{proof}
The rational polyhedral cone generated by the columns of $D$ meets
the open negative orthant if and only if $Dn<0$ for some
$n\in\mathbb N^t$, by rational approximation and homogeneity.
If they are disjoint, separation gives a nonzero functional
nonnegative on this cone and nonpositive on the open negative
orthant. Its coordinates are nonnegative, and the rational
inequalities admit such a functional with integral coordinates.
The two alternatives exclude one another. For an integral $n$
with $d=Dn<0$, multiply the product by $x^{-\mathbf1-d}$
to arrange its mutable exponent precisely $-\mathbf1$.
This operation preserves the minimizing chamber.

Under the value semigroup hypothesis the extended functional
$(\lambda,0)$ is nonnegative at the least exponent of every
regular Laurent polynomial. It is nonnegative on the mutable
coordinate exponents and zero on the units. Strict negativity in
every mutable coordinate would give a negative pairing, a
contradiction. The initial-ideal test of
Lemma~\ref{lem:reciprocal-initial-ideal} is one way to verify the
semigroup hypothesis from a positively graded presentation.
\end{proof}

\begin{corollary}\label{cor:negative-vertex-semigroup}
Fix a weight $w$ separating the points of the source character lattice
$M=\mathbb Z^{r+c}$. The \emph{initial algebra} of $A$ at $w$
is the subalgebra of $k[M]$ generated by the least terms of its
nonzero elements. The following are equivalent.
\begin{enumerate}
\item A regular function has a $w$-minimizing exponent strictly
negative in every mutable coordinate.
\item The $w$-minimizing exponents of nonzero elements of $A$
form the entire lattice $M$.
\item The initial algebra of $A$ in its marked Laurent ring is
$k[M]$, or equivalently every initial mutable character is a unit
in that initial algebra.
\end{enumerate}
If they fail, some nonzero $\lambda\in\mathbb R_{\ge0}^r$ pairs
nonnegatively with the mutable projection of the $w$-minimizing
exponent of every regular function.
\end{corollary}
\begin{proof}
The exponent semigroup contains the nonnegative mutable coordinate
vectors and both signs of every unit coordinate vector. If it contains
a vector strictly negative in the mutable coordinates, adding these
known elements gives $-e_i$ for each mutable $i$. This proves
the equivalence of the first two conditions. Since $w$ separates
the lattice, the initial algebra is the semigroup algebra of these
exponents, giving the third condition. In the absence of a negative
exponent even their real conical hull misses the open negative
mutable orthant. Indeed, a conical combination in this open set could be
approximated by a rational one and its denominators cleared.
Separation as in Proposition~\ref{prop:seed-negative-partners}
gives the stated functional without any finite generation assumption.
\end{proof}

\subsection{Gauss valuations of the forward basis}

\begin{theorem}\label{thm:forward-vi-main}
In the Laurent chart obtained from a Newton vertex negative in
all mutable coordinates, the forward theta basis is valuatively
independent for every real Gauss weight over a characteristic-zero
field. Equivalently, the Newton polytope of a nonzero finite linear
combination is the convex hull of the polytopes of its nonzero summands.
\end{theorem}
\begin{proof}
The local functions in Appendix~\ref{sec:one-sided-independence} count
analytic curves with a prescribed incoming direction and final position
in the skeleton. Their tropical paths are the \emph{one-input spines}.
Proposition~\ref{prop:formal-local-chart} identifies these functions
with the marked Laurent expansions. In these
expansions, distinct theta functions have distinct minimizing monomials
at a generic weight, by Theorem~\ref{thm:marked-path-vi}. Thus the
terms of least weight cannot cancel. Continuity of the support minima
gives the assertion for every real weight.
\end{proof}

Suppose the forward mirror is identified with a full upper LP algebra
whose seed torus has logarithmic volume form proportional to that of
the marked chart. Then every regular Laurent monomial in a reachable
mirror seed is a scalar multiple of a forward theta function. Indeed,
mutation preserves the logarithmic volume form up to sign, so
Proposition~\ref{prop:volume-transfer} gives valuative independence in
that seed. The monomial has a one-point Newton polytope there. As in
Corollary~\ref{cor:reachable-monomials-minimum-basis}, its theta
expansion therefore has a single term.

Under the hypotheses of Proposition~\ref{prop:monic-integral-descent},
a regular seed monomial with coefficient one equals $\pm\theta_p$.
Indeed, the first assertion of that proposition, applied with this
seed as the initial one, makes the monomial integral in every chart.
Its coefficient in the integral theta basis is therefore an integer.
The theta function is integral in the seed where the given function
is a monomial, so the reciprocal coefficient is also an integer.
The coefficient is consequently $1$ or $-1$. It is $1$ when the
monomial's least coefficient in the marked initial chart is positive.

Scaling a path with initial label $p$ gives one with initial label $jp$.
In Lemma~\ref{lem:local-support-dilation} we use this observation to prove
$\nu^z_a(\Theta_{jp})=j\nu^z_a(\Theta_p)$ for every positive integer $j$.
This is the homogeneity needed to construct $P_a$ in
Theorem~\ref{thm:reciprocal-polytope}.

\subsection{Change of initial seed for the polynomial factors}

We keep an intrinsic theta function fixed and change the seed in which
we express it. For the classical initial-seed mutation formula for
$F$-polynomials, see \cite[(6.21)]{FZIV}
and \cite[Proposition~3.4]{FuGy}.

Let $t'$ be obtained from $t$ by mutation at $k$, and let
$E=\widehat F_k$ be the exchange Laurent polynomial in $x_kx'_k=E$.
Suppose both seeds give first-star charts of the same $U$ as in
Section~\ref{subsec:counting-hypotheses}. In each marking choose a
forward chart given by Theorem~\ref{thm:corner} or
Proposition~\ref{prop:negative-vertex-chart}. Write $z,z'$ for their
coordinates. By Lemma~\ref{lem:boundary-gauss-change}, a skeletal label
with coordinates $p$ in $t$ has coordinates $p'=T(p)$ in $t'$,
where $p'_j=p_j$ for $j\ne k$ and $p'_k=\tau_E(p)-p_k$.

\begin{proposition}\label{prop:theta-factor-mutation}
Let $W_j$ be the expression for $z'_j$ in the $z$-coordinates.
Then
\[ A^{t'}_{p'}(W_1,\ldots,W_r)=\frac{z^p A_p^t(z_1,\ldots,z_r)}{\prod_{j=1}^n W_j^{p'_j}} . \]
If the chart of $t'$ is its polynomial-corner chart, then
$W_j=\Theta^t_{T^{-1}(e'_j)}$, where $e'_j$ is the $j$-th
standard cocharacter in $t'$. In general, choose an integral $b'$
in the new straight cone sufficiently far inside that $b'+e'_j$ also
lies in it for every $j$. Then
\[ W_j=\frac{\Theta^t_{T^{-1}(b'+e'_j)}}{\Theta^t_{T^{-1}(b')}} . \]

Suppose both charts are polynomial-corner charts and $E$ itself is
a polynomial with nonzero constant term. Set $H_k=A^t_{-e_k}$ and
$d=\deg_{z_k} A^t_p=-\tau_E(p)$. Then $H_k$ is independent of
$z_k$. Upon identifying $z_j=z'_j$ for $j\ne k$, we have
\begin{equation}\label{eq:theta-factor-mutation}
A^{t'}_{p'}(z'_1,\ldots,z'_r)=(z'_k)^d H_k^{p_k}
   A^t_p(z'_1,\ldots,z'_{k-1},H_k/z'_k,z'_{k+1},\ldots,z'_r).
\end{equation}
In particular, if $A^t_p=\sum_{j=0}^d B_j z_k^j$, then the
coefficient of $(z'_k)^{d-j}$ in $A^{t'}_{p'}$ is
$H_k^{p_k+j}B_j$. Whenever $p_k+j<0$, $H_k^{-p_k-j}$ divides
$B_j$ in $\mathbb Z[z_i:1\le i\le r,\ i\ne k]$.
\end{proposition}

\begin{proof}
Both expressions for the theta function represent the same element of
the augmented algebra, so $z^p A_p^t=W^{p'} A^{t'}_{p'}(W_1,\ldots,W_r)$.
The $W_j$ are algebraically independent, hence this identity
determines $A^{t'}_{p'}$. In the new straight cone the theta functions
are characters. Taking the quotient for $b'+e'_j$ and $b'$ gives
$z'_j$, and at a polynomial corner we can take $b'=0$.

When $E$ is polynomial with nonzero constant term,
$\tau_E(e_j)=0$ for every $j$. Thus $T^{-1}(e'_j)=e_j$ for
$j\ne k$, while $T^{-1}(e'_k)=-e_k$. In the old polynomial-corner
chart this gives $z'_j=z_j$ for $j\ne k$ and
$z'_k=z_k^{-1}H_k$. The degree formula in Proposition~\ref{prop:axis}
applies with $E$. The ordinary exchange has no mutable monomial
factor, so polynomiality of $E$ forces its mutable normalizing
factors to be trivial. Thus $E$ and $G_k$ differ by a coefficient
unit. Both have minimum exponent zero in each coordinate, so this
unit is a scalar.
The degree formula gives $\deg_{z_k}A^t_p=-\tau_E(p)$ and
$\deg_{z_k}H_k=-\tau_E(-e_k)=0$. Substitute $z_k=H_k/z'_k$
in $z^p A_p^t$ and divide by $(z')^{p'}$. Since $p'_k=-p_k-d$,
the term $B_j z_k^j$ contributes $H_k^{p_k+j}B_j(z'_k)^{d-j}$.
The divisibility follows because the new factor is a polynomial with
integer coefficients.
\end{proof}

For example, take the coefficient-free seed of type $A_2$ with
exchanges $xx'=1+y$ and $yy'=1+x$, and mutate at $x$.
For $p=(0,-1)$, Proposition~\ref{prop:axis} gives
$A_p^t=1+z_x$ and $H_x=1+z_y$. The new coordinates of the label
are $(-1,-1)$. Substitution in~\eqref{eq:theta-factor-mutation} gives
$A^{t'}_{(-1,-1)}=z'_x(1+(1+z'_y)/z'_x)=1+z'_x+z'_y$.

\section{Reciprocity and return}\label{sec:reciprocity}

\subsection{Reciprocity from a marked inclusion of algebras}

We call the theta functions obtained by taking the mirror a second time
\emph{reverse theta functions}. Suppose they are regular on $U$ in
its marked torus coordinates. We ask when they span $k[U]$ and when
the Gauss valuation of a reverse theta function at a forward label
agrees with that of the forward theta function at the reverse label.

Let $M$ and $N$ be dual lattices, and give the field $k$ the
trivial valuation, which is zero on $k^\times$. Write $x^q$ for the character of exponent
$q\in M$ and $z^p$ for the character of exponent $p\in N$.
We write $\nu^x_p$ for the Gauss valuation of weight $p$ on
$k[M]$, and $\nu^z_q$ for its counterpart on $k[N]$.
Both use the minimum convention.

In the geometric setting, $A$ is the source coordinate ring, $B$
its augmented mirror algebra and $C$ the algebra constructed from
$\Spec B$. For the criterion we take $k$-algebras $A,B,C$ with
bases $(\theta_p)$ of $B$ and $(\varphi_q)$ of $C$ satisfying
\[ A\subset k[M],\qquad B=\bigoplus_{p\in N}k\theta_p\subset k[N],\qquad C=\bigoplus_{q\in M}k\varphi_q\subseteq A . \]

Suppose there are full-dimensional rational polyhedral cones
$\sigma\subset N_{\mathbb R}$ and $\tau\subset M_{\mathbb R}$
such that
\[ \theta_p=z^p\quad(p\in\sigma\cap N),\qquad \varphi_q=x^q\quad(q\in\tau\cap M). \]
For every $p,q$, assume the coefficient of $z^p$ in $\theta_p$
and that of $x^q$ in $\varphi_q$ are both one. The remaining
exponents lie in $p+(\tau^\vee\setminus\{0\})$ and
$q+(\sigma^\vee\setminus\{0\})$, respectively. Here the dual cone
consists of linear forms nonnegative on the given cone.

For real weights $d(a)$ on a basis $(b_a)$, the \emph{split filtration}
has $F^\lambda=\operatorname{span}_k\{b_a:d(a)\ge\lambda\}$.
Its order is $v_d(\sum c_ab_a)=\min_{c_a\ne0}d(a)$, with
$v_d(0)=+\infty$. It is multiplicative if
$v_d(uv)\ge v_d(u)+v_d(v)$. We require both of the following
filtrations to be multiplicative.
\begin{enumerate}
\item Each $0\ne f\in A$ gives weights
  $p\mapsto\nu^x_p(f)$ on the basis of $B$.
\item Each $0\ne g\in B$ gives weights
  $q\mapsto\nu^z_q(g)$ on the basis of $C$.
\end{enumerate}
In the geometric setting the regular-function inequality~\eqref{eq:regular-function-inequality} says that
the order attached to any output label in a nonzero theta product is
at least the sum of the orders attached to its inputs. Applied on
each side, it gives these two multiplicative filtrations once the
reverse functions are regular in the chosen source chart.

\begin{lemma}\label{lem:rec-straight}
Let $L$ be a lattice and $D\subset k[L]$ a subalgebra with a
basis $(b_a)_{a\in L}$.
Assume $b_a=z^a$ on the lattice points of a full-dimensional rational
cone $\eta$. If a multiplicative split filtration on $D$ has weights
$d(a)=\langle w,a\rangle$ on $\eta$ for a real weight $w\in L^*_{\mathbb R}$,
then $v_d(g)\le\nu_w(g)$ for all nonzero $g\in D$, where
$\nu_w$ is the Gauss valuation in the $z$-coordinates.
\end{lemma}
\begin{proof}
Choose $b\in L\cap\operatorname{int}\eta$. For sufficiently large
$n$, the exponents of $z^{nb}g$ all lie in $\eta$.
Its Laurent expansion is therefore its basis expansion. Multiplicativity
gives
\[ n\langle w,b\rangle+v_d(g)\le v_d(z^{nb}g) =n\langle w,b\rangle+\nu_w(g). \]
\end{proof}

\begin{theorem}\label{thm:one-sided-return}
Under the preceding hypotheses the two bases satisfy
\[ \nu^z_q(\theta_p)=\nu^x_p(\varphi_q)\qquad (p\in N,\ q\in M). \]
Moreover $C=A$. If $0\ne f\in A$ has Newton polytope $P$,
its expansion in the basis $(\varphi_q)$ takes at most
$\lvert P\cap M\rvert$ successive subtractions. Every basis element
used has Newton polytope contained in $P$.
\end{theorem}
\begin{proof}
For $f=\varphi_q$, the first filtration has weight
$\langle p,q\rangle$ on $\sigma$, by the support condition.
Lemma~\ref{lem:rec-straight} applied to $\theta_p$ yields
$\nu^x_p(\varphi_q)\le\nu^z_q(\theta_p)$. Apply the same lemma
on $\tau$ to the second filtration, taking $g=\theta_p$.
This gives the reverse inequality.

Choose $\omega\in\operatorname{int}\sigma$ whose pairing with $M$
is injective, and let $q_0$ be the exponent of $f$ with least
$\omega$-weight. The cone
\[ \eta=\{p\in\sigma:\langle p,e-q_0\rangle\ge0 \text{ for every }e\in\Supp (f)\} \]
is full dimensional, since all inequalities with $e\ne q_0$ are
strict at $\omega$. On $\eta$ the weights supplied by $f$ are
linear with slope $q_0$. Lemma~\ref{lem:rec-straight} and the
reciprocity already proved give, for every $p\in N$,
\[ \nu^x_p(f)\le\nu^z_{q_0}(\theta_p)=\nu^x_p(\varphi_{q_0}). \]
Homogeneity and continuity give the same inequality for all real $p$.
A polytope is the intersection of the closed half-spaces determined by
its minimum support function. Hence
$\Newt (\varphi_{q_0})\subseteq P$.

The monomial $x^{q_0}$ has coefficient one in $\varphi_{q_0}$,
and all its other terms have strictly greater $\omega$-weight.
Subtract from $f$ its $x^{q_0}$-coefficient times $\varphi_{q_0}$.
The remainder belongs to $A$. If it is nonzero, its Newton polytope
lies in $P$ and its least weight strictly increases.
Iteration terminates since $P\cap M$ is finite. This proves both the
bound and $A=C$.
\end{proof}

\begin{corollary}\label{cor:return-independence}
The basis $(\varphi_q)$ of $A$ is valuatively independent for all real Gauss
valuations on the marked torus. If $k$ has characteristic zero and
$\varphi_q\in\mathbb Z[M]$ for every $q$, then
$A\cap\mathbb Z[M]=\bigoplus_q\mathbb Z\varphi_q$.
\end{corollary}
\begin{proof}
For a nonzero sum $f=\sum c_q\varphi_q$, uniqueness of expansion
and the theorem give $\Newt(\varphi_q)\subseteq\Newt(f)$ whenever $c_q\ne0$. The opposite inclusion
in the convex hull of these polytopes follows by addition. Taking a
minimum gives $\nu^x_p(f)=\min_{c_q\ne0}\nu^x_p(\varphi_q)$.
For $f\in\mathbb Z[M]$ every subtraction in the proof has integral
coefficient and leaves an integral Laurent polynomial.
\end{proof}

\subsection{Reciprocal Newton polytopes from the forward mirror}
\label{subsec:one-sided-polytopes}

Let $L$ be the source cocharacter lattice and $L^*$ its dual. In
the chart of Proposition~\ref{prop:negative-vertex-chart}, write
$C\subset L_{\mathbb R}$ for its full-dimensional straight cone.
Use a lattice-separating weight in $\operatorname{int}C$ to define
the smallest source polytopes $Q_a$. For $a\in L^*$, the forward
basis determines the integer function
$h_a(p)=\nu^z_a(\Theta_p)$ on $L$. A Laurent polynomial on the torus of $U$
with reciprocal values at every $p$ would have $h_a$ as the
minimum support function of its Newton polytope.

\begin{theorem}\label{thm:reciprocal-polytope}
There is a unique lattice polytope $P_a\subset L^*_{\mathbb R}$
with minimum support function $h_a$. It is computable from finitely
many forward theta functions. For $p$ in the interior of $C$, the
unique minimizing exponent of $P_a$ is $a$.

Suppose a nonzero $f\in k[U]$ has $a$ as its unique minimizing
exponent at some $\omega\in\operatorname{int}C$. Then $P_a\subseteq\Newt(f)$. A regular function on $U$ with Newton
polytope $P_a$ and coefficient one at $a$ exists precisely when
$P_a=Q_a$, where $Q_a$ is the smallest polytope of
Theorem~\ref{thm:minimum-source-newton-polytope}. Its possible
coefficients are the solutions of the first-star divisibility
equations on $P_a\cap L^*$, with coefficient one prescribed at $a$.
\end{theorem}
\begin{proof}
The coefficient of $\theta_{p+q}$ in $\theta_p\theta_q$ is one.
Positivity of all Laurent coefficients and structure constants gives
$h_a(p+q)\ge h_a(p)+h_a(q)$. Lemma~\ref{lem:local-support-dilation} gives
$h_a(jp)=j h_a(p)$ for every positive integer $j$.
Lemma~\ref{lem:integral-support}
gives the asserted lattice polytope. On $C$ the equality $\Theta_p=z^p$ gives $h_a(p)=\langle a,p\rangle$,
so each interior point exposes $a$.

All vertices of $P_a$ lie in the box with coordinate bounds
$h_a(e_i)\le t_i\le-h_a(-e_i)$. As in the proof of
Proposition~\ref{prop:source-polar-polytope}, one evaluation in each
chamber cut out by differences of lattice points of this box recovers
all vertices. We compute the evaluations by the multiplication rule of
Theorem~\ref{thm:corner-free-product} and translation into $C$ as in
Proposition~\ref{prop:negative-vertex-chart}.

On a full-dimensional cone about $\omega$ inside $C$, the minimum
support function $\tau_f$ has slope $a$.
The regular-function inequality makes the split filtration with
weights $\tau_f(p)$ multiplicative. Apply
Lemma~\ref{lem:rec-straight} on that cone to $\Theta_p$.
It gives $\tau_f(p)\le\nu^z_a(\Theta_p)=h_a(p)$ for every $p$.
This is precisely $P_a\subseteq\Newt (f)$.
A solution supported in $P_a$ with a nonzero coefficient at $a$
has unique minimizing exponent $a$ at every interior point of $C$.
Its Newton polytope therefore contains $Q_a$. Applying the
inclusion to a normalized section on $Q_a$ gives $P_a\subseteq Q_a$,
so a solution forces equality. Conversely, if $P_a=Q_a$, that section
is a solution. The initial and adjacent charts cover $U$ in
codimension one, so normality and Lemma~\ref{lem:recognition-rank-one}
identify regularity with the stated divisibilities. The quotient bounds in the proof of
Theorem~\ref{thm:minimum-source-newton-polytope} make the linear system finite.
\end{proof}

Thus for every source label $a\in L^*$ and forward label $p\in L$,
\begin{equation}\label{eq:forward-newton-bound}
P_a\subseteq Q_a,\qquad \nu^z_a(\Theta_p)=h_a(p)=\min_{b\in P_a}\langle b,p\rangle\ge\min_{b\in Q_a}\langle b,p\rangle.
\end{equation}
By~\eqref{eq:dual-polytope-minima}, the last minimum is
$\min_{e\in T_p}\langle a,e\rangle$. As $a$ varies, the inequality
therefore gives $\Newt(\Theta_p)\subseteq T_p$. Equality $P_a=Q_a$
for every $a$ is equivalent to $\Newt(\Theta_p)=T_p$ for every $p$.

We call $P_a$ the \emph{reciprocal Newton polytope} at $a$.
A \emph{normalized reciprocal function} at $a$ is a regular function
supported in $P_a$ whose coefficient at $x^a$ is one. We say
that the label $a$ is \emph{realized} if such a function exists.
The theorem identifies realization with $P_a=Q_a$.
A Laurent polynomial is \emph{pointed at} $a$ in this chart if $x^a$
has coefficient one and is its unique term of least weight throughout
$\operatorname{int}C$. In particular, every normalized reciprocal
function at $a$ is pointed there.

\begin{corollary}\label{cor:reciprocal-from-windows}
If the test of Theorem~\ref{thm:reciprocal-polytope} has a solution
$f_a$ for every $a$, any choice of one solution for each $a$
is a basis of $k[U]$ satisfying
$\tau_{f_a}(p)=\nu^z_a(\Theta_p)$. It is valuatively independent
for every real Gauss weight of the marked source torus.
\end{corollary}
\begin{proof}
The solutions have Newton polytopes $P_a=Q_a$ and least terms
$x^a$ at a separating weight in $\operatorname{int}C$.
Theorem~\ref{thm:minimum-source-newton-polytope} gives the basis and
its valuative independence. Taking the minimum on $P_a$ gives
$\tau_{f_a}(p)=h_a(p)$.
\end{proof}
At a polynomial corner the resulting source basis exists under the
hypotheses of Theorem~\ref{thm:polynomial-corner-reciprocity}.
Identifying one of its elements with a prescribed reverse theta
function also requires the comparison in the marked chart.

\section{Reciprocity at a polynomial corner}\label{sec:polynomial-corner-reciprocity}

Let $U$ be the smooth spectrum of a full upper LP algebra over
$k[u_1^{\pm1},\ldots,u_c^{\pm1}]$, where $k$ is algebraically
closed of characteristic zero. Assume that the marked seed has a
polynomial corner. By Section~\ref{subsec:counting-hypotheses}, $U$
then satisfies the geometric hypotheses of Theorem~\ref{thm:incidence}.
We normalize the ordinary exchanges by $G_i(0)=1$. Write $M$ for
the character lattice of the source torus and $L=M^*$ for its
cocharacter lattice. We write $x^a$ on the source and $z^p$ in
the Laurent expansion of the forward mirror. Fix a weight $\mu\in L_{\mathbb R}$
whose coordinates are strictly positive and rationally independent. Use its
order on $M$ for the smallest source polytopes $Q_a$ and their dual
polytopes $T_p$ from Proposition~\ref{prop:source-polar-polytope}. Let
$P_a$ be the polytope of Theorem~\ref{thm:reciprocal-polytope} for the
forward chart at the polynomial corner.

\begin{theorem}\label{thm:polynomial-corner-reciprocity}
For every $p\in L$ and $a\in M$,
\begin{equation}\label{eq:corner-reciprocity}
\Newt(\Theta_p)=T_p,\qquad P_a=Q_a,\qquad \nu^z_a(\Theta_p)=\min_{b\in Q_a}\langle b,p\rangle .
\end{equation}
In particular, any choice of normalized source sections $f_a$
on the smallest Newton polytopes $Q_a$ satisfies
$\nu^x_p(f_a)=\nu^z_a(\Theta_p)$.
\end{theorem}

The exchanges may depend cyclically on one another, and the pole orders
are arbitrary. We prove that each vertex of $T_p$ occurs in
$\Theta_p$. Lemma~\ref{lem:polynomial-contact-coefficient} gives a
sufficient condition for a theta coefficient to be positive, in terms
of simple zeros of the exchange polynomials. Proposition~\ref{prop:corner-disk}
produces these zeros on a disk in a curve over a function field. We
apply the proposition after perturbing $p$, then recover the desired
vertex by continuity.

\subsection{A criterion for positivity of a theta coefficient}

Write $r$ for the number of mutable coordinates and $n=\operatorname{rank}M$.
Let $w,Y_1,\ldots,Y_n$ be indeterminates and write $Y=(Y_1,\ldots,Y_n)$.
For $p\in L$ and $1\le i\le r$, define
\[ D_i=-\min_{\alpha\in\Supp G_i}\langle\alpha,p\rangle,\qquad H_i(Y,w)=w^{D_i}G_i((w^{p_j}Y_j)_{1\le j\le n}). \]
Each $H_i$ is a polynomial, with $H_i(0,w)=w^{D_i}$.
We extend mutable vectors by zero on the designated coordinates.
For a real additive valuation $\operatorname{val}$ we use the absolute
value $|h|=\exp(-\operatorname{val}(h))$, with $|0|=0$.

\begin{lemma}\label{lem:polynomial-contact-coefficient}
Let $\Omega$ be a complete algebraically closed valued extension of
the trivially valued field $k$, with real-valued absolute value.
For each mutable $i$, let $A_i\in\Omega[w]$ be a monic polynomial
of degree $d_i$ with $0\le d_i\le D_i$. Suppose the roots of
$\prod_i A_i$ are simple, nonzero and of norm less than one.
Take $Y_i=A_i C_i$, where $C_i$ is a polynomial whose coefficients
have norm less than one. Take the designated $Y_j$ to be constants
of norm less than one. If $H_i(Y(w),w)$ has a simple zero at each
root of $A_i$, then $[z^{p+d}]\Theta_p>0$.
\end{lemma}

\begin{proof}
If every $d_i$ is zero, the coefficient is one by
Theorem~\ref{thm:corner}. Otherwise choose $m_i\ge0$ with
$\deg C_i\le m_i$, allowing $C_i=0$, and increase one bound to
at least one. The generic marked line below will then have at least
three boundary contacts.
Vary all coefficients of the $C_i$ and the designated constants
independently, with parameter list $c$. Let $u$ denote the lower
coefficients of the monic $A_i$. Over $R=k[[c]]$ consider
\[ \mathcal B=k[[c,u]]/ \bigl(\text{coefficients of }H_i(Y,w)\bmod A_i(w),\ 1\le i\le r\bigr). \]
At $c=0$ the conditions are $A_i\mid w^{D_i}$. Thus the closed
fiber is supported at $A_i=w^{d_i}$. In the regular local ring
$k[[c,u]]$, the remainder coefficients followed by the $c$'s
form a system of parameters, hence a regular sequence. A power of
the maximal ideal of $\mathcal B$ lies in $(c)$, so completeness
lifts generators of the closed fiber to a finite set of $R$-module
generators. The regular sequence $c$ on $\mathcal B$ then makes
this finite module free over the regular local ring $R$.

For $d_i>0$, take the product of the constant term and discriminant
of $A_i$ and its resultant with $H_i(Y,w)/A_i$. Multiply these
products by the pairwise resultants of the $A_i$. Call the result
$\Delta\in\mathcal B$. The assumed solution sends $\Delta$ to a
nonzero value. This substitution converges since all formal parameters
have norm less than one. Freeness implies that $\Delta$ is still
nonnilpotent on the generic fiber. Over any algebraically closed
extension of $\Frac R$, we can therefore choose a solution with
$\Delta\ne0$.

To specify the marked points of the domain, take independent
indeterminates $h_{ij}$ and set $L_i(w)=1+\sum_{j=1}^{m_i}h_{ij}w^j$.
Give $K_0=\overline{k(h_{ij})}$ the trivial valuation. On the
rational function field $K_0(b_1,\ldots,b_n,b_*)$ in algebraically
independent variables take the monomial Gauss valuation with positive,
rationally independent weights. Pass to a complete algebraically
closed extension of this valued field. Substitute $C_i=b_iL_i$ and $Y_j=b_j$ for $j>r$.
These substitutions embed $R$. In fact distinct multidegrees in
the $b_i$ have distinct values. In any one multidegree the
exponents of the $h_{ij}$ distinguish the coefficient monomials,
and every nonempty set of occurring multidegrees has a least value.
Choose the generic solution with $\Delta\ne0$ over this embedding.
Its coefficients $u$ have norm at most one by integrality over $R$,
and less than one by nilpotence on the closed fiber. All roots of the
$A_i$ now have norm less than one.

Take the product $U\times\mathbb G_m$ with one further designated
coordinate. Writing $z_*$ for the dual character,
Lemma~\ref{lem:intrinsic-torus-product} gives the marked expansion
$\Theta^{U\times\mathbb G_m}_{(q,h)}=\Theta^U_q z_*^h$.
This extra coordinate ensures that both end contacts and the slope
at which we shall cut are nonzero.

Write $m=(m_1,\ldots,m_r,0,\ldots,0)$ and
$d=(d_1,\ldots,d_r,0,\ldots,0)$. On the marked projective line take coordinate functions
$w^{p_i}Y_i(w)$ and $b_*w$, and evaluate at $w=1$.
The contacts at $0,\infty$ are $(p,1)$ and
$(-p-d-m,-1)$. Each root of $L_i$ gives a further contact
$(e_i,0)$. These roots have distinct residues away from $0,1,\infty$.
Once ordered, they give a generic marked domain. Indeed their
elementary symmetric functions recover the independent coefficients
of the $L_i$, and adjoining the roots is a finite extension.
At such a root the unchanged coordinates have norm less than one,
so $G_i$ has residue one. At a root of $A_i$, the condition
$\Delta\ne0$ makes the ordinary partner regular and nonzero.
The two end contacts lie off the positive mutable axes. This proves
all the open conditions in Theorem~\ref{thm:incidence}.

Let $\xi$ be the vector of valuations of all $b_j$, including
$b_*$. We identify $\xi$ with the Gauss point of the evaluation
torus, the multiplicative seminorm given by its monomial Gauss valuation.
At $w=1$, the factors $A_i(1)$ and $L_i(1)$ have
norm one. Over the field of the marked configuration the valuation
is trivial. Distinct Laurent monomials in the evaluation coordinates
have distinct values, so a polynomial relation over the field of
the marks would have a unique term of least value. Thus the evaluation
coordinates are algebraically independent over that field. The domain
and evaluation give the canonical point of the product skeleton over the
star-shaped tree of the marked line (one leg for each mark) and the Gauss
point $\xi$, as in \cite[Proposition~8.8 and \S14.1]{KY}. The fiber is in the
generic finite \'etale locus. Choose a compactification satisfying
\cite[Assumption~2.4]{KY} and take a common toroidal refinement for the
prescribed contact rays and the two toric tail rays at the cut below,
generated by $\pm(p+d,1)$, as in the proof of
Lemma~\ref{lem:local-zero-trace}. Let $\beta$ be the numerical
class of the map on this model. Partition the fiber according to
the spine of each map. The degree contributed by maps with a given
spine is its \emph{spine degree}, counting residue field degrees.
The spine of our map therefore has positive degree
\cite[Construction~9.2, Lemmas~8.20--8.22]{KY}.

Choose $A$ large enough for every effective summand of $\beta$ in
the sense of \cite[Definition~4.28]{KY}, and let $\operatorname{Wall}_A$
be the corresponding union of walls.
On each arm, a rational shear fixing the wall pointwise takes the
incoming slope at a bend to the outgoing slope. Since the coordinates of $\xi$ are
linearly independent over $\mathbb Q$, the transversality argument
in the proof of Lemma~\ref{lem:local-zero-trace} applies to every arm
and to its straight extension at a generic cut.

The evaluation leg is contracted, and the hull of the boundary legs
is stable since their number is $2+\sum_i m_i\ge3$.
Write $S^B$ for this hull, with the induced map. Its central vertex
still maps to $\xi$. By \cite[Theorem~9.14]{KY}, its spine degree
equals the positive fixed-evaluation degree of the original marked spine.

Toward $w=0$, the initial slope is $(p+d,1)$. Indeed, on the
Gauss segment each factor $w-a$ of $A_i$ has valuation
$\min(\operatorname{val}(w),\operatorname{val}(a))$, whose initial
derivative is one. The factors of $L_i$ have derivative zero there.
The same calculation on the other arms gives slopes $(e_i,0)$
at the roots of $L_i$ and $(-p-d-m,-1)$ toward infinity.
The walls in $\operatorname{Wall}_A$ and the codimension-one cones of
the toric fan lie in finitely many rational hyperplanes. Since $\xi$
lies on none of these hyperplanes, a sufficiently small ball about $\xi$
misses both $\operatorname{Wall}_A$ and the codimension-one skeleton
of the fan. Cut the arm of $S^B$
toward $w=0$ close enough to $\xi$ that the segment from $\xi$
to the cut is straight and lies in this ball. The opposite
derivatives at the cut meet the hypotheses of the concatenation
formula \cite[Theorem~12.4]{KY}. It expresses the positive spine
degree as a sum of products of nonnegative degrees over effective
class decompositions. The arm to $0$ thus has positive degree
in some class. The short straight segment has degree one in class
zero and degree zero in every other class \cite[Lemma~12.3]{KY}.
Concatenating it back to $\xi$ therefore preserves the degree
\cite[Theorem~12.4]{KY}.
Continuing the finite end by a straight toric tail as in
\cite[Construction~9.3, Lemma~9.12 and Definition~20.20]{KY}
gives a contribution to the
coefficient with input $(p,1)$ and exponent $(p+d,1)$.
The point $\xi$ lies in the interior of the straight positive
orthant. Proposition~\ref{prop:formal-local-chart} identifies the
augmented local expansion there with the marked Laurent polynomial.
The torus product identity now gives the asserted coefficient.
\end{proof}

\subsection{Gauss valuations and a slope criterion for a disk}

Write $A=k[U]$ and $K=\Frac A$. For a nonzero
Laurent polynomial $h$, let $\ell(h)$ be its unique exponent of
least $\mu$-weight. The product $\prod_{i=1}^r G_i/x_i$ is regular
and its least exponent is negative in every mutable coordinate.
Thus the smallest Newton polytopes $Q_a$ exist for all $a\in M$, by
Theorem~\ref{thm:minimum-source-newton-polytope}. For $v\in T_p$
and $0\ne h\in A$ their minimality gives
\begin{equation}\label{eq:corner-order-bound}
 \nu^x_p(h)\le\min_{b\in Q_{\ell(h)}}\langle b,p\rangle
                    \le\langle\ell(h),v\rangle .
\end{equation}

To obtain a curve over a function field, we look for regular functions
on which both inequalities in \eqref{eq:corner-order-bound} are
equalities. When their least exponents are linearly independent,
$\nu^x_p$ restricts to a Gauss valuation in these functions.

\begin{lemma}\label{lem:corner-gauss-subfield}
Suppose $s_1,\ldots,s_n\in A$ have linearly independent least
exponents $a_j=\ell(s_j)$, and $\nu^x_p(s_j)=J_j=\langle a_j,v\rangle$
for $v\in T_p$. Then the $s_j$ are algebraically independent
and $\nu^x_p$ on $k(s_1,\ldots,s_n)$ is their Gauss valuation
of weight $J=(J_1,\ldots,J_n)$.
\end{lemma}
\begin{proof}
Distinct monomials in the $s_j$ have distinct least exponents.
In any polynomial in these functions, let $h$ be the sum of the
terms of smallest $J$-weight $b$. Its least exponent comes from
one of these terms, so $\langle\ell(h),v\rangle=b$.
Equation~\eqref{eq:corner-order-bound} gives $\nu^x_p(h)\le b$,
and the valuation inequality gives the reverse. All remaining terms
have strictly greater value. The formula extends to fractions.
\end{proof}

On the generic curve, $\nu^x_p$ will be a positive multiple of the
order at a point $P$. The comparison between this order and the
slope on an annulus gives a criterion for the annular inverse of
a rational parameter $t$ to extend to a disk with center mapping to $P$.
We use Berkovich analytification, denoted by a superscript
$\mathrm{an}$. In an annular coordinate $t$, the \emph{Gauss point}
of radius $s$ is the seminorm
$\lvert\sum_j c_jt^j\rvert_s=\max_j\lvert c_j\rvert s^j$.
These points form the \emph{skeleton} of the annulus. Over an
algebraically closed field $\Omega$ they have type~II for
$s\in|\Omega^\times|$ and type~III otherwise
\cite[\S2]{BPR}.

\begin{lemma}[A slope criterion for a disk]\label{lem:annular-slope-disk}
Let $\Omega$ be a complete algebraically closed non-archimedean
field with a nontrivial real-valued absolute value. Let $\bar C$ be a smooth connected projective curve over
$\Omega$, and let $C$ be the complement of a finite set of points.
Suppose $t\in\mathcal O(C)^\times$ has a simple zero at
$P\in\bar C\setminus C$. Suppose an open annulus in $C^{\rm an}$
is mapped isomorphically by $t$ onto $r_1<|t|<r_2$. For
$r_1<\rho<r_2$ with $\rho\notin|\Omega^\times|$, let $b$
correspond to the Gauss point of radius $\rho$. For a nonzero rational
function $h$, let $\sigma_b(h)$ be the slope of $-\log|h|$ on
the skeleton of this annulus, with respect to the length coordinate
$-\log|t|$ increasing toward smaller radii.

Suppose $\ord_P(h)\le\sigma_b(h)$ for every
$0\ne h\in\mathcal O(C)$. Then the component of
$\bar C^{\rm an}\setminus\{b\}$ toward smaller radii contains $P$
and maps isomorphically under $t$ to the open disk $|t|<\rho$.
The inverse agrees with the given annular section near $b$ on this side.
\end{lemma}
\begin{proof}
Fix $D\in\bar C\setminus C$ different from $P$. Such a point
exists because $t$ has a pole. Choose a skeleton of the punctured
curve with legs toward $D$ and $P$, and with $b$ inside the
given annular edge. A \emph{semistable vertex set} is a finite set of
type~II points whose complement is a union of open disks and finitely
many open annuli, allowing punctured disks at the missing points of the
curve. We may include two nearby type~II points on either
side of $b$ in such a set
\cite[Proposition~3.13(3) and Corollary~5.10]{BPR}.
Truncate the legs to obtain a finite connected metric graph.

Write $\gamma$ for the genus of $\bar C$. For each $N>0$,
Riemann--Roch gives a nonzero $h_N$ with a zero of order at least
$N$ at $P$ and poles only at $D$, of total order at most
$N+\gamma$. Its divisor is $N(P-D)$ plus at most $\gamma$
differences of pairs of points, counted with multiplicity. For each
pair consider the unit current on the graph directed from the pole
toward the zero, both retracted to the graph. It is the gradient
of a continuous piecewise affine potential, harmonic away from the
two points, with one unit entering at the pole and leaving at the
zero. The slope formula \cite[Theorem~5.15]{BPR} identifies
$\sigma_b(h_N)$ with the sum of these currents across $b$ in
the chosen direction, with multiplicities. Indeed, one can first
adjoin the legs toward all zeros and poles, then sum the balancing
equations along the added trees. This gives the currents on the
same fixed graph for all $N$.

For each pair of points, the absolute current on any edge is at most
one. Indeed, if the current on an edge is nonzero, choose a regular
level strictly between the potentials at its two ends. The currents
across that level, oriented toward increasing potential, sum to one.
The residual pairs in the divisor of $h_N$ therefore contribute at
most $\gamma$ to the current across $b$, independently of $N$.
Write $I$ for the signed unit current across $b$ from $D$ to $P$.
Since $h_N$ is regular on $C$, the hypothesis gives
$N\le\ord_P(h_N)\le\sigma_b(h_N)\le NI+\gamma$.
Letting $N$ tend to infinity gives $I=1$.

The edge through $b$ must separate $D$ from $P$. Otherwise
there would either be no current through that edge, or a path joining
its two ends while avoiding its interior. Such a path would cross an
intermediate regular potential level and contribute positive current
in addition to the unit already crossing the edge. The sign of $I$
places $P$ on the side of smaller $|t|$. Each complementary open
disk attaches to the skeleton at a single point, and the type III
point $b$ has exactly two tangent directions
\cite[Lemma~5.12(2)]{BPR}. Thus $b$ also separates $P$ from $D$
on the analytic curve.

Apply this to every boundary point $D\ne P$. On the component $B$
toward smaller radii, $t$ has no pole and its only zero is $P$.
In fact $|t|<\rho$ throughout $B$. Cut the annular edge again
at a nearby type III point of radius $\rho'<\rho$. On the closure
of the component beyond this second cut, the maximum principle gives
$|t|\le\rho'$, since its only boundary point has radius $\rho'$.
On the part between the cuts the annular coordinate gives the strict
bound directly. The boundary of $B$ in $\bar C^{\rm an}$ is
$\{b\}$, so $B$ is open and closed in the inverse image of the
disk $|t|<\rho$. The finite flat map
$\bar C^{\rm an}\to(\mathbb P^1_\Omega)^{\rm an}$ therefore
restricts to a finite flat map from $B$ to that disk. Its degree
is the length of the fiber at zero, which is one since $P$ is a
simple zero. The restriction is consequently an isomorphism. Its
inverse agrees on $r_1<|t|<\rho$ with the original section, whose
image there lies in $B$.
\end{proof}

\subsection{Constructing the disk on the generic curve}

We construct a disk on which the zeros of the mutable coordinates give
the polynomials required by Lemma~\ref{lem:polynomial-contact-coefficient}.
The weight $p$ selects a point on the completion of the generic curve.
The least exponents at $\mu$ give a section over an annulus, which
we extend to a disk containing that point.

\begin{proposition}\label{prop:corner-disk}
Assume $0\ne p\in L$, $v\in T_p$, and $s_1,\ldots,s_n$ satisfy
Lemma~\ref{lem:corner-gauss-subfield}. Suppose their least exponents
$a_j$ form a lattice basis, each $s_j$ has a unique $p$-least
exponent, and these latter exponents also form a lattice basis.
If each designated unit is a Laurent monomial in the $s_j$ up to
a scalar, then $v\in L$ and $[z^v]\Theta_p>0$.
\end{proposition}

\begin{proof}
The case of a torus follows from $\Theta_p=z^p$. We may thus
assume $n\ge2$, since in dimension one an exchange independent
of the only coordinate would be constant. Since $p\ne0$ and the
$p$-least exponents form a basis, some $J_j$ is nonzero. Let
$g=\gcd(J_1,\ldots,J_n)>0$. The equalities
$\langle a_j,v\rangle=J_j$ give $v/g\in L$, since the $a_j$
form a basis of $M$. The degree of $s^b$, for $b\in\mathbb Z^n$,
is $J(b)=\sum_j b_jJ_j$. Take monomial
coordinates $y_1,\ldots,y_{n-1},t$
in the $s_j$ corresponding to a basis of $\ker J$ and a vector
of degree $g$. Write $E=k(y_1,\ldots,y_{n-1})$. By the lemma,
$E$ has trivial $p$-valuation and $k(s_1,\ldots,s_n)=E(t)$.
The $p$-least exponents form a lattice basis. Hence degree-zero
monomials in the $s_j$ have residues generating all toric residue
characters, and reduction identifies $E$ with the residue field
$\kappa(\nu^x_p)$.
The same basis gives $\nu^x_p(K^\times)=g\mathbb Z$.

Invert in $A$ the nonzero homogeneous polynomials in the $s_j$
and denote the result by $S$. This is a regular curve ring of
finite type over $E$. Indeed the corresponding localization of
$k[s_1,\ldots,s_n]$ is $E[t^{\pm1}]$, and $K/E$ has
transcendence degree one. Any subfield of $K$ algebraic over $E$
has trivial $p$-valuation, so reduction embeds it in $E$ over
$E$. Thus $E$ is algebraically closed in $K$. In characteristic
zero the curve $C=\Spec S$ is smooth and geometrically integral.
The $p$-valuation gives an $E$-rational point $P$ on the
smooth projective completion $\bar C$, with
$\nu^x_p=g\ord_P$ and $\ord_P(t)=1$. Since $t$ is invertible
on $C$, this point is in $\bar C\setminus C$.

We now give $E$ the valuation induced by $\mu$, and set
$\rho=|t|_\mu$. The $\mu$-least exponents of $y_1,\ldots,y_{n-1},t$
form a lattice basis, so no nonzero power of $\rho$ lies in
$|E^\times|$. There is a section of $C\to\mathbb P^1_t$
on an annulus about radius $\rho$ over the completion $\widehat E$.
To construct it, replace the $s_j(x)$ by their $\mu$-least
monomials, including their coefficients. Invert this monomial change
of coordinates to obtain $x^0_i(y,t)$. Write $x_i=x^0_i U_i$
in the original equations for $s_j$, dividing each equation by
its prescribed monomial in $y,t$. At radius $\rho$ the equations
reduce to $U^{a_j}=1$, and every omitted coefficient has norm
strictly less than one. The same strict bound holds on a small
closed annulus around $\rho$. At $U=(1,\ldots,1)$ the Jacobian
of the reduced equations has rows $a_j$ and determinant $\pm1$.
The Jacobian of the full equations in $U_1,\ldots,U_n$ differs
from it by a matrix of norm less than one, so its determinant is
a unit in the annular algebra. Newton iteration in this algebra gives
a section with $\lVert U_i-1\rVert<1$. The inverse function
theorem makes its image open near the central circle. Real annular
radii allow the coefficients to stay in $\widehat E$.

Let $\Omega$ be the completed algebraic closure of $\widehat E$.
The Gauss point $b$ of the section at radius $\rho$ is of type III,
since $|\Omega^\times|$ is the divisible hull of $|E^\times|$.
Orient the annulus toward decreasing $|t|$. For a rational function
$h$ let $\sigma(h)$ denote the outgoing slope of $-\log|h|$
at $b$, with respect to $-\log|t|$. The $t$-exponent of $x^0_i$
is $v_i/g$, since $\langle a_j,v\rangle=J_j$. The unique
$\mu$-least Laurent monomial remains dominant on a short annulus.
Consequently $\sigma(h)=\langle\ell(h),v\rangle/g$ for $h\in A\setminus\{0\}$.
For every homogeneous denominator in $S$ the bound
\eqref{eq:corner-order-bound} is an equality, by
Lemma~\ref{lem:corner-gauss-subfield}. It follows that
\begin{equation}\label{eq:corner-slope-bound}
 \ord_P(h)\le\sigma(h)\qquad(0\ne h\in S).
\end{equation}

We extend the slope bound first to the completion $\widehat E$, then
to its completed algebraic closure $\Omega$. It suffices to work
with a finite-dimensional $E$-subspace $V\subset S$, since each
element of $S\otimes_E\Omega$ belongs to some $V\otimes_E\Omega$.
For an integer $j$, its subspace $V_j$ of functions of order
at least $j$ at $P$ is the kernel of finitely many Laurent
coefficient functionals at this rational smooth point. Thus $V_j$
commutes with base change. In $V_j\otimes_E\widehat E$ each nonzero
function still has slope at least $j$, by density. Indeed a
nonzero annular series has a unique dominant power at the type III
radius $\rho$, and a sufficiently small perturbation in Gauss norm
preserves that power. Restrictions to the annulus are injective,
since it is open in a geometrically integral analytic curve.
Choose a nonzero vector in $V_j\otimes_E\widehat E$ and pass to
the kernel of the coefficient of its dominant power. Repeating gives
a basis with distinct dominant powers $t^{k_i}$, all $k_i\ge j$.
After any further scalar multiplication by elements of $\Omega^\times$,
the Gauss norms of these basis vectors belong to distinct cosets
$\rho^{k_i}|\Omega^\times|$. In any nonzero linear combination the
largest term therefore survives, with slope at least $j$.
This proves \eqref{eq:corner-slope-bound} after base change.

Write $r_1<|t|<r_2$ for the domain of the annular section over
$\Omega$, with $r_1<\rho<r_2$. The function $t$ is a unit on
$C_\Omega$ and has a simple zero at $P$. Since
\eqref{eq:corner-slope-bound} holds for every nonzero regular function on
$C_\Omega$, Lemma~\ref{lem:annular-slope-disk} gives an inverse to $t$
on $|t|<\rho$. This inverse agrees with the annular section on
$r_1<|t|<\rho$. The two sections therefore define a section
$\{|t|<r_2\}\longrightarrow\bar C_\Omega^{\mathrm{an}}$.
It sends $0$ to $P$ and every other point into $C_\Omega^{\mathrm{an}}$.

Extend the valued field to a complete algebraically closed field
containing $t_0$ of norm $\rho$, and substitute $t=t_0 w^g$.
The section is now defined on $|w|<(r_2/\rho)^{1/g}$, a disk of
radius greater than one. Since $g\ord_P(x_i)=p_i$, the function
$w^{-p_i}x_i(w)$ extends analytically across $w=0$ with nonzero
value there. At the Gauss point of the closed unit disk its strictly
dominant power is $w^{v_i-p_i}$
and its norm is $|x_i|_\mu<1$.
The dominant exponent is nonnegative because the function is analytic
at zero. Weierstrass preparation gives exactly $d_i=v_i-p_i$
roots in the open unit disk, counted with multiplicity. All are nonzero.
For a designated unit the normalized
function is constant, since the unit is a monomial in $y,t$.
To see that the roots are simple and distinct, consider the open set
where all $s_j$ are invertible. There the algebraically independent
$y_j$ define a dominant map to $\mathbb G_m^{n-1}$ with generic
fiber $C$. The locus in this open set outside the first star has
dimension at most $n-2$, so its image misses the generic point.
On the first star a zero of $x_i$ lies in the $i$-th adjacent
torus, where the other mutable coordinates and $G_i/x_i$ are units.
The divisor $x_i=0$ is reduced by Proposition~\ref{prop:upper-factorial},
and its restriction to $C$ is a localization. Its closed points
have separable residue fields over $E$, so the divisor stays reduced
after base change. The section of $t$ is locally an isomorphism,
and the power map $w\mapsto t_0w^g$ is unramified away from zero.

Let $A_i$ be the monic product of the resulting factors for $x_i$.
The analytic quotient $C_i=w^{-p_i}x_i(w)/A_i$ has Gauss norm
less than one, since $\lVert A_i\rVert=1$.
If $A_*=\prod_i A_i$ has positive degree, replace each mutable
quotient by its polynomial remainder modulo $A_*^2$. Division by
this monic polynomial of Gauss norm one does not increase the norm
of remainders, first for polynomials and then by continuity on the
closed disk. The replacement agrees to first order at every root.
At each root of $A_i$ the ordinary partner is regular and nonzero,
so the substituted exchange has a simple zero there. Finally,
$d_i=v_i-p_i$ lies between zero and $D_i$. Indeed $Q_{e_i}=\{e_i\}$
and $G_i/x_i$ is regular of least exponent $-e_i$. Testing the
identity~\eqref{eq:dual-polytope-minima} at
$e_i,-e_i$ gives the two bounds. Both signs of a designated
coordinate are regular units. Lemma~\ref{lem:polynomial-contact-coefficient}
now gives $[z^v]\Theta_p>0$. When $A_*=1$, the coefficient is
the straight coefficient one of Theorem~\ref{thm:corner}.
\end{proof}

\subsection{Realizing the vertices of the dual polytope}

\begin{proof}[Proof of Theorem~\ref{thm:polynomial-corner-reciprocity}]
By~\eqref{eq:dual-polytope-minima} and~\eqref{eq:forward-newton-bound},
$\Newt(\Theta_p)\subseteq T_p$
and the two support functions that we must compare are
\[ \min_{e\in T_p}\langle a,e\rangle=\min_{b\in Q_a}\langle b,p\rangle, \qquad \nu^z_a(\Theta_p)=\min_{b\in P_a}\langle b,p\rangle . \]
For fixed $a$, the right sides are continuous and homogeneous
in $p$. For $p=0$ both polytopes in the first assertion of the
theorem are $\{0\}$. Fix $p\ne0$ and a vertex $v$ of $T_p$.

Choose a basis $q_1,\ldots,q_n$ of $L$. Order exponents
lexicographically by their pairings with $p,q_1,\ldots,q_n$.
Let $m(a)$ be the least exponent of $Q_a$ in this order.
Every finite set of these comparisons persists after replacing the
order by the weight $p_\epsilon=p+\epsilon q_1+\cdots+\epsilon^n q_n$
for small positive rational $\epsilon$. Choose normalized sections
$f_a$ as in Theorem~\ref{thm:minimum-source-newton-polytope}.
The map $a\mapsto m(a)$ is injective. Otherwise choose
$p_\epsilon$ selecting a common least exponent of distinct
$f_a,f_b$, and cancel its coefficient in a linear combination.
The increase in valuation would contradict valuative independence.
It follows that the lexicographically least exponent of every
nonzero regular function equals $m(b)$ for some label $b$ in
its basis expansion. Also $m(ka)=km(a)$ for every integer $k>0$.
Multiplication by a designated unit $x^c$ and its inverse gives
$Q_{a+c}=Q_a+c$, hence $m(a+c)=m(a)+c$.

We first find a full-dimensional open convex cone $\mathcal V$
inside the minimizing normal cone at $v$, on which $m(a)=Ba$
for an integral linear transformation $B$. We can work in the
span of the mutable labels. Designated translations belong to
the \emph{lineality space} of the normal cone, its largest vector
subspace, by the support formula
and $Q_{a+c}=Q_a+c$. In the interior of that normal cone the
primary minimum on $Q_a$ is the linear function $\langle a,v\rangle$.
We refine this cone successively for $q_1,\ldots,q_n$.
Define the polytopes at rational labels by homogeneity.
On a cone where all preceding minima are linear, the corresponding
minimizing faces satisfy the same Minkowski inclusion as the
original polytopes. Indeed $Q_{a+b}\subseteq Q_a+Q_b$ passes to
minimizing faces whenever the two minima agree, and the minimizing
face of a sum is the sum of the faces. Repeating this gives the
inclusion through all preceding weights.

For each rational label $a$ in this cone, let $F_a$ be the face
selected by the preceding weights and set
$h(a)=\min_{b\in F_a}\langle b,q_j\rangle$.
The inclusion for the faces makes $h$ superadditive. It is
homogeneous and takes integer values on lattice points. Its
absolute value is locally bounded since $Q_a$ is contained in the
sum of the coordinate polytopes $Q_{\pm e_i}$ with coefficients
$|a_i|$ of the corresponding signs. The function extends by
concavity to a locally Lipschitz function on the open cone.
As in Lemma~\ref{lem:integral-support}, its gradients at
differentiability points are integral. Indeed, at such a point $u$
choose mutable lattice points $b_k$ within bounded distance of
$ku$. For large $k$, both $b_k$ and $b_k+e_i$ lie in the cone. The integers
$h(b_k+e_i)-h(b_k)$ tend to $\partial_i h(u)$ by homogeneity
and differentiability. Only finitely many gradients occur on a
smaller interior ball. Concavity and density of differentiability
points express $h$ there as the minimum of their supporting
linear forms. It follows that $h$ is linear on a smaller open ball
and, by homogeneity, on its open conical span.
After the last refinement the $q_j$ determine $m(a)$, giving
an integral matrix on the mutable labels. Extend it by the identity
on designated labels and add their full real span to the cone.
This gives $B$ and $\mathcal V$. If all coordinates are designated,
the assertion follows directly on the torus.

In fact $B$ is unimodular. Injectivity of $m$ forces $B$ to
have full rank, since a lattice point sufficiently far inside
$\mathcal V$ and its translate by a fixed integer kernel vector
would have the same image. Let $D=|\det B|$ and choose an interior
lattice point $a$ of $\mathcal V$. For any coordinate character
$e_j$, the least exponent of the regular function $x_j f_{ka}$
is $kBa+e_j=m(b)$ for some $b$. For $k$ sufficiently large
this lies in $B\mathcal V$. Then $c=DB^{-1}m(b)$ belongs to
$\mathcal V\cap M$, and $m(c)=Bc=Dm(b)=m(Db)$.
Injectivity gives $c=Db$, hence $b\in\mathcal V$ and $e_j=B(b-ka)$.
Thus $BM=M$. On $\mathcal V$ the primary minimum now gives
$\langle a,v\rangle=\langle Ba,p\rangle$, so $v=B^*p$, where
$B^*:L\to L$ is the dual lattice map.

Choose a basis $a_1,\ldots,a_r$ of the mutable character lattice
inside the last mutable cone. To find such a basis, extend a
primitive interior vector to a lattice basis and add large multiples
of it to each other vector. Complete this by the designated coordinate
characters $a_{r+1},\ldots,a_n$. Take $s_j=f_{a_j}$ for $j\le r$
and $s_j=x^{a_j}$ for $j>r$. Set $a_0=\sum_{j=1}^n a_j$.
It lies in $\mathcal V$ and exposes $v$ on $T_p$. For all
$0\le j\le n$ one has $m(a_j)=Ba_j$, using designated translation
for the coordinate characters.

For every sufficiently small positive rational $\epsilon$, the unique
$p_\epsilon$-least exponent of $Q_{a_j}$ is $Ba_j$ for all
$0\le j\le n$. Choose $N_\epsilon\in\mathbb Z_{>0}$ such that
$p'=N_\epsilon p_\epsilon\in L\setminus\{0\}$. The minimum at
$a_j$ on $T_{p'}$ is $\langle Ba_j,p'\rangle$. Since
$a_0=\sum_{j=1}^n a_j$ and $B$ is linear, the minimum at $a_0$
equals the sum of the $n$ individual minima. Any minimizer at $a_0$
thus attains all $n$ minima, so it is unique and equals $v'=B^*p'$.
The least exponents of the $s_j$ form the basis $a_j$ at $\mu$
and the basis $Ba_j$ at $p'$. We also have
$\nu^x_{p'}(s_j)=\langle a_j,v'\rangle$.
Proposition~\ref{prop:corner-disk} now gives $[z^{v'}]\Theta_{p'}>0$.

Since $v'$ minimizes $a_0$ on $T_{p'}$ and occurs in $\Theta_{p'}$,
the containment $\Newt(\Theta_{p'})\subseteq T_{p'}$ gives equality
of the two minima at $a_0$. By the support identities at the start
of the proof, $P_{a_0}$ and $Q_{a_0}$ therefore have the same
minimum at $p'$, and by homogeneity at $p_\epsilon$. Keeping
$a_0$ fixed and letting $\epsilon$ tend to zero gives equality
at $p$. Thus $\Newt(\Theta_p)$ meets the face of $T_p$ on
which $a_0$ is minimal. That face is $\{v\}$. This proves
$\Newt(\Theta_p)=T_p$, since $v$ was an arbitrary vertex.
The same support identities now give $P_a=Q_a$ for every $a$
and the asserted equality of valuations.
\end{proof}

Reciprocity also determines how the least exponent of a section on
$Q_a$ changes with the initial seed. Let $t,t'$ be two reachable
seeds of $U$ with polynomial corners. Write $x,x'$ for their
coordinates and $z,z'$ for the corresponding forward coordinates.
Write $M_t,M_{t'}$ for the character lattices of the two seed tori
and $L_t,L_{t'}$ for their duals.
Let $T:L_{t,\mathbb R}\to L_{t',\mathbb R}$ be the change of Gauss
weights from $t$ to $t'$ of Lemma~\ref{lem:boundary-gauss-change}.
In each seed choose a strictly positive separating weight.

\begin{corollary}\label{cor:corner-least-exponent-change}
Let $f_a$ be a normalized section on the smallest Newton polytope
$Q_a$ in $t$. Its least exponent $a'$ in $t'$ is given by
\[ a'_j=\nu^z_a\bigl(\Theta^t_{T^{-1}(e'_j)}\bigr), \]
where $e'_j$ are the coordinate cocharacters of $t'$ and the
theta functions are expanded in the forward chart for $t$. The map
$a\mapsto a'$ extends to an integral piecewise linear homeomorphism
$M_{t,\mathbb R}\to M_{t',\mathbb R}$, linear on the cones of a
finite rational fan.

If $t'$ is obtained by mutation at $i$ and the exchange Laurent
polynomial $\widehat F_i$ is polynomial with nonzero constant term,
let $H_i=z_i\Theta^t_{-e_i}$. Then
\[ a'_j=a_j\quad(j\ne i),\qquad a'_i=-a_i+\tau_{H_i}(a). \]
\end{corollary}
\begin{proof}
By Section~\ref{subsec:smallest-newton-polytope}, $f_a$ is supported
on the smallest Newton polytope at $a'$ in $t'$. Multiplying by
a scalar normalizes its least coefficient there and preserves its
Gauss valuations. In the forward chart for $t'$, the theta function of label $e'_j$
is $z'_j$ by Theorem~\ref{thm:corner}. In $t'$, reciprocity~\eqref{eq:corner-reciprocity} gives
$a'_j=\nu^{x'}_{e'_j}(f_a)$. Transporting the valuation to $t$ and
applying reciprocity there gives the first formula. Each coordinate is
the minimum of finitely many integral linear forms. Interchanging the
two seeds gives the inverse map on lattice points, hence on the real
spaces by homogeneity and continuity. For the stated mutation,
$T^{-1}(e'_j)=e_j$ for $j\ne i$ and $T^{-1}(e'_i)=-e_i$.
Substitution gives the second formula.
\end{proof}

For example, take the type $A_2$ exchanges $xx'=1+y$ and
$yy'=1+x$. The least exponent of $y'$ in $(x,y)$ is $(0,-1)$.
In $(x',y)$ the same function is $(x'+1+y)/(x'y)$, with least
exponent $(-1,-1)$. Its first coordinate thus changes from $0$
to $-1$. For cluster $g$-vectors, defined using principal
coefficients, changing the initial seed negates the mutated coordinate
\cite[Proposition~6.8, (6.15)]{FZIV}.

\section{Boundary valuations and polynomial frozen variables}
\label{sec:boundary-valuations}

Let $D_1,\ldots,D_b$ be distinct prime divisors on a normal variety
$X$ over $k$, and let $U=X\setminus\bigcup_iD_i$. A function
in $\mathcal O(U)$ extends to $X$ precisely when its order along
each $D_i$ is nonnegative. Write $\ord_{D_i}$ for the order of
vanishing at the generic point of $D_i$, with poles counted negatively.
For a $k$-basis $(\vartheta_m)_{m\in M}$ of $\mathcal O(U)$,
indexed by a set $M$, we ask when the basis is valuatively independent
for each boundary order. This means that the order of any
nonzero finite linear combination is the minimum of the orders of its
nonzero terms. It holds precisely when the initial forms of the terms
of least order cannot cancel. For a polynomial
frozen variable, we must consider each prime component of its zero divisor.

\subsection{Computing an order in a Laurent chart}

In torus coordinates $x=(x_1,\ldots,x_d)$, let $v_{x,w}$ denote the
Gauss valuation of weight $w$. We first give a condition under which
the terms of least weight in a basis expansion cannot cancel.

\begin{lemma}\label{lem:boundary-leading-independence}
Let $(b_m)_{m\in M}$ be nonzero Laurent polynomials in $x$. Suppose
that their expansions have the form
\[ b_m=x^{p_m}\left(c_m+\sum_{u\ne0} c_{m,u}x^u\right), \qquad c_m\in k^\times, \]
with distinct $p_m$. Let $w$ be a real weight. Assume that all
corrections satisfy $\langle w,u\rangle\ge0$, and that there is one
real linear functional $a$ such that $a(u)>0$ for every correction
with $\langle w,u\rangle=0$. Then, for any nonzero finite sum,
\[ v_{x,w}\!\left(\sum_m d_m b_m\right) =\min_{d_m\ne0}\langle w,p_m\rangle . \]
\end{lemma}

\begin{proof}
Order the exponents occurring in the sum first by $w$, then by
$a$, and finally lexicographically in the given lattice coordinates.
This order is compatible with translation. The first exponent in the
support of $b_m$ is $p_m$. The least of the $p_m$ with
$d_m\ne0$ therefore occurs only once in the sum, with nonzero
coefficient. Its $w$-weight is the asserted minimum.
\end{proof}

A common pointed rational polyhedral cone containing the correction
vectors supplies the functional $a$ whenever $w$ is nonnegative on
the cone. One can also test independence directly. For a discrete
valuation $v$ trivial on $k$, let $F^j=\{g:v(g)\ge j\}$,
with $v(0)=+\infty$. A family satisfies the minimum formula for all
finite linear combinations if and only if the classes of the elements of
value $j$ are linearly independent in $F^j/F^{j+1}$ for each $j$.
Indeed, cancellation of the terms with least value is exactly a linear
relation in this quotient.

We can characterize independence for all real Gauss valuations by
comparing pairs of Newton polytopes. For a vertex $a$ of a polytope
$P$, write $N_P(a)$ for the cone of weights attaining their
minimum on $P$ at $a$. Its interior in the full weight space
consists of the weights for which $a$ is the unique minimizer. This
convention applies to polytopes of every dimension.

\begin{proposition}\label{prop:boundary-vertex-criterion}
Let $(b_m)_{m\in M}$ be a linearly independent family of Laurent
polynomials in one chart, and let $\Delta_m$ be the Newton polytope of
$b_m$. The following are equivalent.
\begin{enumerate}
\item For each real Gauss weight $w$, every nonzero finite sum satisfies
      \[ v_{x,w}\Bigl(\sum_m d_m b_m\Bigr)=\min_{d_m\ne0} v_{x,w}(b_m). \]
\item For every two distinct labels $m,n$ and every common vertex
      $a$ of $\Delta_m$ and $\Delta_n$, the interiors of
      $N_{\Delta_m}(a)$ and $N_{\Delta_n}(a)$ are disjoint.
\end{enumerate}
Under these conditions, the Newton polytope of any such combination is
the convex hull of the union of the polytopes of its nonzero summands.
\end{proposition}

\begin{proof}
If a weight lies in both interiors, each of the two initial forms is a
nonzero scalar multiple of $x^a$. Choose the two coefficients to
cancel this monomial. Independence ensures the combination is nonzero,
and its value is strictly greater than the common minimum.

Conversely, fix a finite combination. For a dense open set of weights,
every polytope involved has a unique minimizing vertex. These vertices
are distinct by condition (2). The initial form of the combination is
therefore nonzero at the least of their weights. Both sides of the
minimum formula are continuous in $w$, as minima of finitely many
linear functions, so the formula holds everywhere. The final assertion
follows because a polytope is determined by its minimum in every
direction.
\end{proof}

\begin{corollary}\label{cor:boundary-dense-charts}
Suppose a family of Laurent polynomials is Laurent in each of a
collection of charts obtained from the initial chart by reciprocal
exchanges and integral changes of monomials. In each such chart assume
the expansions have distinct distinguished exponents and that all their
corrections belong to a common pointed rational polyhedral cone $C$
in its exponent space. The cone may vary with the chart. Transport
weights from the initial chart by Lemma~\ref{lem:boundary-gauss-change}.
If the union of the inverse images of the dual cones $C^\vee$ is
dense in the initial real weight space, the family satisfies
Proposition~\ref{prop:boundary-vertex-criterion}(1) in the initial
chart and in every reachable chart in which the family is Laurent.
\end{corollary}

\begin{proof}
For any weight in the indicated union, use the corresponding chart.
There the weight is nonnegative on every correction and a linear
functional strictly positive on $C\setminus\{0\}$ supplies the
secondary order. Lemmas~\ref{lem:boundary-leading-independence}
and~\ref{lem:boundary-gauss-change} give the minimum formula in the
initial chart. Continuity for each finite combination extends it to
every real weight. Applying Lemma~\ref{lem:boundary-gauss-change} along
the inverse word gives the assertion in the other charts.
\end{proof}

In a single chart the pointed expansions give the minimum formula
throughout $C^\vee$.

\subsection{Extension across the boundary}

For an integral Weil divisor $D$ on $X$, the nonzero global
sections of the divisorial sheaf $\mathcal O_X(D)$ are the rational
functions $g$ with $\operatorname{div}_X(g)+D\ge0$.

\begin{proposition}\label{prop:boundary-extension}
Write $h_i(m)=\ord_{D_i}(\vartheta_m)$. The minimum formula
for the prescribed basis at each $\ord_{D_i}$ holds if and only if,
for every tuple of integers $(b_i)$,
\[ \Gamma\!\left(X,\mathcal O_X\!\left(-\sum_i b_iD_i\right)\right) =\bigoplus_{\substack{m\in M\\ h_i(m)\ge b_i\ \text{for every }i}} k\vartheta_m . \]
In particular the labels with all $h_i(m)\ge0$ give a basis of
$\mathcal O(X)$.
\end{proposition}
\begin{proof}
By normality the sections on the left are the functions in
$\mathcal O(U)$ whose specified orders meet the thresholds.
Taking the minimum of the orders of their summands gives the identity.
Conversely, if cancellation raises the order of a nonzero sum at
$D_i$, take $b_i$ equal to that raised order. Choose the other
thresholds below the orders of every summand there. The sum belongs to
the left side although one of its unique basis summands is excluded
on the right.
\end{proof}

For each $i$, suppose a torus chart identifies $\ord_{D_i}$ with
a Gauss weight $w_i$ and satisfies Lemma~\ref{lem:boundary-leading-independence}
with distinguished exponents $p_i(m)$. Then the minimum formula
holds and $h_i(m)=\langle w_i,p_i(m)\rangle$. Each divisor can be
treated in its own chart. Suppose an open
neighborhood of the generic point of $D_i$ is
$\Spec  k[t,y_2^{\pm1},\ldots,y_d^{\pm1}]$, with
$D_i$ cut out there by $t=0$. The normalized order along $D_i$
is the Gauss valuation of weight $(1,0,\ldots,0)$ in these
coordinates. Indeed, the coefficient of the lowest power of $t$ is a
nonzero rational function in the residue field
$k(y_2,\ldots,y_d)$. Lemma~\ref{lem:boundary-gauss-change} then
computes the same order in any chart connected to this one by reciprocal
exchanges and integral changes of monomials.

Suppose the labels in a fixed chart form a lattice and a marked mirror
function $W_i=\sum_{u\in E_i}c_u z^u$ has been identified. Its
tropicalization is $W_i^{\mathrm{trop}}(m)=
\min_{u\in E_i}\langle u,m\rangle$. Its Newton inequalities describe the order in
Proposition~\ref{prop:boundary-extension} when
$\langle w_i,p_i(m)\rangle=W_i^{\mathrm{trop}}(m)$.
When the exponent map $p_i$ is given by a finite homogeneous
piecewise-linear map on the real label space, this equality can be checked
on generators of the cones of a common refinement with the domains of
linearity of $W_i^{\mathrm{trop}}$, including generators of their
lineality spaces. Both sides are linear on each cone. Under this
identification the regular labels satisfy the individual inequalities
$\langle u,m\rangle\ge0$ for all $u\in E_i$ and all $i$.

\emph{Rees filtrations.} Suppose $\mathcal O(X)$ is the Rees algebra
$R=\bigoplus_{\ell\ge0}z^\ell F_\ell B$ of an exhaustive increasing
multiplicative filtration on a $k$-algebra $B$, with $1\in F_0B$
and $\operatorname{div}_X(z)=D_1$. Assume the Rees grading extends
to $\mathcal O(U)$, the basis is homogeneous and the minimum
formula holds at the boundary divisors. Since $R[z^{-1}]=B[z^{\pm1}]$,
Proposition~\ref{prop:boundary-extension} gives a basis of $B$
consisting of the degree-zero $\vartheta_m$ with $h_i(m)\ge0$ for
$i>1$. Those with $-h_1(m)\le\ell$ form a basis of $F_\ell B$, since
$\vartheta_m\in F_\ell B$ is equivalent to $z^\ell\vartheta_m\in R$. The least such
$\ell$ is $-h_1(m)$. Indeed, if $h_1(m)>0$, then
$z^{-1}\vartheta_m$ would be regular of Rees degree $-1$.
If $h_1=W_1^{\mathrm{trop}}$, this threshold equals
$\max_{u\in\Supp(W_1)}\langle-u,m\rangle$. For a compatible auxiliary
grading with homogeneous basis and $\deg z=0$, counting labels of
threshold $\ell$ in a given degree gives the dimension of that
component of $F_\ell B/F_{\ell-1}B$, with $F_{-1}B=0$.

\subsection{Multiplication by frozen variables}

A regular function whose divisor is supported on the removed boundary
is a unit on $U$. Under positivity one can check directly when this
unit permutes the theta basis.

\begin{lemma}\label{lem:positive-boundary-unit}
Suppose $k$ has characteristic zero and the structure constants of a
basis of a $k$-algebra are nonnegative integers. If both $f$
and $f^{-1}$ belong to the additive monoid of nonnegative integral
combinations of basis elements, multiplication by $f$ permutes the
basis. If, in one chart, the basis elements have distinct least
exponents for a common translation-invariant order and $f=x^v$,
the permutation adds $v$ to those exponents.
\end{lemma}

\begin{proof}
The two multiplication maps are inverse automorphisms of the free
commutative additive monoid on the basis. Its indecomposable nonzero
elements are exactly the individual basis vectors. The final assertion
follows by multiplying a Laurent expansion by $x^v$ and taking its
least exponent.
\end{proof}

\begin{proposition}\label{prop:polynomial-frozen-basis}
Partition the boundary primes into nonempty sets $I_1,\ldots,I_s$.
Suppose the boundary orders satisfy the minimum formula and
there are regular functions with
$\operatorname{div}_X(f_i)=\sum_{D\in I_i}e_DD$, $e_D>0$.
If multiplication by each $f_i$ permutes the basis lines on $U$,
that is, the one-dimensional spaces $k\vartheta_m$,
set
\[ h_i(m)=\min_{D\in I_i} \left\lfloor\frac{\ord_D(\vartheta_m)}{e_D}\right\rfloor, \qquad M_0=\{m:h_i(m)=0\text{ for every }i\}. \]
The $f_i$ are algebraically independent and
\[ \mathcal O(X)=\bigoplus_{m\in M_0}k[f_1,\ldots,f_s]\,\vartheta_m . \]
After any base change from $k[f_1,\ldots,f_s]$, the images of the
$\vartheta_m$ for $m\in M_0$ form a basis over the new coefficient ring.
\end{proposition}
\begin{proof}
The regular labels are the ones with $h_i\ge0$ for every $i$, by
Proposition~\ref{prop:boundary-extension}. The commuting permutations induced
by the $f_j$ add $\delta_{ij}$ to $h_i$. Each orbit thus has
a unique line with all $h_i=0$, and the regular lines in that orbit
are precisely its nonnegative translates. The resulting vectors are
nonzero multiples of distinct basis vectors. In any one orbit this
proves algebraic independence of the $f_i$. Taking all orbits
gives the direct sum, and tensoring it proves base change.
\end{proof}

\begin{boundarycorollary}\label{cor:coordinate-boundary-reduced}
In Proposition~\ref{prop:polynomial-frozen-basis}, suppose the labels
form a lattice $M$ and there is a Laurent chart in which
$\vartheta_m$ has least exponent $m$ for a fixed translation-invariant
total order on $M$. Assume $h_i(nm)=nh_i(m)$ for every $i$,
$m\in M$ and integer $n\ge1$. For any $J\subseteq\{1,\ldots,s\}$,
the quotient $\mathcal O(X)/(f_j:j\in J)$ is geometrically reduced.
The images of the regular basis elements with $h_j(m)=0$ for every
$j\in J$ form a basis of this quotient.
\end{boundarycorollary}
\begin{proof}
Multiplication by $f_j$ sends the regular basis lines onto those
with $h_j(m)>0$, by the translation of orders in
Proposition~\ref{prop:polynomial-frozen-basis}. The ideal is therefore
spanned by the regular basis vectors with $h_j(m)>0$ for some
$j\in J$. This gives the stated quotient basis over every extension
of the ground field.

Over such an extension, lift a nonzero class to a sum of the
surviving basis vectors, and let $m$ be its least exponent.
For $n\ge1$ the $n$-th power of the lift has least exponent $nm$.
Since the least exponents of the basis vectors are distinct, its
basis expansion has a nonzero coefficient at $\vartheta_{nm}$.
The homogeneity of the $h_i$ makes this label regular and gives
$h_j(nm)=0$ for $j\in J$. Thus every power of the class is nonzero.
\end{proof}

\subsection{Boundary orders from the first star}

For an upper LP algebra with polynomial frozen variables, the initial
seed and its first mutations determine all the required boundary
orders. Here the Gauss weights include the frozen coordinates.

\begin{theorem}\label{thm:polynomial-star-descent}
Let $k$ be a field and let $\mathcal U_S$ be the full upper
algebra of an LP seed over
$S=k[f_1,\ldots,f_s,u_1^{\pm1},\ldots,u_c^{\pm1}]$.
Let $F=\prod_i f_i$, and suppose
$B=\mathcal U_S[F^{-1}]$ has a prescribed $k$-basis
$(\vartheta_m)_{m\in M}$. For each first-star chart $t$, let
$v_{t,i}$ assign weight one to $f_i$ and zero to every other
coordinate. In the initial torus choose a weight $w$ separating
lattice points. Assume the minimum formula for the basis at
$w$, the positive and negative coordinate weights, and the
$v_{t,i}$. Define
\[ h_i(m)=\min_t v_{t,i}(\vartheta_m),\qquad M_0=\{m:h_i(m)=0\text{ for all }i\}. \]
Then
\[ \mathcal U_S=\bigoplus_{m\in M_0}k[f_1,\ldots,f_s]\,\vartheta_m . \]
Write $E(m)$ for the least exponent of $\vartheta_m$ at $w$.
The set $E(M_0)$ is invariant under translations by the exponent
lattice of the $u_j$. Choosing one representative from each orbit
gives an $S$-basis of $\mathcal U_S$. For every $S$-algebra $S'$,
its images form an $S'$-basis of $S'\otimes_S\mathcal U_S$.

These boundary orders are explicit. Let $J_i$ be the mutable
indices $j$ for which the ordinary exchange $F_j$ is an
associate of $f_i$. With $\Delta_m=\Newt(\vartheta_m)$ one has
\[ h_i(m)=\min\bigl(\{\tau_{\Delta_m}(e_i)\}\cup \{\tau_{\Delta_m}(e_i+e_{x_j}):j\in J_i\}\bigr), \qquad \tau_\Delta(v)=\min_{a\in\Delta}\langle v,a\rangle , \]
where $e_i$ is the pure $f_i$-weight and $e_{x_j}$ the
$j$-th mutable weight of the initial torus.
\end{theorem}
\begin{proof}
Let $L_t$ be the Laurent rings in the first star over $S$.
The LP starfish theorem gives $\mathcal U_S=\bigcap_t L_t$.
Clearing a common denominator gives $B=\bigcap_t L_t[F^{-1}]$.
An element of $B$ lies in $L_t$ precisely when each
$v_{t,i}$ is nonnegative. At the $j$-th neighbor this
valuation has initial weight $e_i+d_{ij}e_{x_j}$ by
Lemma~\ref{lem:boundary-gauss-change}, where $d_{ij}$ is the
minimum $f_i$-exponent of $\widehat F_j$. Its normalizing
monomial uses only mutable variables, and $F_j$ is irreducible.
Thus $d_{ij}=1$ for $j\in J_i$ and zero otherwise. The
minimum formula now identifies the regular subspace with the span of
the labels satisfying every $h_i\ge0$.

Distinct labels have distinct least exponents at $w$, since otherwise
their leading terms could be canceled. These exponents are exactly the
least exponents of nonzero elements of $B$. For each coefficient
character $z^a$ generated by the $f_i,u_j$ and their inverses,
multiplication by $z^a$ translates this set bijectively. Write
$\sigma_a m$ for the unique label with $E(\sigma_a m)=a+E(m)$.
In the expansion of $z^a\vartheta_m$ the label $\sigma_a m$
occurs, and all other least exponents have strictly larger
$w$-weight. At each test valuation $v$, the minimum formula gives
$v(\vartheta_{\sigma_a m})\ge v(z^a)+v(\vartheta_m)$.
The inverse translation gives equality. In particular translation
by $f_j$ adds $\delta_{ij}$ to $h_i$, and translation by a
$u_j$ preserves every $h_i$.

In each orbit under the $f_i$, take the unique label of orders zero.
Its translates by nonnegative powers index precisely the regular
labels in that orbit. If the $u_j$ are also used, keep a single
zero-order label modulo their free translations. In either case let
$b$ be the chosen representative for a regular label $m$ and
$z^a$ the unique translating monomial, so $m=\sigma_a b$.
Set $g_m=z^a\vartheta_b$. These functions are regular, and their
distinct least exponents prove independence.

For spanning, let the coordinate box of a Laurent polynomial be the
product of the smallest closed intervals containing its coordinate
exponents. Equality at the two coordinate weights in each direction
gives $g_m$ and $\vartheta_m$ the same box. Each constituent of
$g_m$ has its box inside this box by the minimum formula.
Its expansion in the regular labels has a nonzero coefficient at
$\vartheta_m$ and only strictly later least exponents elsewhere.
Inside any fixed box there are only finitely many possible labels,
since $E$ is injective. Triangular substitution thus expresses
every $\vartheta_m$ as a finite sum of the $g_n$, proving both
module bases. Tensoring a free module gives the stated base change.
\end{proof}

In particular the frozen divisor has explicit reduced components.
In the factorization argument of Proposition~\ref{prop:upper-factorial},
apply the prime normalization to $f_i$ in the initial Laurent ring.
Its order at the adjacent prime for $x_j$ is one for $j\in J_i$
and zero otherwise. Hence
$f_i=p_i\prod_{j\in J_i}x_j$
where $p_i$ and the indicated $x_j$ are pairwise nonassociate primes.
Their orders in the initial chart are respectively $e_i$ and
$e_i+e_{x_j}$, by the same first-star comparison.

When the structure constants are nonnegative integers in characteristic
zero and $1$ belongs to the basis, the coefficient characters also
permute the basis lines themselves. The two coordinate tests confine
every constituent of a character and of its inverse to its singleton
box. Thus each is a scalar multiple of a basis element. Multiplication
by these two elements gives nonnegative integral matrices inverse up to
a positive scalar. They preserve the extreme rays of the real cone of
finite nonnegative basis expansions. If the basis is monic for a common
translation-invariant order, each character therefore translates the
basis vectors with this normalization.

For example, the rank-one seed with exchange $s$ over $k[s]$ has
upper algebra $k[x,y]$ with $xy=s$. Its two boundary orders on
$x^m s^n$ are $n$ and $m+n$. Specializing the free module of
Theorem~\ref{thm:polynomial-star-descent} at $s=0$ gives a basis
of $k[x,y]/(xy)$. The specialized exchange is zero. In general,
identifying the tensor product with an upper algebra requires an LP
seed in the fiber and a comparison of its Laurent rings.

\section{Explicit reciprocal functions and criteria for realization}\label{sec:explicit-reciprocal-functions}

At a polynomial corner Theorem~\ref{thm:polynomial-corner-reciprocity}
determines the Newton polytope at every label. Each reachable seed
determines a cone of weights generated by the slopes of its mutable
coordinates and both signs of the slopes of its designated units
(Corollary~\ref{cor:source-seed-slopes}). On each such \emph{source
seed cone} we identify the unique normalized function with that polytope. For
some labels outside the seed cones we construct a function directly
from the exchanges. In a chart defined by a negative Newton vertex
we give criteria for equality in $P_a\subseteq Q_a$.

\subsection{Principal opens and source seed cones}

Keel, White and Yu identify the associated graded algebra of the
filtration defined by a regular function with the mirror algebra of
its nonvanishing locus \cite[Corollary~15.3, Remark~15.4 and
Proposition~21.11(1)]{KWY}. We use this identification after
curve-class augmentation to compute Gauss valuations. Let $U$ be a smooth
affine log Calabi--Yau variety over an algebraically closed field of
characteristic zero with a volume-compatible torus, and
write $V=U_f=\{f\ne0\}$ for the \emph{principal open subset}
defined by $0\ne f\in k[U]$. Suppose $V$ contains a
volume-compatible torus $T$. Its cocharacter lattice identifies
the integral skeleta of $U$ and $V$. Since $f|_T$ is invertible,
it is a scalar times a character. Thus $\tau_f$ is integral and
linear in this chart.

For $j\in\mathbb Z$, let $F^j$ be the $\mathbb Z$-span of the
$\theta_p$ with $\tau_f(p)\ge j$. Write
$\operatorname{gr}_f B_{\mathbb Z}(U)=\bigoplus_j F^j/F^{j+1}$
for the associated graded algebra and $\operatorname{in}_f(\theta_p)$
for the image of $\theta_p$ in degree $\tau_f(p)$. We denote the
theta function of $V$ at the same label by $\theta_p^V$.

\begin{proposition}\label{prop:principal-open}
The filtration $F^\bullet$ of $B_{\mathbb Z}(U)$ is multiplicative, and
there is an isomorphism preserving the theta bases
\[ \operatorname{gr}_f B_{\mathbb Z}(U)\ \simeq\ B_{\mathbb Z}(V),\qquad \operatorname{in}_f(\theta_p)\longmapsto\theta_p^V . \]
When $V=T$ the graded algebra is the group algebra of the
cocharacter lattice of $T$.
\end{proposition}
\begin{proof}
For each nonzero product coefficient the regular-function inequality
gives $\tau_f(r)\ge\sum_i\tau_f(p_i)$. We identify the terms for which
equality holds. Choose a smooth projective compactification $Y$ of
$U$ with the boundary refined so that
the closure of $U\setminus T$ contains no essential boundary stratum,
as in the paragraph following \cite[Lemma~18.3]{KY}. Write
$\operatorname{div}_Y(f)=Z+D_f$, where $Z$ is the closure, with
multiplicities, of the zero divisor in $U$, and $D_f$ is supported
on the boundary. On each boundary valuation representing a label the
order supplied by $D_f$ is $\tau_f$. Keel and Yu's boundary intersection
identity \cite[Lemma~16.1]{KY} gives
\[ Z\cdot\beta=-D_f\cdot\beta =\tau_f(r)-\sum_i\tau_f(p_i) \]
for each contributing class $\beta$.
By the open restriction formula of Keel and Yu
\cite[Proposition~18.1(1),(4)]{KY}, the terms with $Z\cdot\beta=0$
count the structure disks in $V$, with their multiplicities.
Augmentation sums their coefficients across boundary refinements
\cite[Proposition~17.3(1) and Remark~17.8]{KY}.
These are exactly the terms surviving in the graded product. The theta
basis splits the filtration, proving the asserted isomorphism also on
the underlying modules. For a torus, multiplication adds the two
labels with coefficient one \cite[Lemma~9.4 and \S18]{KY}.
\end{proof}

\subsubsection{Gauss valuations from principal opens}
\label{subsec:source-probe}

A regular function that becomes invertible on a volume-compatible
torus in its nonvanishing locus defines a filtration of the forward
mirror algebra. We first give a criterion for its order to agree
with a Gauss valuation. Let $L$ be a lattice and let $B\subset k[L]$ have a
basis $(\theta_p)_{p\in L}$. Suppose $\theta_p=z^p$ on the lattice
points of a full-dimensional rational polyhedral cone $\sigma$.
Let $\rho:L\to\mathbb Z$ give a multiplicative decreasing filtration
with split order $v_\rho(\sum c_p\theta_p)=\min_{c_p\ne0}\rho(p)$.
Write $\operatorname{in}_\rho$ for the symbol in the associated
graded ring.

\begin{lemma}\label{lem:probe-graded-test}
Suppose $\rho(p)=\langle a,p\rangle$ on a full-dimensional rational
cone $\eta\subseteq\sigma$, and choose $b\in L\cap\operatorname{int}\eta$.
The following conditions are equivalent:
\begin{enumerate}
\item $\operatorname{in}_\rho(\theta_b)$ is a nonzerodivisor.
\item The split order on $B$ equals the Gauss valuation of weight $a$
      in the given chart.
\item $\operatorname{gr}_\rho B$ is a domain.
\end{enumerate}
\end{lemma}
\begin{proof}
Under the first condition multiplication by $\theta_b$ increases
the order of every nonzero element by exactly $\rho(b)$, by the
definition of multiplication of initial forms. For $g\in B\setminus
\{0\}$ and sufficiently large $n$, the support of $z^{nb}g$
lies in $\eta$. There the monomial expansion is the basis expansion.
Consequently
\[ n\langle a,b\rangle+v_\rho(g) =v_\rho(z^{nb}g)=n\langle a,b\rangle+\nu^z_a(g). \]
This proves the second condition. The associated graded ring of a
valuation on a domain is a domain, and the third condition implies
the first.
\end{proof}

\begin{theorem}\label{thm:one-sided-probe}
Let $U$ be a smooth affine log Calabi--Yau variety over an
algebraically closed field $k$ of characteristic zero, with
$\Pic (U)=0$ and a marked volume-compatible torus.
Suppose $B_k(U)$ has a Laurent realization as above, straight
on $\sigma$, with the cocharacter lattice of the marked source torus
as its label lattice $L$.
For $0\ne f\in k[U]$, suppose $U_f$ contains a volume-compatible
torus. Then $\tau_f$ is linear on $\sigma$, with integral slope
$a_f\in L^*$, and for every nonzero finite sum
\begin{equation}\label{eq:single-function-valuation}
\nu^z_{a_f}\Bigl(\sum_p c_p\theta_p\Bigr) =\min_{c_p\ne0}\tau_f(p).
\end{equation}
In particular $\nu^z_{a_f}(\theta_p)=\tau_f(p)$ for every label.
More generally, if $f_1,\ldots,f_l$ are regular functions invertible
on one volume-compatible source torus, then for any
$\lambda_j\in\mathbb R_{\ge0}$ the Gauss weight
$a=\sum_j\lambda_j a_{f_j}$ satisfies
\[ \nu^z_a\Bigl(\sum_p c_p\theta_p\Bigr) =\min_{c_p\ne0}\sum_j\lambda_j\tau_{f_j}(p). \]
\end{theorem}
\begin{proof}
Set $V=U_f$ and choose a volume-compatible torus $T'\subset V$.
The ring $k[U]$ is factorial since $U$ is smooth with trivial
Picard group, and $k[V]$ is its localization. Choose equations of
the finitely many prime divisors in $V\setminus T'$, and let $g$
be their product, with $g=1$ for the empty boundary. Then $T'$
is open in $V_g$ with complement of codimension at least two.
Normality gives $k[V_g]=k[T']$. Both opens are affine, so $V_g=T'$.
Two applications of Proposition~\ref{prop:principal-open} give
\[ \operatorname{gr}_f B_k(U)=B_k(V),\qquad \operatorname{gr}_g B_k(V)=B_k(T') , \]
where the subscript $k$ denotes extension of the augmented algebra
to $k$. The second graded algebra is a group algebra of a lattice,
hence a domain. The split order takes a finite value on every nonzero
element. Multiplying initial forms now shows that $B_k(V)$ is a
domain as well. Thus the split $f$-order on $B_k(U)$ is a valuation.

On $\sigma$ we have $\theta_p\theta_q=\theta_{p+q}$, so
$\tau_f(p+q)=\tau_f(p)+\tau_f(q)$ for $p,q\in\sigma\cap L$.
This monoid generates $L$ as a group. The additive function on it
therefore extends to a unique $a_f\in L^*$, and continuity and
homogeneity of $\tau_f$ give linearity on $\sigma$.
Lemma~\ref{lem:probe-graded-test} with $\eta=\sigma$ proves the
first formula. For integral $\lambda_j\ge0$, apply the formula to
$\prod_j f_j^{\lambda_j}$, which is invertible on the same torus.
Its value at $p$ is $\sum_j\lambda_j\tau_{f_j}(p)$.
Homogeneity gives rational coefficients. For each finite sum both
sides are continuous in $(\lambda_1,\ldots,\lambda_l)$, which gives
the assertion over $\mathbb R$.
\end{proof}

More generally, the proof of~\eqref{eq:single-function-valuation} requires only that
$U_f$ contain a volume-compatible torus and that $B_k(U_f)$
be a domain. In particular it applies, regardless of the Picard
group, whenever such a torus is a principal open of $U_f$.
For a common torus that is already principal open in $U$, this
observation applies simultaneously to every product of regular
functions invertible on that torus.

\begin{boundarycorollary}\label{cor:initial-form-coefficient}
Under the hypotheses of Theorem~\ref{thm:one-sided-probe}, suppose
$U_f$ is a split torus. Write a theta function in the marked chart as
$\theta_p=\sum_{a\in E}c_a z^a$, where $E=\Supp(\theta_p)$.
If $a_f$ attains its minimum on $E$ at a unique exponent $a_0$,
then $c_{a_0}=1$.
\end{boundarycorollary}
\begin{proof}
Proposition~\ref{prop:principal-open} identifies the associated graded
algebra with the group algebra of the cocharacter lattice of $U_f$,
and sends the initial form of each theta function to its labelled
character with coefficient one. The same normalization holds for
every character $z^a$ in the fraction field. Indeed, choose
$b\in L\cap\operatorname{int}\sigma$ and $n$ sufficiently large
that $nb,nb+a\in\sigma$. Then
$z^a=\theta_{nb+a}/\theta_{nb}$, so its initial form is a quotient
of two such characters. By~\eqref{eq:single-function-valuation},
the initial form of $\theta_p$ is
$c_{a_0}\operatorname{in}_f(z^{a_0})$.
Comparing the two coefficient-one characters gives $c_{a_0}=1$.
\end{proof}

Theorem~\ref{thm:one-sided-probe} also yields a finitely generated
cone of Gauss weights for each volume-compatible torus $T'\subset U$.
Indeed, every regular function
invertible on $T'$ is a unit of $k[U]$ times a product of the prime
equations of the divisorial components of $U\setminus T'$.
The units modulo scalars form a subgroup of the character lattice of
$T'$ and so are finitely generated. The cone is generated by the
slopes of the prime equations and by both signs of the slopes of a
basis of the units. The union of these cones over a choice of tori
consists of weights for which the basis is valuatively independent.
If the union is dense, continuity gives valuative independence at
every real Gauss weight.

\begin{corollary}\label{cor:source-seed-slopes}
Suppose $k[U]$ is a full upper LP algebra with Laurent coefficients,
and use the straight mirror chart of Theorem~\ref{thm:one-sided-probe}.
For every reachable source seed $(y_1,\ldots,y_r,u_1,\ldots,u_c)$,
the functions $\tau_{y_i}$ and $\tau_{u_j}$ are linear on all of
$\sigma$. Their slopes form an integral basis of $L^*$. Thus this
seed gives a full-dimensional cone of mirror Gauss weights on which
valuative independence holds, generated by the mutable slopes and by
both signs of the frozen slopes.
\end{corollary}
\begin{proof}
The variables of the seed are regular by the Laurent phenomenon.
Inverting the mutable variables gives its Laurent ring, so its torus
is an open subset of $U$. Each normalized exchange preserves the
logarithmic volume form up to sign. The functions in the assertion
are therefore units on the same volume-compatible torus, and
Theorem~\ref{thm:one-sided-probe} gives the linearity and the
minimum formula. The map taking source Gauss weights to the values
of the new coordinates is a composition of the tropicalizations of
the normalized exchanges, by Lemma~\ref{lem:boundary-gauss-change}.
Each exchange sends the changed entry $p_i$ to
$\tau_{\widehat F_i}(p)-p_i$ and fixes the others. Since
$\widehat F_i$ is independent of $x_i$, this is an integral
unimodular map on each linearity cone. One full-dimensional cone of
the composite meets $\operatorname{int}\sigma$ in an open set.
The slope vectors are the rows of this unimodular matrix.
\end{proof}

\begin{corollary}\label{cor:one-negative-probe}
Under the polynomial corner hypotheses of Theorem~\ref{thm:corner},
the prescribed functions $\Theta_m$ are valuatively independent at
every weight with at most one negative mutable coordinate, with
arbitrary frozen coordinates. If $w_i\le0$ and the remaining mutable
weights are nonnegative, then
\[ \nu^z_w(\Theta_m)=\langle w,m\rangle+(-w_i)\tau_{G_i}(m). \]
\end{corollary}
\begin{proof}
By Proposition~\ref{prop:upper-factorial}, $\Pic(U)=0$.
The initial torus gives the weights with nonnegative
mutable coordinates by taking $x_j$ and the units $u_a^{\pm1}$.
On the $i$-th adjacent seed torus the function $t_i=G_i/x_i$
is a unit. Indeed, the proof of Theorem~\ref{thm:corner} expresses it
as the actual mutant times a monomial in the unchanged coordinates.
It is regular on $U$. On the straight nonnegative orthant the
constant term of $G_i$ gives $a_{t_i}=-e_i$. Every $x_j$ with
$j\ne i$ and every $u_a^{\pm1}$ is invertible on the same torus,
with slopes $e_j$ and $\pm e_{r+a}$. The multi-function formula
in Theorem~\ref{thm:one-sided-probe}, with coefficient $-w_i$
for $t_i$, gives the assertion.
\end{proof}

\subsubsection{Uniqueness on a source seed cone}

On a source seed cone, every normalized reciprocal function is a
seed monomial. This also determines the function on a face shared by
several seed cones.

\begin{theorem}\label{thm:seed-windows}
Let $y_1,\ldots,y_r,u_1,\ldots,u_c$ be a reachable seed of the
source, and let $K\subset L^*_{\mathbb R}$ be the closed cone generated
by the mutable slopes of Corollary~\ref{cor:source-seed-slopes} and both
signs of its unit slopes. For every $a\in K\cap L^*$, the space of regular functions supported
in $P_a$ is one-dimensional. Its element with coefficient one at $x^a$ is
the normalized seed monomial of slope $a$, and its Newton polytope is
$P_a$. Thus seeds whose cones contain $a$ give the same normalized
function. Moreover, for any $b\in L^*$, the inclusion
$P_b\subseteq P_a$ implies $b=a$.
\end{theorem}

\begin{proof}
Write $d=r+c$, list all the seed coordinates as
$Y=(y_1,\ldots,y_r,u_1,\ldots,u_c)$, and let $g_j$ be the slope
of $\tau_{Y_j}$ on the straight cone $C=\sigma$. The map
$G:L\longrightarrow\mathbb Z^d$ with components
$p\mapsto\langle g_j,p\rangle$ is unimodular. The tropical change of seed
$T(p)=(\tau_{Y_1}(p),\ldots,\tau_{Y_d}(p))$ is a piecewise
unimodular bijection by
Lemma~\ref{lem:boundary-gauss-change} applied along the mutation word.
In particular $E=G^{-1}T$ maps $L$ onto $L$.

For fixed $p\in L$, Theorem~\ref{thm:one-sided-probe} gives
$h_a(p)=\langle a,E(p)\rangle$ for all real $a\in K$.
Indeed $a=G^t\lambda$ with $\lambda_j\geq0$ at the mutable
indices, and $\langle\lambda,T(p)\rangle$ is the value supplied
by that theorem. Since $K$ is full dimensional, a generic interior
weight exposes a single exponent of $\Theta_p$. Its value agrees
with $\langle a,E(p)\rangle$ on an open set, so that exponent is
$E(p)$. Consequently $E(p)$ still attains the minimum at every
boundary weight of $K$.

Now suppose $a\in K\cap L^*$ and $P_b\subseteq P_a$. Comparing
minimum support functions gives, for all $p\in L$,
\[ \langle a,E(p)\rangle=h_a(p)\leq h_b(p) \leq\langle b,E(p)\rangle . \]
Surjectivity of $E$ gives $b=a$. If a nonzero regular function
$f$ is supported in $P_a$, choose an interior weight of $C$
separating its exponents and let $b$ be its minimizing exponent.
Theorem~\ref{thm:reciprocal-polytope} gives
$P_b\subseteq\Newt(f)\subseteq P_a$. Hence $b=a$ and
$\Newt(f)=P_a$. In particular the coefficient at $x^a$ is
nonzero for every such $f$, which proves that the space has dimension
at most one.

Finally, $\lambda=(G^t)^{-1}a$ is integral and has nonnegative
mutable entries. The seed monomial $Y^\lambda$ is therefore
regular on $U$. Its minimum at $p$ is
$\langle\lambda,T(p)\rangle=h_a(p)$, so its Newton polytope is
$P_a$. Dividing by its coefficient at $x^a$ gives the asserted
normalization.
\end{proof}

\subsection{Explicit formulas at a polynomial corner}

For exchanges that are linear in the variables with negative exponents,
principal minors give candidates for the reciprocal functions. For
nonlinear exchanges, orders of vanishing can instead determine the
Newton polytope directly.

\subsubsection{Principal minors}

Let $A=k[U]$ be a smooth affine full upper LP algebra with a
polynomial corner. Theorem~\ref{thm:polynomial-corner-reciprocity}
identifies the reciprocal Newton polytopes at every label. For the
exchange polynomials below we can also give a reciprocal function
as a product of principal minors. Choose polynomial representatives
$G_i$ of the ordinary exchanges with $G_i(0)=1$, including the
designated units among the corner coordinates. Write $t_i=G_i/x_i$
for the corresponding regular functions.
Fix a set $S$ of mutable indices with an ordered partition
$S=S_1\sqcup\cdots\sqcup S_\ell$. Suppose that for $i\in S_h$
one has
\[ G_i=b_i+\sum_{j\in S_h}c_{ij}x_j,\qquad c_{ii}=0 , \]
where $b_i,c_{ij}$ are polynomials in the coordinates outside $S$
and the variables of the earlier parts $S_1,\ldots,S_{h-1}$.
There is no restriction on their degrees in those coordinates.
For $J\subseteq S_h$ define the regular principal minor
\[ d_J=\det\bigl(\delta_{ij}t_i-c_{ij}\bigr)_{i,j\in J}, \qquad d_\varnothing=1. \]
For graph LP algebras these are the determinants of Lam and
Pylyavskyy \cite[\S4.2, Proposition~4.5]{LPlinear}. Daisey
\cite[\S4]{Daisey} uses the same principal minors for general
linear exchanges. Here the minors give coordinates on tori on which
the valuations of the reciprocal functions can be computed.
Multiplying $d_J$ by $\prod_{i\in J}x_i$ gives a polynomial
with constant term one, since the resulting matrix has constant
part equal to the identity. Thus $d_J$ is pointed at
$-\mathbf1_J$ on the interior of the positive orthant, where
$\mathbf1_J=\sum_{i\in J}e_i$.

\begin{theorem}\label{thm:principal-minor-reciprocity}
Let $a\in M$, and let $S$ consist of the mutable indices at which
$a_i<0$. Suppose the exchanges indexed by $S$ have the ordered partition
above. In each part $S_h$, order its indices
$i_{h1},\ldots,i_{hm_h}$ by decreasing
$\lambda_{hj}=-a_{i_{hj}}$, where $m_h=|S_h|$. Let $J_{hj}$ consist of the first
$j$ indices, with $J_{h0}=\varnothing$, and set $\lambda_{h,m_h+1}=0$. Define
\[ \eta_{hj}=d_{J_{hj}}/d_{J_{h,j-1}},\qquad R=k[\eta_{hj}^{\pm1},x_i^{\pm1}\ (i\notin S),u^{\pm1}]. \]
Here and below the outside $x_i$ range only over mutable indices.
Let $s$ be the product of those outside variables and all the
nonempty prefix minors. Then $R$ is a Laurent polynomial ring.
For each $i\in S$, Cramer's rule gives $X_i\in R$
for the original coordinate $x_i$. Then $A[s^{-1}]=R$ if and
only if the $X_i$ have no pair with a common irreducible factor
in $R$. When these equivalent conditions hold, a normalized
reciprocal function at $a$ is given by
\[ f_a=\prod_{k=1}^c u_k^{a_{r+k}}\prod_{i\notin S}x_i^{a_i}\prod_{h=1}^\ell\prod_{j=1}^{m_h} d_{J_{hj}}^{\lambda_{hj}-\lambda_{h,j+1}}. \]
It satisfies $\nu^x_p(f_a)=\nu^z_a(\Theta_p)$ for every $p\in L$.
\end{theorem}

\begin{proof}
The case $S=\varnothing$ follows from the initial torus. Otherwise
write $D_h=\operatorname{diag}(t_i)_{i\in S_h}-(c_{ij})_{i,j\in S_h}$.
In part $h$ the exchange equations become $D_hx_{S_h}=b_{S_h}$.
All $\eta_{hj}$ and their inverses belong to $A[s^{-1}]$.
Proceed in the prescribed order of the parts. Once the earlier
variables have been recovered, all coefficients of the system for
part $h$ are known Laurent polynomials. Expanding the prefix minor
at its last diagonal entry gives
$d_{J_{hj}}=t_{i_{hj}}d_{J_{h,j-1}}+H_{hj}$, where $H_{hj}$
depends only on already recovered variables and earlier diagonal
entries in this part. Subtracting $H_{hj}/d_{J_{h,j-1}}$ from
$\eta_{hj}$ gives $t_{i_{hj}}$. After treating all entries of
the part, Cramer's rule with denominator $d_{S_h}$ recovers
$x_{S_h}$ as Laurent polynomials $X_i$ in the indicated
coordinates. Their calculation uses only the ordinary exchanges and
the prefix order. The coordinates generate the entire rational
function field because they recover all $x_i,u_j$. Their number
equals its transcendence degree, so they are algebraically independent.
Thus $R\subseteq A[s^{-1}]$ is the asserted Laurent ring, and
all factors of $s$ are units in $R$, since
$d_{J_{hj}}=\prod_{l=1}^j\eta_{hl}$.

If $R=A[s^{-1}]$, the distinct nonassociate prime elements
$x_i$ of Proposition~\ref{prop:upper-factorial} remain distinct
prime elements or become units after localization. This proves the
gcd condition. Conversely consider a height-one prime of $R$.
When no $X_i$ vanishes there, every $f\in A$ belongs to its
local ring by the initial Laurent expression. If just $X_i$
vanishes there, the unchanged coordinates are units there, and
$t_i=G_i/x_i$ belongs to $R$ by the same recursion.
Proposition~\ref{prop:upper-factorial} shows that $f$ belongs to the rank-one
ring generated by $x_i,t_i$ over the unchanged Laurent ring.
Thus $f$ is again regular there. The gcd condition exhausts all
height-one primes, so normality gives $f\in R$. Since $s$ is
a unit in $R$, this proves $A[s^{-1}]=R$.

For fixed outside coordinates the Jacobian of $t_S$ in $x_S$
is lower triangular by parts. Its $h$-th diagonal matrix is
$-\operatorname{diag}(x_i^{-1})_{i\in S_h}\,D_h$.
The change from $\eta$ to $t_S$ is triangular also within each
part, with diagonal entries one, and
$\prod_{h,j}\eta_{hj}=\prod_h d_{S_h}$. Wedging with the
logarithmic differentials of the outside coordinates now identifies
the logarithmic volume forms of the two tori up to sign. All the
$d_{J_{hj}}$, the outside variables, and both signs of the
designated unit characters are invertible on the constructed torus.
Their product $f_a$ is regular and has coefficient one at its
least exponent $a$ by the determinant calculation. Apply
Theorem~\ref{thm:one-sided-probe} to $f_a$ on this torus. Its
slope on the straight orthant is $a$, giving the asserted Gauss
values.
\end{proof}

An alternative sufficient condition for the gcd test is generation
after inverting the outside mutable variables by the initial
coordinates, their designated inverses, the outside inverses, and
the ordinary partners $t_i$ for $i\in S$. Indeed, the recursion
places all these generators in $R$, so $A[s^{-1}]=R$.
Regard the outside mutable coordinates as units. The first-star
intersection for $S$ equals the localization of $A$ at those
coordinates: localization commutes with a finite intersection, and
an outside neighboring ring after localization contains the initial
Laurent ring. For the remaining directions the normalizers can
change by monomials in unchanged coordinates, which do not change
the adjacent Laurent rings. If for each $J\subseteq S$ of size
at least two the restricted ordinary exchanges in its stratum have
unit gcd, Theorem~\ref{thm:local-exact-exchange-strata} therefore
gives the required generation. Nonvanishing of the stratum ideals
follows from the constant terms of the corner representatives.

\begin{corollary}\label{cor:component-linear-reciprocity}
Draw an edge $j\to i$ if $G_i$ depends on $x_j$. The
strongly connected components are the maximal sets of indices in
which each index is reachable from every other along directed edges. Suppose each
$G_i$ is affine-linear jointly in the mutable variables of the
strongly connected component containing $i$. The principal-minor
formula gives a normalized reciprocal function at every label if and
only if the initial mutable coordinates and their ordinary partners
generate the full upper algebra over the designated Laurent ring.
In that case the functions form a valuatively independent reciprocal
basis. The dependence of an exchange on variables of
other components is arbitrary.
\end{corollary}
\begin{proof}
Order the strongly connected components topologically, so every
variable outside the component of an exchange on which it depends
lies in an earlier component. For each label intersect this partition
with the negative mutable indices, discarding empty parts. Variables
outside the negative set are held fixed throughout. After inverting
them, their partners can be eliminated from the ordinary presentation
using the exchange relations. The sufficient condition just proved therefore applies to every
label. Corollary~\ref{cor:reciprocal-from-windows} gives the basis
and its valuative independence. Conversely, all coefficients
$c_{ij}$ are polynomials in the initial coordinates. Every
principal minor, and therefore each of these products, belongs to the
ordinary exchange subalgebra. If those products realize every label,
the same corollary makes them span $A$, forcing ordinary
generation. In this polynomial corner each restricted stratum
ideal contains a polynomial with nonzero constant term, so the
ordinary presentation injects into the function field by
Theorem~\ref{thm:local-exact-exchange-strata}.
\end{proof}

\subsubsection{Newton polytopes for nonlinear exchanges}

Work in the polynomial corner of Theorem~\ref{thm:corner}, and scale
the ordinary exchanges so that $G_i(0)=1$. For $a\in M$,
let $S=\{i:a_i<0,\ i\text{ mutable}\}$ and
$\lambda_i=-a_i$ for $i\in S$. Assume $S\ne\varnothing$.
For distinct $i,j\in S$ define
\[ d_{ij}=\operatorname{ord}_{x_i}(G_j-1),\qquad K_a=a+\operatorname{conv}\!\bigcup_{i\in S} \lambda_i\Newt(G_i). \]
Here $d_{ij}$ is the least $x_i$-exponent among the terms of
$G_j-1$. Each $G_j$ is nonconstant and has constant term one,
so its Newton polytope contains the origin.

\begin{theorem}\label{thm:nonlinear-simplex-reciprocity}
Under the hypotheses of Theorem~\ref{thm:reciprocal-polytope},
suppose that for every pair of distinct indices in $S$ the orders
$d_{ij},d_{ji}$ are positive and
\begin{equation}\label{eq:pair-pole-bound}
\frac{\lambda_i}{d_{ij}\lambda_j} +\frac{\lambda_j}{d_{ji}\lambda_i}\le 1.
\end{equation}
Then $P_a=Q_a=K_a$. In particular the forward valuations in all
directions are
\begin{equation}\label{eq:nonlinear-theta-valuation}
\nu^z_a(\Theta_p)=\langle a,p\rangle+ \min_{i\in S}\lambda_i\tau_{G_i}(p).
\end{equation}
The same conclusions hold under the alternative hypothesis
\begin{equation}\label{eq:pair-pole-congruence}
G_i^{\lambda_i}\equiv G_j^{\lambda_j} \pmod{(x_i^{\lambda_i},x_j^{\lambda_j})} \qquad(i,j\in S),
\end{equation}
and on every positive integral multiple of any label satisfying it.
In particular the congruences hold if $d_{ij}\ge\lambda_i$ for
all distinct $i,j\in S$.
No degree bound is imposed on the terms beyond the indicated orders.
\end{theorem}

\begin{proof}
The axis coefficient at shift $D_i(p)e_i$ is one by
Proposition~\ref{prop:axis}, where $D_i(p)=-\tau_{G_i}(p)$. Thus
$h_a(p)\le\langle a,p\rangle+\lambda_i\tau_{G_i}(p)$ for every
$p$ and $i\in S$. Taking minimum support functions gives
$K_a\subseteq P_a$. It suffices to construct a regular function
with coefficient one at $x^a$ and all its exponents in $K_a$,
since $P_a\subseteq Q_a$ and every strictly positive weight has
its unique minimum on $K_a$ at $a$.

Take indeterminates $U_i$, $i\in S$, and the finite set
$\Gamma=\{k\in\mathbb N^S:\textstyle\sum_i k_i/\lambda_i\le1\}$.
We construct $H=\sum_{k\in\Gamma} h_k U^k$ with $h_0=1$.
For fixed exponents $k_j$, $j\ne i$, occurring in $\Gamma$,
the possible values of $k_i$ are the integers from zero to
$N_i=\left\lfloor\lambda_i\left(1- \sum_{j\ne i}k_j/\lambda_j\right)\right\rfloor$.
Require the corresponding polynomial in $U_i$ to be divisible by
$(1+U_i)^{r_i}$, with
$r_i=[\lambda_i-\sum_{j\ne i}d_{ij}k_j]_+$, where $[t]_+=\max(t,0)$.
Write $c_{ij}=d_{ij}-\lambda_i/\lambda_j$. The pair hypothesis gives
$c_{ij}>0$ and $c_{ij}c_{ji}\ge1$. In particular $r_i\le N_i$.
When $r_i>0$, the last $r_i$ coefficients are uniquely determined
by the coefficients with $k_i\le N_i-r_i$. Indeed their consecutive
powers of $U_i$ form a basis modulo $(1+U_i)^{r_i}$, since $U_i$
is invertible in the quotient. An exponent in this last portion satisfies
\[ k_i>N_i-r_i=\left\lfloor\sum_{j\ne i} c_{ij}k_j\right\rfloor, \qquad\text{hence}\quad k_i>\sum_{j\ne i}c_{ij}k_j. \]
It cannot also be in the last portion for direction $j$, since
the two strict inequalities would imply $k_i>c_{ij}c_{ji}k_i$.
The origin lies in none of these portions. Set the unprescribed
coefficients to zero except for $h_0=1$, and compute each remaining
coefficient in increasing total degree. Its determining coefficients
all have smaller total degree. This satisfies every divisibility
condition simultaneously.

Substitute $u_i'=G_i-1$ and let $f=x^aH((u_i')_{i\in S})$. In direction
$i$, the terms with the fixed outside powers $k_j$, $j\ne i$,
contain a factor $x_i^sG_i^{r_i}$ before multiplication by $x^a$,
where $s=\sum_{j\ne i}d_{ij}k_j$. Since $s+r_i\ge\lambda_i$,
the only expression with an $x_i$-denominator needed for this part
of $f$ is a power of the regular ordinary partner $G_i/x_i$.
Consequently $f$ belongs to the $i$-th adjacent Laurent ring.
At directions outside $S$ it already has nonnegative mutable
exponents, so the first-star intersection gives $f\in k[U]$.
Every $u_i'$ has zero constant term, which gives coefficient one
at $x^a$. For $k\in\Gamma$ we have
\[ \sum_i k_i\Newt(G_i)\ \subseteq\ \operatorname{conv}\!\bigcup_i\lambda_i\Newt(G_i), \]
by taking convex coefficients $k_i/\lambda_i$ and assigning their
unused sum to zero. Thus $\Newt(f)\subseteq K_a$, as required.

Under the alternative hypothesis form a polynomial $H_0$ as follows.
For any nonnegative exponent $b$ with $b_i<\lambda_i$ for some
$i\in S$, take its coefficient from $G_i^{\lambda_i}$. Give all
other exponents coefficient zero. Two prescriptions agree since the
congruence identifies exactly the coefficients in the common range
$b_i<\lambda_i, b_j<\lambda_j$. The support of $H_0$ is contained
in the finite union of the supports of the prescribed powers, and
$H_0\equiv G_i^{\lambda_i}\bmod x_i^{\lambda_i}$ for every $i$.
Thus $x^aH_0$ is regular and has the required support and least
coefficient. If $d_{ij}\ge\lambda_i$ for every ordered pair,
both sides of every pair congruence equal one in the quotient.
Taking powers gives the assertion for multiples by homogeneity of $P_a$.
\end{proof}

The inequalities~\eqref{eq:pair-pole-bound} hold throughout $\max_i\lambda_i\le
2\min_i\lambda_i$ whenever all $d_{ij}\ge3$. This gives a
full-dimensional cone for smooth sources of unbounded internal degree
in every rank. Take a distinct designated unit $u_j$ for each
exchange and arbitrary nonzero
polynomials $B_j(x_{\ne j})$. The exchanges
$G_j=1+u_j B_j\prod_{i\ne j}x_i^3$ are irreducible, all double-zero
strata are empty, and $\partial G_j/\partial u_j$ is nonzero on
$G_j=0$. Thus Theorem~\ref{thm:local-exact-exchange-strata} gives the
full upper presentation and the one-exchange hypersurface opens give
smoothness.

\begin{example}\label{ex:cyclic-unequal-poles}
Consider a cyclic seed with two mutable coordinates $x,y$ and
three invertible coordinates $u,v,w$. Take
$G_x=1+(u+v)y^3$ and $G_y=1+wx^3$. Both exchanges are irreducible
by linearity in $u$ and $w$, respectively, and primitivity over the
remaining coordinates. The locus $x=y=0$ is empty in the ordinary
exchange presentation. Theorem~\ref{thm:local-exact-exchange-strata}
therefore identifies it with the full upper algebra. The two one-exchange
opens are smooth since along $G_x=0$ we have
$\partial G_x/\partial u=y^3\ne0$, and along $G_y=0$ we have
$\partial G_y/\partial w=x^3\ne0$.

Take $a=(-1,-m,0,0,0)$ with $m\ge4$. The inequality~\eqref{eq:pair-pole-bound} fails because
$m/3+1/(3m)>1$. The congruence~\eqref{eq:pair-pole-congruence} also fails, since
$G_x-G_y^m\equiv (u+v)y^3\bmod (x,y^m)$.
Mutating first at $y$, then at $x$, gives $y'=G_y/y$ and
$x''=((y')^3+u+v)/x$. The seed monomial $x''(y')^{m-3}$ has least
exponent $a$, so Theorem~\ref{thm:seed-windows} identifies it as the
normalized reciprocal function. Its Laurent expansion is
\[ x''(y')^{m-3}=\frac{(1+wx^3)^m+(u+v)y^3(1+wx^3)^{m-3}}{xy^m}. \]
Let $q=3e_x+e_w$ and $J=[3e_y+e_u,3e_y+e_v]$.
The two numerator terms have Newton polytopes $[0,mq]$ and
$J+[0,(m-3)q]$, respectively. Their supports have different
$y$-degrees. Hence $Q_a=P_a$ is $a$ plus the convex hull of their
union. At $p=(0,0,-m,-m,-1)$, its minimum is $3-2m$, whereas the
minimum on $K_a=a+\operatorname{conv}([0,mq]\cup J)$ is $-m$.
Thus the formula using only the two dilated exchange polytopes misses
the minimum by $m-3$, even though reciprocity holds.
\end{example}

The valuation formula~\eqref{eq:nonlinear-theta-valuation} determines a projection of the forward Newton
polytope. Fix $p$, choose a set $S$ of mutable indices with
$D_i(p)>0$ for $i\in S$, and take
$\lambda_i$ proportional to $1/D_i(p)$, clearing denominators.
If the pair inequalities hold, the projection of the Newton polytope
of the polynomial factor $A_p$ to these coordinates is exactly
$\operatorname{conv}\bigl(0,\ D_i(p)e_i:\ i\in S\bigr)$.
Indeed take $a_i=-\lambda_i$ on $S$ and zero elsewhere.
The valuation formula gives $\sum_{i\in S}e_i'/D_i(p)\le1$ for
every correction exponent $e'$, and both the origin and each $D_i(p)e_i$ occur in the projection.

\subsection{Realization criteria in a marked Laurent chart}

Work over an algebraically closed field $k$ of characteristic zero
under the hypotheses of Theorem~\ref{thm:reciprocal-polytope}.
Write $M=L^*$. By Theorem~\ref{thm:reciprocal-polytope},
a label $a$ is realized precisely when $P_a=Q_a$.
At a polynomial corner every label is realized
(Theorem~\ref{thm:polynomial-corner-reciprocity}).

To test realization in a general marked chart, it suffices to
approximate $P_a$ by rescaling the supports of regular sections.
Fix a Euclidean norm on $M_{\mathbb R}$ with closed unit ball $D$
and give $L_{\mathbb R}$ the dual norm.

\begin{proposition}\label{prop:cartier-closure}
For $a\in M$, the following conditions are equivalent.
\begin{enumerate}
\item The label $a$ is realized.
\item For every $\epsilon>0$ there are an integer $n\ge1$,
      a label $b\in M$ with $\|b/n-a\|\le\epsilon$, and a
      regular Laurent polynomial pointed at $b$ whose support lies in
      $n(P_a+\epsilon D)$.
\end{enumerate}
In particular, $a$ is realized if and only if $na$ is realized
for any specified integer $n\ge1$. If $b_j$ are realized labels and
$b_j/n_j\to a$ for positive integers $n_j$, then $a$ is
realized.
In particular, every integral point in the closure of the union of the source seed
cones of Corollary~\ref{cor:source-seed-slopes} is realized.
\end{proposition}

\begin{proof}
For $K_a=Q_a$ and for $K_a=P_a$, the assignment $a\mapsto K_a$
is homogeneous under positive integral dilation and satisfies
$K_{a+b}\subseteq K_a+K_b$. For $Q_a$ this is
Theorem~\ref{thm:minimum-source-newton-polytope}. For $P_a$
it follows by taking minima on $\Newt(\Theta_p)$ for each $p$.
The finite collection $K_{\pm e_i}$ thus bounds $K_d$ in a
ball of radius $c\|d\|_1$. Using the inclusions for both
$a=b+(a-b)$ and $b=a+(b-a)$ gives the same Lipschitz bound
in Hausdorff distance. Homogeneity gives the bound also at rational
points.

If $H$ has least exponent $b$ at any separating
$\omega\in\operatorname{int}C$, its polytope contains $Q_b$.
Applying this containment to the approximating sections and taking
the Hausdorff limit gives $Q_a\subseteq P_a$. The opposite inclusion is
Theorem~\ref{thm:reciprocal-polytope}. Equality gives realization
by taking a representative on $Q_a$, whose powers give the
approximations in the other direction. The same continuity shows
that equality persists at all limits of rescaled realized labels.
Finally apply Theorem~\ref{thm:seed-windows} on each source seed
cone and approximate real points by rational ones.
\end{proof}

For the original polynomial coefficient ring the same limit
assertion holds if all the approximating functions are regular
there. Indeed their supports satisfy the closed homogeneous frozen
inequalities of Corollary~\ref{cor:polynomial-polyhedral-newton-sections}.
In the limit the entire $Q_a$ satisfies these inequalities, since
each approximate support contains its corresponding smallest
polytope in the localized algebra. Every representative on $Q_a$
therefore descends to the original ring.

For a single reciprocal function, descent to the polynomial
coefficient ring can be tested directly on its Newton polytope.

\begin{corollary}\label{cor:polynomial-star-reciprocal}
Keep the polynomial seed, its marked torus, and the original first-star
weights $v_{t,i}$ of Theorem~\ref{thm:polynomial-star-descent}.
For each frozen variable $f_i$, let $J_i$ consist of the mutable
indices $j$ for which $F_j$ is an associate of $f_i$.
Suppose a marked forward counting construction on
$B=\mathcal U_S[F^{-1}]$ gives the reciprocal polytope $P_a$ for
one label $a$, in this same source character lattice. Assume that
every normalized solution in $B$ supported in $P_a$ has Newton
polytope $P_a$, as in Theorem~\ref{thm:reciprocal-polytope}.
Then there is a normalized solution in $\mathcal U_S$ if and only if
there is one in $B$ and
\[ \tau_{P_a}(e_i)\ge0,\qquad \tau_{P_a}(e_i+e_{x_j})\ge0\quad(j\in J_i),\qquad 1\le i\le s. \]
When these inequalities hold, every normalized solution in $B$
supported in $P_a$ belongs to $\mathcal U_S$.
\end{corollary}

\begin{proof}
For every nonzero $g\in B$, the first-star calculation above gives
$g\in\mathcal U_S$ exactly when $v_{t,i}(g)\ge0$ for every
$t,i$. For any normalized solution supported in $P_a$ these
values depend only on $P_a$ by the hypothesis on its Newton
polytope. This proves both implications. The first-star calculation
uses the localizations of the original Laurent rings, so it also
applies to exchanges that become units after inverting $F$.
\end{proof}

In the marked counting chart this test treats a whole coefficient orbit.
Let $\epsilon_i$ denote the exponent of the now invertible
character $f_i$, and write $a(n)=a+\sum_i n_i\epsilon_i$ for
$n\in\mathbb Z^s$. The exponent corrections in each forward
$\Theta_p$ have zero component at every unit by
Proposition~\ref{prop:local-bend-cone}. Hence
\[ \nu_{a(n)}(\Theta_p)=\nu_a(\Theta_p) +\Big\langle\sum_i n_i\epsilon_i,p\Big\rangle, \qquad P_{a(n)}=P_a+\sum_i n_i\epsilon_i . \]
Multiplication by $\prod_i f_i^{n_i}$ identifies the normalized
solution spaces in the localized ring. If that space is nonempty,
write $b_i=\min_t\tau_{P_a}(v_{t,i})$, using the initial weights for
these orders. The labels of this orbit realized over the original
polynomial ring are exactly $a(n)$ with $n_i\ge-b_i$ for all
$i$, since translation adds $n_i$ to every $i$-th boundary
order. Thus the least polynomial translate is computed from $P_a$.

\subsubsection{Finite tests on a rational cone}
\label{subsec:finite-reciprocal-cones}

Work in the marked Laurent setting of
Theorem~\ref{thm:reciprocal-polytope}. Fix a lattice-separating
weight $\omega$ in the interior of the straight cone. For
$a\in M=L^*$, the smallest source polytope at $\omega$-least
exponent $a$ is $Q_a$, and the polytope supplied by the forward
counts is $P_a\subseteq Q_a$. Thus $a$ is realized precisely when
$P_a=Q_a$. Both families extend continuously and positively
homogeneously to $M_{\mathbb R}$ by
Proposition~\ref{prop:cartier-closure} and its proof. We use the equality
of their real extensions when discussing loci in $M_{\mathbb R}$.
All the tests below involve only lattice polytopes at specified
integral arguments. For the algorithmic assertions, we assume that
arithmetic and equality are computable in a field of definition.

Let $K=\sum_{i=1}^m\mathbb R_{\ge0}v_i$ for $v_i\in M$, and let
$s=\sum_i v_i$. The cone may have lineality or smaller dimension.
We write $\operatorname{relint}K$ for the interior relative to the
linear span of $K$. For either family $S=P,Q$, say that $S$ is \emph{Minkowski additive} on
$K$ if $S_{a+b}=S_a+S_b$ for all $a,b\in K$.

\begin{proposition}\label{prop:single-sum-realization-cone}
For $S=P$ or $S=Q$, Minkowski additivity on $K$ is equivalent to the
single finite test
\[ S_s=\sum_i S_{v_i}. \]
Under this condition $S_{\sum_i t_i v_i}=\sum_i t_iS_{v_i}$ for
all $t_i\ge0$. The following two conclusions apply separately.
\begin{enumerate}
\item If $Q$ is Minkowski additive on $K$, then
$\{a\in K:P_a=Q_a\}$ is a union of closed faces of $K$. On each
face $F$ it suffices to test one integral point of
$\operatorname{relint}F$ to decide whether all of $F$ is realized.
In particular realization at one integral point of
$\operatorname{relint}K$ is equivalent to realization at every
integral point of $K$.
\item If $P$ is Minkowski additive on $K$, then
$\{a\in K:P_a=Q_a\}$ is one computable rational polyhedral cone.
It equals $K$ precisely when $P_{v_i}=Q_{v_i}$ for every $i$.
If both families are Minkowski additive, the equality locus in $K$
is a single face of $K$.
\end{enumerate}
\end{proposition}
\begin{proof}
Both families satisfy $S_{a+b}\subseteq S_a+S_b$. For $Q$ this is
Theorem~\ref{thm:minimum-source-newton-polytope}, and for $P$ it
follows by minimizing on each $\Newt(\Theta_p)$ and using
Theorem~\ref{thm:reciprocal-polytope}. Fix any $p\in L$ and write
$k(a)=\min_{b\in S_a}\langle b,p\rangle$. This is homogeneous and
superadditive. Suppose the single test holds. For
$a=\sum_i t_i v_i$ choose an integer $N\ge\max_i t_i$. Then
\[ \sum_i t_i k(v_i)\le k(a) \le N k(s)-\sum_i(N-t_i)k(v_i)=\sum_i t_i k(v_i). \]
The middle inequality uses $Ns=a+\sum_i(N-t_i)v_i$.
Testing all $p$ gives the asserted Minkowski formula. The converse
follows by addition. In particular each such $k$ is the restriction
of a linear form to $K$.

Suppose $Q$ is additive. For fixed $p$ the difference between the
minimum on $P_a$ and the minimum on $Q_a$ is a nonnegative concave
function of $a\in K$. If it vanishes at $b\in\operatorname{relint}F$,
it vanishes throughout $F$. Indeed for any $a\in F$, a sufficiently
small $\epsilon>0$ has $b+\epsilon(b-a)\in F$. Express $b$ as a
strict convex combination of $a$ and that point, and use concavity
and nonnegativity. Apply this to every $p\in L$. Each face admits an
integral relative-interior point, for example the sum of the given
generators belonging to the face, with empty sum zero. This proves
the first assertion and the finite face test.

Now suppose $P$ is additive. Choose a finite complete rational fan
in $L_{\mathbb R}$ refining the normal fans of the $P_{v_i}$. One
can take the normal fan of their Minkowski sum with any
full-dimensional lattice polytope. Let $p_1,\ldots,p_h$ be the
primitive rays. Every $P_a$ for $a\in K$ is cut out by its support
inequalities at these rays, since its minimum support function is
linear in $p$ on every cone of the fan. This applies also to
lower-dimensional $P_a$. Set
$\ell_j(a)=\min_{b\in P_a}\langle b,p_j\rangle$. On $K$ this is
a rational linear form, computable from its values at the $v_i$.
By Proposition~\ref{prop:source-polar-polytope}, the equality locus
is exactly the following finite intersection:
\[ \left\{a\in K: \langle a,e\rangle\ge\ell_j(a)\quad (1\le j\le h,\ e\in\operatorname{Vert}(T_{p_j}))\right\}. \]
Indeed the inequalities say that $Q_a\subseteq P_a$, and the other
inclusion always holds. The condition on generators follows by
convexity. If both families are additive, the same locus is convex
and is a union of faces by the first assertion. A nonempty convex
union of faces of a polyhedral cone is a face: a relative-interior
point of the union belongs to a face in the union, and the face
property implies that the whole union lies in that face. Nonemptiness follows
from $P_0=Q_0=\{0\}$.
\end{proof}

\begin{corollary}\label{cor:one-product-realization-cone}
Let $f_i\in k[U]$ have $\omega$-least exponents $v_i$ with
coefficient one. Choose positive integers $n_i$ and set
$a_0=\sum_i n_i v_i$. The equalities
\[ \Newt\!\left(\prod_i f_i^{m_i}\right)=P_{\sum_i m_i v_i} \qquad(m_i\in\mathbb N) \]
are equivalent to the one equality
$\sum_i n_i\Newt(f_i)=P_{a_0}$.
If it holds, both families are Minkowski additive on $K$, and every
lattice label in $K$ is realized. Each such product gives a
normalized solution at its monoid label. Labels outside the monoid
obtain their solutions from the smallest source polytopes.
\end{corollary}
\begin{proof}
Always
\[ P_{a_0}\subseteq\sum_i n_iP_{v_i} \subseteq\sum_i n_iQ_{v_i}\subseteq\sum_i n_i\Newt(f_i), \qquad P_{a_0}\subseteq Q_{a_0}\subseteq\sum_i n_iQ_{v_i}. \]
The proposed equality forces equality at every inclusion. At each
weight $p$, the minimum on each $P_{v_i}$ is at least the minimum
on $Q_{v_i}$, in turn at least the minimum on $\Newt(f_i)$.
The weighted sums agree, and all $n_i>0$, so these three minima
agree individually. Thus $P_{v_i}=Q_{v_i}=\Newt(f_i)$.
Apply Proposition~\ref{prop:single-sum-realization-cone} with generators
$n_i v_i$, using homogeneity. Both families are additive and agree
on $K$. Products have the asserted Newton polytopes by Minkowski
addition and have monic least term by multiplication. At every
integral $a\in K$, equality $P_a=Q_a$ supplies a regular function
on $P_a$ even if $a$ is outside the monoid of the $v_i$. Conversely
the all-products formula includes $m_i=n_i$.
\end{proof}

For nonlinear exchanges the corollary determines whether every product
of the prescribed functions realizes the reciprocal polytope at its least exponent. The
functions, their least exponents, and the finitely many polytopes
in the test suffice. Other representatives at the same labels are
governed by Theorem~\ref{thm:reciprocal-polytope}.

\begin{corollary}\label{cor:finite-mixed-realization-cover}
Suppose a rational polyhedral cone $K$ has a finite cover by rational
polyhedral cones on each of which at least one of $P,Q$ is Minkowski
additive. Then the realized locus in $K$ is a computable finite union
of rational polyhedral cones, even if the cover is initially unknown.
In particular this decides uniform realization when the hypothesis
holds on $M_{\mathbb R}$.
\end{corollary}
\begin{proof}
Enumerate finite lists of rational polyhedral subcones of $K$.
Whether a list covers $K$ is a finite rational polyhedral computation,
for example by subdivision along the defining hyperplanes, including
all relative faces. On each cone apply the single-sum test to both
families. Under the hypothesis the search finds a cover on which
every cone passes at least one test. Proposition~\ref{prop:single-sum-realization-cone}
then computes the locus on each cone. The calculation requires exact arithmetic in the coefficient field.
\end{proof}

Several products can also establish realization at an individual
label without an additivity assumption on a cone.
\begin{proposition}\label{prop:product-window-integer-hull}
Keep the regular functions $f_i$ above, and set $R_i=\Newt(f_i)$.
For $a\in K\cap M$ form the rational polyhedron
$\Lambda_a=\{\lambda\in\mathbb R_{\ge0}^m:\sum_i\lambda_iv_i=a\}$
and the rational polytope
\[ H_a=\bigcap_{\lambda\in\operatorname{Vert}(\Lambda_a)} \sum_i\lambda_iR_i . \]
Then $Q_a\subseteq\operatorname{conv}(H_a\cap M)$. If
$H_a\cap M\subseteq P_a$, the label $a$ is realized, in every
positive integral multiple as well. The same bound $H_a$ results if
the intersection uses all nonnegative real representations of $a$.
\end{proposition}
\begin{proof}
The nonempty polyhedron $\Lambda_a$ lies in the nonnegative orthant
and has no lines, so it has finitely many vertices and is the sum
of their convex hull and its recession cone. At every vertex the
coefficients $\lambda_i$ are rational. Clear denominators and
use homogeneity and subadditivity of $Q$ to obtain
$Q_a\subseteq\sum_i\lambda_i Q_{v_i}\subseteq\sum_i\lambda_i R_i$.
Its vertices are lattice points, proving the first assertion.
The proposed containment then gives $Q_a=P_a$, with equality
preserved under dilation. For a recession vector $r$ of $\Lambda_a$,
the polytope $\sum_i r_i R_i$ contains $\sum_i r_i v_i=0$.
A point in all the vertex bounds is therefore in every convex
combination of those bounds and remains there upon addition of a
recession contribution. This proves the last assertion.
\end{proof}

\subsubsection{A numerical criterion using divisors}
\label{subsec:exchange-divisor-realization}

We continue under the hypotheses of Proposition~\ref{prop:cartier-closure}
over the Laurent coefficient ring. Write $M=L^*$ for the character
lattice of the marked torus. Fix $a\in M$ and set $\Delta=P_a$.
If $M$ has rank zero, its only label is realized by $1$.
Assume $M$ has positive rank and choose a primitive $q\in L$
in the interior of the straight cone $C$.
Then $a$ is the unique minimizer of $q$ on $\Delta$.

Take a smooth projective toric compactification $Z$ of the marked
torus whose fan refines the normal fan of $\Delta$ and contains the
rays of $q$ and the mutable coordinate weights $e_i$.
One may first add a full-dimensional lattice polytope to $\Delta$
and refine the normal fan of the sum. Toric subdivisions then insert
the desired rays and give a smooth projective refinement. This works
also when $\Delta$ has smaller dimension. Write
\[ H_\Delta=-\sum_\rho \min_{b\in\Delta}\langle \rho,b\rangle\,T_\rho , \qquad \lambda_i=\min_{b\in\Delta}\langle e_i,b\rangle , \]
where $\rho$ denotes the primitive generator of a ray and $T_\rho$
the corresponding prime divisor. Sections of $nH_\Delta$ restrict
to precisely the Laurent polynomials supported in $n\Delta$.

In the open orbit of $T_{e_i}$, the ordinary exchange $F_i$
defines an irreducible hypersurface. Blow up its closure and take the
exceptional prime dominating that hypersurface. Denote its valuation
by $v_i$. Choose a common smooth projective model
$\pi:X\longrightarrow Z$, isomorphic over the torus, on which these
valuations occur as prime divisors $E_i$. For example, blow up the
closures successively and resolve, preserving their smooth generic
points. Distinct centers have distinct generic points on the toric
boundary. Let $T$ be the strict transform of $T_q$ and define
\[ D_\Delta=\pi^*H_\Delta-\sum_i\max(0,-\lambda_i)E_i. \]

A divisor is \emph{pseudoeffective} if its numerical class lies in the
closure of the cone generated by effective divisors. It is \emph{nef}
if it has nonnegative degree on every integral curve.
For a pseudoeffective divisor $D$, write $\sigma_T(D)$ for the
asymptotic order along $T$ after an arbitrarily small ample
perturbation. For an integral ample $H$ it is the limit, as positive
rational $\epsilon$ tends to zero, of the infimum of
$\operatorname{mult}_T(B)/n$ for effective divisors $B$ linearly
equivalent to $n(D+\epsilon H)$ with $n\epsilon$ integral. Here
$\operatorname{mult}_T(B)$ is the coefficient of $T$ in $B$. This number
depends only on the numerical class of $D$ and is independent of
$H$ \cite[III, \S1]{Nakayama}.

\begin{theorem}\label{thm:numerical-realization}
The label $a$ is realized if and only if $D_\Delta$
is pseudoeffective and $\sigma_T(D_\Delta)=0$.
In particular, $a$ is realized if $D_\Delta$ is nef.
\end{theorem}

\begin{proof}
We first identify the sections of the divisor. At the generic exchange
center the parameters normal to the center are $x_i,F_i$. On the
exceptional divisor $F_i/x_i$ is a transcendental residue coordinate
over the center. The normalizing factor $\widehat F_i/F_i$ is a
unit in the complementary Laurent ring. Thus $v_i$ is the order
at $F_i$ in the $i$-th adjacent chart of
Proposition~\ref{prop:upper-factorial}, and $v_i(x_j)=\delta_{ij}$.
The coefficient of $E_i$ in $\pi^*H_\Delta$ is $-\lambda_i$.
For a Laurent polynomial $f$ supported in $n\Delta$, each
monomial has $v_i$-value at least $n\lambda_i$.
Membership in $H^0(X,nD_\Delta)$ imposes on $f$ the additional
conditions $v_i(f)\ge0$ when $\lambda_i<0$. When
$\lambda_i\ge0$ these conditions already hold. In fact
$\operatorname{div}(f)+n\pi^*H_\Delta$ is effective, being the
pullback of an effective Cartier divisor on $Z$. Conversely a
section of $nD_\Delta$ is a section of $n\pi^*H_\Delta$, hence
of $nH_\Delta$ on $Z$. The first-star description now gives
\[ H^0(X,nD_\Delta)= \{f\in k[U]:\Supp(f)\subseteq n\Delta\}. \]
For a nonzero such section its order of vanishing along $T$ is
$\min_{b\in\Supp(f)}\langle q,b\rangle-n\langle q,a\rangle$.
If $a$ is realized, the corresponding section in degree one avoids
$T$ at its generic point. It gives an effective representative
of $D_\Delta$. Adding small ample perturbations avoiding $T$
gives $\sigma_T(D_\Delta)=0$.

Conversely, assume pseudoeffectivity and $\sigma_T(D_\Delta)=0$.
Represent an ample Cartier divisor $H$ by a divisor supported on
the complement of the torus, using $\Pic((\mathbb G_m)^{\dim Z})=0$.
There exist positive integers $n_j,m_j$ with $n_j\to\infty$,
$m_j/n_j\to0$, and sections $f_j$ of
$n_jD_\Delta+m_jH$ whose vanishing orders along $T$ are $o(n_j)$.
Indeed take successively smaller rational perturbations in the
definition and clear denominators, taking powers if necessary.
These sections are Laurent polynomials on the torus.
At each strict toric boundary divisor their exponents $b$ satisfy
\[ \langle \rho,b/n_j\rangle \ge \min_{e\in\Delta}\langle\rho,e\rangle-O(m_j/n_j). \]
By completeness of the fan, these finitely many inequalities have
recession cone zero. As $m_j/n_j\to0$, they confine the rescaled
supports to $\Delta+\epsilon_jD$ for some $\epsilon_j>0$ tending
to zero, where $D$ is the fixed unit ball.
At $T$ the vanishing order gives
$\min_{b\in\Supp(f_j)}\langle q,b/n_j\rangle\to\langle q,a\rangle$.
Since $q$ exposes $a$, a sequence of minimizing exponents divided
by $n_j$ converges to $a$.

At $E_i$ the section inequality gives $v_i(f_j)\ge-O(m_j)$,
since the coefficient of $E_i$ in $D_\Delta$ is nonpositive.
Multiplication by $\prod_i x_i^{c m_j}$ for a fixed sufficiently
large integer $c$ clears all these poles, translating the rescaled
supports by $o(1)$. The resulting Laurent polynomials $G_j$ belong to
$k[U]$. The same support estimate holds after enlarging
$\epsilon_j$ by a quantity tending to zero, and the translated
supports still contain points converging to $a$.
Fix a lattice-separating $\omega\in\operatorname{int}C$.
Since $\omega$ also exposes $a$, the
$\omega$-least exponents $b_j$ of $G_j$ satisfy $b_j/n_j\to a$.
Normalize the least coefficients to one.

Since $Q_{b_j}\subseteq\Newt(G_j)$, the continuity argument in
Proposition~\ref{prop:cartier-closure} gives $Q_a\subseteq\Delta$.
This proves realization. Finally a nef class becomes ample
under every positive ample perturbation. A sufficiently divisible
multiple has a section avoiding $T$, which proves the last assertion.
\end{proof}

The same numerical criterion tests whether a regular Laurent polynomial
supported in a lattice polytope $\Delta$ can have a nonzero
coefficient at a given vertex $a$. We require the open normal cone
at $a$ to meet $\operatorname{int}C$ and choose $q$ and $\omega$
in the intersection. The proof then gives $Q_a\subseteq\Delta$
under the stated numerical condition.

\begin{corollary}\label{cor:uniform-divisor-realization}
Let $\Omega\subseteq M_{\mathbb R}$ be a rational polyhedral cone.
Suppose a given finite list $E\subset L$ detects every strict
inclusion $P_a\subsetneq Q_a$ for $a\in\Omega\cap M$, in the
sense that their support minima differ at some $p\in E$.
Then the realized labels in $\Omega$ are the integral points of a
computable finite union of rational polyhedral cones. A list of
supporting normals defining all the $P_a$ has this detecting
property, as do the rays of any common finite complete pointed refinement
of their normal fans.
\end{corollary}
\begin{proof}
At $p\in E$ the two minima are the values at $a$ of the minimum
support functions of $\Newt(\Theta_p)$ and $T_p$, respectively,
by Theorem~\ref{thm:reciprocal-polytope} and
Proposition~\ref{prop:source-polar-polytope}. Intersect $\Omega$
with a common refinement of the normal fans of these finitely many
computable polytopes. Within each cone their equalities are homogeneous
rational linear equations and detect precisely the realized lattice
labels. For a defining list of supporting normals any larger polytope
violates at least one of the defining inequalities.
\end{proof}

On the projective model of Theorem~\ref{thm:exchange-cox-model}
the \emph{Cox ring}, the ring of sections graded by divisor class,
is finitely generated under the same negative-vertex
hypothesis. Corollary~\ref{cor:cox-reciprocal-face} gives the
corresponding test in its divisor classes.

\begin{samepage}
\begin{corollary}\label{cor:simultaneous-exchange-orders}
Let $v_i$ be the exchange valuations of Theorem~\ref{thm:numerical-realization}.
Fix a lattice polytope $\Delta$ and a vertex $a$ whose open normal
cone meets $\operatorname{int}C$. Let $r$ be the number of mutable
coordinates, and set $\lambda=(\lambda_i)_{i=1}^r$, where $\lambda_i=\min_{b\in\Delta} b_i$.
There is a rational polytope $J\subseteq\lambda+\mathbb R_{\ge0}^r$
with the following property. For every integer $n\ge1$, the lattice
points of $nJ$ are exactly the vectors $t=(t_i)\in\mathbb Z^r$
with $t_i\ge n\lambda_i$ for which some Laurent polynomial $f$
supported in $n\Delta$ has coefficient one at $x^{na}$ and
satisfies $v_i(f)\ge t_i$ for all $i$. In particular, integral
bounds $t\ge\lambda$ are attainable in degree one whenever $nt$
is attainable in some degree $n\ge1$.
\end{corollary}
\end{samepage}
\begin{proof}
Choose a separating weight in the intersection of the two open cones,
and extend $t$ by zeros in the frozen coordinates to obtain $\bar t\in M$.
Since $v_i(x^{\bar t})=t_i$, the first-star criterion says that the
order bounds hold exactly when $x^{-\bar t}f$ is regular. Its least
exponent is then $na-\bar t$. Conversely, we may normalize a regular
function with this least exponent and multiply it by $x^{\bar t}$.
Thus such an $f$ exists exactly when
$Q_{na-\bar t}\subseteq n\Delta-\bar t$.
Choose defining inequalities $\langle b,p\rangle\ge c$ for $\Delta$
with integral normals, using both signs of each equation of its affine hull.
For each inequality, Proposition~\ref{prop:source-polar-polytope}
expresses the corresponding support bound as
\[ \langle na-\bar t,e\rangle+\langle\bar t,p\rangle\ge nc \qquad(e\in\operatorname{Vert}(T_p)). \]
Together with $t_i\ge n\lambda_i$, these inequalities define the
$n$-th dilate of a rational polyhedron $J$.

To see that $J$ is bounded, consider degree one. In the finite-dimensional
space of Laurent polynomials supported in $\Delta$, the subspaces of
$v_i$-order at least $j$ decrease to zero as $j$ increases. Thus
each attainable $t_i$ is bounded above. The monomial $x^a$ shows
that $\lambda$ is an integral point of $J$. If $J$ were unbounded,
a nonzero integral recession direction would give an unbounded sequence
of integral points starting at $\lambda$.
\end{proof}

\section{The second mirror in the marked torus}
\label{sec:second-mirror}

To count curves on a smooth augmented forward mirror, we first show that
the logarithmic volume form on its marked torus extends across the
complement. We then compare the reverse theta functions with the source
functions in their common Laurent field.
Theorem~\ref{thm:first-reverse-partner-return} reduces equality of the
two algebras to equality at the first mutation partners.

\subsection{The volume form and smoothness of the forward mirror}
\label{subsec:reverse-projective-traces}

\begin{proposition}\label{prop:reverse-volume-automatic}
Let $U$ be a smooth affine log Calabi--Yau variety over an
algebraically closed field $k_0$ of characteristic zero. Suppose $U$
contains a dense split torus and $\Pic(U)=0$. For an algebraically
closed field $k$ of characteristic zero, set
$X=\Spec(k\otimes_{\mathbb Z}B_{\mathbb Z}(U))$.
Then $X$ admits a normal projective compactification $\widetilde X$
whose boundary $\widetilde X\setminus X$ is the support of a reduced
divisor $\Delta$, with $(\widetilde X,\Delta)$ log canonical and
$K_{\widetilde X}+\Delta\sim0$.
On every dense open split torus $T\subset X$,
the standard logarithmic volume form extends to a generator of the
dualizing sheaf of $X$, which on the smooth locus is the sheaf of
top differential forms. In particular, if $X$ is smooth and the
marked Laurent realization of Proposition~\ref{prop:negative-vertex-chart}
exists, then Keel and Yu's construction applies to $X$ with this marked torus.
\end{proposition}
\begin{proof}
Choose a compactification of $U$ whose boundary supports an ample
divisor, as in \cite[Lemma~19.2]{KY}. The augmented algebra is
independent of this choice. It is normal by
Corollary~\ref{cor:augmentation-normal}. The compactified augmented
fiber $(\overline X,Z)$ of \cite[Proposition~19.1]{KY} is a
projective semi-log-canonical pair with $K_{\overline X}+Z$ trivial
and complement $X$. Its affine complement is Gorenstein by part~(2)
of that proposition. We apply these assertions to the augmented fiber
over $\mathbb Q$ and extend scalars to $k$. Since $X$ is integral by
Corollary~\ref{cor:augmentation-normal}, it lies in a unique irreducible
component of $\overline X$. On the normalization of this component,
let $\Delta$ be the sum of the strict transform of
$Z$ and the conductor divisor, the reduced divisor above the double
locus. By the normalization formula for a semi-log-canonical pair,
the resulting pair $(\widetilde X,\Delta)$ is log canonical, with
$K_{\widetilde X}+\Delta$ trivial. Normalization is finite, so
$\widetilde X$ is projective. It is an isomorphism over $X$, and
$\Delta$ is the reduced complement of $X$.

A generator gives a rational top form $\Omega$ with
$\operatorname{div}_{\widetilde X}(\Omega)=-\Delta$. At every divisorial
valuation of $k(X)$ the order of $\Omega$ is at least $-1$, by
properness of $\widetilde X$ and log canonicity. In particular this holds at
every toric divisorial valuation of any dense open split torus
$T\subset X$. Write $z$ for coordinates
on $T$ and $\Omega_T=\bigwedge_j d\log z_j$. The quotient
$\Omega/\Omega_T$ is a unit on $T$, hence equals $c z^h$.
At the valuation of any primitive cocharacter $v$, its contribution
to the order of $\Omega$ gives
$\operatorname{ord}_v(\Omega)=\langle h,v\rangle-1\ge-1$.
Applying this inequality to $v$ and $-v$ gives $h=0$.
The restriction of $\Omega$ generates the dualizing sheaf over the
smooth locus of $X$, hence over the normal Gorenstein variety $X$
by reflexivity. Its poles on a log resolution are at most logarithmic by log canonicity.
When the marked Laurent realization of
Proposition~\ref{prop:negative-vertex-chart} exists, it identifies
the marked torus with a principal open subset of $X$.
\end{proof}

\subsubsection{A smoothness criterion from the first mutation partners}

In the polynomial-corner chart write $Z_1,\ldots,Z_n$ for the
standard forward coordinates and $T_i=\Theta_{-e_i}$. By
Theorem~\ref{thm:corner}, $H_i=Z_iT_i$ is a polynomial in the
other coordinates. Indeed the $i$-th correction degree is zero
if $i$ is mutable, and corrections have zero designated-unit
exponents. We may invert any $Z_i$ whose inverse is already regular. Write $J$ for the remaining indices and $S$ for the Laurent ring
in the inverted coordinates.

\begin{proposition}\label{prop:forward-stratum-smoothness}
Suppose every forward theta function is Laurent in
$(T_i,Z_j:j\ne i)$ for each $i\in J$.
For every nonempty $I\subseteq J$, restrict the $H_i$, $i\in I$,
to $Z_i=0$ for $i\in I$, with all other coordinates invertible.
Suppose their common zero locus in that torus has codimension at
least one when $|I|=1$, and at least two when $|I|\ge2$, with
the empty locus allowed. Then the forward ring is $S[Z_i,T_i:i\in J]/(Z_iT_i-H_i:i\in J)$.
It is smooth if, on each nonempty such locus, the differentials
of the restricted $H_i$, $i\in I$, are independent.
More generally, it is singular if and only if for some $I$
there are a common zero $z$ and scalars $\lambda_i$, $i\in I$,
not all zero, such that the differential
$\sum_{i\in I}\lambda_i\,dH_i(z)$ vanishes on every coordinate
direction other than those $Z_j$ with $j\in I$ and
$\lambda_j\ne0$.
\end{proposition}
\begin{proof}
On the stratum where $Z_i=0$ exactly for $i\in I$, the equations
determine $T_j$ for $j\notin I$. The coordinates $T_i$ for
$i\in I$ are free. Thus the stratum has dimension $n$ minus the
codimension of the common zero locus in the complementary torus. The stratum with no zeros is the
initial torus. All other strata have dimension less than $n$.
Every component has dimension at least $n$ by the number of
equations. Thus the torus is dense and the equations form a complete intersection.
The quotient is Cohen--Macaulay and generically reduced. It is therefore
reduced, and density of the torus makes it a domain. The strata with at least two zeros have
codimension at least two. On a one-zero stratum $T_i$ is a free
coordinate, so its vanishing adds another codimension. The initial
torus and the tori obtained by replacing $Z_i$ by $T_i$ therefore
cover every height-one point. Their Laurent rings intersect in the
proposed presentation by Serre's condition $S_2$, which says that each
local ring has depth at least the smaller of two and its dimension. The actual forward
ring lies in that intersection and contains all the generators,
which proves equality.

At a point of the stratum the equations for $i\in I$ have
differentials $T_i\,dZ_i-dH_i(z)$. The other equations eliminate
their $T_j$. A dependence with coefficients $\lambda_i$ first
forces the vanishing in the statement. Conversely, when that
vanishing holds, set $T_j=(\sum_i\lambda_i\,\partial_jH_i(z))/
\lambda_j$ for each nonzero $\lambda_j$. These are free stratum
coordinates, so the resulting point is singular. Independence of
the differentials on the complementary torus excludes such a
dependence. This proves both smoothness assertions.
\end{proof}

To verify the Laurent inclusions one may use a second straight chart,
as in Corollary~\ref{cor:generic-linear-exchange-singularity}, or apply
Lemma~\ref{lem:recognition-rank-one} to the finite generators of
Lemma~\ref{lem:coercive-theta-generators}. These tests also apply to
reducible $H_i$.

\subsubsection{A smoothability obstruction from a principal open}

Let $U$ satisfy the hypotheses of Theorem~\ref{thm:corner-free-product}
over an algebraically closed field $k$ of characteristic zero. For
$0\ne f\in k[U]$, suppose $V=U_f$ contains a volume-compatible
torus. Write $B=B_k(U)$ and $B_0=B_k(V)$ for the respective augmented
algebras. Let $F^jB$ be the span of the $\theta_p$ with
$\tau_f(p)\ge j$. For a finitely generated $k$-algebra $A$, write
$T^1(A/k)=\operatorname{Ext}^1_A(\mathbb L_{A/k},A)$, where
$\mathbb L_{A/k}$ is the cotangent complex. This vector space
classifies first-order deformations of $A$ over $k$.

\begin{proposition}\label{prop:principal-open-smoothability}
The algebra $\mathcal B=\sum_{j\in\mathbb Z}t^{-j}F^jB\subset
B[t^{\pm1}]$ is a finitely generated flat $k[t]$-algebra. Its
fiber over zero is $B_0$, and away from zero its family is constant
with algebra $B$. If $\Spec B$ is smooth and $\Spec B_0$ singular,
this family is a \emph{smoothing} of $\Spec B_0$, a flat family with
that special fiber and smooth general fiber. If $B_0$
is singular and $T^1(B_0/k)=0$, then $\Spec B$ is singular.
\end{proposition}
\begin{proof}
Choose a coercive regular function $g$ and set $\rho=-\tau_g$
as in Theorem~\ref{thm:corner-free-product}. Take
a common rational subdivision of the cones on which $\rho$ and
$\tau_f$ are linear. Each cone is pointed because $\rho$ is linear there and positive
away from zero. For a label $p$ write $b_p=t^{-\tau_f(p)}\theta_p$.
These elements form a free $k[t]$-basis of $\mathcal B$.
We claim that the $b_h$ for the finitely many Hilbert generators of
the cones generate the algebra. Write $p=\sum h_i$ within one cone.
The product of the corresponding $\theta_{h_i}$ has coefficient
one at $\theta_p$. Any other label $s$ differs from $p$ by a
nonzero nonnegative mutable shift, has $\rho(s)\le\rho(p)$, and
has $\tau_f(s)\ge\tau_f(p)$, by convexity and the two linearities.
Thus the product of the $b_{h_i}$ has coefficient one at $b_p$,
and every other coefficient belongs to $k[t]$. The induction
of Lemma~\ref{lem:coercive-theta-generators}, ordered within a finite
$\rho$-level by decreasing sum of mutable coordinates, proves the
claim. Proposition~\ref{prop:principal-open} identifies the special
fiber, and $t$-inversion gives $B[t^{\pm1}]$.

For the last assertion, suppose $\Spec B$ is smooth of dimension
$n$. The Rees algebra $\mathcal B$ is a finitely generated domain
and $t$ cuts out its special fiber, so every component of that fiber
has dimension $n$. Let $I=\operatorname{Fitt}_n(\Omega_{\mathcal B/k[t]})$
be the $n$-th Fitting ideal of the module of relative K\"ahler
differentials. For a presentation with $q$ generators this ideal
is generated by the minors of size $q-n$ of the relation matrix,
with the unit ideal understood when $q\le n$. It becomes the
unit ideal after inverting $t$. Thus $t^N\in I$ for some $N$.

If $T^1(B_0/k)=0$, we can trivialize the family modulo each $t^j$
by induction. Given a trivialization modulo $t^j$, flatness identifies
the kernel at the next order with $B_0\otimes_k(t^j)/(t^{j+1})$.
Since the family is flat, the cotangent complex commutes with passage
to the special fiber \cite[Lemma~92.6.2, Tag~08QQ]{Stacks}.
The square-zero deformation theorem \cite[Lemma~92.16.1, Tag~08SP]{Stacks}
then makes the set of isomorphism classes of lifts a torsor under
$T^1(B_0/k)\otimes_k(t^j)/(t^{j+1})$. This group vanishes, so the
given lift is isomorphic to the trivial lift and the trivialization extends.
At order $N+1$, base change for differentials and Fitting ideals yields
$t^N\in\operatorname{Fitt}_n(\Omega_{B_0/k})B_0[t]/(t^{N+1})$.
Comparison of the $t^N$-coefficient forces
$\operatorname{Fitt}_n(\Omega_{B_0/k})=B_0$. At a singular closed
point of this equidimensional fiber the Zariski tangent space has
dimension greater than $n$, a contradiction.
\end{proof}

\subsection{The reverse chart of a smooth forward mirror}
\label{subsec:automatic-reverse-chart}

Let $U$ satisfy the hypotheses of Proposition~\ref{prop:negative-vertex-chart}
and Corollary~\ref{cor:augmentation-normal}.
Write $L$ for the forward label lattice, $M=L^*$, and
$B=B_k(U)\subset k[L]$ for the augmented forward algebra in that chart.
Let $\sigma\subset L_{\mathbb R}$ be the chosen full-dimensional
straight cone. Write $D_+=\mathbb N^r\oplus\{0\}^c$.
Each forward expansion has the form $\Theta_p=z^p A_p$, with
$\Supp(A_p-1)\subseteq D_+\setminus\{0\}$. Let
$K=D_+^\vee=\{q\in M_{\mathbb R}:\langle q,d\rangle\ge0
\text{ for every }d\in D_+\}$. Assume that
$Y=\Spec B$ is smooth. Proposition~\ref{prop:reverse-volume-automatic}
makes its marked torus $T_z=\Spec k[L]$ volume-compatible. Write $R=B_k(Y)$ for
its augmented counting algebra with basis $(\phi_q)_{q\in M}$
in the marking given by $T_z$.

\begin{theorem}\label{thm:automatic-reverse-chart}
There is a unique marked Laurent realization $R\hookrightarrow k[M]$
that sends $\phi_u$ to $x^u$ for every $u\in K\cap M$.
In this realization $\Phi_q$, the image of $\phi_q$, has
coefficient one at $x^q$, nonnegative integral coefficients, and
\[ \Supp(\Phi_q)\subseteq P_q\cap (q+\sigma^\vee)\cap M. \]
Here $P_q$ is the reciprocal polytope of Theorem~\ref{thm:reciprocal-polytope}.
\end{theorem}
\begin{proof}
Choose $v\in L\cap\operatorname{int}\sigma$. The forward chart
gives $Y_{z^v}=T_z$ and makes $z^t$ regular on $Y$ for every
$t\in\sigma\cap L$. If $p+q+d$ occurs in $\phi_p\phi_q$,
the regular-function inequality~\eqref{eq:regular-function-inequality} for these characters gives
$\langle t,d\rangle\ge0$, so $d\in\sigma^\vee$.
Since $v$ is interior to $\sigma$, one has $\langle v,d\rangle>0$ for every nonzero such $d$.
Proposition~\ref{prop:principal-open}, applied to $z^v$ on $Y$,
identifies its based associated graded with the torus counting algebra.
Hence the term of shift zero occurs with coefficient one.

Consider the regular function $G=\Theta_{-v}$ on $Y$.
Its coefficient at $-v$ is one and its support lies in $-v+D_+$. Thus $\tau_G(q)=-\langle v,q\rangle$
for $q\in K$. For each label $p+q+d$ occurring in $\phi_p\phi_q$ with
$q\in K$, convexity gives
$0\le\langle v,d\rangle \le -\langle v,p\rangle-\tau_G(p)$.
Indeed the upper bound follows by bounding $\tau_G(p+q+d)$ above
by the value of its term of exponent $-v$. If also $p\in K$,
only $d=0$ is possible. In general the bound cuts a compact
subset out of $\sigma^\vee$, so the possible integral shifts
are finite independently of $q$. Choose a sufficiently large
multiple $q$ of an interior lattice point of $K$ so that all possible
outputs lie in $K$. The product $\phi_p\phi_q$ then belongs
to its straight semigroup algebra. Multiplication by $\phi_q$ is injective. For any additive total
order increasing on $\sigma^\vee\setminus\{0\}$, the least label in a
nonzero sum contributes with unchanged coefficient to the least output.
Localizing the straight algebra gives the asserted embedding.
Division by $x^q$ proves the support and coefficient assertions
except for the inclusion in $P_p$. It also proves uniqueness, since
all translated basis functions are determined on the straight cone.

For a forward label $t\in L$, the function $\Theta_t$ is regular
on $Y$. The filtration it defines on $R$ gives
$\phi_p$ weight $\nu_p^z(\Theta_t)$, equal to
$\langle t,p\rangle$ on $K$. By Lemma~\ref{lem:rec-straight},
$\nu_p^z(\Theta_t)\le\nu_t^x(\Phi_p)$ for every $p\in M$.
The left side is the minimum support function of $P_p$ at $t$.
Comparison for all lattice directions gives
$\Newt(\Phi_p)\subseteq P_p$.
\end{proof}

Write $\mathsf T_Y(q_1,\ldots,q_s)$ for the coefficient of the unit
in the product $\phi_{q_1}\cdots\phi_{q_s}$ in $B_k(Y)$.
\begin{corollary}\label{cor:reverse-finite-window-traces}
For $q\in M$, set $E=P_q\cap(q+\sigma^\vee)\cap M$ and choose
$b\in K\cap M$ with $b+E\subset K$. The coefficients in
$\Phi_q=\sum_{e\in E} c_e x^e$ are the unique solution of
\[ \sum_{e\in E}c_e\,\mathsf T_Y(b+e,-b-s) =\mathsf T_Y(b,q,-b-s)\qquad(s\in E). \]
\end{corollary}
\begin{proof}
Multiplication by $\phi_b$ in the marked chart translates the whole
support into $K$. Hence
$\phi_b\phi_q=\sum_{e\in E}c_e\phi_{b+e}$ in the theta basis of
$R$. Taking the trace against $\phi_{-b-s}$ gives the equations.
A matrix entry can be nonzero only when $s-e\in\sigma^\vee$,
and its diagonal entries equal one by the shift-zero rule in the
theorem. Ordering the indices by a linear functional positive on the nonzero shifts makes the matrix unitriangular.
\end{proof}

Proposition~\ref{prop:laurent-generator-traces} computes these traces
directly from Laurent generators of the marked smooth forward mirror.
If $Y$ itself has an LP seed satisfying the hypotheses of Section~\ref{sec:incidence}, Theorem~\ref{thm:incidence}
computes its traces in the seed coordinates. To use these traces in
the marking of $T_z$, transform the contact orders by the change of
tropical coordinates.

\subsection{Reverse functions on seed cones and the criterion for marked return}
\label{subsec:reverse-seed-cones}

Suppose now that $k[U]$ is the full upper LP algebra of the marked
seed. Identify the reverse algebra $R$ with its image in $k[M]$
under Theorem~\ref{thm:automatic-reverse-chart}, and write
$h_a(p)=\nu_a^z(\Theta_p)$.

\begin{proposition}\label{prop:all-reverse-seed-cones}
Let $t$ be a reachable source seed to which
Theorem~\ref{thm:seed-windows} applies, and let $K_t\subset M_{\mathbb R}$
be its cone of slopes. Then
$\Phi_a\Phi_b=\Phi_{a+b}$ for all $a,b\in K_t\cap M$.
Write $g_1,\ldots,g_{r+c}$ for the slopes of the seed coordinates,
including its designated units. The $\Phi_{g_j}$ are algebraically
independent, and every reverse function belongs to
$k[\Phi_{g_1}^{\pm1},\ldots,\Phi_{g_{r+c}}^{\pm1}]$ inside the
fixed field $k(M)$.
\end{proposition}
\begin{proof}
The proof of Theorem~\ref{thm:seed-windows} supplies a surjection
$E_t:L\to L$ such that $E_t(p)\in\Supp(\Theta_p)$ and
$h_b(p)=\langle b,E_t(p)\rangle$ for $b\in K_t$.
Let $s$ occur in $\phi_a\phi_b$, with $b\in K_t\cap M$.
The regular-function inequality for $\Theta_p$ on $Y$ gives
\[ \langle s-b,E_t(p)\rangle \ge h_s(p)-h_b(p)\ge h_a(p)\qquad(p\in L). \]
If $a\in K_t$, surjectivity implies $s=a+b$, by testing both
signs of a lattice basis. The coefficient is one by the shift-zero
rule of Theorem~\ref{thm:automatic-reverse-chart}.
Thus the straight subalgebra is the semigroup algebra of
$K_t\cap M$, since its elements at distinct labels are independent.
The $g_j$ are a lattice basis. The mutable coordinates give
nonnegative exponents and the unit coordinates allow both signs.

For arbitrary $a$, choose preimages under $E_t$ of the
positive and negative vectors of a fixed lattice basis of $L$.
The same inequalities bound all coordinates of $s-b$
independently of $b$. Translation by a sufficiently large interior multiple
$b$ then sends every output into $K_t$. Hence $\Phi_b\Phi_a$
lies in the straight subalgebra, and $\Phi_b$ is a monomial in
the indicated generators. Division proves Laurent membership.
\end{proof}

For the return criterion we use Proposition~\ref{prop:all-reverse-seed-cones}
only at the initial seed and its immediate neighbors. Assume that
Theorem~\ref{thm:seed-windows} applies at these seeds.
Write $x'_i=\widehat F_i/x_i$ for the first mutation partner,
$a_i$ for its slope on the forward straight cone, and
$c_i\in k^\times$ for its coefficient at $x^{a_i}$.
Set $C_i=k[x_j^{\pm1}:j\ne i,\ j\text{ mutable}][u^{\pm1}]$.
\begin{theorem}\label{thm:first-reverse-partner-return}
There is a nonzero $H_i^*\in C_i$ with
$\Phi_{a_i}=H_i^*/x_i$. After inverting the other mutable
coordinates, one has
\[ R[x_j^{-1}:j\ne i]=C_i[x_i,H_i^*/x_i]. \]
In particular the following conditions are equivalent.
\begin{enumerate}
\item The marked reverse algebra equals $k[U]$.
\item $\Phi_{a_i}=c_i^{-1}x'_i$ for every mutable $i$.
\item $\widehat F_i\mid H_i^*$ in $C_i$ for every mutable $i$.
\end{enumerate}
Thus the traces at the first-mutant labels decide whether the
augmented reverse algebra recovers the source in the given torus
marking.
\end{theorem}
\begin{proof}
Theorem~\ref{thm:seed-windows} gives
$P_{a_i}=\Newt(x'_i)$. This polytope has constant $x_i$-degree
$-1$. Theorem~\ref{thm:automatic-reverse-chart} gives the
asserted form of $\Phi_{a_i}$ and
$\Newt(H_i^*)\subseteq\Newt(\widehat F_i)$.

The initial reverse chart supplies $x_i$, all the other mutable
coordinates, and the unit characters with their inverses.
After inverting the other mutables the left side of the stated
equality contains $C_i[x_i,H_i^*/x_i]$. The two source seeds
on the $i$-th edge have the same unchanged coordinate slopes.
By Proposition~\ref{prop:all-reverse-seed-cones} the localization
is contained in both $C_i[x_i^{\pm1}]$ and
$C_i[(H_i^*/x_i)^{\pm1}]$. Their intersection is precisely
$C_i[x_i,H_i^*/x_i]$, by comparison of the coefficients of
each power of $x_i$. This also covers a unit $H_i^*$.

If $R=k[U]$, Theorem~\ref{thm:seed-windows} identifies the
unique normalized regular polynomial supported in $P_{a_i}$
with $c_i^{-1}x'_i$. This proves the second condition, which
implies the third. Conversely, if $H_i^*=\widehat F_i Q_i$
with $Q_i\in C_i$, the Newton polytope of a product is the
Minkowski sum. Thus $\Newt(\widehat F_i)+\Newt(Q_i)\subseteq
\Newt(\widehat F_i)$. Boundedness forces $\Newt(Q_i)=\{0\}$,
so $Q_i$ is a scalar. The coefficient at $x^{a_i}$ determines
this scalar as $c_i^{-1}$. Consequently the reverse functions
are Laurent on the initial torus and its immediate neighbors, by
Proposition~\ref{prop:all-reverse-seed-cones}. Their union contains
every codimension-one point of the normal variety $U$.
Thus $R\subseteq k[U]$. Theorem~\ref{thm:one-sided-return} applies
with the compatible cones and filtrations of
Theorem~\ref{thm:automatic-reverse-chart} and proves equality.
\end{proof}

For the rank-one source $xx'=1+u+v$ of the introduction, we can
compute every reverse function. Write $X,U_1,U_2$ for the torus
characters dual to $x,u,v$. For $p=(i,j,\ell)$, set
$d=\max(0,-j,-\ell)$. Proposition~\ref{prop:axis} gives
\[ \Theta_p=X^i U_1^j U_2^\ell(1+X)^d =X^i U_1^{j+d}U_2^{\ell+d}W^d,\qquad W=\frac{1+X}{U_1U_2}. \]
The generators $X^{\pm1},U_1,U_2,W$ all occur among these theta
functions. Hence the forward mirror $Y$ has coordinate ring
$k[X^{\pm1},U_1,U_2,W]/(U_1U_2W-1-X)$. Indeed, eliminating $X$
identifies this ring with $k[U_1,U_2,W,(U_1U_2W-1)^{-1}]$, which
embeds in the marked Laurent field by the formula for $W$.
In particular $Y$ is smooth.

Use $U_1,U_2$ as mutable coordinates and $X$ as a designated unit.
Both ordinary exchanges are $1+X$. Their exchange Laurent polynomials
are $(1+X)/U_2$ and $(1+X)/U_1$, respectively, so mutation at
either index gives $W$. The three seed tori together cover the locus
where at least two of $U_1,U_2,W$ are nonzero. Their complement has
codimension two, and inverting either mutable coordinate gives the
required one-exchange model. The torus volume form extends across $Y$
by Proposition~\ref{prop:reverse-volume-automatic}. We may therefore
apply Corollary~\ref{cor:independent-face} to the two equal exchanges.
In the original coordinates this gives
\[ \Phi_a=x^{a_x}u^{a_u}v^{a_v}(1+u+v)^{\max(0,-a_x)} \qquad\bigl(a=(a_x,a_u,a_v)\in M\bigr). \]
The two Laurent markings agree on the positive orthant, where both
send each theta function to the character of its label. In particular
$\Phi_{-e_x}=(1+u+v)/x=x'$, so
Theorem~\ref{thm:first-reverse-partner-return} gives marked return.

The exchange $G=1+u+v+2uv$ gives a smooth rank-one source whose
forward mirror is singular. Indeed, $G$ is primitive and linear in
$v$, and $\partial G/\partial v=1+2u$ is nonzero on $G=0$.
Its Newton polygon is the unit square, so Proposition~\ref{prop:axis}
gives $\Theta_{(i,j,\ell)}=X^iU_1^jU_2^\ell(1+X)^d$ with
$d=\max(0,-j)+\max(0,-\ell)$.
It follows that $k[Y]$ is generated by $X^{\pm1},U_1,U_2$ and
$W_h=(1+X)/U_h$ for $h=1,2$. Eliminating $X$ gives
\[ k[Y]=k[U_1,U_2,W_1,W_2,(U_1W_1-1)^{-1}]/(U_1W_1-U_2W_2). \]
Indeed, the ring on the right is a domain and inverting $U_1U_2$
recovers the marked Laurent ring. The origin is an ordinary double
point of $Y$. Theorem~\ref{thm:polynomial-corner-reciprocity}
still gives $P_a=Q_a$ for every $a\in M$.

\section{Recovering effective numerical curve classes}\label{sec:coefficient-recovery}

\subsection{Recovering the curve class from the boundary grading}

Augmenting a theta product sets every curve parameter equal to one.
To recover the class of a contributing curve, we use its intersection
numbers with the boundary. Keel and Yu's grading computes these numbers
from the input and output labels.

Let $U$ be a smooth affine log Calabi--Yau variety containing a dense
split torus over an algebraically closed field $k$ of characteristic
zero. Choose a smooth projective compactification $\overline U$ of $U$ and write its
simple normal crossings boundary as $D_1+\cdots+D_s$, where the
$D_j$ are its prime components. In the group
$G=N_1(\overline U,\mathbb Z)$ of integral one-cycles modulo
\emph{numerical equivalence}, two cycles have the same class if they
have the same intersection number with every Cartier divisor. Classes
represented by nonnegative sums of curves form the \emph{effective
monoid} $Q\subset G$. In its monoid ring $R=\mathbb Z[Q]$ the basis
elements multiply by $t^\beta t^{\beta'}=t^{\beta+\beta'}$.

Write $H_R=\bigoplus_m R\vartheta_m$ for the Keel--Yu algebra of
the pair $(\overline U,\sum_j D_j)$. For a skeletal label $m$,
write $D_j^{\mathrm{trop}}(m)$ for the value of the corresponding
valuation on a local equation for $D_j$ at its center on $\overline U$.
This value is independent of the equation since local units have value
zero. Set $h(m)=(D_j^{\mathrm{trop}}(m))_{j=1}^s$.
By the \emph{boundary grading}
of Keel and Yu~\cite[Lemma~16.1 and Theorem~16.2]{KY}, each term
$t^\beta\vartheta_r$ in $\prod_i\vartheta_{p_i}$ satisfies
\begin{equation}\label{eq:curve-boundary-grading}
(D_j\cdot\beta)_j=\sum_i h(p_i)-h(r).
\end{equation}
Thus the labels determine
$\beta$ whenever the boundary intersections distinguish numerical
curve classes. We show this holds integrally when $\Pic(U)=0$.

\begin{lemma}\label{lem:integral-boundary-pairing}
If $\Pic(U)=0$, the intersection map
$w:G\to\mathbb Z^s,\ \beta\mapsto(D_j\cdot\beta)_j$, is a split
injection with cokernel
$\operatorname{Hom}(\mathcal O(U)^\times/k^\times,\mathbb Z)$.
\end{lemma}
\begin{proof}
The variety $\overline U$ is rational. We first show that its Picard group is
free of finite rank and that intersection identifies $G$ with the
integral dual of $\Pic(\overline U)$. Injectivity follows from numerical
equivalence. The remaining assertions hold for projective space,
and weak factorization reduces them to invariance under blowups
of smooth centers \cite[Theorem~0.1.1]{AKMW}.

Let $\pi:\widetilde X\to X$ be such a blowup, with exceptional
divisors $E_i$ over the components of codimension at least two.
The divisor sequence gives
$\Pic(\widetilde X)=\pi^*\Pic(X)\oplus\bigoplus_i\mathbb Z[E_i]$.
Indeed, subtracting the pullback of a divisor's pushforward leaves an
exceptional divisor, and pushing forward a principal relation gives
the corresponding principal relation on $X$.
Every integral curve on $X$ lifts to a one-cycle of degree one over
it. For a curve contained in the center, take the closure of a rational
point of the exceptional projective space over its function field.
Otherwise take the strict transform. By the projection formula these
lifts give the prescribed values on $\pi^*\Pic(X)$. A line in a
fiber of $E_i$ has degree $-1$ on $E_i$ and zero on the other
summands, so the values on all exceptional classes can be prescribed
independently. Conversely, extend a functional on $\Pic(X)$ by
zero on the exceptional summands, realize it on $\widetilde X$
and push the resulting cycle forward. This proves the invariance.

Since $\overline U$ is proper and $\Pic(U)=0$, divisor localization gives
\[ 0\longrightarrow \mathcal O(U)^\times/k^\times\longrightarrow \mathbb Z^s\xrightarrow{\,e_j\mapsto[D_j]\,}\Pic(\overline U)\longrightarrow0 . \]
The sequence splits since $\Pic(\overline U)$ is free. Taking its integral
dual gives the assertion.
\end{proof}

Write $H_1=\bigoplus_m\mathbb Z\psi_m$ for the augmentation of
$H_R$ at $t^\beta=1$. Inverting the curve monomials gives the
\emph{coefficient torus} $\Spec\mathbb Z[G]$ over $\mathbb Z$.
Keel, White and Yu observe that the fibers over this torus
are isomorphic in their example with trivial Picard group
\cite[\S21, before Conjecture~21.9]{KWY}.
The following integral formula recovers the coefficients in the
effective monoid ring.

\begin{corollary}\label{cor:full-integral-coefficient-gauge}
Assume $\Pic(U)=0$. Choose an integral retraction
$\rho:\mathbb Z^s\to G$ of $w$, and write $d(m)=\rho h(m)$.
If $\prod_i\psi_{p_i}=\sum_r n^r\psi_r$, then
\[ \prod_i\vartheta_{p_i}=\sum_{r:\,n^r\ne0} n^r t^{\,\sum_i d(p_i)-d(r)}\vartheta_r . \]
Every exponent in the sum belongs to $Q$.
The homomorphism $\vartheta_m\mapsto t^{d(m)}\psi_m$ identifies
$H_R$ with $\bigoplus_m R\,t^{d(m)}\psi_m$ in
$\mathbb Z[G]\otimes_{\mathbb Z}H_1$.
\end{corollary}
\begin{proof}
The boundary grading~\eqref{eq:curve-boundary-grading} and
Lemma~\ref{lem:integral-boundary-pairing} show that at most one
effective class contributes to each output label $r$. Applying $\rho$
identifies it as $\sum_i d(p_i)-d(r)$. Augmentation preserves its
coefficient, giving the product formula and multiplicativity of the
stated map. Independence of the $\psi_m$ proves injectivity and
identifies the image.
\end{proof}

Some $d(m)$ can lie outside $Q$, which accounts for the use of
the group ring in the embedding. The augmentation specializes only
the curve variables. The designated units of an LP seed still have
their specified values.

\subsection{Gauss valuations on the curve and character lattices}

Choose the compactification and local Laurent chart of
Proposition~\ref{prop:formal-local-chart}, assuming $\Pic(U)=0$.
Let $L$ be the cocharacter lattice of the source torus. Choose a
generic point $\xi$ in the interior of the normal cone $C$ at the
negative vertex of Proposition~\ref{prop:local-support}.
Write $V_m=V_{\xi,m}\in\mathbb Z[Q][L]$ for the expansion with
curve classes of the theta function of label $m$ in the chart at
$\xi$. Its monomials have the form $t^\beta z^e$, with curve class
$\beta\in Q$ and character exponent $e\in L$. Write
$\Theta_m=\sum_e c_{m,e}z^e$ for its augmentation.
For $p\in C\cap L$ let $\delta_p$ be the numerical class of
the straight term of $V_p$.

\begin{proposition}\label{prop:curve-character-gauss}
There is a unique homomorphism $\ell:L\to G$ such that
$\ell(p)=\delta_p-d(p)$ on $C\cap L$. One has
\[ V_m=\sum_{e\in\Supp\Theta_m} c_{m,e}\,t^{d(m)+\ell(e)}z^e . \]
Thus each exponent of $\Theta_m$ occurs in one effective numerical
class. The Newton polytope of $V_m$ is the image of
$\Newt(\Theta_m)$ under $e\mapsto(d(m)+\ell(e),e)$.

Let $K$ be a field containing $\mathbb Q$, with a real valuation
$v$ that is zero on $\mathbb Q^\times$.
For a homomorphism $\chi:G\to K^\times$, specialize $t^\beta$
to $\chi(\beta)$, obtaining $V_m^\chi$. Write
$\gamma=v\circ\chi$. The Gauss extension of $v$ with weight
$b\in L^*_{\mathbb R}$ is
$\nu_b(\sum_e u_e z^e)=\min_{u_e\ne0}(v(u_e)+\langle b,e\rangle)$.
For any nonzero finite
sum with coefficients $a_m\in K$,
\[ \nu_b\!\left(\sum_m a_mV_m^\chi\right)= \min_{a_m\ne0}\bigl(v(a_m)+\gamma(d(m))+h_{b+\gamma\circ\ell}(m)\bigr), \qquad h_w(m)=\min_{e\in\Supp\Theta_m}\langle w,e\rangle . \]
\end{proposition}
\begin{proof}
For $p\in C\cap L$, Proposition~\ref{prop:local-support} gives
$\Theta_p=z^p$. Positivity of the effective coefficients gives
$V_p=t^{\delta_p}z^p$ for a unique $\delta_p\in Q$.
After extending scalars to $\mathbb Z[G]$, compose the inverse
of the isomorphism in Corollary~\ref{cor:full-integral-coefficient-gauge}
with the local expansion. The resulting homomorphism sends
$\psi_m$ to $t^{-d(m)}V_m$. In the faithful augmented chart,
$\psi_p\psi_{p'}=\psi_{p+p'}$ for $p,p'\in C\cap L$. Hence
$p\mapsto\delta_p-d(p)$ is additive there. The monoid $C\cap L$
generates $L$, since any lattice point translates into $C$ by
a sufficiently large multiple of an interior lattice point. This
proves existence and uniqueness of $\ell$.

The augmented expansion followed by $z^e\mapsto t^{\ell(e)}z^e$
agrees with this homomorphism on the $\psi_p$ for $p\in C\cap L$.
For each $m$, multiply $\Theta_m$ by $z^{Nb_0}$ for a fixed
interior lattice point $b_0$ and $N$ large enough that all its
translated exponents lie in $C$. Now $\psi_{b_0}^N\psi_m$
is a linear combination of straight elements. The two homomorphisms
agree on this product, and the image of $\psi_{b_0}$ is a unit
in $\mathbb Z[G][L]$. This proves the expansion. Since
$V_m\in\mathbb Z[Q][L]$, all classes occurring are effective.
Taking convex hulls gives the Newton polytope.

Nonzero integer coefficients have valuation zero. At a generic
translated weight $b+\gamma\circ\ell$, each $\Theta_m$ in the
sum has a unique minimizing exponent, and the exponents are distinct
by Theorem~\ref{thm:marked-path-vi}. Otherwise one could cancel the
common least term in a rational linear combination. Multiplication
by $a_m$ preserves that exponent. The terms of minimum valuation therefore have distinct characters
and cannot cancel. For fixed coefficients both sides are minima of
finitely many affine functions of $b$. Continuity gives the formula
at every weight.
\end{proof}

For example, take $K=\Frac(k[G])$, with the Gauss valuation for any
weight on $G$, and $\chi(\beta)=t^\beta$. The proposition gives
the minimum formula for all Gauss weights on $G\oplus L$. For a
nonzero finite sum $\sum_m a_m V_m$, with $a_m\in k[G]$, its
Newton polytope is therefore the convex hull of the Newton polytopes
of its nonzero summands. This applies in particular to $a_m\in k[Q]$.

\begin{corollary}\label{cor:augmentation-normal}
Let $U$ be a smooth affine log Calabi--Yau variety over an
algebraically closed field of characteristic zero. If $U$
contains a dense split torus and $\Pic(U)=0$, then for every
algebraically closed field $k$ of characteristic zero the augmented
mirror algebra $B_k(U)=k\otimes_{\mathbb Z}B_{\mathbb Z}(U)$ is a
finitely generated normal domain.
\end{corollary}
\begin{proof}
Specializing the finitely generated $R$-algebra $H_R$ at
$t^\beta=1$ gives the finitely generated $k$-algebra
$B_k(U)=k\otimes_{\mathbb Z}H_1$ \cite[Theorem~16.8]{KY}.
Since the augmentation is independent of the compactification, we may
choose a smooth projective compactification whose boundary is the
support of an effective ample Cartier divisor $F$
\cite[Proposition~17.3 and Lemma~19.2]{KY}.
Over the coefficient torus Corollary~\ref{cor:full-integral-coefficient-gauge}
identifies the mirror family with the scalar extension of the augmented
algebra. Write $k(G)=\Frac(k[G])$ and
$\mathbb Q(G)=\Frac(\mathbb Q[G])$. The generic fiber over $\mathbb Q(G)$ is normal by
\cite[Proposition~19.1]{KY}, hence geometrically normal in characteristic
zero. Thus $k(G)\otimes_k B_k(U)$ is normal, and normality descends to
$B_k(U)$.

It remains to show that the spectrum is connected. Filter $B_k(U)$
by assigning $\psi_p$ degree $F^{\mathrm{trop}}(p)$, which is
positive for $p\ne0$ \cite[Lemma~15.2]{KY}. For a term
$t^\beta\vartheta_r$ in $\vartheta_p\vartheta_q$, the
boundary grading gives
$\deg\psi_p+\deg\psi_q-\deg\psi_r=F\cdot\beta$.
By ampleness this difference vanishes exactly in class zero. The
associated graded algebra therefore has the class-zero multiplication
of \cite[Proposition~15.12]{KY} and is a reduced Stanley--Reisner
ring over $k$. An idempotent of positive filtered degree would have
a nonzero leading form whose square is zero. Hence every idempotent
lies in degree zero, which is $k$. A connected Noetherian normal
ring is a domain.
\end{proof}

\section{Variation of the exchange coefficients}
\label{sec:coefficient-variation}

\subsection{Generic constancy and spaces of source sections}

Let $\mathcal P$ be an irreducible complex variety. In fixed marked mutable
coordinates $x$ and unit coordinates $u$, suppose the ordinary
exchanges have polynomial-corner representatives
\[ G_{i,t}=\sum_{a\in A_i} c_{ia}(t) (x,u)^a \quad(t\in \mathcal P),\qquad 0\in A_i\subset\mathbb N^{r+c}, \qquad c_{ia}\in\mathcal O(\mathcal P)^\times . \]
Assume that each fiber defines a smooth full upper LP source $U_t$
satisfying the hypotheses of Theorem~\ref{thm:corner}. We compare
the augmented mirrors of the different fibers in their common marked
Laurent coordinates.

\begin{theorem}\label{thm:generic-coefficient-obstruction}
On a nonempty Zariski open subset of $\mathcal P$, the marked forward
algebra and all the forward, reciprocal and smallest source Newton
polytopes are constant. Outside a countable union of proper closed
subsets, the forward Laurent polynomials themselves agree at every label.

Suppose there are integers $n_{ia}$ with
\[ \sum_{a\in A_i} n_{ia}=0\quad\text{for each }i, \qquad \sum_{i,a} n_{ia}a=0, \qquad \prod_{i,a} c_{ia}^{n_{ia}}\text{ nonconstant on }\mathcal P . \]
Then there exist two very general fibers with identical marked forward
theta bases whose source upper algebras cannot be identified by any
diagonal torus translation.
\end{theorem}

We first study how the traces vary in a family.

\begin{proposition}\label{prop:trace-drop-locus}
In an algebraic coefficient family satisfying the hypotheses of
Section~\ref{sec:incidence} on each smooth full upper fiber, fix the
coordinates and designated units. The trace for each fixed
list of contact orders is lower semicontinuous. For each positive integer
$d$ its locus of values below $d$ is closed and can be computed
by elimination over an effective ground field.
\end{proposition}
\begin{proof}
Write $\mathcal P$ for the parameter base and $M$ for the irreducible
space of distinct marked positions and torus evaluation points.
The incidence equations of Theorem~\ref{thm:incidence} and its
open conditions define a relative scheme over $\mathcal P\times M$.
Denominators can be cleared uniformly using the union of the
exchange supports, since the extra clearing factors are units
at the roots on the specified open locus. Restrict to the open
subscheme where the Jacobian of the remainder equations in the
unknown monic coefficients is invertible. Its projection is
\'etale by the square presentation. Over the generic point of
$M$ in each geometric fiber, this open subscheme contains all
the counted solutions. Indeed, in the proof of
Theorem~\ref{thm:incidence} the fixed-domain mapping space is
smooth at these maps and evaluation is \'etale, which makes each
root of the fixed-evaluation equations simple.

Form the $d$-fold fiber product of this \'etale subscheme over
$\mathcal P\times M$ and remove the pairwise diagonals. Its image in
$\mathcal P\times M$ is open. Its projection to $\mathcal P$ is also open,
since $\mathcal P\times M\to \mathcal P$ is smooth. A geometric point $t$ is
in this projection precisely when the fiber over $t$ of the
first open image is nonempty. By irreducibility of $M$, this
is equivalent to having $d$ distinct simple roots at the generic
marks and evaluation of that fiber, hence to trace at least $d$.
Both images are constructible projections given by polynomial
equations and inequations, so effective elimination computes them
and their complements. The trace is zero if a mutable defect is negative or the orders of a
unit fail to balance.
\end{proof}

\begin{corollary}\label{cor:constructible-reciprocal-locus}
Fix the marked coordinates, designated units and polynomial-corner
supports $A_i$. Among the smooth full upper LP sources of
Theorem~\ref{thm:polynomial-corner-reciprocity} with these supports,
only finitely many integral augmented forward algebras occur as
subrings of the marked Laurent ring $\mathbb Z[L]$. Their number
has a computable upper bound depending only on the $A_i$.
For any finite-type family of these varieties, the base has a finite
stratification by locally closed sets on each of which the marked
algebra, all the forward Newton polytopes, and all the polytopes
$P_a=Q_a$ are constant.
\end{corollary}
\begin{proof}
Let $\Delta$ be the convex hull of the supports of $1,x_i,G_i/x_i,u_j,u_j^{-1}$
for all $i,j$. These functions are regular for every choice of the
exchange coefficients under consideration. Since each $G_i$ has a
nonzero constant term, $\Delta$ contains both signs of each coordinate
vector. A general linear combination of the functions has Newton polytope $\Delta$
and is therefore coercive. By Lemma~\ref{lem:coercive-theta-generators},
the normal fan of $\Delta$ gives a finite set of theta generator
labels that depends only on the exchange supports.

For a fixed label $m$, the support box in
Corollary~\ref{cor:trace-reconstruction} depends only on $m$ and
the $A_i$. Thus a single finite set of contact lists $I$
suffices to compute the Laurent polynomials of the
generators. Lists with a negative mutable defect or unbalanced unit
orders have trace zero, and we assign $B_I=0$ to them.
For each remaining list $I$, let $N_I$
be the sum of its mutable defects, hence the number of unknown
monic coefficients and of scalar remainder equations in
Theorem~\ref{thm:incidence}. The exchange supports give upper bounds
$d_{I1},\ldots,d_{IN_I}$ for their total degrees in these unknowns.
Then
the affine B\'ezout bound gives
$B_I=\prod_{j=1}^{N_I}\max(1,d_{Ij})$ as an upper bound for
their isolated solutions. The empty product is one.
The counted solutions are isolated by Theorem~\ref{thm:incidence},
so there are at most
$\prod_I(1+B_I)$ possible tuples of trace values, each of which
determines the marked algebra. In a family their level sets are
locally closed by Proposition~\ref{prop:trace-drop-locus}.

Consider two fibers with the same marked ring. For a strictly
positive separating weight, the least exponent of each theta
polynomial is its label and its least coefficient is one.
Expand $\Theta_p$ from the first fiber in the theta basis of the
second. The least label in this expansion is $p$, with coefficient
one. By valuative independence (Theorem~\ref{thm:forward-vi-main}),
the Newton polytope of $\Theta_p$ in the second fiber is contained
in its polytope in the first. Exchanging the fibers proves equality.
The forward polytopes determine every $P_a$ by
Theorem~\ref{thm:reciprocal-polytope}, and
Theorem~\ref{thm:polynomial-corner-reciprocity} gives $P_a=Q_a$.
\end{proof}

\begin{proof}[Proof of Theorem~\ref{thm:generic-coefficient-obstruction}]
Since $\mathcal P$ is irreducible, one stratum of
Corollary~\ref{cor:constructible-reciprocal-locus} is dense and open.
This proves the assertions about the marked algebra and the polytopes.
For each label $m$, the finitely many traces determining $\Theta_m$
are bounded and lower semicontinuous, so each attains its maximum on
a nonempty open subset. Intersecting these opens over the countable
lattice of labels gives the asserted agreement of the theta bases.

The product $\prod_{i,a}c_{ia}^{n_{ia}}$ is invariant under diagonal
translation and under rescaling each exchange by a scalar. By
Proposition~\ref{prop:marked-exchange-rigidity}, a translation
identifying two marked source rings identifies their exchanges up to
coefficient units. For exchanges with the same finite support, such a
unit is a scalar. The product must therefore take the same value on
the two fibers. Over $\mathbb C$ the complement of a countable
union of proper closed subsets is dense in $\mathcal P$. Since the
product is nonconstant, it takes at least two values there.
\end{proof}

\begin{corollary}\label{cor:forward-smoothness-strata}
Let the ordinary exchanges vary algebraically, with fixed coordinates
and designated units, over a finite-type characteristic-zero
parameter variety. Suppose every geometric fiber, including those
over generic points, is the smooth spectrum of a full upper LP algebra
satisfying the hypotheses of
Section~\ref{sec:incidence}.
The parameter variety has a finite stratification by locally closed
sets on each of which the augmented forward algebra has the same
integral presentation in the same finite set of theta generators.
In particular the locus of smooth forward mirrors is constructible.
For defining data over a field with effective polynomial elimination
and primary decomposition the strata and their smoothness can be
computed.
\end{corollary}
\begin{proof}
Over the function field of an irreducible component, the first-star
equations give a coercive function. To see why existence over the
algebraic closure suffices, fix a support bound. The corresponding
space of sections is the kernel of a linear system over the function
field. Extending scalars preserves the set of exponents that occur in
this space, and a general section over the function field has this
entire support. Spread the function and its
Laurent expressions on the first star to a dense open, shrinking it further to preserve the supports of the exchanges and of the coercive
function. In the notation of Theorem~\ref{thm:finite-forward-presentation},
we can therefore use the same weight $w$, the same finite set $H$ of
generator labels, and the same bound $D$ in every fiber of this open
set. In particular the set $S_D=\{p:w(p)\le D\}$ is fixed.

Theorem~\ref{thm:incidence} computes the generic value of each of
the finitely many traces. By Proposition~\ref{prop:trace-drop-locus},
each trace attains its maximum on a nonempty open subset of the
irreducible base. On the intersection of these opens the relations of
Theorem~\ref{thm:finite-forward-presentation} are constant. Apply the
same argument to each irreducible component of the complement.
Induction on dimension gives a finite stratification. Under the stated
effectiveness hypotheses, elimination computes each stratum.
\end{proof}

\subsubsection{Sections with prescribed Newton polytope bounds}

Let $\mathcal P^\circ$ be a nonempty open subset as in
Theorem~\ref{thm:generic-coefficient-obstruction} on which the $P_a$
are constant. Theorem~\ref{thm:polynomial-corner-reciprocity}
identifies them with the smallest source polytopes in each fiber. Write
$V_\Delta(t)=\{f\in k[U_t]:\Supp(f)\subseteq\Delta\}$ for a lattice
polytope $\Delta\subset M_{\mathbb R}$, and let
$I(\Delta)=\{a\in M:P_a\subseteq\Delta\}$. This set is finite because
$a\in P_a$.

\begin{proposition}\label{prop:closed-realization-locus}
For each lattice polytope $\Delta$, the spaces $V_\Delta(t)$
are the fibers of a vector bundle of rank $\lvert I(\Delta)\rvert$
on $\mathcal P^\circ$. Its formation commutes with base change.
\end{proposition}
\begin{proof}
In one fiber choose a normalized regular function $f_a$ supported
in $P_a=Q_a$ for each $a$. These functions form the basis of
Theorem~\ref{thm:minimum-source-newton-polytope}. Valuative independence
implies that the Newton polytope of a linear combination is the convex
hull of the polytopes of its nonzero terms. Thus $V_\Delta(t)$ has
basis $(f_a)_{a\in I(\Delta)}$ and the stated dimension.

To construct the bundle, take as unknowns the coefficients of the
monomials with exponents in $\Delta\cap M$. In a negative mutable slice
of exponent $-h$, its Laurent coefficient must be $G_{i,t}^h$
times a quotient in the unchanged coordinates. The coordinate minima
and maxima of the slice minus the corresponding extrema of
$h\Newt(G_{i,t})$ bound the support of the quotient. These bounds
are uniform because the exchange supports are fixed. Introducing the
quotient coefficients as further unknowns gives a matrix $D_\Delta$
of regular functions. Each quotient is unique, so projection on the
section coefficients identifies $\ker D_\Delta(t)$ with $V_\Delta(t)$.
The matrix has constant rank on the reduced base. Locally a minor
of maximal rank is invertible and all larger minors vanish. Row and
column operations therefore split the matrix into an identity matrix
and a zero matrix. Its kernel is locally a direct summand and commutes
with base change.
\end{proof}

The local splitting also gives the tangent condition for constant
rank. For a tangent vector $v\in T_t\mathcal P^\circ$, differentiate
the matrix entries in direction $v$ to obtain $\partial_vD_\Delta$.
The induced map $h\mapsto[(\partial_vD_\Delta)h]$ from the kernel
of $D_\Delta(t)$ to its cokernel
vanishes. More generally,
on a reduced irreducible base where $P_b\subseteq Q_b(t)$ and the
bounds defining $D_\Delta$ are fixed, one has
$\dim V_\Delta(t)\le |I(\Delta)|$. All labels are realized precisely
when equality holds for every $\Delta=P_a$. Indeed the possible least
exponents in $V_\Delta(t)$ belong to $I(\Delta)$, and equality forces
each of them to occur. These are closed conditions on the ranks of the
matrices, since equality is the largest possible kernel dimension. Noetherianity reduces their intersection to a finite one.
At a point of simultaneous realization, a nonzero kernel-to-cokernel
derivative for $D_{P_a}$ implies that one of its next larger minors
is nonzero generically. In that case some label in $I(P_a)$ fails to
be realized on a dense open subset. This last criterion applies to
families for which only the inclusion is known.

\begin{corollary}\label{cor:relative-reciprocal-basis}
On any reduced affine subvariety $Z\subseteq \mathcal P^\circ$, normalized
reciprocal functions can be chosen to give a free
$\mathcal O(Z)$-basis of the Laurent polynomials satisfying the
first-star divisibility conditions. Its members indexed by $I(\Delta)$
form a basis of the sections supported in $\Delta$. Both assertions
commute with specialization to every residue field of $Z$.
\end{corollary}
\begin{proof}
On the bundle for $P_b$, taking the coefficient at $b$ is
surjective in every fiber. Affineness gives a global section mapping
to one. Its Laurent polynomial $f_b$ is normalized and supported in
$P_b$. If $b\in I(\Delta)$, then $f_b$ belongs to the bundle for
$\Delta$. The quotient coefficients outside the prescribed bounds vanish
on each fiber and hence on the reduced base. The resulting family
$(f_b)_{b\in I(\Delta)}$ is a basis on each fiber, so it gives a
trivialization of that bundle. Every Laurent polynomial has bounded
support, and the same bounds apply to its quotients over $\mathcal O(Z)$
by reducedness. Taking the union over lattice boxes proves the assertion
for the algebra.
\end{proof}

\subsection{Common refinements of the source normal fans}

For a fixed source, can a single finite fan describe the minimizing
faces of all its smallest Newton polytopes as the Gauss weight varies?
Fix the initial character lattice $M$, its dual $L$, and a
lattice-separating weight $\omega$. Assume the full upper algebra
with Laurent coefficients has a section whose $\omega$-least vertex
is negative in every mutable coordinate. Use the smallest Newton
polytopes $Q_a$ and the dual family $T_p$ from
Theorem~\ref{thm:minimum-source-newton-polytope} and
Proposition~\ref{prop:source-polar-polytope}. A common refinement of
the normal fans of all $Q_a$ means a finite complete pointed rational
fan on which the minimum support function of every $Q_a$ is linear
on every cone. Individual polytopes may have smaller dimension.

\begin{proposition}\label{prop:discover-source-fan}
Let $p_1,\ldots,p_l$ be the primitive ray generators of a nonzero
pointed rational cone $\sigma\subset L_{\mathbb R}$, and set
$p_*=\sum_jp_j$. The minimum support functions of all $Q_a$ are
linear on $\sigma$ if and only if
\[ T_{p_*}=\sum_j T_{p_j}. \]
Failure of the equality gives an integral label $a$ for which
$Q_a$ subdivides $\sigma$. There is a search using only finite
first-star systems which produces a common refinement whenever one
exists, without a bound or a prior choice for its rays.
\end{proposition}
\begin{proof}
The two polytopes have the same minimum at a lattice label $a$ precisely when
$\min_{b\in Q_a}\langle b,p_*\rangle =\sum_j\min_{b\in Q_a}\langle b,p_j\rangle$.
Equality holds precisely when a single point of $Q_a$ attains all
the minima on the right. This is equivalent to linearity of its
support function on $\sigma$. Equality of minima at all lattice
weights characterizes equality of rational polytopes. If the
inclusion $T_{p_*}\subseteq\sum_jT_{p_j}$ is strict, strict
separation of a vertex outside the smaller polytope supplies an
integral $a$ at which the equality fails.

Start with any complete pointed rational fan. Test its maximal cones
and, at a failure, refine the entire fan by the normal fan of the
computed $Q_a$. Each such step strictly subdivides a cone, and the
polytopes at each step are computable from the first-star divisibility conditions.
If a finite common refinement exists, intersect it with the starting
fan. Every successive fan is a coarsening of that one finite fan,
which has only finitely many coarsenings because every cone of a
coarsening is a union of its cones. Strict refinement must therefore terminate.
When the procedure terminates, the equality holds on every cone,
so the resulting fan is a common refinement.
\end{proof}

The finite-refinement property is symmetric in the two families
$(Q_a)$ and $(T_p)$. On a common refinement the values of every
$Q_a$ at the rays determine its support function by linearity.
Taking a common refinement of the normal fans of the finitely many
corresponding $T_p$ then makes the support minima of \emph{all}
$T_p$ linear there, since on each original cone these minima are
nonnegative linear combinations of the ray values. Interchanging
the two lattices proves the converse. In particular the same test
can be run in label space using only $Q_{a_*}$ and the ray polytopes
$Q_{a_j}$. This recovers a finite fan whenever one exists.

In the setting of Theorem~\ref{thm:reciprocal-polytope}, use the
mirror chart determined by the vertex exposed by $\omega$.
On any such common refinement, all reciprocal polytopes are realized
on the source if and only if
$\Newt(\Theta_p)=T_p$
at the ray generators and at the sum of those generators in each
maximal cone. Indeed the minimum at $a$ on the counted polytope
is the minimum of $P_a$ at $p$. Equality on the rays and their
sum forces these minima of $P_a$ to be attained simultaneously,
by the same argument as for $Q_a$. Both support functions then
agree throughout the cone and, over the full fan, $P_a=Q_a$.
Conversely realization gives these equalities at every $p$.
In the dual search on label space we test
$P_{a_*}=\sum_jP_{a_j}$ for generators of a cone and their sum. If equality fails,
strict separation supplies a forward polytope $\Newt(\Theta_p)$
by which to refine. When this search terminates, the $P_a$ are
Minkowski additive on the label cones. Realization at their rays
suffices by taking products and applying the dilation identity for
$Q_a$. Below we give examples for which the refinement continues indefinitely.

\begin{theorem}\label{thm:source-family-direction-strata}
Let ordinary Laurent-coefficient LP exchanges vary algebraically over
a finite-type parameter variety with fixed marked coordinates. Use
the same separating $\omega$ and suppose all geometric fibers,
including those over generic points, admit a section whose
$\omega$-least exponent is negative in every mutable coordinate.
\begin{enumerate}
\item For any designated finite set $E\subset L$, there is a finite
locally closed stratification on which
$\min_{Q_a}\langle -,p\rangle$ is constant for every $a\in M$
simultaneously and every $p\in E$.
\item If a reference fiber has a finite common refinement for all its
$Q_a$, the locus of fibers with precisely the same $Q_a$ for every
label is constructible. A determining list consists of the rays and
the ray sums of one such fan of the reference fiber, tested on the
dual polytopes $T_p$. Equivalently, one can test the $Q_a$ at
the rays and ray sums of the dual fan constructed above.
\item If the geometric generic fiber of every irreducible closed
subvariety admits a finite common refinement, there is a finite
locally closed stratification on which $Q_a$ is constant for every
label $a$. The construction finds a fan on each successive generic fiber. The finite computations are effective over an effective field
with elimination and primary decomposition, given effective comparisons
of the weights of lattice points.
\end{enumerate}
\end{theorem}
\begin{proof}
Over the function field of each irreducible component there is a section
whose least exponent is negative in every mutable coordinate. Indeed, take
the finite support of such a section over the geometric function field, prescribe coefficient one
at its least label and descend the resulting first-star linear
system. Spread the section and its expressions in the finite star
to a dense open, shrinking to keep the Newton polytope and its
negative least exponent $v$ fixed. This also keeps fixed the finite
support bounds for the corrections $R_a=x^{a-lv}g^l$, with $l$ chosen
as in Theorem~\ref{thm:minimum-source-newton-polytope}. At any finite
list of labels, the bounded linear systems determining $Q_a$ can be
checked by rank tests. For these tests one may use the ordinary $F_i$,
whose ratios to the associated exchanges are Laurent units in the
unchanged variables. Shrink so the answers to all tests used in the
generic calculation agree in every geometric fiber. The generic
computation of $T_p$ uses only finitely many such labels, first for
its coordinate box and then for its separating weights. Applying
this to $E$ makes all corresponding $T_p$ constant on a dense
open. Repeating on the proper closed complement proves (1) since
the minimum of $T_p$ at $a$ is the required value at every $a$.

For (2), apply (1) to the prescribed rays and their sums and select
the strata with the reference values of these polytopes. On any
such stratum the same fan passes the test in
Proposition~\ref{prop:discover-source-fan}. Every $Q_a$ then has
the same support function as the reference on the rays and hence
on all weights. The converse follows from the defining duality.
For (3), the hypothesis on the generic fibers ensures that the
procedure terminates at the generic point. Its finitely many
determining polytopes are constant on a dense open by the preceding
argument. Induction on dimension then applies to each irreducible
component of the complement. Over an effective field we can find a
section with the required negative least exponent by solving the
first-star linear equations on successively larger lattice boxes.
\end{proof}

The same generic comparison provides a uniform reduction statement.
For a seed over $\mathbb Q$, spread its finite first star and the
section chosen above over a localization of $\mathbb Z$.
For fixed finite $E$, the generic computation uses finitely many
matrix ranks, hence survives after deleting finitely many primes.
At every remaining prime of good reduction its values agree for all
labels at all weights in $E$. If the characteristic-zero fiber
has a finite common refinement, use its determining list to obtain
simultaneous agreement of all $Q_a$ in these characteristics.

\begin{corollary}\label{cor:source-family-symmetry-probes}
The reference-fiber test also applies to a finite list of cones whose
images cover all directions under specified symmetries. More precisely,
suppose each of the fibers being compared has upper-algebra
automorphisms $\phi_g$ satisfying
$\nu_p(\phi_g h)=\nu_{g(p)}(h)\ (h\ne0)$
for the same maps $g:L_{\mathbb R}\to L_{\mathbb R}$. Suppose $g$
is given around $\omega$ by a unimodular integral linear map
$D_g$, and $g(\omega)$ lies inside the negative vertex cone
used on that fiber. Let $\sigma_1,\ldots,\sigma_b$ be pointed
rational cones on which the reference fiber passes the test of
Proposition~\ref{prop:discover-source-fan}. If the sets
$g(\sigma_i)$ cover $L_{\mathbb R}$, comparison at the original
cones' rays and ray sums determines equality of all the $Q_a$.
The collection of images may be infinite.
\end{corollary}
\begin{proof}
On each fiber use a normalized source basis $(f_a)$ with polytopes
$Q_a$. Each $f_a$ has unique least vertex $a$ on the interior
of the chosen negative vertex cone. Thus the least exponent of
$\phi_g f_a$ at $\omega$ is $D_g^t a$, by taking support
minima on a neighborhood of $\omega$. The transported basis is
independent for every Gauss valuation by the valuation identity in the statement. Consequently $\Newt(\phi_g f_a)=Q_{D_g^t a}$. For the
reverse of the minimality inclusion, expand a normalized section
with polytope $Q_{D_g^t a}$ in the transported basis. The distinct
least exponents show that $\phi_g f_a$ occurs, and independence at all
weights gives the opposite polytope inclusion. It follows that
$\min_{Q_a}\langle -,g(p)\rangle =\min_{Q_{D_g^t a}}\langle -,p\rangle$.
The finite comparisons first give equality at every label for $p$
in the listed cones, as in the theorem. This identity then gives
equality on every image. For a sequence of seed mutations inducing an
automorphism of the marked upper algebra, the valuation identity follows
from Lemma~\ref{lem:boundary-gauss-change}.
\end{proof}

On an irreducible parameter variety over an uncountable algebraically
closed field, all $Q_a$ agree with the geometric generic ones at
very general points even without a common refinement. Indeed (1)
gives a dense open with the generic value of each $T_p$, and
there are countably many integral $p$. In (3) over such a field
it suffices to assume a finite common refinement exists at every
closed point. For each fixed rational fan its locus is constructible
by (1) and Proposition~\ref{prop:discover-source-fan}. There are
countably many rational fans, and no irreducible subvariety is covered
on its closed points by countably many proper closed subsets. Thus
on each such subvariety at least one fan locus contains a dense open,
giving the hypothesis on its geometric generic fiber. The same
argument with the countably many bounded supports in the first adjacent charts descends
the existence of such a section from closed points. There is also exact
constancy along every orbit of diagonal translations of the marked
torus because multiplying monomial coefficients by nonzero character
values changes neither their exponents nor the order at $\omega$,
and the translation identifies the two full upper algebras.

Even without this finiteness hypothesis one can still compare the dimensions of all spaces of sections. Fix a reference fiber $s$ and write $V_\Delta(t)$ for
the regular functions on fiber $t$ supported in the lattice polytope $\Delta$.
Restrict to a locally closed parameter set on which the ordinary
exchange supports are fixed. The locus
\[ \{t:\dim V_\Delta(t)\ge\dim V_\Delta(s) \text{ for every lattice polytope }\Delta\} \]
is closed. For each fixed polytope, write the slice coefficients as
the exchange powers times bounded auxiliary quotients. The support
bounds are uniform on this parameter set. Multiplication by each
nonzero exchange power is injective, so the auxiliaries are unique
for a given section. Thus the dimension is the kernel dimension of
a finite matrix, and the intersection of the rank loci is closed. By
Noetherianity finitely many of its rank conditions suffice, without
an a priori bound on the polytopes required. Equality of dimensions for
every polytope is equivalent to equality of every labeled $Q_a$.
Indeed when $a$ is the $\omega$-least vertex of $\Delta$,
subtracting the dimension for the hull of
$(\Delta\cap M)\setminus\{a\}$ tests the existence of a section
with nonzero coefficient at $a$. We use dimension zero for the empty support. Apply this with the smallest polytope on each of the two
fibers to get both inclusions. The reverse dimension inequalities
are open for a single polytope. Their infinite intersection does
not acquire a finite determining set from Noetherianity alone.

\subsubsection{Infinitely many normal directions}

By Corollary~\ref{cor:reachable-monomials-minimum-basis}, every
reachable seed variable has marked Newton polytope $Q_a$ for its
least exponent $a$. If the $Q_a$
had a common finite refinement, all reachable coordinate valuations
would be linear on its cones. The same would hold in any fixed
reachable marking, since the tropical map from that marking to the
initial one and its inverse are piecewise integral linear on finite
fans by Lemma~\ref{lem:boundary-gauss-change}.

\begin{theorem}\label{thm:mutual-degrees-infinite-fan}
At a reachable seed of a full upper LP algebra with Laurent
coefficients over a field, choose two mutable variables $x,y$.
Let $d>0$ be the degree in $y$ of the ordinary exchange at $x$,
and $e>0$ the degree in $x$ of the exchange at $y$. If
$de\ge4$, there is no finite common refinement of the normal fans
of the smallest source polytopes $Q_a$, at any separating weight.
In the marked counting setting of Theorem~\ref{thm:seed-windows},
the reciprocal polytopes $P_a$ likewise have no finite common
refinement or finite common list of defining normals. In this
counting setting there is also no finite cover of the full label space by rational cones on
each of which at least one of the two families $P,Q$ is
Minkowski additive.
\end{theorem}
\begin{proof}
Keep the other seed coordinates fixed. Over their rational function
field, the two ordinary exchanges are univariate polynomials with
nonzero constant and leading coefficients. The exchange at $x$
has no LP normalizing factor in $y$. Such a factor would require
the exchange at $y$ to divide the constant coefficient in $y$
of the exchange at $x$, by reducing the LP normalization test
modulo the exchange at $y$ \cite[Lemma~2.7 and \S2.2]{LP}.
A positive-degree polynomial in
$x$ with nonzero constant term cannot divide a nonzero scalar in
the Laurent polynomial ring. The same applies in the other
direction. On mutating $y$, the LP rule replaces the opposite
exchange $F(y)$, up to a nonzero scalar of this field, by
$z^d F(c/z)$, where $z$ is the new coordinate and $c\ne0$
is the associated exchange at $y$ evaluated at $x=0$.
All factors removed in this substitution divide the scalar $c$,
and the factor $z^d$ gives a nonzero constant term. The exchange at
the mutated index stays the same. Thus the degrees are preserved and
the leading and constant coefficients stay nonzero along alternating mutations.

On the plane where the other coordinate weights vanish, write
$h_0,h_1$ for the initial two weights and $h_j$ for the
successive coordinate values in this sequence. Laurent Gauss
transport gives
\[ h_{j+1}=b_j\min(0,h_j)-h_{j-1},\qquad b_j\in\{d,e\}\text{ alternating}. \]
Only the leading and constant coefficients enter this minimum. Each step
is a piecewise unimodular homeomorphism of the plane. For $k\ge2$,
take the inverse image of the ray generated by
$(h_{k-1},h_k)=(-1,0)$. At that point all previous $h_j$,
$0\le j<k$, are strictly negative. To see this, start with
$q_0=0,q_1=1$, extend the alternating coefficients periodically, and set
$q_{l+1}=b_{k-l}q_l-q_{l-1}$, so that while the terms are
positive they are exactly $-h_{k-l}$. The even and odd
subsequences satisfy
\[ q_{l+4}=(de-2)q_{l+2}-q_l, \qquad q_2=b_{k-1}>0,\quad q_3=de-1>1. \]
Both subsequences are strictly increasing because $de-2\ge2$.
This proves the signs, including the case $de=4$.

Near this ray the map from $(h_0,h_1)$ to $(h_{k-1},h_k)$
is consequently linear and invertible. The linear function $h_k$
changes sign across the ray, producing distinct linear expressions
for $h_{k+1}$ on its two sides. These rays are distinct as $k$
varies, since on the $k$-th ray the first zero occurs at $h_k$.
Restricting a finite fan to the plane gives only finitely many rays,
which contradicts the necessary condition above. In the counting
setting the Newton polytopes of all these variables are also the
corresponding $P_a$ by Theorem~\ref{thm:seed-windows}. A fixed
finite set of defining normals permits only finitely many normal
cones, each generated by the subset active on a face, and hence
would give a finite common refinement.

Finally suppose such a finite cover of label space existed. Choose
finitely many integral generators for each cone and take the common
refinement of the normal fans of the $P$- or $Q$-polytopes at
its generators, using the additive family on each cone. At every
label of a reachable seed variable the two polytopes agree.
Each such polytope is therefore a nonnegative Minkowski combination
of the chosen polytopes on one cone of the cover, by
Proposition~\ref{prop:single-sum-realization-cone}. Its support function
is linear on every cone of the common refinement. The infinitely
many distinct bends just proved exclude this possibility.
\end{proof}

The smooth family following Theorem~\ref{thm:nonlinear-simplex-reciprocity}
contains examples with both mutual degrees at least three in every
mutable rank at least two. Thus equality $P_a=Q_a$ for all labels
is compatible with infinitely many distinct normal directions.

\subsection{Arithmetic obstructions and marked return}

Work over an algebraically closed field $k$ of characteristic zero,
and regard $\mathbb Q_{>0}$ as a subset of $k^\times$.
Let $U$ be a marked source satisfying the hypotheses
of Proposition~\ref{prop:negative-vertex-chart}. Recall that marked
return holds if the forward mirror $Y$ is smooth and its augmented
counting algebra equals $k[U]$ in the Laurent chart of
Theorem~\ref{thm:automatic-reverse-chart}. When comparing markings we
call return in the given marking \emph{strict return}.

\begin{proposition}\label{prop:arithmetic-marked-return}
Suppose $k[U]$ is the full upper LP algebra of the marked seed with
Laurent coefficients. Let $h$ be a monomial
in a reachable seed with nonnegative mutable exponents and arbitrary
integral powers of the designated units. Let $a$ be the slope of
$\tau_h$ on the forward straight cone and $c_h$ the coefficient of
$x^a$ in the initial chart. If marked return holds, then
$c_h^{-1}h=\Phi_a$, so the normalized monomial has nonnegative
integral Laurent coefficients in the initial chart.

Write the original ordinary exchanges as
$F_i=\sum_{\alpha\in A_i}c_{i\alpha}x^\alpha$, using the whole
character lattice $M$ of the initial torus.
If marked return holds after some translation of the initial torus,
then every integral family $(n_{i\alpha})$ satisfying
\[ \sum_{\alpha\in A_i} n_{i\alpha}=0\quad\text{for each }i, \qquad \sum_{i,\alpha} n_{i\alpha}\alpha=0 \]
satisfies
$\prod_{i,\alpha} c_{i\alpha}^{\,n_{i\alpha}} \in\mathbb Q_{>0}$.
\end{proposition}
\begin{proof}
By Corollary~\ref{cor:source-seed-slopes}, $a$ lies in the closed
cone of the chosen seed. Theorem~\ref{thm:seed-windows} identifies
$c_h^{-1}h$ as the unique regular function supported in $P_a$ with
coefficient one at $x^a$. Theorem~\ref{thm:automatic-reverse-chart}
gives the same support bound and coefficient for $\Phi_a$, and marked
return makes $\Phi_a$ regular on $U$. Hence $c_h^{-1}h=\Phi_a$,
whose Laurent coefficients are nonnegative integers.

The first mutation partner $\widehat F_i/x_i$ is a scalar times a
Laurent monomial times $F_i$, so the ratios of coefficients within
one exchange are the same for both. Translating $x$ by a
torus point $t$ multiplies $c_{i\alpha}$ by $t^\alpha$, up to
a common factor for that exchange. The two relations on $n$ cancel
all these factors. In a translated marking with return, the product
can therefore be evaluated using the nonzero integral coefficients
of the normalized first partners, all of which are positive.
\end{proof}

In an irreducible complex parameter variety with fixed supports and
nonvanishing coefficients, suppose the sources satisfy the stated
forward hypotheses on a nonempty open set. If one of the indicated
products is a nonconstant rational function on the parameter variety,
the points of this open set admitting return after torus translation
lie in a countable union of proper closed subsets. Indeed, each
positive rational value cuts out a proper closed subset on the locus
of definition. For fixed polynomial-corner supports,
Theorem~\ref{thm:finite-return-coefficients} determines the finite
set of joint values of the coefficient invariants for which translated
return occurs.

\subsubsection{Smooth examples with no marked return}

\begin{theorem}\label{thm:smooth-both-sides-nonreturn}
In every mutable rank $r\ge2$ there are smooth affine full upper
LP sources with a polynomial corner whose augmented forward mirror is
smooth and for which every reciprocal system has a unique normalized
solution, but marked return fails after every diagonal translation of
the initial torus. The obstruction to marked return persists on taking
products, in the product marking, with any other polynomial-corner source
satisfying the hypotheses of Theorem~\ref{thm:polynomial-corner-reciprocity}
whose forward mirror is smooth and has the required volume form.
\end{theorem}
\begin{proof}
Work over $\mathbb C$. Take a designated unit $u$, let
$A=1+u$ and $B=1+\lambda u$, and choose $\lambda>0$
irrational. The two ordinary exchanges at $x,y$ are
$AB+y$ and $A+x$. They are irreducible, being monic and
linear in the opposite variable. The consecutive variables
$x,y,t,w,v$ have the five relations
\begin{equation}\label{eq:five-source-exchanges}
xt=y+AB,\qquad yw=t+B,\qquad tv=1+Aw,\qquad wx=1+Bv,\qquad vy=x+A .
\end{equation}
For instance the LP substitution after the first mutation gives
$A+AB/t$, whose common factor $A$ is removed, giving $t+B$.
At the succeeding steps, the other exchange of the pair is obtained
by reciprocating its linear polynomial at the constant term of the
mutated exchange and removing common factors. This gives, in order,
$1+Aw,1+Bv,x+A,y+AB$. Thus each of the five consecutive pairs is a seed with a polynomial corner.

We show that these relations present the full upper LP algebra and that
its spectrum is smooth.
The five coordinates have no common zero, since $tv=1+Aw$. On inverting one coordinate
the relations eliminate its two nonneighbors and give the sole equation
that the product of its two neighbors equals its exchange expression
above. More explicitly the equations on inverting $x,y,t,w,v$,
respectively, are
\[ yv=x+A,\quad xt=y+AB,\quad yw=t+B,\quad tv=1+Aw,\quad wx=1+Bv . \]
Substitution verifies the other relations in each localization.
Each open is integral and smooth. In the first three equations the
coefficient of the inverted coordinate is one. In the last two it is
nonzero wherever the right-hand side vanishes. The initial torus is
dense in each open. The locus $x=y=0$ has $A=0$, $t=-B\ne0$,
$v=-1/B$ and $w$ free, so has codimension two. On $x\ne0$
or $y\ne0$ the initial and adjacent tori miss only the locus
where both opposite variables in the corresponding conic vanish,
also of codimension two. The generators are Laurent in the full
mutation class. Normality and the first-star intersection give the
claim, and the torus volume extends by the same covering.

We compute the forward ring. Denote the five cyclic source
coordinates by $(s_1,\ldots,s_5)=(x,y,t,w,v)$, taking subscripts modulo five.
In each consecutive seed take the standard dual coordinates
from Theorem~\ref{thm:corner}. A shared coordinate has the same
dual theta label on both sides of a mutation. Indeed its positive
coordinate ray and the positive $u$-ray have zero tropical
exchange value there, since each of the linear exchanges has
nonzero constant term at the polynomial corner. Write $C_i$
for the successive mutable dual coordinates and $Z_+$ for the
common positive unit-coordinate theta. After mutation the new
positive mutable ray is the old negative axis. There the only
possible correction is in the unchanged mutable direction, with
degree one. With indices taken cyclically, the axis formula gives
\begin{equation}\label{eq:five-mirror-exchanges}
C_i C_{i+2}=1+C_{i+1}.
\end{equation}

At any nonnegative order of $u$, both $A$ and $B$ have
value zero. At that order each tropical exchange is
$(p,q)\mapsto(q,\min(0,q)-p)$, with the pairs in cyclic order.
Thus the five nonnegative mutable cones cover $\mathbb R^2$. Their successive rays in the original label marking are
$(1,0),(0,1),(-1,0),(-1,-1),(0,-1)$.
Every label with last entry $m\ge0$ is in the nonnegative
orthant of one of these five markings. Its theta function is
therefore a monomial in the corresponding pair of $C_i$ times
$Z_+^m$. The fixed-unit-exponent rule grades the whole counting
algebra by $m$. Its nonnegative part is $D[Z_+]$, where
$D=\mathbb C[C_1,\ldots,C_5]$ in the initial two mutable Laurent
coordinates.

Use instead $u'=u^{-1}$ and replace $(x,y,t,w,v)$ by
$(x/u,y/u^2,t/u,uw,uv)$. This is an invertible change in each
marked source torus. The exchange relations~\eqref{eq:five-source-exchanges} have the same form with
$A'=u'+1$, $B'=u'+\lambda$. Each ordinary exchange still has
a polynomial corner. Let $Z_-$ be the theta function on the positive
$u'$-coordinate ray of this marking. The same argument gives the
nonpositive part $D[Z_-]$, with the same degree-zero coordinates, since at unit
order zero the two mutable markings agree. The positive $u'$-ray
of the new $(x/u,y/u^2,u')$-chart has old label $(-1,-2,-1)$.
In the $(w,v,u)$-chart this label is $(1,1,-1)$. Equivalently the rescaled $(uw,uv,u')$-coordinates have
orders $(0,0,1)$. The mutable exchanges at $w,v$ are
$1+Bv$ and $1+Aw$. Both have tropical value zero at
$(1,1,-1)$. The corner support bound thus gives exactly
$Z_-=C_4 C_5/Z_+$. The forward ring is
$D[Z_+,Z_-]/(Z_+Z_--C_4C_5)$.
There are no additional relations, by reduction to the separate
positive and negative powers in the Laurent field.

The mirror relations~\eqref{eq:five-mirror-exchanges} rule out a common zero of the $C_j$. Where
$C_j\ne0$, they identify the coordinate ring with the polynomial
ring in the two neighbors $C_{j-1},C_{j+1}$ localized at
$C_{j-1}C_{j+1}-1=C_j$. The neighbors are algebraically independent
by the birational recurrence, so $\Spec D$ is smooth. On $C_j=0$
the relations give $C_{j-1}=C_{j+1}=-1$. The whole divisor
therefore lies in the chart $C_{j+1}\ne0$, where it is the affine
line with equation $C_j=0$. In particular $C_4=0$ and $C_5=0$
are disjoint prime divisors. Their
product cuts out a smooth reduced divisor, so the forward
threefold is smooth, including where $Z_+=Z_-=0$. The logarithmic
area form on the surface is regular and nowhere zero by the torus
charts in codimension one (consecutive zeros are impossible).
Wedge it with $d\log Z_+$ on $Z_+\ne0$. On $Z_-\ne0$
the resulting form equals its wedge with $-d\log Z_-$. These opens cover in
codimension one. Thus the forward mirror is itself a smooth counting
source with the prescribed torus volume.

In the two original exchanges the constant coefficients are one and
the coefficients of $u$ are $1+\lambda$ and $1$, respectively.
The ratio of the two $u$-coefficients, each divided by its constant
coefficient, is invariant under torus translation and exchange
rescaling. Its value $1+\lambda$ is irrational. By
Proposition~\ref{prop:arithmetic-marked-return}, no diagonal
translation gives marked return.

Finally, disjoint seeds in independent variables give the product
source and the tensor product of the augmented forward algebras,
with product basis. For the source assertion, the generators of both
coordinate rings are Laurent on the reachable product tori, while
the initial and first adjacent product charts cover in codimension
one. For the basis assertion, project the contact orders to the
two factors and use the same generic marked line in both. The
ordinary incidence equations for their monic polynomials separate,
with independent evaluation points. No roots in different factors
coincide generically. After fixing one factor, its finitely many
roots are unmarked points of the fixed domain for the other. Evaluation
there is dominant on every counting component by the fixed-contact
deformation argument in Theorem~\ref{thm:incidence}. The same
argument at a contact gives generic residues in its divisor orbit,
so a projected positive mutable axis satisfies its exchange
nonvanishing condition even when the combined contact is not an
axis. Projected zero contacts give interior marks. Thus the
degrees of traces of two or three theta functions multiply over the
common parameter field. The tensor-product multiplication consequently
satisfies the trace equations of the product, whose unitriangular
pairing determines every coefficient. The marked realizations agree
by translation into the product of the straight cones.
The coefficient invariant already used involves only the first two
exchanges and persists under products. To see that every reciprocal
system in the example has a unique normalized solution, compute the
slopes of the original regular functions on the strictly positive
orthant. There $A,B$ both have order zero, so the slopes of
$x,y,t,w,v$ have zero unit entries and successive mutable entries
\[ (1,0),\quad(0,1),\quad(-1,0),\quad(-1,-1),\quad(0,-1). \]
Consecutive pairs, together with both signs of the unit slope, thus
cover the full space of weights. Theorem~\ref{thm:seed-windows}
identifies the unique normalized solution at each lattice point with
the corresponding normalized seed monomial.

To raise the mutable rank
by one, take the independent factor $ss'=1+\eta$ with $\eta$ a designated unit.
It is smooth, and its forward algebra by the axis formula is
$k[S^{\pm1},V,W]/(VW-1-S)$, also smooth and volume-compatible.
Its two seeds have opposite mutable slopes on the positive orthant
and both allow arbitrary unit slopes. The cones of product seeds
are Cartesian products, so these factors preserve the complete
cover by seed cones. Iterating gives all the stated ranks.
\end{proof}

The forward threefold $Y$ in the rank-two example also has
$\Pic(Y)\simeq\mathbb Z$. To see this, first invert $C_1C_2$
in $D$. The cyclic relations give
$D[(C_1C_2)^{-1}]=\mathbb C[C_1^{\pm1},C_2^{\pm1}]$.
We have seen that both deleted divisors are prime and reduced. They
are principal, so the localization sequence gives $\Cl(D)=0$.

Now $\mathbb C[Y][Z_+^{-1}]=D[Z_+^{\pm1}]$ is factorial. The deleted prime
divisors are $P_j=V(Z_+,C_j)$ for $j=4,5$, each of multiplicity
one in $\operatorname{div}(Z_+)$. A rational function with divisor
supported on $P_4\cup P_5$ is a unit of $D[Z_+^{\pm1}]$, hence
has the form $dZ_+^n$ with $d\in D^\times$. Thus $P_4+P_5$
generates all relations. Smoothness gives
\[ \Pic(Y)=\Cl(Y)\simeq (\mathbb Z P_4\oplus\mathbb Z P_5)/\mathbb Z(P_4+P_5) \simeq\mathbb Z. \]
Consequently $Y$ admits no presentation as the spectrum of a full
upper LP algebra over a Laurent polynomial coefficient ring, by
Proposition~\ref{prop:upper-factorial}.

\subsubsection{The translation criterion}

Keep the smooth forward mirror and notation of
Theorem~\ref{thm:first-reverse-partner-return}. Use the actual
associated exchange $\widehat F_i$ and the reverse numerator
$H_i^*=x_i\Phi_{a_i}$ in the same Laurent ring. Thus the monomial
normalizer of the original LP exchange is included in the comparison.
For $t\in\operatorname{Hom}(M,k^\times)$, let $t$ act by
$x^\alpha\mapsto t^\alpha x^\alpha$, where $t^\alpha=t(\alpha)$.

\begin{proposition}\label{prop:torus-translate-return}
There is a torus translation of the source in whose marking return
holds if and only if the following two conditions hold.
\begin{enumerate}
\item For every $i$, the supports of $H_i^*$ and $\widehat F_i$
coincide, say with common support $A_i$.
\item Write $h_{i\alpha}$ and $c_{i\alpha}$ for their respective
coefficients on $A_i$. For a $\mathbb Z$-basis of the lattice
\[ \mathcal L=\left\{n\in\bigoplus_i\mathbb Z^{A_i}: \sum_{\alpha\in A_i} n_{i\alpha}=0\ \text{for all }i,\quad \sum_{i,\alpha}n_{i\alpha}\alpha=0\right\} \]
one has $\prod_{i,\alpha}h_{i\alpha}^{\,n_{i\alpha}}
=\prod_{i,\alpha}c_{i\alpha}^{\,n_{i\alpha}}$.
\end{enumerate}
When nonempty, the set of translations giving return is a coset of
$\operatorname{Hom}(M/D,k^\times)$, where $D\subset M$ is the
lattice generated by differences of two exponents in the same $A_i$.
This coset has codimension $\operatorname{rank}D$ in the translation torus.
\end{proposition}

\begin{proof}
A translation identifies the source with its translate, preserves
the logarithmic volume, and fixes the cocharacter labels. Transport
of the generic evaluation point and boundary model identifies its
labeled counts, as in the discussion following
Corollary~\ref{cor:diagonal-exchange-stabilizer}. Thus the forward
algebra in its theta marking is the same for each translate.
The reverse counts therefore come from the same marked forward
mirror. Theorem~\ref{thm:automatic-reverse-chart} determines the
Laurent realization by its characters on the straight cone, so each
$H_i^*$ has the same Laurent expansion for every translate.

By Theorem~\ref{thm:first-reverse-partner-return}, equality in the
translated marking is equivalent to
$H_i^*=\gamma_i\,t(\widehat F_i)\ \text{for all }i$
for nonzero scalars $\gamma_i$. In fact the leading coefficient
of each reverse partner is one, which fixes the normalization.
Translation of the denominator $x_i$ contributes only a scalar.
The coefficient equalities in the statement are therefore necessary.

Conversely the ratios $h_{i\alpha}/c_{i\alpha}$ give a character
of $\bigoplus_i\mathbb Z^{A_i}$ trivial on $\mathcal L$. It
descends to the image of the map to $\mathbb Z^r\oplus M$
which sends the generator at $(i,\alpha)$ to $(e_i,\alpha)$.
Over an algebraically closed field a character of this image extends
to the whole lattice, by Smith normal form and existence of roots
in $k^\times$. Its values supply the required $\gamma_i$ and
$t$. Two choices of $t$ differ by a character trivial on $D$,
and every such character preserves the proportionalities. This
proves the assertion about the coset. If specified characters must
remain fixed, use their quotient lattice in place of $M$ in the
character test, which amounts to requiring $t$ to be one on
their exponents.
\end{proof}

In positive mutable rank this stabilizer is proper, since an ordinary
exchange is a nonunit in the Laurent ring of unchanged coordinates
and therefore has two distinct exponents. Realization of a label is
invariant along the translation orbit. The forward polytopes are fixed,
and translation preserves the support of every regular function on $U$.
Dividing the translated function by its leading coefficient restores
the normalization. In particular, if one marking
returns, every translate realizes every label, but a general translate
fails to return in its prescribed marking. Smoothness of the source and its forward mirror is preserved along the
orbit.

The rationality condition of Proposition~\ref{prop:arithmetic-marked-return}
also characterizes when the first partners can simultaneously have
positive integral normalized coefficients after unrestricted torus
translation. Write $A_i=\Supp(\widehat F_i)$ and $c_{i\alpha}$
for its coefficient at $\alpha\in A_i$. Choose reference exponents
$\alpha_i^0$ at one common interior weight of the forward straight
cone. Let $E$ consist of the pairs $e=(i,\alpha)$ with
$\alpha\in A_i\setminus\{\alpha_i^0\}$. Write
$r_e=c_{i\alpha}/c_{i\alpha_i^0}$ and $\beta_e=\alpha-\alpha_i^0$.
On the kernel of the map $\mathbb Z^E\to M$ sending the basis
element indexed by $e$ to $\beta_e$, the
character defined by the $r_e$ takes values in $\mathbb Q_{>0}$.
This kernel is a direct summand since its quotient is a lattice.
Extend its character to $\mathbb Z^E$ with positive rational
values $s_e$ on the generators. All $\beta_e$ pair positively
with some integral covector $v$, since the chosen reference exponents
are strict common minima. Multiplying $s_e$ by
$N^{\langle v,\beta_e\rangle}$ for a common sufficiently divisible
positive integer $N$ makes every $s_e$ integral and preserves the
character on the kernel. Now $s_e/r_e$ descends to the span of
the $\beta_e$ and extends to a torus point, as above. This gives
the required translation. Actual return further requires the
specified reverse coefficients of Proposition~\ref{prop:torus-translate-return}.

\subsubsection{The finite set of returning coefficient tuples}

Fix polynomial-corner supports $A_i$ as above. Divide each
representative $G_i$ by its nonzero constant coefficient. Write
$c=(c_{i\alpha})_{\alpha\in A_i\setminus\{0\}}$ for the remaining
coefficients and $\pi(c)$ for the complete tuple of their character
invariants for diagonal translation, given by a lattice basis as in
Proposition~\ref{prop:torus-translate-return}.

\begin{theorem}\label{thm:finite-return-coefficients}
For these supports there is a computable finite set $S$ of
coefficient tuples for the normalized ordinary exchanges. Every entry
of every tuple in $S$ is a positive integer.
Let $U$ be a source with these supports satisfying the hypotheses of
Theorem~\ref{thm:polynomial-corner-reciprocity}, and let $Y$ be its
augmented forward mirror. Then
\begin{enumerate}
\item $c\in S$ if and only if $Y$ is smooth and strict marked return holds.
\item $\pi(c)\in\pi(S)$ if and only if $Y$ is smooth and marked return
holds after a diagonal translation.
\end{enumerate}
Every tuple in $S$ occurs for some such $U$.

In an algebraic family satisfying the same hypotheses on each fiber,
the fibers with smooth forward mirror and strict marked return
form a closed subset of the base. The same holds for marked return
after a diagonal translation. If a coefficient invariant is
nonconstant on an irreducible base, the latter subset is proper.
\end{theorem}
\begin{proof}
First suppose $Y$ is smooth. Write
$E_i=\operatorname{conv}(A_i)\cap M$ and
$\Delta_i=\operatorname{conv}(A_i)-e_i$. The ordinary partner
$G_i/x_i$ is the first mutation partner multiplied by a Laurent
unit in the designated coordinates and a nonnegative monomial in
unchanged mutable variables, up to a scalar. It is thus a monomial
in the adjacent seed, normalized with least exponent $-e_i$.
Theorems~\ref{thm:seed-windows} and
\ref{thm:automatic-reverse-chart} give
$\Supp\Phi_{-e_i}\subseteq P_{-e_i}=\Delta_i$.
By Proposition~\ref{prop:all-reverse-seed-cones}, the same
multiplication by the unchanged characters carries the normalized
actual reverse partner to $\Phi_{-e_i}$. Hence
Theorem~\ref{thm:first-reverse-partner-return} gives
\[ \text{strict return}\quad\Longleftrightarrow\quad G_i=x_i\Phi_{-e_i}\quad(1\le i\le r). \]
Set $c_{i0}=1$ and $c_{i\alpha}=0$ for $\alpha\in E_i\setminus A_i$.
The normalized coefficients give strict return if and only if
\[ \sum_{\alpha\in E_i}c_{i\alpha}\,\mathsf T_Y(\alpha,-\beta) =\mathsf T_Y(e_i,-e_i,-\beta) \qquad (\beta\in E_i,\ 1\le i\le r). \]
Indeed translating $\Delta_i$ by $e_i$ sends every exponent into
the reverse cone $K=D_+^\vee$ of Theorem~\ref{thm:automatic-reverse-chart}.
The proof of Corollary~\ref{cor:reverse-finite-window-traces}
gives this unitriangular system on $E_i$. Its unique solution gives
the nonnegative integers in $x_i\Phi_{-e_i}$, with constant one.
The contact lists and matrix sizes depend only on the supports.
Proposition~\ref{prop:laurent-generator-traces} computes their entries
from the marked generators of $Y$. For the equation indexed by $\beta$, the monic
polynomials needed have degrees at most $\beta_j$ for coordinate $j$,
using the regular basis $z_1,\ldots,z_n$.

To enumerate candidate forward rings, use the contact lists and
uniform trace bounds in the proof of
Corollary~\ref{cor:constructible-reciprocal-locus}. For each integral
trace tuple in these ranges, compute the Laurent generators by
Corollary~\ref{cor:trace-reconstruction} and adjoin the marked
characters $z_j$. Discard inconsistent systems. The resulting
list contains every marked forward ring with the given supports.

Elimination gives a presentation of each candidate. We need the
smooth candidates on which the torus volume form extends without
zeros or poles. Inverting the $z_j$ gives the torus. On each smooth
candidate, test whether
$dz_1\wedge\cdots\wedge dz_n$ is divisible by $\prod_j z_j$
in the module of top differentials and whether the quotient generates
that module. On each such ring compute the degrees of the reverse traces
and solve the specified finite systems. Discard a candidate unless
these systems have unique solutions giving normalized, positive
integral polynomials of exactly the prescribed supports.

Finally check each surviving tuple as an ordinary seed. Test
irreducibility over $\overline{\mathbb Q}[u^{\pm1}]$ in the
mutable variables. The regular function $\prod_i(G_i/x_i)$ has a
Newton vertex negative in every mutable coordinate.
Proposition~\ref{prop:source-polar-polytope} computes the dual
polytopes needed for a lattice box containing zero in its interior.
Theorem~\ref{thm:polyhedral-newton-sections} then gives the rational
cone of its section semigroup. Compute its Hilbert basis, lift those
generators by the finite first-star equations, and evaluate the degree
variable at one. This gives generators of the full upper algebra
before any smoothness test. Compute its presentation and test
smoothness. On a smooth
full upper source the first-star tori cover in codimension one by
Proposition~\ref{prop:upper-factorial}. Their logarithmic forms agree
up to sign, so the logarithmic volume form in the initial chart
trivializes the canonical bundle of the source.
Compute the forward Laurent generators of this source and check
smoothness and the partner equalities on its forward ring. Let $S$
consist of the tuples that pass these tests. The computations use
exact polynomial arithmetic over $\overline{\mathbb Q}$, and the
tests and counts persist over every algebraically closed extension
of characteristic zero.

The construction includes the tuple of every $U$ with smooth
forward mirror and strict marked return. Conversely, the tests produce
a variety with these properties for each tuple in $S$.
The normalized exchanges determine its marked full upper algebra,
so every $U$ with that tuple has the same properties.
This proves (1).

By Proposition~\ref{prop:torus-translate-return} and its lattice
calculation, two normalized coefficient tuples have the same
$\pi$-value precisely when they differ by a diagonal translation.
Thus $\pi(c)\in\pi(S)$ precisely when some translate of $c$
belongs to $S$. Translation preserves the marked forward algebra
and hence its smoothness. Applying (1) to the translated source
proves (2).

In a family the two loci are the inverse images of the finite sets
$S$ and $\pi(S)$ under the regular maps $c$ and $\pi\circ c$.
Both are closed. If a coefficient invariant is nonconstant on an
irreducible base, its image is infinite, so the translated return
locus is proper.
\end{proof}

\appendix

\section{Valuative independence in the marked Laurent chart}
\label{sec:one-sided-independence}

For a tropicalized marked curve, its \emph{spine} is the restriction
of its tropicalization to the subtree spanned by the marked ends
\cite[Proposition~4.24]{KY}. Keel and Yu's local functions count
curves with such a spine and prescribed final position and tangent.
We cut a complete spine near its central vertex and express
the intrinsic trace as the constant term of a product of local functions
(Lemma~\ref{lem:local-zero-trace}). This identity gives a homomorphism
from the intrinsic algebra to a ring of formal series. A vertex of a
regular function's Newton polytope that is negative in every mutable
coordinate gives a cone on which these series are the marked Laurent
expansions (Proposition~\ref{prop:formal-local-chart}). We use broken
paths to prove homogeneity in the label and valuative independence of
these expansions.

\subsection{Regular functions bound local spine expansions}
\label{subsec:local-spine-support}

Let $U$ be a smooth affine log Calabi--Yau variety over a field $k$
of characteristic zero with a volume-compatible dense split torus $T$.
We give $k$ the trivial absolute value when taking analytifications.
Choose a smooth projective compactification with simple normal crossings
boundary as in \cite[Assumption~2.4]{KY}. The closure of $U\setminus T$
contains no boundary stratum, and the cones of the essential skeleton
form a smooth toric fan. Write
$L$ for the cocharacter lattice of $T$ and use the torus to
identify $L_{\mathbb R}$ with the skeleton. Suppose that the
coordinates $x_1,\ldots,x_r,u_1,\ldots,u_c$ of $T$ extend to
regular functions on $U$ and that the $u_j$ are units there.
For $g\in k[U]\setminus\{0\}$, write
$\tau_g(y)=\min_{a\in\Supp(g|_T)}\langle a,y\rangle$.
Write $Q$ for the monoid of effective numerical curve classes of
the compactification. For $m\in L\setminus\{0\}$ and a generic
point $\xi\in L_{\mathbb R}$, let $V_{\xi,m}$ be the
\emph{local one-input series} of \cite[Definition~20.20]{KY}.
The coefficient of $q^\beta z^e$ counts curves of class
$\beta\in Q$ whose spines have incoming derivative $-m$ and
final derivative $-e$ at $\xi$, where $e\in L$. Set $V_{\xi,0}=1$.
We may choose $\xi$ with coordinates linearly independent over
$\mathbb Q$. An \emph{infinitesimal cylinder} here is counted using
a sufficiently short interval crossing a rational hyperplane at a
general point, with prescribed derivatives on its two sides
\cite[Definition~20.2]{KY}.

\begin{lemma}\label{lem:local-regular-slope}
Suppose that $q^\beta z^e$ occurs in $V_{\xi,m}$, and that
$\tau_g$ is linear near $\xi$, with gradient $a$. Then $\langle a,e\rangle\geq\tau_g(m)$.
For an infinitesimal cylinder of nonzero count with incoming and
outgoing derivatives $-v,-w$, if its two ends lie in the same open
linearity domain of $\tau_g$ with gradient $a$, then
$\langle a,w-v\rangle\geq0$.
\end{lemma}

\begin{proof}
The assertion for $m=0$ is immediate. For $m\ne0$ the last
derivative of an interval ending at a generic point is nonzero. Otherwise the last change to zero slope would be a bend on a wall
at $\xi$. Both end directions in Definition~20.2 are also nonzero.
Use the torically extended representative of the counted body from
\cite[Construction~9.3, Definition~9.6 and Theorem~11.9]{KY}. Its two nonzero contacts lie at the ends of the extended interval.
Any additional evaluation marks have order zero. One can compute the nonzero count over a
skeletal parameter on a smooth marked rational curve by
\cite[Theorem~8.18 and Lemma~8.21]{KY}. On the interval between its
two contact marks, the skeletal map takes values on the target
skeleton. Thus the valuation of the pullback of $g$ restricts to
$\tau_g\circ h$ on that interval, where $h$ is the spine map.
The pullback is a nonzero rational function, since its value at a
finite skeletal point has finite valuation.

Orient the interval from the incoming end to the outgoing end.
Both contact marks lie on this path. Every branch off the path
therefore misses the compactification boundary, so the pullback of $g$
has no poles on such a branch. After passing
to an algebraically closed valued extension, the slope of its
valuation toward an off-path branch is the number of zeros on that
branch, with multiplicity. Indeed one may factor the rational
function on the smooth rational curve into linear factors and apply
the slope rule separately to each zero and pole. Balancing at each
point of the path shows that the derivative along the oriented path
can only decrease. Equivalently this is the principal-divisor slope
calculation of \cite[Lemma~15.6]{KY}, with divisor
$-\operatorname{div}(g)$.

Far enough toward the infinite incoming end its derivative is
$-\tau_g(m)$. Indeed, along $h(t)=h_0-tm$ with $t\to-\infty$, only
the exponents minimizing $\langle a',m\rangle$ can attain the
minimum. Near the finite endpoint it is $-\langle a,e\rangle$.
For a short cylinder the two derivatives are instead
$-\langle a,v\rangle$ and $-\langle a,w\rangle$. Comparing
them proves both assertions.
\end{proof}

\begin{proposition}\label{prop:local-bend-cone}
For a cylinder with derivatives $-v,-w$ as above, call $w-v$
its \emph{bend}. All bends of infinitesimal cylinders of nonzero count belong to
$P_0=\mathbb N^r\times\{0\}^c$. The exponents of $V_{\xi,m}$
belong to $m+P_0$. In particular the pointed bend condition of
\cite[Assumption~20.27]{KY} holds, and in every bounded total mutable
degree only finitely many effective classes contribute and their characters
$q^\beta$ can be specialized to one.
\end{proposition}
\begin{proof}
Apply Lemma~\ref{lem:local-regular-slope} to $x_i$ and to both
$u_j$ and $u_j^{-1}$. All these functions have linear
tropicalization on the entire skeleton. The positive orthant in $L$ gives the full-dimensional monoid in
that assumption. On $P_0$ its monomial-adic degree is the total mutable degree.
A \emph{twig} is the restriction of the tropicalized curve to the
closure of a connected component off its spine
\cite[Definition~4.22]{KY}. We bound the degrees of the twigs using
the characters of $T$.
On the open set where all $x_i$ are invertible, the coordinates
define a morphism to $T$. Its composite with $T\hookrightarrow U$
agrees with the inclusion on the dense torus, hence on the whole open set. It follows that $U\setminus T=V(\prod_i x_i)$.
For a counted extended two-contact map $f$, work on its smooth domain
curve over an algebraically closed valued extension. At a point $b$ away
from the two contact marks the vector of orders of the pulled-back
characters is
\[ \ell_b=(\operatorname{ord}_b f^*x_1,\ldots, \operatorname{ord}_b f^*x_r,0,\ldots,0)\in P_0. \]
It is nonzero exactly when $f(b)\in U\setminus T$. These vectors
are the outgoing weights of the additional tropical ends, since the
slope of each character toward such a puncture is its order there.
The two contact orders are $v$ and $-w$, with $v=m,w=e$ for
a one-input spine. Degree zero of the principal divisors gives
$\sum_b\ell_b=w-v$. On the pointed cone of $P_0$ the
piecewise-linear norm of \cite[Definition~4.11]{KY} is bounded by
a constant times total mutable degree, by compactness of its unit
simplex. Thus the sum of the degrees of all twigs in
\cite[Definition~4.12]{KY} is bounded by the same constant times
the total degree of $w-v$. In particular a zero correction has no twigs. With this degree bound,
the finiteness and specialization statements follow from
\cite[Lemmas~20.29--20.31 and Definition~20.32]{KY}.
\end{proof}

At a general point of a rational hyperplane the \emph{wall function}
records infinitesimal cylinder counts by their bends
\cite[Proposition~20.13 and Definition~20.32]{KY}. To specialize
such a series we set its curve monomials $q^\beta$ equal to one.

\begin{proposition}\label{prop:local-support}
Suppose $a$ is a vertex of the Newton polytope of a regular function
$g$ and $a_i<0$ for all mutable coordinates. Let $C$ be its
closed minimizing normal cone and choose a generic $\xi$ in its
interior. The local series $V_{\xi,m}$ is a finite Laurent polynomial in the characters with coefficients in the
effective curve ring. Its torus exponents
belong to
\[ m+\left\{(d,0):d\in\mathbb N^r,\quad \sum_i(-a_i)d_i\leq\langle a,m\rangle-\tau_g(m)\right\}. \]
After specialization the coefficients are nonnegative integers and
the coefficient at $z^m$ is one. For $m\in C$ the specialization
is $z^m$. Every specialized wall function in the interior of $C$
is one.
\end{proposition}
\begin{proof}
Write $e=m+(d,0)$ using Proposition~\ref{prop:local-bend-cone}.
Lemma~\ref{lem:local-regular-slope} gives
$\tau_g(m)\leq\langle a,e\rangle$, which is exactly the asserted
bound. The possible $d$ form a finite set, so
\cite[Lemma~20.31]{KY} gives finiteness of the whole polynomial with
its curve variables. Positivity follows from the enumerative counts,
and \cite[Lemma~20.34]{KY} gives the specialized coefficient at
$z^m$. For $m\in C$ the bound vanishes. For a general point of a wall in the interior of $C$, choose a
cylinder with both ends in that interior. Its possible bend $(d,0)$ has
$0\leq\langle a,(d,0)\rangle$, which forces $d=0$. The
normalized constant term of the specialized wall function is one
by \cite[Proposition~20.33]{KY}.
\end{proof}

At a polynomial corner each ordinary mutation partner bounds one coordinate of the exponent.
Suppose $G_i$ is a polynomial in $(x,u)$, independent of
$x_i$, with a nonzero constant term, and $G_i/x_i$ extends
regularly to $U$. For $\xi$ in the strict positive orthant of
all marked coordinates the gradient of $\tau_{G_i/x_i}$ is the
negative $i$-th coordinate functional. Consequently $0\le d_i\le-\tau_{G_i}(m)$ for $1\le i\le r$.
These bounds give finitely many exponents, and every specialized wall
function in the interior of the positive orthant is one. Both conclusions
also follow by applying Proposition~\ref{prop:local-support} to
$g=\prod_i(G_i/x_i)$.

\subsection{Constant terms of the local spine functions}
\label{subsec:local-zero-trace}

We use the hypotheses and compactification of
Section~\ref{subsec:local-spine-support}. In particular,
$U\setminus T=V(\prod_i x_i)$. For a generic $\xi\in L_{\mathbb R}$,
let ${\rm SP}_{\xi,m,e}$ be the set of one-input spines ending at
$\xi$ with initial and terminal derivatives $-m,-e$. Write $N(S,\beta)$
for the count of spine $S$ in class $\beta$ and abbreviate
$V_{\xi,m}$ to $V_m$. Thus
\[ V_m=\sum_{\beta,e} n_{\xi}(m,e,\beta)q^\beta z^e,\qquad n_{\xi}(m,e,\beta)=\sum_{S\in{\rm SP}_{\xi,m,e}}N(S,\beta) \]
by \cite[Definition~20.20]{KY}. Write $\operatorname{CT}$ for the
coefficient of $z^0$, keeping all curve coefficients.

\begin{lemma}\label{lem:local-zero-trace}
For $s=2,3$ and nonzero labels $m_1,\ldots,m_s\in L$, we have
\begin{equation}\label{eq:local-constant-term-trace}
\operatorname{CT}(V_{m_1}\cdots V_{m_s}) =\sum_{\beta\in Q}\eta(m_1,\ldots,m_s,\beta)q^\beta,
\end{equation}
where $\eta$ is the algebraic curve count of \cite[Definition~1.1]{KY}.
The identity holds in each curve class. A single choice of $\xi$ outside
a countable union of proper rational linear subspaces works for all
the labels and classes.
\end{lemma}

\begin{proof}
Fix the inputs and the class $\beta$. All decompositions of $\beta$
into effective summands form a finite set. Indeed, an ample divisor bounds
the degrees of the summands. Adding and subtracting a basis of numerical
divisors from a sufficiently large multiple of this ample divisor bounds
each of the integral numerical coordinates of every summand. For each
summand and each input the local series at $\xi$ has only finitely many
terminal exponents in that class by \cite[Proposition~20.21]{KY}.
In the completion used there, one embeds $Q$ in a strictly convex
toric monoid and completes at the ideal generated by monomials with
nonzero curve exponent. Reduction modulo any power of this ideal gives a
Laurent polynomial in the torus characters. For fixed contact orders only finitely many
combinatorial types of complete tropical curves and their spines
occur, even as the curve class varies \cite[Lemma~4.29]{KY}.
We may therefore choose a common toroidal subdivision containing the
input rays and the toric-tail rays at the cuts.
To compare classes, take a common further subdivision of two models
containing the prescribed rays. Each extended smooth map in
\cite[Construction~9.3]{KY} meets the boundary only at its contact
marks, in the open strata of their divisors. Thus it misses subsequent
blowup centers in boundary strata of codimension at least two. The
incidence conditions and evaluations agree there, and the toric tails
use the same torus. The extended curve class pushes forward to the
original class. The same holds for the toric tail classes by
\cite[Definition~7.9 and Proposition~7.11]{KY}, and hence also for
the body class obtained by subtracting them. We can therefore work
on the subdivision and group by pushforward to the original model,
as in \cite[Definition~1.1 and Proposition~17.3(1)]{KY}.

Choose an integer $A$ big with respect to every effective summand
of $\beta$ as in \cite[Definition~4.28]{KY}, so that $A$ bounds
the twig degrees required by concatenation \cite[Theorem~12.4]{KY}.
Let $\operatorname{Wall}_A$ be the finite union of polyhedral walls
constructed for this bound \cite[Constructions~4.15--4.16]{KY}. Take a ball $B$
about $\xi$ disjoint from
$\operatorname{Wall}_A$ and the codimension-one skeleton of the toric fan.
Use the tropical domain with a single vertex and $s+1$ infinite marked
legs. For $s=3$ this is the central vertex of the tripod
$\operatorname{Sk}(\mathcal M_{0,4})$, the skeleton of the moduli
space of four-pointed smooth rational curves. At this vertex the bounded edge
length is zero. This point belongs to the skeleton of the moduli space of smooth
marked curves. All arms meet at $\xi$.

Choose the coordinates of $\xi$ linearly independent over $\mathbb Q$.
We check transversality on this fixed domain stratum.
The central vertex maps to $\xi$, since the evaluation leg is contracted
by \cite[Lemma~8.22]{KY}. The supporting line of an arm through $\xi$
has the form $\xi+\mathbb R a$, $a$ rational. It misses every rational
subspace of codimension two and is contained in no rational hyperplane.
Trace the arm outward until a bend. The bend lies on just one rational
hyperplane $n^\perp$ and is crossed with nonzero normal slope. The
bending condition in \cite[Definition~4.17]{KY} makes the difference of
the two slopes tangent to the hyperplane. The two slopes are therefore
related by an invertible rational shear fixing $n^\perp$ pointwise.
Induction over the finitely many bends expresses each subsequent supporting
line as $B_0\xi+\mathbb R v$, where $B_0$ is invertible over
$\mathbb Q$ and $v$ is rational. It has the same avoidance properties.
At infinity the excluded end loci of
\cite[Construction~4.15]{KY} have inverse images that are finite unions
of rational cones of positive codimension. If a terminal ray ended
in an excluded locus, its supporting line would lie in the linear
span of one of these cones, contradicting the same avoidance property.
Subdivision at a point away from the fan faces
and walls now gives all the conditions in \cite[Definition~4.20]{KY},
including transversality of the straight cut extensions.

Take the canonical point of the product skeleton over this domain
point and $\xi$ \cite[Proposition~8.8 and \S14.1]{KY}. It lies in
every dense Zariski open of the irreducible parameter space
\cite[Lemma~8.2]{KY}, hence in the finite \'etale locus of the
evaluation map \cite[Definition~1.1 and Proposition~3.13]{KY}.
The analytic fiber has degree $\eta(m_1,\ldots,m_s,\beta)$.
Partition it according to the spine of each map, counting residue
field degrees \cite[Construction~9.2 and Lemmas~8.20--8.21]{KY}.
There is no toric tail to remove since all original marks are infinite.

Let $i$ be the interior mark of a complete spine $S$. Denote by
$N_i(S,\beta)$ the degree contributed by maps of class $\beta$
with spine $S$ to the fiber of evaluation at $i$.
To compute the count $N(S,\beta)$ used in concatenation, attach
a second contracted leg $w$ at the central vertex and call the
resulting spine $S^+$. The spine $S$ is stable and transverse, so
$N(S,\beta)=N_w(S^+,\beta)$ by
\cite[Proposition~9.5 and Definition~9.6]{KY}.

Evaluate at $w$ and forget $i$. The smooth $(s+2)$-pointed line
remains stable after forgetting $i$, also when $s=2$. Renaming
$w$ as $i$ recovers the original marked tree and evaluation point
$\xi$. By the Cartesian square of \cite[Lemma~9.13]{KY}, the
evaluation fiber for the additional mark is obtained from the original
fiber by extension of the completed residue field. This comparison
holds for each prescribed spine. Indeed, on every skeletal lift the
leg at the forgotten mark is contracted \cite[Lemma~8.22]{KY}, and
the fixed marked tree places its attachment at the central vertex.
Forgetting the mark preserves the curve class. Finite \'etale degree
is unchanged by extension of the residue field, so
$N_w(S^+,\beta)=N_i(S,\beta)$.

Orient every incoming arm outwards from the central vertex. Its slope
there is $e_j$ by the convention for $V_{m_j}$. Choose a positive
cut distance $\delta$ so small that the initial segment of every arm
lies inside $B$, for all the fixed-class complete maps and for all the
tuples occurring in the local product. The finite bound on the directions gives such a $\delta$, and
$B$ contains no bends. An arm with
nonzero input has nonzero final slope, since a final bend to slope zero
would have to occur on a wall containing $\xi$. Its cut position is
$y_j=\xi+\delta e_j$. Choose $\delta$ outside the countably many
rational affine conditions for a cut to lie on a rational hyperplane.

Balancing at the central vertex gives $\sum_j e_j=0$. Conversely,
any tuple with this sum zero determines the entire central map. It is
constant on the evaluation leg and sends the $j$-th finite segment
to $t\mapsto \xi+t e_j$, $0\le t\le\delta$. This spine has no
bending vertices and its image is contained in $B$, so its count
is one in class zero and zero in every other class by
\cite[Lemma~12.3]{KY}.

Construction~9.3 of \cite{KY} counts the cut arm by extending it with
a toric tail. Its extended class is the sum of the body class and the
class of this tail. Lemma~9.12 of \cite{KY} relates the corresponding
bounds on the two classes. Attach the straight segment from $y_j$
to $\xi$, with derivative $-e_j$. Every arm at $\xi$ with
terminal slope $e_j$ contains this segment, since $B$ meets no wall.
Removing or attaching the segment thus gives inverse operations on the
spines. The segment has count one in class zero, so concatenation
preserves the count in each body class
\cite[Lemma~12.3 and Theorem~12.4]{KY}. Summing the cut-arm degrees
with given input, slope and class gives $n_\xi(m_j,e_j,\beta_j)$.

Now attach a transverse arm at each $y_j$ to the balanced central
spine. The outward derivatives on the two sides of each cut are
opposite. The metric tree and the complete spine are recovered
by concatenation. Successive applications of \cite[Theorem~12.4]{KY}
give the product of the arm counts, summed over the effective splittings
of $\beta$. The central term has class zero. Together with the
partition of the complete fiber, this gives
\[ \eta(m_1,\ldots,m_s,\beta)= \sum_{\substack{e_1+\cdots+e_s=0\\ \beta_1+\cdots+\beta_s=\beta}} \prod_{j=1}^s n_{\xi}(m_j,e_j,\beta_j), \]
which proves~\eqref{eq:local-constant-term-trace} coefficient by coefficient.
\end{proof}

\begin{proposition}\label{prop:formal-local-chart}
Assume $U,T,x_i,u_j$ and the compactification satisfy the hypotheses
of Section~\ref{subsec:local-spine-support}. Write $R=\mathbb Z[Q]$.
For generic $\xi$, sending the intrinsic basis element of label $m$
to $V_{\xi,m}$ defines an injective $R$-algebra homomorphism to
\[ R[[z_1,\ldots,z_r]][z_1^{-1},\ldots,z_r^{-1},z_{r+1}^{\pm1},\ldots,z_{r+c}^{\pm1}]. \]
After augmentation $q^\beta\mapsto1$, the assignment
$\theta_m\mapsto V_{\xi,m}|_{q=1}$ is again injective, with the same
target over $\mathbb Z$. If a regular function has a Newton
vertex negative in every mutable coordinate and $\xi$ lies in the
interior of its minimizing normal cone $C$, the images are Laurent
polynomials. For any $b\in L\cap\operatorname{int}C$, the augmented
map identifies $B_{\mathbb Z}(U)[\theta_b^{-1}]$ with $\mathbb Z[L]$.
Under the hypotheses of Theorem~\ref{thm:incidence} it agrees with
the marked realization of Proposition~\ref{prop:negative-vertex-chart}.
\end{proposition}
\begin{proof}
Write $V_m=V_{\xi,m}$ for this proof.
Proposition~\ref{prop:local-bend-cone} and \cite[Lemmas~20.31 and~20.34]{KY}
show that at exponent $m$ the coefficient of $V_{\xi,m}$ is
a single effective monomial $q^{\delta_m}$, since the coefficient
is a finite sum of nonnegative counts and its augmentation is one. All other
exponents have strictly larger total mutable degree. The unit coordinates of the exponent remain fixed. Thus each coefficient in
a finite product involves only finitely many decompositions of its
mutable displacement and, for each of them, finitely many curve classes.
For $m\ne0$ a one-input spine ending at the generic $\xi$ cannot
have final slope zero. If it did, its last nonzero bend would lie on a
wall containing $\xi$. Hence $\operatorname{CT}(V_m)=0$.
Together with $V_0=1$, this includes lists with zero labels in
Lemma~\ref{lem:local-zero-trace}.

Let $V$ be the $R$-linear map defined by the local functions.
For two elements $a,b$ of the intrinsic algebra, write
$D=V(a)V(b)-V(ab)$. Applying~\eqref{eq:local-constant-term-trace} to the two products and using the
Frobenius structure theorem gives $\operatorname{CT}(D V_m)=0$ in $R$ for every $m$. If $D\ne0$, choose an exponent $u$ of $D$ of smallest total
mutable degree. Such a degree exists because the formal supports are bounded below in
total mutable degree. Then $\operatorname{CT}(D V_{-u})=q^{\delta_{-u}}[z^u]D\ne0$,
since $R$ embeds in the group ring of the numerical curve lattice.
This proves multiplicativity. The same least-degree argument gives
injectivity on finite linear combinations of distinct basis labels,
both over $R$ and after augmentation.

For the stated vertex, Proposition~\ref{prop:local-support} gives
finite Laurent expansions, and after augmentation $\theta_p$ maps
to $z^p$ for every $p\in C\cap L$. Fix $b\in L\cap\operatorname{int}C$.
For any $e\in L$, choose $N\ge0$ with $e+Nb\in C$. The
localized map sends $\theta_{e+Nb}\theta_b^{-N}$ to $z^e$.
Its image thus equals $\mathbb Z[L]$, and injectivity follows from
the injectivity of the augmented map.
When Theorem~\ref{thm:incidence} applies, the realization of
Proposition~\ref{prop:negative-vertex-chart} sends these same
generators to $z^e$, so the two maps agree.
\end{proof}

\subsection{Broken paths and Gauss valuations}
\label{subsec:path-vi}

On the same model we have $P_0=\mathbb N^r\oplus\{0\}^c$ as the
monoid of possible corrections by Proposition~\ref{prop:local-bend-cone}. Write $V_{\xi,m}$ for the
local spine series of \cite[Definition~20.20]{KY}. Choose $\xi$ with
coordinates linearly independent over $\mathbb Q$. At a general
point $y$ of a rational hyperplane $n^\perp$, write $F_y$ for
the specialized wall function of \cite[Definition~20.32]{KY}.
It has constant term one and nonnegative integral coefficients, with
support in $P_0\cap n^\perp$
\cite[Lemma~20.3 and Propositions~20.13, 20.33]{KY}.
We use one function at the point, including all contributions on the
same supporting hyperplane.

An \emph{integral broken path} from $m\ne0$ to $\xi$ is a piecewise affine
map $(-\infty,0]\to L_{\mathbb R}$ with endpoint $\xi$ and initial
derivative $-m$. Its successive derivatives are $-a$, with
$a\in L$. At each of its finitely many bends choose the primitive
normal $n$ so that $k=\langle n,a\rangle>0$. The next exponent
$a'$ lies in the support of $z^aF_y^k$. Crossing with unchanged exponent
requires no mark.

\begin{lemma}\label{lem:path-support}
The terminal exponents of these integral paths are precisely the
support of $V_{\xi,m}|_{q=1}$, with the specialization taken at each
finite correction degree. The support for $m=0$ is $\{0\}$.
\end{lemma}
\begin{proof}
The transversality argument in Lemma~\ref{lem:local-zero-trace}
applies to each segment, also for rational exponents. To see this
explicitly, write
\[ T_\delta a=a+\langle n,a\rangle\delta,\qquad \delta=(a'-a)/\langle n,a\rangle\in n^\perp_{\mathbb Q} \]
at a bend from $a$ to $a'$. The shear $T_\delta$ fixes the wall
pointwise and has inverse $T_{-\delta}$. Let $A_j$ be the product
of the first $j$ shears and $N$ the total number of bends.
Undoing these shears sends the path to a straight line with velocity
$-m$. Thus the supporting line of segment $j$ is
$A_jA_N^{-1}\xi+\mathbb R A_jm$. As in
Lemma~\ref{lem:local-zero-trace}, each such line avoids all rational
subspaces of codimension two and the exceptional loci at infinity.
For a counted spine, trace the successive segments backwards from
$\xi$. Its terminal derivative is nonzero, and the bending condition
\cite[Definition~4.17]{KY} gives the same rational shear at each
crossing. On each line the intersections with rational hyperplanes
form a countable set. Choosing the cuts outside it makes the straight
extensions transverse at every finite order
\cite[Construction~4.15 and Definition~4.20]{KY}.

We record the contribution of a straight interval. At any generic
endpoint $b$, \cite[Lemma~20.34]{KY} gives augmented coefficient
one at the uncorrected exponent $v$ of $V_{b,v}$.
Every individual bend belongs to $P_0$, as we see by applying
the slope calculation
of Lemma~\ref{lem:local-regular-slope} on its two sides, for each mutable coordinate
and for both signs of the units. Thus a counted path from $v$ to $v$
must be the straight ray. Its coefficient one, nonnegativity and
class finiteness \cite[Lemma~20.31]{KY} imply that this ray contributes
in one class $\gamma_b$ with count one. Let $a,b$ be generic
points on a line in this order with derivative $-v$. Concatenating
the straight ray at $a$ with the straight segment $[a,b]$ gives
the ray at $b$. For every class $\beta$, the concatenation formula
of \cite[Theorem~12.4]{KY} therefore reads
\[ \mathbf 1_{\beta=\gamma_b} =\sum_{\gamma_a+\eta=\beta} N([a,b],\eta). \]
The sum is over effective $\eta$. Hence the segment has exactly one
contributing class, of multiplicity one. This also applies to segments crossing fan faces, where the contributing
class can be nonzero. As in Lemma~\ref{lem:local-zero-trace}, an effective
numerical class has only finitely many effective summands.
Consequently one can always take $A$ big for all summands, as
required in the cited concatenation formula.

Now fix a counted spine and one class $\beta$ in which it contributes.
Cut out disjoint small symmetric intervals about its nonzero bends,
with generic cut endpoints. Shrink them until the cylinder counts
stabilize for every effective summand of $\beta$, using
\cite[Definition~20.2]{KY}. Concatenation at the cuts and positivity
give a nonzero stabilized cylinder count at each bend. By
\cite[Proposition~20.13 and Definition~20.32]{KY} this is an allowed
choice from $z^aF_y^k$. Conversely, for an integral broken path
choose at each bend one effective class of nonzero infinitesimal
count realizing the specified term. Such a class exists by positivity,
and only finitely many classes occur for this fixed choice of
exponents \cite[Lemma~20.29]{KY}. Take generic, disjoint sufficiently
short representatives for those chosen classes. The complementary
straight pieces each contribute with multiplicity one by the preceding
argument. Take $\beta$ to be the sum of the chosen cylinder and
straight classes. Enlarge $A$ for all its effective summands.
The chosen representatives remain transverse, and the derivatives
match at every cut. Repeated concatenation in this class has a
strictly positive summand, giving the required contribution to
$V_{\xi,m}$. We cut and glue before specializing the curve
variables. For each fixed terminal exponent there are only
finitely many contributions \cite[Proposition~20.21 and Lemma~20.31]{KY}. For
$m=0$ use $V_{\xi,0}=1$.
\end{proof}

Assume for the remainder of this subsection that the marked realization
sends $\theta_m$ to the finite Laurent polynomial
$\Theta_m=V_{\xi,m}|_{q=1}$ for every $m$. By
Proposition~\ref{prop:formal-local-chart}, this holds in the interior
of the normal cone of a Newton vertex negative in all mutable
coordinates. Over a field of characteristic zero write $\nu_b$ for the
Gauss valuation in this Laurent chart of weight $b\in L^*_{\mathbb R}$,
trivial on the coefficient field.

\begin{lemma}\label{lem:local-support-dilation}
For every $m\in L$ and positive integer $s$,
\[ s\Supp(\Theta_m)\subseteq\Supp(\Theta_{sm})\subseteq s\Newt(\Theta_m). \]
In particular $\Newt(\Theta_{sm})=s\Newt(\Theta_m)$ and
$\nu_b(\Theta_{sm})=s\nu_b(\Theta_m)$ for every real weight $b$.
\end{lemma}
\begin{proof}
The assertion for $m=0$ follows from $\Theta_0=1$.
For $m\ne0$, take a path $\gamma$ with terminal exponent $e$.
The path $t\mapsto\gamma(st)$ has the same image and its exponents
are multiplied by $s$. At each crossing the $s$-th power of the
chosen monomial occurs in the corresponding wall power by positivity.
Thus $se\in\Supp(\Theta_{sm})$ by Lemma~\ref{lem:path-support}.

For the other inclusion, expand $\theta_m^s$ in the intrinsic basis.
Every label $l$ that occurs belongs to $sm+P_0$, since $z^l$
occurs in $\Theta_l$ and all coefficients are nonnegative.
The coefficient of $z^{sm}$ in $\Theta_m^s$ is one, and only the
label $sm$ can contribute to it. Hence $\theta_{sm}$ occurs with
coefficient one. Its Laurent support is contained in the support of
$\Theta_m^s$, whose Newton polytope is $s\Newt(\Theta_m)$.
Taking convex hulls and then minima proves the last assertions.
\end{proof}

\begin{theorem}\label{thm:marked-path-vi}
Over a field $k$ of characteristic zero, for every real weight $b$
and every finitely supported nonzero family $(c_m)$ in $k$ we have
\[ \nu_b\Bigl(\sum_m c_m\Theta_m\Bigr)=\min_{m:\,c_m\ne0}\nu_b(\Theta_m),\qquad \Newt \Bigl(\sum_m c_m\Theta_m\Bigr) =\operatorname{conv}\!\bigcup_{m:\,c_m\ne0}\Newt (\Theta_m). \]
\end{theorem}
\begin{proof}
We adapt the minimizing-path argument of \cite[\S3]{CMMM} to wall
functions whose exponents may lie in several independent tangent
directions. Write $h_m(b)=\nu_b(\Theta_m)$.

Allow rational exponents, and write $\operatorname{conv}_{\mathbb Q}$
for the set of rational convex combinations. Each crossing shear
$T_\delta$ is defined by $\delta\in\operatorname{conv}_{\mathbb Q}
\Supp(F_y)$. For a path with integral input, multiply
its exponents by one positive integer $s$ so that, at every
crossing, $s\langle n,a\rangle$ times each rational barycentric
coefficient is integral. Each scaled choice is then a term of the
integral power of $F_y$ by positivity. By Lemma~\ref{lem:local-support-dilation}, no
rational path for $m$ improves on $h_m(b)$. An integral path
attaining it exists since $\Theta_m$ is a finite Laurent polynomial.

For such a minimizer, at any crossing pull $b$ back by the transposes
of all later shears, obtaining $b_+$. Its chosen $\delta$ minimizes
$\langle b_+,\delta\rangle$ on the rational convex hull of the
wall support. Indeed, otherwise move $\delta$ a small rational
fraction of the way to a better finite convex combination. Keep the other shears and the initial exponent fixed. The final pairing decreases,
since its change is the change in $\langle b_+,\delta\rangle$
multiplied by $\langle n,a\rangle>0$. A crossing at which the
original path goes straight is included by marking it with $\delta=0$.

For each support monomial in these finitely many choices, take an
effective class $\beta$ with nonzero coefficient. Let $I$ be the
ideal generated by the monomials with nonzero curve exponent in the
completion. Choose $\ell$ large enough that each chosen $q^\beta$
is nonzero modulo $I^\ell$. In the resulting finite diagram
\cite[Proposition~20.17]{KY}, their coefficients are constant on the
relative interiors of the wall cells. Undoing the crossing shears gives
a straight line ending at $A_N^{-1}\xi$. On this line each chosen crossing solves a linear equation
with nonzero normal derivative. Its position varies continuously
with the shear. For a sufficiently small change, the crossings
therefore occur in the same order, with the same normal signs and
in the same relative interiors. For a rational change the lines
supporting the perturbed segments avoid the codimension-two strata of the diagram
and the exceptional end loci at every order, by
the proof of Lemma~\ref{lem:path-support}.

The finite diagram thus preserves the individual wall monomials.
Once we clear the new denominators, their products occur in the full
wall powers by positivity, regardless of their curve degree. At
every other crossing we use the constant term. Only the finitely many
selected crossings change the exponent. We have constructed
a rational path with smaller final pairing, contrary to minimality.

Take $b$ with rationally independent coordinates. Each transported
covector is likewise rationally independent, since the transport is
rational and invertible. At a minimizing crossing all the support
points in the chosen finite convex combination must minimize
$\langle b_+,\cdot\rangle$, and at most one lattice point can do so.
Thus the chosen $\delta$ is determined by the wall and $b_+$,
also when the wall has several independent tangent directions. Two
minimizing paths with the same terminal exponent and endpoint agree
backwards. At a first possible disagreement in their finite union of
nonzero bend events the wall, $b_+$, and outgoing exponent $a'$
coincide. Tangency fixes $k=\langle n,a'\rangle>0$, and the unique
choice gives the same incoming exponent $a'-k\delta$ and the
same preceding covector. The paths therefore have the same initial
label. Only the zero input has terminal exponent zero. Consequently distinct labels have distinct minimizing
characters. Their coefficients are nonzero in $k$, so the least
terms in a finite linear combination cannot cancel at this weight.
Continuity of the finite minima extends the equality to all real
weights. Equality of minimum support functions gives the Newton
polytope assertion.
\end{proof}

\begin{corollary}\label{cor:effective-marked-vi}
On the prepared model, under the negative-vertex hypothesis of
Proposition~\ref{prop:local-support}, the finite local Laurent functions
over $R=\mathbb Z[Q]$ satisfy the minimum formula for every
real character weight trivial on $R$. The same holds for their
Newton polytopes in the character lattice.
\end{corollary}
\begin{proof}
Positivity and finiteness show that the support in the character lattice
of each local function is unchanged by setting the curve monomials
equal to one. The distinct minimizing characters of Theorem~\ref{thm:marked-path-vi}
therefore remain distinct before augmentation. Multiplication by a nonzero element of the
domain $R$ preserves the corresponding nonzero coefficient.
Continuity gives the assertion at every weight.
\end{proof}

\begin{proposition}\label{prop:split-vi}
In the setting of Theorem~\ref{thm:split-coordinate}, the theta basis before splitting is valuatively independent for all real Gauss weights in the marked chart if and only if the theta basis after splitting has this property.
\end{proposition}
\begin{proof}
Let $w$ be the old coordinate and $u_1,\ldots,u_d$ the new ones.
By Theorem~\ref{thm:split-coordinate}, the old function
$z^pw^s A_{p,s}(z,w)$ gives $z^p u^m A_{p,s}(z,\sum_i u_i)$
whenever $\sum_i m_i=s$. For generic weights $(c,b_1,\ldots,b_d)$,
choose $j$ with $b_j=\min_i b_i$. If the old minimizing exponent
at $(c,b_j)$ is $(h,t)$, the new minimizing exponent is
$(h,m+(t-s)e_j)$. Indeed, the least term in
$(\sum_i u_i)^{t-s}$ is $u_j^{t-s}$ with coefficient one.

For a generic weight, valuative independence is equivalent to
distinct basis elements having distinct minimizing exponents. In
the forward direction the new exponent recovers $t$ by summing
its last $d$ entries. Thus it determines the old label $(p,s)$
by independence of the old basis, and then $m$ by the formula.
Conversely, lift every old label $(p,s)$ to $(p,se_j)$. A common
minimizing exponent $(h,t)$ for two old functions would give the
same exponent $(h,te_j)$ for their lifts. Every generic old weight
admits a generic lift with $b_j$ least. Continuity for each finite
linear combination extends both implications to all real weights.
\end{proof}

\subsection{Atomicity in the intrinsic local charts}
\label{subsec:intrinsic-atomicity}

Work on the prepared compactification of
Proposition~\ref{prop:formal-local-chart}. Write $Q$ for its monoid of
effective numerical curve classes, $R=\mathbb Z[Q]$, and
$R_+=\mathbb N[Q]$. Denote the intrinsic algebra over $R$ by
$B_R$, its basis by $(\vartheta_m)_{m\in L}$, and its local map at
a generic endpoint $\xi$ by $V_\xi$. Thus
$V_\xi(\vartheta_m)=V_{\xi,m}$. We call an element of $B_R$
\emph{positive in the intrinsic local charts} if every coefficient of
$V_\xi(f)$ belongs to $R_+$ at every generic endpoint of the
marked source torus.

A nonzero positive element is \emph{additively indecomposable} if it
cannot be written as a sum of two nonzero positive elements. The
atomicity assertion below identifies all such elements.
For cluster varieties, Mandel proves atomicity in the scattering atlas
by reading a theta coefficient near its label
\cite[Theorem~1]{MandelAtomic}. In the intrinsic local charts the
same argument recovers the coefficient with its curve class.

\begin{theorem}\label{thm:intrinsic-atomicity}
Let $E\subset L$ be finite and $m\in E$. There is a neighborhood
of $m$ in $L_{\mathbb R}$ such that at every generic endpoint
$\xi$ in this neighborhood,
\[ [z^m]V_{\xi,l}=\delta_{ml}\quad(l\in E) \]
in $R$. Consequently the elements positive in the intrinsic local
charts form $\bigoplus_{m\in L}R_+\vartheta_m$. Its nonzero
additively indecomposable elements are precisely
$q^\beta\vartheta_m$ for $\beta\in Q$ and $m\in L$.
After curve augmentation the corresponding elements are the
individual $\theta_m$.
\end{theorem}

\begin{proof}
For $m=0$ the assertion follows from $V_{\xi,0}=1$ and
$\operatorname{CT}(V_{\xi,l})=0$ for $l\ne0$, as in the proof
of Proposition~\ref{prop:formal-local-chart}. Suppose $m\ne0$.
Write $v\preceq w$ if $w-v\in\mathbb N^r\oplus\{0\}^c$.
If a counted path from $l\in E$ has terminal exponent $m$, its
successive exponents increase for this order by
Proposition~\ref{prop:local-bend-cone} and Lemma~\ref{lem:path-support}.
If it bends, the exponent $v$ immediately before its last nonzero
bend therefore belongs to one of the finite intervals
$\{v:l\preceq v\preceq m\}$, and $v\ne m$.

Such a bend is supported on a rational hyperplane carrying a nonzero
infinitesimal cylinder count from $v$ to $m$, by the cutting
argument of Lemma~\ref{lem:path-support}. There are only finitely many
possible hyperplanes, independently of $\xi$. Indeed, for each
of the finitely many $v$, \cite[Lemma~20.29]{KY} bounds the set
of contributing classes for the fixed displacement $m-v$.
All these classes survive at some common finite order of the
curve-monomial filtration, since the ambient toric monoid is strictly
convex. Now \cite[Lemma~20.6]{KY} gives finitely many primitive wall
normals at that order.

Enlarge this collection of hyperplanes to contain all the cells of
$\operatorname{Wall}_A$ for one $A$ big with respect to the zero
class, and all codimension-one cones of the toric fan on the prepared
model. These are finite rational conical sets
\cite[Constructions~4.15--4.16]{KY}. For each resulting hyperplane
$n^\perp$ with $\langle n,m\rangle\ne0$, take a neighborhood
of $m$ on which $\langle n,\cdot\rangle$ has the same sign as
$\langle n,m\rangle$. At a generic $\xi$ in their intersection
the ray $\xi+t m$, $t\ge0$, misses every hyperplane in the
collection. For the hyperplanes containing $m$ this follows from
$\langle n,\xi+t m\rangle=\langle n,\xi\rangle\ne0$.

Following the terminal segment backwards gives exactly this ray up
to the last bend. Hence a contributing path with terminal exponent
$m$ is straight and its input is $m$. The straight ray avoids
$\operatorname{Wall}_A$ and the codimension-one skeleton of the
toric fan. Its count in class zero is one by \cite[Lemma~12.3]{KY}.
The augmented coefficient of $z^m$ in $V_{\xi,m}$ is one by
\cite[Lemma~20.34]{KY}. All the counts are nonnegative and the
coefficient at this exponent is a finite sum over classes. Thus its
entire coefficient in $R$ is one.

Now write $f=\sum_{l\in E} c_l\vartheta_l$ with $c_l\in R$.
For $\xi$ as above we have $[z^m]V_\xi(f)=c_m$.
Positivity in the local charts therefore gives $c_m\in R_+$
for each $m$. Conversely every local coefficient of an individual
$\vartheta_m$ belongs to $R_+$ by its counting definition.
These elements therefore generate the positive additive monoid freely.
The coefficient identity also holds after $q^\beta\mapsto1$, so
the same argument applies to the augmented algebra.
\end{proof}

Thus one can test positivity on any nonempty set of generic points
whose closure meets each ray $\mathbb R_{>0}m$,
$0\ne m\in L$. Indeed the sign conditions used in the proof hold
in an open neighborhood of every point of that ray. For the zero label
any generic point suffices.

\subsection{Comparing Laurent charts}\label{subsec:based-scattering}

Agreement on a basis of characters determines the marked Laurent
expansions. To transfer valuative independence between two charts,
it suffices that the change of coordinates preserve the logarithmic
volume form up to a scalar.

\subsubsection{Identifying the theta functions in a scattering chart}

Let $k$ be a field of characteristic zero and $M$ the full character
lattice of a marked torus. Write $k(M)=\Frac (k[M])$.
Suppose $A$ and $B$ have $k$-bases $(\theta_m)$ and
$(\psi_m)$ with the same labels, and have faithful Laurent realizations
$\phi:A\hookrightarrow k[M]$ and $\iota:B\hookrightarrow k[M]$.
A \emph{based isomorphism} $F:A\to B$ sends $\theta_m$ to
$\psi_m$ for every $m$.

\begin{lemma}\label{lem:scattering-chart-marking}
Suppose $F$ is a based isomorphism and $u_1,\dots,u_n$ are a
$\mathbb Z$-basis of $M$. If there are labels $m_j$ such that
\[ \phi(\theta_{m_j})=\iota(\psi_{m_j})=z^{u_j} \qquad(1\le j\le n), \]
then $\phi(\theta_m)=\iota(\psi_m)$ for every $m$.
\end{lemma}
\begin{proof}
The characters $z^{u_j}$ generate $k(M)$, so their preimages
$\theta_{m_j}$ generate $\Frac(A)$. The two embeddings $\phi$
and $\iota F$ agree on these generators and hence on every element of $A$.
\end{proof}

A based isomorphism can itself be obtained from trace identities.
Let the traces on both algebras be the coefficient of the distinguished
unit in the respective bases. Suppose the traces at label zero on every
pair and every triple of labels agree and the bilinear pairing
$(b,c)\mapsto\Tr _B(bc)$ is nondegenerate. The pair
and triple identities for the label-preserving linear bijection $F$ give
$\Tr _B((F(ab)-F(a)F(b))F(c))=0$ for every $c\in A$.
Nondegeneracy proves multiplicativity.
This is the uniqueness argument of the Frobenius structure theorem
\cite[Theorem~1.2]{KY}. Gross and Siebert define multiplication in their
mirror algebra by punctured Gromov--Witten invariants~\cite{GS}.
For a smooth affine log Calabi--Yau variety containing a dense torus,
Johnston identifies the Gross--Siebert and Keel--Yu algebras when
the compactifying boundary is the support of an effective nef divisor
\cite[Theorem~1.1 and Corollary~1.2]{Johnston}. To identify the
Laurent expansions one must also identify the characters in
Lemma~\ref{lem:scattering-chart-marking}.

In particular, suppose $\iota(\psi_m)$ are finite broken-line theta
functions at a generic endpoint of a positive scattering diagram in
the sense of \cite[\S2.1]{CMMM}, with the same marked character lattice.
In that definition the nonzero character exponents on each wall are
positive multiples of a single vector.
Under the hypotheses of Lemma~\ref{lem:scattering-chart-marking},
the intrinsic functions $\phi(\theta_m)$ then satisfy
\[ v_w\!\left(\sum_m c_m\phi(\theta_m)\right) =\min_{c_m\ne0}v_w\bigl(\phi(\theta_m)\bigr) \]
for every real Gauss weight $w$ and every nonzero finite sum with
$c_m\in k$. Before specializing the scattering coefficients,
\cite[Theorems~3.7 and~3.8]{CMMM} give distinct minimizing exponents
for distinct labels at a generic weight. Continuity gives the formula at
all weights. For specialized scattering coefficients one may instead apply
\cite[Theorem~3.10]{CMMM} to the functions indexed by the remaining
character lattice, provided its Assumption~3.9 holds, the functions are
finite before specialization, and every nonzero integer multiple of a
coefficient monomial has nonzero image in $k$.

For an intrinsic chart obtained from a negative Newton vertex,
Theorem~\ref{thm:marked-path-vi} establishes the required independence
directly. The general intrinsic independence conjecture for skeletal
valuations is formulated in \cite[Conjecture~21.7]{KWY}.

\subsubsection{Changing torus coordinates}

\begin{proposition}\label{prop:volume-transfer}
Let $K=k(x_1,\ldots,x_n)=k(y_1,\ldots,y_n)$ over a field of
characteristic zero. Suppose the logarithmic volume forms
$\Omega_x=\bigwedge_i d\log x_i$ and
$\Omega_y=\bigwedge_i d\log y_i$ are proportional over $k$.
A family of functions in $K$ which is Laurent in both charts
is valuatively independent for all real Gauss valuations in one
chart if and only if it is so in the other.
\end{proposition}

\begin{proof}
It suffices first to consider a primitive integral weight in the
$x$-chart. Its Gauss valuation $v$ is divisorial, and on a
smooth toric model containing its ray the order of $\Omega_x$
at $v$ is $-1$. The same is therefore true of $\Omega_y$.
Some $y_i$ has nonzero value. Otherwise all the $y_i$ would
be units at the generic point of the divisor, and $\Omega_y$
would be regular there. After a unimodular change of the $y$-characters,
the first has positive value $d$ and all the others have value zero.
Since the pole of $\Omega_y$ is simple, its nonzero residue is,
up to sign, $d$ times the wedge of the logarithmic differentials
of the last $n-1$ residues. These residues are therefore algebraically
independent in characteristic zero. In a Laurent expansion by powers
of the first character, each nonzero coefficient in the remaining
characters has value zero. The powers have distinct values. This is
the Gauss rule, and normalization of the divisorial valuation gives
$d=1$.

The minimum formula in the $y$-chart thus proves it for
all primitive weights in the $x$-chart. Positive rescaling gives
the formula at rational weights. For a fixed nonzero finite sum its
Gauss value in the $x$-chart is the minimum of finitely many
linear functions, hence is continuous in the weight. The same holds
for the minimum of the values of its nonzero summands. Density
gives the assertion for all real weights. Interchanging the charts
gives the converse.
\end{proof}

In particular a based isomorphism with a scattering algebra transfers
valuative independence as soon as its two Laurent realizations induce
a volume-preserving birational change of torus coordinates. The
individual character identities of
Lemma~\ref{lem:scattering-chart-marking} give the stronger conclusion
that the Laurent polynomials themselves agree.

\begin{corollary}\label{cor:open-torus-volume-transfer}
Suppose an integral variety $X$ over a field $k$ of characteristic
zero has a rational top form
$\Omega$ regular and nowhere vanishing on its smooth locus. If on
one dense open split torus it is the standard logarithmic volume up
to a scalar, then the same holds on every dense open split torus
of $X$. Hence a family of regular functions on $X$ satisfying
all-weight valuative independence in one such chart satisfies it
in every such chart.
\end{corollary}

\begin{proof}
On a second torus with coordinates $y_i$, the quotient
$\Omega/\Omega_y$ is a regular unit. It equals $c y^h$ for
some $c\in k^\times$ and $h\in\mathbb Z^n$. If $h_j\ne0$,
differentiating
$y^h\bigwedge_{i\ne j}d\log y_i$ and dividing by
$(-1)^{j-1}h_j$ would make $y^h\Omega_y$ exact. Yet $\Omega_x$ is not exact in $k(X)$.
Indeed, expand rational coefficients in the iterated Laurent series
field in the $x_i$. The iterated residue of an exact top form vanishes, since a partial
derivative in $x_i$ has zero coefficient at $x_i^{-1}$. The
iterated residue of $\Omega_x$ is one.
Thus $h=0$, and Proposition~\ref{prop:volume-transfer} applies.
\end{proof}

\section{Recognizing the basis of an upper LP algebra}
\label{sec:star-recognition}

Suppose an intrinsic algebra and an LP seed are given in the same
Laurent ring. To identify the theta functions as a basis of the upper
LP algebra, we must show that they extend to the adjacent seed charts
and span the intersection of their Laurent rings with the initial one.

\subsection{Extension and spanning along mutations}

We check extension on algebra generators. To prove spanning, we
compare least exponents along a finite tree of mutations.

\subsubsection{Regularity in the first star}

Let $R$ be a UFD, let $M=\mathbb Z^n$, and write
$\mathcal L=R[M]=R[z_1^{\pm1},\ldots,z_n^{\pm1}]$. Designated Laurent coefficient
characters may be included among the $z_j$. Polynomial coefficients are
kept in $R$. Fix a target LP seed with this initial Laurent ring, and let
$I$ be its set of mutable indices. For $i\in I$ write
\[ \mathcal L_{\ne i}=R[z_j^{\pm1}:j\ne i],\qquad z_i'=H_i/z_i,\qquad \mathcal L_i=\mathcal L_{\ne i}[(z_i')^{\pm1}]. \]
Here $H_i\in \mathcal L_{\ne i}$ is the exchange Laurent polynomial,
including its normalizing monomial. Set $W_i=\mathcal L\cap \mathcal L_i$.
The $z_i$-\emph{slice} of degree $u$ of a Laurent polynomial is
its coefficient at $z_i^u$, regarded as an element of $\mathcal L_{\ne i}$.
By the LP starfish theorem \cite[Theorem~1.3]{Branching},
the full upper LP algebra of the seed is
$\mathcal U=\mathcal L\cap\bigcap_{i\in I}\mathcal L_i$.

\begin{lemma}\label{lem:recognition-rank-one}
Fix $i\in I$, and view $\mathbb Z^{n-1}$ as the sublattice of $M$
with zero $i$-th coordinate. An element
$f=\sum_{u\in\mathbb Z}z_i^u f_u$ of $\mathcal L$ lies in $W_i$ precisely
when $H_i^{-u}$ divides $f_u$ in $\mathcal L_{\ne i}$ for every $u<0$.
The functions
\[ E_{u,b}=z^b z_i^{[u]_+}(z_i')^{[-u]_+}, \qquad u\in\mathbb Z,\quad b\in\mathbb Z^{n-1}, \]
form an $R$-basis of $W_i$, where $[u]_+=\max(u,0)$.
\end{lemma}
\begin{proof}
In the neighboring torus the coefficient of $(z_i')^{-u}$ is
$H_i^u f_u$. This gives the divisibility test. In the initial torus
$E_{u,b}=z_i^uH_i^{[-u]_+}z^b$, so division of each admissible slice by
$H_i^{[-u]_+}$ gives exactly its coordinates in the asserted basis.
\end{proof}

\begin{corollary}\label{cor:finite-star-generator-test}
Suppose $A\subset \mathcal L$ is an $R$-algebra and $g_1,\ldots,g_t$ are
known $R$-algebra generators of $A$. Then $A\subseteq\mathcal U$
if and only if the negative $z_i$-slices of the $g_j$ satisfy the
divisibilities of Lemma~\ref{lem:recognition-rank-one} for every $i\in I$.
\end{corollary}
\begin{proof}
Each $W_i$ is a subring. Thus the generators lie in every $W_i$
if and only if $A$ does.
\end{proof}

Suppose $A$ has an exhaustive increasing multiplicative filtration
indexed by the nonnegative integers, and the initial forms of
$g_1,\ldots,g_t$ generate $\operatorname{gr}A$ over $R$.
Lifting a homogeneous expression for the leading form and subtracting
it lowers the degree. Induction proves that the $g_j$ generate $A$.
Thus Hilbert-basis generators of a finite toric-face degeneration
lift to algebra generators, whose Laurent expansions suffice for
Corollary~\ref{cor:finite-star-generator-test}.

\subsubsection{Spanning along a finite mutation tree}

Let $\Lambda$ be a lattice of rank $n$ and let
$(\vartheta_m)_{m\in\Lambda}$ be the prescribed basis of an $R$-algebra
$A\subset \mathcal L$. Here the basis and its Laurent realization are given
independently of $\mathcal U$. A \emph{marking of the basis} consists of a continuous,
positively homogeneous piecewise integral-linear homeomorphism
$\mu:\Lambda_{\mathbb R}\longrightarrow M_{\mathbb R}$ inducing a
bijection of lattices, together with a pointed rational polyhedral cone
$C\subset M_{\mathbb R}$, so $C\cap(-C)=\{0\}$. We say that the basis is \emph{pointed in
this marking} if, in the initial Laurent expansion of every $\vartheta_m$,
the coefficient of $z^{\mu(m)}$ is $1$ and
$\Supp(\vartheta_m)\subseteq\mu(m)+C$.

We can recover the marking from the leading terms of products.
In the application to a counting source, the cone $P$ below records
nonnegative corrections at regular coordinates and zero corrections
at units.

\begin{lemma}\label{lem:straight-contact-pointing}
Suppose $G:M\longrightarrow\Lambda$ is a lattice isomorphism and
$\vartheta_{Ga}=z^a$ for every $a\in\mathbb N^n$. Let
$\tau:\Lambda_{\mathbb R}\longrightarrow V_{\mathbb R}$ be a continuous
positively homogeneous piecewise unimodular homeomorphism for another
lattice $V$. Suppose there is a pointed rational polyhedral cone
$P\subset V_{\mathbb R}$ such that in every product
$\vartheta_p\vartheta_q$ the term with $\tau$-coordinate
$\tau(p)+\tau(q)$ has coefficient $1$, and every other term has
$\tau$-coordinate in
$\tau(p)+\tau(q)+(P\setminus\{0\})$.
Let $G_\tau:M\to V$ have columns $\tau(Ge_j)$.
Then $G_\tau$ is unimodular, and the basis is pointed in the marking
$(G_\tau^{-1}\tau,G_\tau^{-1}P)$.
\end{lemma}
\begin{proof}
Each product of the monomials in the hypothesis is a single basis element.
Their prescribed leading outputs imply $\tau(Ga)=G_\tau a$ for all
$a\in\mathbb N^n$. Homogeneity and continuity extend the equality
to the real positive orthant. One full-dimensional unimodular piece of
$\tau G$ meets its interior in an open set, where it agrees with $G_\tau$.
Thus $G_\tau$ is unimodular.

Write $\vartheta_m=\sum_a c_a z^a$. Choose a positive integral vector
$b$ sufficiently far inside the orthant that every $a+b$ in this
finite sum and $G_\tau^{-1}\tau(m)+b$ belong to $\mathbb N^n$.
Then $\vartheta_{Gb}\vartheta_m=\sum_a c_a\vartheta_{G(a+b)}$.
The term on the right indexed by $a$ has $\tau$-coordinate
$G_\tau(a+b)$. Comparing with the prescribed leading term and the
cone for the other outputs gives the assertion.
\end{proof}

For a vertex $v$ of the Newton polytope of $H_i$, define the map
$\lambda_{i,v}(u,b)=u e_i+b+[-u]_+v$.
It is a bijection from $\mathbb Z\times\mathbb Z^{n-1}$ onto $M$.
At any pointed marking $(\mu,C)$ with $A\subseteq W_i$, there is
a unique vertex $v$ of $\Newt(H_i)$ such that
$\Newt(H_i)-v\subset C$. Its coefficient is a unit of $R$.
Indeed, take the label with $\mu(m)=-e_i$. The $-1$-slice of
$\vartheta_m$, stripped of $z_i^{-1}$, has constant coefficient one
and support in $C$. By Lemma~\ref{lem:recognition-rank-one} it is
$H_i q$. Take a generic weight strictly positive on $C\setminus\{0\}$.
The least exponents $v,b$ of the two factors have sum zero and their
coefficients have product one. Newton additivity gives
$\Newt(H_i)+b\subseteq\Newt(H_iq)\subset C$.
Two vertices with this property would have both differences in the
pointed cone, so they coincide. At this vertex the exponent
$\lambda_{i,v}(u,b)$ of $E_{u,b}$ has unit coefficient and all
corrections lie in $C$.

\begin{theorem}\label{thm:finite-tree-recognition}
Let a finite tree carry a pointed marking $(\mu_s,C_s)$ at every
vertex and a mutable index $i$ on each edge. Write $J$ for the set
of indices used and $W_J=\bigcap_{i\in J}W_i$, with $W_\varnothing=\mathcal L$.
Assume $A\subseteq W_J$. On an edge $s,t$ of index $i$, use the
vertices $v_s,v_t$ determined above for $H_i$. Suppose
\[ \lambda_{i,v_s}^{-1}\mu_s=\lambda_{i,v_t}^{-1}\mu_t \quad\text{on each edge},\qquad \bigcap_s\mu_s^{-1}(C_s)=\{0\}. \]
Then $A=W_J$ with the prescribed $R$-basis. In particular, if
$A\subseteq\mathcal U$, then $A=\mathcal U$.
\end{theorem}

The edge equality is equivalent to
$\mu_t(m)=\mu_s(m)+[-\mu_s(m)_i]_+(v_t-v_s)$.
For finite piecewise-linear markings both conditions can be tested
on finite rational polyhedral subdivisions.

\begin{proof}
Fix a vertex $s$. An infinite sum in exponents of $M$ will be used
at $s$ only when its support lies in a finite union of translates of
$C_s$. Choose an integral linear functional strictly positive on the
nonzero vectors of $C_s$. On each such translate the lattice points of
bounded functional value form a finite set. In particular, a matrix which
sends each exponent to a finite sum supported in that exponent plus
$C_s\setminus\{0\}$ strictly raises the functional by at least one on
every nonzero entry. For such a matrix $D$, the series
$\sum_{k\ge0}(-D)^k$ defines a coefficientwise inverse to $1+D$.
For a fixed exponent only finitely many starting exponents and powers
of $D$ can contribute, by the bound on the functional. Each column
of $D$ has finite support. The finite geometric-sum identity thus
gives the inverse at every exponent. When $C_s=\{0\}$, $D=0$.

Index the monomial basis at $s$ by $m\mapsto z^{\mu_s(m)}$.
The matrix sending this indexed monomial basis to the Laurent expansions
of the correspondingly indexed $\vartheta_m$ has the form $1+D$
just described. Thus every $f\in \mathcal L$ has a unique formal expansion in
the $\vartheta_m$ at $s$ within this class of sums. Its coefficients
$c^s_m(f)$ satisfy
\[ c^s_m(f)\ne0\quad\Longrightarrow\quad \mu_s(m)\in\Supp(f)+C_s . \]

We prove equality of these coefficients across one edge when $f\in W_i$.
The edge condition identifies the same basis element $E_{u,b}$ with
a given label $m$ at both ends. Denote it on this edge by $E_m$.
The expansions $\vartheta_m=\sum_n T_{nm}E_n$ in this $R$-basis
of $W_i$ give a matrix with finite columns.
Fix an end $j\in\{s,t\}$ and a nonzero coefficient $T_{nm}$.
If $E_n=E_{u,b}$, the $u$-slice of $\vartheta_m$, after removing
$z_i^u$, is
\[ H_i^{[-u]_+}\sum_{b'}d_{b'}z^{b'}, \qquad d_b=T_{nm}\ne0 . \]
Over the domain $R$ the initial forms of two nonzero factors have
nonzero product for every real weight. This proves that their Newton
polytopes add under multiplication. Consequently the Newton polytope of
this slice, with its $i$-coordinate restored, contains
$u e_i+b+[-u]_+v_j=\mu_j(n)$, even if the monomial itself cancels.
Pointedness of $\vartheta_m$ gives
$\mu_j(n)\in\mu_j(m)+C_j$.

If $E_n$ contributes to the coefficient of $z^{\mu_j(m)}$, the
tangent-cone containment also gives $\mu_j(m)\in\mu_j(n)+C_j$.
Pointedness of $C_j$ and injectivity of $\mu_j$ force $n=m$.
The coefficients of $E_m$ and $\vartheta_m$ there are units,
so $T_{mm}$ is a unit. Write $T=(1+D)\Delta$, where $\Delta$
is its diagonal. Every nonzero entry of $D$ advances strictly in both
markings.

Write $b(f)$ for the coefficients of $f\in W_i$ in the $E_m$-basis.
At either end of the edge we compute the same vector
$\Delta^{-1}\sum_{k\ge0}(-D)^k b(f)$. Indeed, the initial vector and the matrix entries agree, and only finitely
many paths contribute to each label. At each end the vector lies in the completion
described above and is sent by $T$ to $b(f)$. Expansion of each
$E_n$ in monomials at that end also preserves the completion, by the
tangent-cone containment. The resulting theta sum therefore represents
$f$ as a Laurent series at that vertex. Uniqueness of the vertex
expansion identifies it with $c^j(f)$. This proves equality across
the edge for arbitrary coefficients in $R$.

Every $f\in W_J$ belongs to each edge ring. Connectedness
gives a single coefficient family $c_m(f)$ throughout the tree. If its
support were unbounded, take a sequence $m_\ell$ in its support with
$\|m_\ell\|\to\infty$ and $m_\ell/\|m_\ell\|\to m_*$ of norm one.
The preceding support inclusion, finiteness of $\Supp(f)$,
and homogeneity and continuity of $\mu_s$ imply
$\mu_s(m_*)\in C_s$ at every vertex. This contradicts the second
hypothesis. The common support is bounded in a lattice and hence finite.
Thus $f$ belongs to $A$, as required.
\end{proof}

\begin{corollary}\label{cor:finite-tree-expansion}
Under Theorem~\ref{thm:finite-tree-recognition}, suppose the markings
are specified on finite rational polyhedral subdivisions. For nonzero
$f\in\mathcal U$, every label in its expansion belongs to the finite set
\[ \Lambda\cap\bigcap_s\mu_s^{-1}\bigl(\Newt (f)+C_s\bigr). \]
The Laurent expansions of the functions indexed by this set determine
the coefficients of $f$ by unitriangular substitution over $R$.
\end{corollary}
\begin{proof}
The support estimate in the proof of the theorem gives the inclusion.
The same limiting argument, applied to the closed sets with
$\Newt (f)$ in place of $\Supp(f)$,
shows that their intersection is bounded. Restriction to a common
subdivision expresses it as a finite union of rational polytopes, so
its lattice points can be enumerated. At one vertex $s$, order these
labels by increasing value of an integral functional strictly positive
on $C_s\setminus\{0\}$, evaluated on $\mu_s(m)$. In this order
read the coefficient of $z^{\mu_s(m)}$ from the current remainder
and subtract that multiple of $\vartheta_m$. Pointedness makes
every off-diagonal contribution strictly later. This procedure recovers
the asserted expansion, including zero coefficients at surplus labels.
\end{proof}

\subsection{Comparing expansions across a common facet}

Let $k$ be a field and let $\mathcal U$ be the full upper LP
algebra in its marked rational function field. Suppose a subalgebra
$A\subseteq\mathcal U$ has a $k$-basis $(\vartheta_m)_{m\in\Lambda}$.
At each of two adjacent seeds $j=s,t$, use its own character
lattice $M_j$ and a marking $(\mu_j,C_j)$ as in
Theorem~\ref{thm:finite-tree-recognition}. Thus the Laurent expansion
in that seed has a nonzero coefficient at $\mu_j(m)$ and is supported
in $\mu_j(m)+C_j$. Write $\Supp_j(f)$ for the support of $f$ in
this chart. Denote the coefficients of the triangular formal
expansion of $f$ at that marking by $c_m^j(f)$.

Consider a nonzero integral weight $\eta_s$ and its transform
$\eta_t$ under the seed exchange, as in
Lemma~\ref{lem:boundary-gauss-change}. They define one discrete
valuation $v$ of the ambient field. Suppose
$\eta_s\in C_s^\vee$ and $\eta_t\in C_t^\vee$, allowing weights
on the boundaries of the dual cones $C_j^\vee$, which consist of
the weights nonnegative on $C_j$. For a subalgebra
$B$ of the field, write $\operatorname{gr}_v B$ for the
associated graded algebra of the filtration by $v\ge l$,
$l\in\mathbb Z$. It embeds in the associated graded ring of the
ambient field, where we compare homogeneous symbols from the two charts.
For $0\ne f\in B$, its \emph{initial form} $\operatorname{in}_v(f)$
is its image in the graded piece of degree $v(f)$.

\begin{proposition}\label{prop:facet-graded-transport}
If $\operatorname{gr}_v A=\operatorname{gr}_v\mathcal U$, then
$c_m^s(f)=c_m^t(f)$ for every $f\in\mathcal U$ and every label
$m$. Conversely, this equality of coefficients implies equality
of the graded algebras if the only $b\in\Lambda_{\mathbb R}$ with
$\mu_s(b)\in C_s,\ \mu_t(b)\in C_t$, and $\langle\eta_s,\mu_s(b)\rangle=0$, is $b=0$.
\end{proposition}
\begin{proof}
For either chart the degree of $\vartheta_m$ is
$d_m=\langle\eta_j,\mu_j(m)\rangle$. Its initial form in that
chart has the same distinguished exponent and corrections in
$C_j\cap\ker\eta_j$. A linear functional strictly positive on
$C_j\setminus\{0\}$, with ties broken by a lattice order, gives
distinct least monomials for distinct labels. The initial forms of
each degree are therefore linearly independent. They form a
$k$-basis of $\operatorname{gr}_v A$, since the terms of lowest
degree in every finite basis expansion survive in the associated
graded algebra.

If the graded algebras agree, express the initial form of a nonzero
$f\in\mathcal U$ in this basis and subtract the corresponding
finite sum of $\vartheta_m$. The remainder belongs to $\mathcal U$
and its degree has increased by at least one unless it is zero. For
a prescribed label $m$, iterate until the remainder $g$ is zero
or $v(g)>d_m$. In the latter case the triangular inverse at
either marking gives $c_m^j(g)=0$. Indeed, a nonzero coefficient
would imply $\mu_j(m)\in\Supp_j(g)+C_j$, hence
$d_m\ge v(g)$. The subtracted finite sums have the same basis
coefficients at both ends, which proves the first assertion.

For the converse, write $c_m$ for the common coefficients of
$f\ne0$. Each nonzero coefficient satisfies both support
inclusions and $d_m\ge v(f)$. Only finitely many such labels can
have $d_m=v(f)$. Otherwise, normalizing an unbounded sequence
of them and taking a subsequence gives a unit vector $b$
satisfying the three conditions, by continuity and homogeneity of
the markings. At the exponents of weight $v(f)$, the triangular
expansion of $f$ now gives
$\operatorname{in}_v(f)=\sum_{d_m=v(f)}c_m
\operatorname{in}_v(\vartheta_m)$, where the nonzero sum is finite.
This proves the graded equality.
\end{proof}

Two finite verifications of the graded hypothesis are available.
Write the exchange as $x'x=H(Z)$, with $Z$ the unchanged
characters, and let $H_0$ be the sum of the terms of $H$
minimizing the weight on $Z$. The associated graded algebra of
each Laurent ring for its Gauss valuation is a Laurent polynomial
ring graded by the coordinate weights. Write $X,Z$ for the
coordinate symbols of the first torus, and evaluate $H_0$ in
these symbols. The symbol of $x'$ is $H_0/X$. In either torus
the initial form of an upper element is Laurent in the symbols of
its coordinates. Their Laurent intersection is therefore an upper
bound for $\operatorname{gr}_v\mathcal U$, and
Lemma~\ref{lem:recognition-rank-one} identifies this intersection
with $k[Z^{\pm1},X,H_0/X]$. It suffices to express these finitely
many generators in initial forms of elements of $A$.

Alternatively, suppose a finite presentation
$\mathcal U=k[g_1,\ldots,g_h]$ with ideal of relations
$I\subset k[T_1,\ldots,T_h]$ is known. Set $w_i=v(g_i)$,
and let $J$ be the ideal of relations of the Laurent initial
forms $\operatorname{in}_v(g_i)$ in one chart. Write
$\operatorname{in}_w I$ for the ideal generated by the minimum-weight
parts of the nonzero polynomials in $I$. If
$\operatorname{in}_w I=J$, it suffices to express these $h$
initial forms in initial forms of elements of $A$. To prove this,
represent $0\ne f\in\mathcal U$ by a polynomial $P(T)$.
The minimum $w$-weight of $P$ is at most $v(f)$. If strictly
smaller, its minimum part belongs to $J$. Express this part as a
sum of multiples of initial forms of relations in $I$, taking
only the multipliers of the required homogeneous weights.
Subtracting the same multiples of the relations from $P$ raises
its minimum weight. The weights are integral and bounded above by
$v(f)$, so after finitely many subtractions the minimum part of
the resulting polynomial evaluates to $\operatorname{in}_v(f)$.
This proves that the $\operatorname{in}_v(g_i)$ generate
$\operatorname{gr}_v\mathcal U$. Over an effective field, elimination
computes $J$ and weight homogenization computes $\operatorname{in}_w I$.

\begin{corollary}\label{cor:facet-graded-spanning}
Let a finite connected graph of seed charts carry such
markings. Suppose that the graded hypothesis of
Proposition~\ref{prop:facet-graded-transport} holds at every edge
and $\bigcap_j\mu_j^{-1}(C_j)=\{0\}$. Then
$\mathcal U=\bigoplus_{m\in\Lambda}k\vartheta_m$. The support of
the basis expansion of $f\in\mathcal U$ is contained in
$\Lambda\cap\bigcap_j\mu_j^{-1}\bigl(\Supp_j(f)+C_j\bigr)$.
\end{corollary}
\begin{proof}
Connectedness gives common coefficients, and the triangular inverse
in each chart gives the stated support inclusion using the Laurent
expansion of $f$ in that chart. An unbounded sequence in this set
would yield, after normalization and passage to a subsequence, a
unit vector in every $\mu_j^{-1}(C_j)$. The set is thus bounded
and contains only finitely many lattice points. The common expansion
is finite and represents $f$.
\end{proof}

\subsection{A divisorial criterion for equality}\label{subsec:divisorial}

For a seed $t$ over $S$, write $x_i$ for its mutable coordinates
and $s_i=\prod_{j\ne i}x_j$, with the product over the other mutable
indices. Write $\mathcal U_S(t)$ for its full upper LP algebra.
In positive mutable rank the union of the principal open sets
$D(s_i)$ is the locus where at most one mutable coordinate vanishes.

\begin{proposition}\label{prop:star-divisors}
Let $S$ be a noetherian UFD, and let $A$ be a noetherian
$S$-domain in the seed field, containing the initial coordinates and
contained in the first-star Laurent rings. Suppose each first mutant
belongs to $A[s_i^{-1}]$ in its direction. Let
$Z=\bigcup_{i<j}V(x_i,x_j)\subset\Spec A$. Then
\[ \mathcal U_S(t)=\Gamma(\Spec A\setminus Z,\mathcal O). \]
Thus $A=\mathcal U_S(t)$ precisely when regular functions on this
open extend to $\Spec A$. For normal $A$ this holds precisely when
no height-one prime contains two distinct mutable coordinates.
\end{proposition}
\begin{proof}
In each direction the two-chart intersection is
$C_i[x_i,F_i/x_i]$, where $C_i=S[x_j^{\pm1}:j\ne i]$ and
$F_i$ is the ordinary exchange polynomial, by Lemma~\ref{lem:recognition-rank-one}.
The hypotheses imply that $A[s_i^{-1}]$ lies in this ring and contains all
its generators. Intersect these localizations and apply the first-star
theorem. At rank zero both rings are $S$.
For normal $A$, codimension at least two of $Z$ gives extension
by the height-one intersection formula. Conversely, equality makes
the $x_i$ pairwise nonassociate primes by
Proposition~\ref{prop:upper-factorial}, which excludes a common
height-one prime.
\end{proof}

The following valuation formula makes the divisor test explicit.
All valuations in the statement are trivial on the ground field and
normalized to have value group $\mathbb Z$. The order of a top form
is its order at the divisor representing the valuation on a smooth model.

\begin{proposition}\label{prop:divisorial-order}
Let $k$ be algebraically closed of characteristic zero, let $T$ be
a split torus with character lattice $M$ and logarithmic volume
$\Omega$, and let $v$ be a divisorial valuation of $k(T)$.
Assume its restriction to characters is $d w_0$, with $d>0$ and
$w_0\in M^\vee$ primitive. Choose a character $t$ with $w_0(t)=1$
and set $S_0=k[\ker w_0]$.
If $\ord _v\Omega=0$, there is a prime Laurent polynomial
$q\in S_0$, unique up to a unit, such that
\[ v\!\left(\sum_jt^j f_j\right) =\min_{j:f_j\ne0}\bigl(dj+\ord _q f_j\bigr), \qquad f_j\in S_0. \]
If $\ord _v\Omega=-1$, then $d=1$ and $v$ is the
monomial valuation of $w_0$.
\end{proposition}
\begin{proof}
All characters in $S_0$ and their inverses have value zero.
The valuation has a center on
$V=\mathbb A^1_t\times\Spec S_0$ in $t=0$.
At the generic point of the center, of codimension $c$, choose
regular parameters $t,h_2,\ldots,h_c$. Complete these by separating
coordinates along the center. Pull their differentials back to a
smooth model representing $v$. Writing each parameter as a power
of a uniformizer times a unit shows that at most one differential in
their wedge can contribute the differential of the uniformizer.
Dividing by $t$, as required for $\Omega$, gives
\[ \ord _v\Omega\ \ge\ \sum_{\ell=2}^c v(h_\ell)-1\ \ge\ c-2 . \]
If $c=1$, the center is the generic point of $t=0$. Its DVR, in
the same fraction field, equals the valuation ring, giving the
monomial valuation and volume order $-1$. This proves the assertion
at a pole. In the remaining case $c=2$, and the center is generically
$t=q=0$ for a unique prime of the factorial ring $S_0$, and
$v(q)=1$.

At the generic point of this hypersurface $q$ is a smooth parameter.
In the chart following the valuation, successive blowups of
$(t,q)$, then $(t/q,q)$, and so on give $t=q^d b$ with $v(b)=0$.
The volume form in this chart is a unit times
$(db/b)\wedge dq\wedge\eta$ for a smooth transverse form $\eta$.
It is regular and nonvanishing near the center, since $b$ is a
unit there. A center of codimension at least two would now give order
at least one by the same differential estimate for a regular volume.
Thus $v$ is the order at the generic point of $q=0$ in this chart.
The residue of $b$ is transcendental over
$\Frac (S_0/(q))$. Substituting $t=q^d b$ in the
Laurent sum, its terms of least order give distinct powers of that
residue with nonzero coefficients. Their sum survives, as required.\qedhere
\end{proof}

In particular a Laurent polynomial is regular at a volume-zero
divisor of this kind precisely when $q^{-dj}$ divides $f_j$ for
every $j<0$. For instance, $t^{-2}q^2+t^{-1}$ passes the test
in degree $-2$ for $d=1$, but its degree $-1$ coefficient fails.
The same slice divisibilities apply to a finite-dimensional space of
Laurent polynomials with prescribed support. After the support bound
and the adjacent Laurentness of a given theta function are proved,
its coefficients satisfy a linear system. If the normalized solution
space consists of a single point, that point gives the theta function.

When the centre equation is linear in one torus coordinate, the
divisorial criterion becomes an explicit condition on coefficients.

\begin{boundarycorollary}\label{cor:linear-centre-valuation}
Use the notation of Proposition~\ref{prop:divisorial-order}, assume
$\ord_v\Omega=0$, and write $w=d w_0$. Denote the prime Laurent
polynomial in that proposition by $q_0$. Choose torus coordinates
$z_1,\ldots,z_n$ and a $w$-homogeneous Laurent polynomial
$q=A z_i+B$ such that $z^{-\alpha}q$ is an associate of $q_0$
in $S_0$ for some $\alpha\in\Supp q$. Here $A=c z^a$ is a
Laurent monomial with $c\in k^\times$, and $A,B$ are independent
of $z_i$, with $B\ne0$. For a nonzero $f=z^\ell P$, where
$P\in k[z_j^{\pm1}:j\ne i][z_i]$ has degree $s$ in $z_i$, write
\[ A^sP\bigl(z_i=(u-B)/A\bigr)=\sum_{j=0}^s u^jR_j(z_{\ne i}). \]
We regard the exponents of $R_j$ as elements of $M$ with zero
$i$-th coordinate. Then
\begin{equation}\label{eq:linear-centre-valuation}
v(f)=\min_{\substack{0\le j\le s\\ e\in\Supp R_j}}
\bigl(\langle\ell+e-sa+j\alpha,w\rangle+j\bigr).
\end{equation}
In particular, $f$ is regular at this divisor if and only if
$\langle\ell+e-sa+j\alpha,w\rangle\ge-j$ for every
$e\in\Supp R_j$.
\end{boundarycorollary}
\begin{proof}
Give $u$ weight $\gamma=\langle\alpha,w\rangle$ and decompose
$\sum_j u^jR_j$ into its nonzero homogeneous parts $H_h$ of weight
$h$. In each part let $j_h$ be the least power of $u$ that occurs.
After division by $u^{j_h}$, the constant coefficient is a nonzero
Laurent polynomial in the coordinates other than $z_i$. It remains
nonzero at the generic point of $q=0$, since there
$z_i=-B/A$ and these other coordinates are algebraically independent.
Consequently $H_h(q)$ has order exactly $j_h$ along $q=0$.

Proposition~\ref{prop:divisorial-order}, applied to the homogeneous
parts after removing a character of their weight, gives
\[ v(f)=\langle\ell-sa,w\rangle+\min_h(h+j_h). \]
For a term $u^jz^e$ in $H_h$ one has
$h=\langle e,w\rangle+j\gamma$. Taking the least $j$ in each
part therefore gives~\eqref{eq:linear-centre-valuation}.
\end{proof}

\subsection{Recovering ordinary exchanges on the mirror}

Proposition~\ref{prop:negative-vertex-chart} supplies a principal
torus on the forward mirror. It admits regular character coordinates:
choose a primitive point $b$ in the interior of the straight cone,
extend it to a lattice basis and add sufficiently large multiples
of $b$ to every other basis vector. The resulting basis lies in the
straight cone, so its characters are forward theta functions.
Inverting these characters gives the whole group algebra because
inverting $\theta_b$ alone already does so. Once a basis and its unit coordinates have been chosen, we can ask
whether the mirror is the upper LP algebra of a seed on this torus.
The criterion below also recovers the ordinary exchange polynomials.

Let $S=k[u_1^{\pm1},\ldots,u_c^{\pm1}]$ and
$\mathcal L=S[z_1^{\pm1},\ldots,z_r^{\pm1}]$. Suppose $B\subset \mathcal L$ is an
affine $k$-domain with known Laurent generators
$h_1,\ldots,h_t$. The proposed mutable coordinates are $z_i$ and
the proposed units are $u_j$. Require first that $z_i,u_j,u_j^{-1}$
belong to $B$. In particular $B[(\prod_i z_i)^{-1}]=\mathcal L$.
For each $i$ define
\[ C_i=S[z_j^{\pm1}:j\ne i],\qquad B_i=B[(\prod_{j\ne i}z_j)^{-1}]. \]
List all nonzero coefficients $f\in C_i$ of
negative powers $z_i^{-s}$ in the generators, recording the
positive integer $s$ beside each coefficient. If the list is empty,
the proposed direction fails. Otherwise factor any one listed $f$
in $C_i$. Choose a nonunit irreducible factor $q$, up to units of
$C_i$, and test the two conditions
\[ \operatorname{ord}_q(f)\ge s\quad\text{for every listed }(f,s), \qquad q\in z_iB_i . \]

\begin{theorem}\label{thm:recover-reverse-exchanges}
There is an LP seed in the proposed coordinates with coefficient ring
$S$ and full upper algebra $B$ if and only if the proposed
coordinates have the required regularity, one factor passes in every
mutable direction, and regular functions on the open where at most one
$z_i$ vanishes extend to $\Spec B$. For normal $B$ the last
condition says that no height-one prime contains two distinct $z_i$. At most one factor up to associates can pass
in each direction. Multiplying it by a Laurent monomial in the other
mutable coordinates to clear all negative powers and remove their
common monomial gives the ordinary exchange polynomial, unique up
to a unit of $S$.
\end{theorem}
\begin{proof}
For a nonunit prime $q\in C_i$, let $W(q)=C_i[z_i,q/z_i]$.
Lemma~\ref{lem:recognition-rank-one} identifies its negative slices:
the coefficient of $z_i^{-s}$ must be divisible by $q^s$.
The first test is therefore equivalent to $B_i\subseteq W(q)$.
Indeed the inverted variables are already units of $C_i$. Since
$C_i[z_i]\subseteq B_i$, the second test is equivalent to the
opposite inclusion. If no negative power occurs, then $B_i$
contains no such power and cannot contain $q/z_i$.
If two factors pass, then their rings $W(q)$ agree. The sets of
coefficients in homogeneous degree $-1$ are their respective
principal ideals in $C_i$, so the factors are associates.

For a seed with full upper algebra $B$,
Proposition~\ref{prop:upper-factorial} gives $B_i=W(F_i)$ and
extension across the simultaneous zeros by
Proposition~\ref{prop:star-divisors}. In particular $F_i$ divides every
listed coefficient, so it is among the factors tested. Conversely,
the normalized representatives of the recovered factors are
irreducible polynomials over $S$, independent of the corresponding
$z_i$ and divisible by no mutable variable. Irreducibility follows
by localizing at the other variables: a polynomial factor which
becomes a unit would be a coefficient unit times a monomial, which
the normalization excludes. These are the ordinary LP seed axioms.
The rank-one rings $W(q)$ agree with the Laurent intersections
for the normalized mutation partners, since the normalizing monomials
are units of $C_i$.

For this seed, $B_i=W(q_i)$ gives containment in the first star
and supplies the first mutant in every indicated localization.
Proposition~\ref{prop:star-divisors} identifies its full upper algebra
with the regular functions on the stated open. The extension
condition now gives equality with $B$.
\end{proof}

Equivalently, the exchange may be read from the image of the divisor
over the torus of unchanged coordinates. Elimination computes the
contraction $K_i=(z_i)B_i\cap C_i$. It is the ideal of the scheme-theoretic image of
$\Spec(B_i/(z_i))\to\Spec C_i$. In $W(q)$ this ideal is
precisely $(q)$: the quotient by $z_i$ is
$(C_i/(q))[Y]$, with $Y=q/z_i$ as a polynomial coordinate.
Thus any seed requires $K_i$ to be a nonzero proper principal
prime. Conversely, if $K_i=(q)$ is such a prime, the inclusion
$q\in z_iB_i$ is automatic and only the negative-slice test
remains at this direction. This description keeps multiplicities in
the image ideal.

The tests use finite Laurent data over a field with effective
arithmetic and polynomial factorization. For example, present $B$
by eliminating the ambient torus variables from the graph of
$(h_1,\ldots,h_t)$, adjoining inverses by equations $vv'-1=0$.
Choose a monomial order eliminating the ambient coordinates. A Laurent
function belongs to $B$ precisely when its normal form modulo the
graph ideal involves only the generator coordinates. Present $B_i$ by adjoining the indicated
inverses to $B$. The second factor test is then ordinary ideal
membership. In the normal case the extension test amounts to computing
the codimension of each $V(z_i,z_j)$. In general the localization
sequence for the complement $Z$ identifies the obstruction with
the first local cohomology group $H^1_Z(B)$ with support in $Z$,
since $B$ is a domain and the open is dense.
For nonempty $Z$ its vanishing is equivalent to depth at least two
along $Z$. The same calculations apply if
a finite presentation with its faithful Laurent realization is
already given.

For $B=B_k(U)$ the inputs may be any verified finite generating
set of forward theta functions and their computed Laurent expansions.
If the criterion holds and $\Spec B$ is smooth, the recovered seed
presents a reverse LP source in this marking. In characteristic zero
one can check smoothness from the presentation. The seed volume
is then regular and nowhere vanishing by the codimension-one covering
of Proposition~\ref{prop:upper-factorial}. A nonzero divisor class
group rules out an upper LP presentation in any marking over these
coefficients, again by that proposition. The recovered seed supplies the data for reverse counts.
Theorem~\ref{thm:one-sided-return} explains how regularity of the
resulting functions in the specified source chart enters marked return.

\subsection{Finite recognition of prescribed theta functions}
\label{subsec:prescribed-recognition}

Equality of two upper algebras leaves open how an isomorphism acts on
their theta functions. We give a finite criterion when enough regular
functions on the counting sources are known.

Let $Q$ be a submonoid of a lattice and let $G=\operatorname{gp}(Q)$
be the subgroup it generates. Set
$R=\mathbb Z[Q]$ and $K=\Frac(R)$. Write $t^\beta$
for the character of $\beta\in G$. Suppose an $R$-algebra has
specified basis $(\theta_p)_{p\in L}$ indexed by a lattice $L$,
with $\theta_0=1$. Fix a pointed rational polyhedral cone
$D\subset L_{\mathbb R}$. We use three properties of the products.
\begin{enumerate}
\item The coefficient of $\theta_{p+q}$ in
  $\theta_p\theta_q$ is a single monomial $t^{\alpha(p,q)}$ with
  coefficient one. Every other term has label $p+q+d$ with
  $d\in (D\cap L)\setminus\{0\}$.
\item On a specified complete finite rational fan $\Sigma_0$, the
  coefficient of $t^0\theta_{p+q}$ is one whenever $p,q$ belong
  to a common cone of $\Sigma_0$.
\item There are nonzero regular source functions $f_j$ with known
  finite Laurent supports $E_j\subset L^*$. Write
  $v_j(p)=\min_{e\in E_j}\langle e,p\rangle$.
  Whenever $\theta_r$ occurs in $\theta_p\theta_q$,
  $v_j(r)\ge v_j(p)+v_j(q)$ for every $j$.
\end{enumerate}
For intrinsic algebras the last inequality is the regular-function
inequality for the counting source. The zero-class rule comes from the
fan of the prepared compactification
\cite[Theorem~15.7 and Proposition~15.12]{KY}.
The first hypothesis concerns coefficients in the effective monoid
ring, before augmentation. Together with the zero-class rule it gives
coefficient $1$ at $\theta_{p+q}$ on every cone of $\Sigma_0$.

Let $\Delta$ be the common refinement of the fans of linearity of
the $v_j$. For each maximal cone $\tau$ of $\Delta$, choose
$a_{j,\tau}\in E_j$ attaining the minimum on $\tau$.
Assume that on every such cone
\[ D\cap\{d:\langle a_{j,\tau},d\rangle\ge0 \text{ for every }j\}=\{0\}. \]

\begin{proposition}\label{prop:prescribed-hilbert-products}
For any finite rational common refinement $\Sigma$ of $\Delta$
and $\Sigma_0$, and any $p_1,\ldots,p_m$ in one of its cones,
$\prod_i\theta_{p_i}=\theta_{\sum_i p_i}$. In particular every
theta function is one product of theta functions indexed by the union
of finite additive generating sets of the monoids $\sigma\cap L$
for $\sigma\in\Sigma$.
\end{proposition}
\begin{proof}
For two inputs $p,q$ choose a maximal $\tau$ containing their
common cone. If $p+q+d$ occurs, then
\[ \langle a_{j,\tau},p+q\rangle =v_j(p)+v_j(q)\le v_j(p+q+d) \le\langle a_{j,\tau},p+q+d\rangle . \]
The intersection hypothesis forces $d=0$, and the first two product
properties give its coefficient. Induction applies since each partial
sum stays in the same cone. The additive generating sets are finite by
Gordan's lemma.
\end{proof}

\begin{proposition}\label{prop:prescribed-finite-comparison}
Let $A_R,B_R$ have respective bases $(\theta_p)$, $(\psi_p)$
with the same marking $L$, each satisfying the hypotheses of
Proposition~\ref{prop:prescribed-hilbert-products}. Take a common
refinement for both algebras, and let $H$ be the union of finite
additive generating sets for the lattice points of its maximal cones.
Set $A_K=K\otimes_R A_R$ and $B_K=K\otimes_R B_R$, and let
$F:A_K\to B_K$ be a $K$-algebra homomorphism. Suppose that
there exist classes $b_h\in G$ satisfying
$F(\theta_h)=t^{b_h}\psi_h\ (h\in H)$.
Then $F$ is an isomorphism, and $F(\theta_p)=t^{b(p)}\psi_p$
for a unique function $b:L\to G$ with $b(h)=b_h$. This function
extends to an integral piecewise linear map on the common refinement.
The map restricts to $A_R\to B_R$ precisely when $b_h\in Q$
for all $h\in H$. Its restriction is an isomorphism onto $B_R$
precisely when $b_h,-b_h\in Q$ for all $h\in H$.
\end{proposition}
\begin{proof}
For $p=\sum_h n_hh$ in a common cone, the proposition applied on
both sides gives
$F(\theta_p)=t^{\sum_h n_h b_h}\psi_p$.
Its exponent is independent of the cone and decomposition since
$\psi_p\ne0$ and $G\to K^\times$ is injective. The resulting
function is additive on each cone. The lattice points in a
full-dimensional rational cone generate $L$ as a group, as is seen
by translating any lattice vector far into its interior. Thus $b$
extends integrally and linearly on each maximal cone. The extensions
agree on each common face because its lattice points span it. The
matrix of $F$ in the two bases is diagonal with invertible entries
over $K$. An entry belongs to $R$ exactly when its exponent
belongs to $Q$. The generator formula gives the stated finite tests
for $F$ and its inverse.
\end{proof}

\subsubsection{Determining the images of the generators}\label{subsec:hilbert-chart-test}

Given the homomorphism $F$, one can verify its effect on the
generators using a support bound and a normalization. Work in a
faithful target Laurent chart over a domain $S$, such as $R$ or $K$. Suppose the prescribed
function $\psi_h$ extends to the specified neighboring charts, is
supported in a known finite set $E$, and satisfies specified linear
conditions on its Laurent coefficients. These may prescribe one
coefficient and impose vanishing at several other exponents.
Write $u=\sum_{e\in E}c_ez^e$.
For each $a<0$, the coefficient $u_a$ of $z_i^a$ must satisfy
$u_a=H_i^{-a}v_a$ in the unchanged Laurent coordinates, as in
Lemma~\ref{lem:recognition-rank-one}. Only finitely many coefficients
of $v_a$ are needed. Indeed, if $P_a$ is the convex hull of the
allowed exponents in that slice, then each exponent $b$ of a
nonzero quotient satisfies
$b+(-a)\Newt (H_i)\subseteq P_a$.
This follows from additivity of Newton polytopes under multiplication
over a domain. Fixing any exponent of $H_i^{-a}$ therefore confines $b$ to a
translate of $P_a$. An empty slice has zero quotient. Coefficient
comparison now gives a finite linear system over $S$.

If the homogeneous system with zero normalization has only the zero
solution, the known function $\psi_h$ is the unique normalized
solution. Consequently an independently computed
$t^{-b_h}F(\theta_h)$ equals $\psi_h$ as soon as its full Laurent
polynomial has support in $E$ and satisfies this system. Over $R$ each slice quotient must have coefficients in $R$.
For the generic comparison one can first work over $K$, then test
the exponents of the generator images for descent to $R$.

\subsection{Integral descent}

We can recover an integral basis from a comparison over a
characteristic-zero field. The exchange polynomials and all coefficient
variables must use the same normalization in the two algebras.

\begin{proposition}\label{prop:monic-integral-descent}
Let $k$ be a field of characteristic zero and set
\[ S_{\mathbb Z}=\mathbb Z[f_1,\ldots,f_s,u_1^{\pm1},\ldots,u_c^{\pm1}],\qquad S_k=k\otimes_{\mathbb Z}S_{\mathbb Z}. \]
Fix an LP mutation pattern over $S_{\mathbb Z}$. For a seed $t$, write
$x_t^{\pm1}$ for the list $x_{t,1}^{\pm1},\ldots,x_{t,r}^{\pm1}$ and set
$\mathcal L_{\mathbb Z,t}=S_{\mathbb Z}[x_t^{\pm1}]$ and
$\mathcal L_{k,t}=S_k[x_t^{\pm1}]$. View these rings in the rational
function field of the pattern over $k$. Assume that
every exchange Laurent polynomial in the pattern is \emph{primitive over
$\mathbb Z$}, meaning that its integer coefficients have greatest common
divisor one. Define
\[ \mathcal U_{\mathbb Z}=\bigcap_t\mathcal L_{\mathbb Z,t},\qquad \mathcal U_k=\bigcap_t\mathcal L_{k,t}. \]
For an initial seed $t_0$, one has
\[ \mathcal U_{\mathbb Z}=\mathcal U_k\cap\mathcal L_{\mathbb Z,t_0}. \]
Suppose $(\Theta_m)_{m\in J}$ is a $k$-basis of $\mathcal U_k$,
each $\Theta_m$ belongs to $\mathcal L_{\mathbb Z,t_0}$, and their
least monomials have coefficient one and distinct exponents for a
single translation-invariant total order on the full exponent lattice
of $f,u,x_{t_0}$. Then
\[ \mathcal U_{\mathbb Z}=\bigoplus_{m\in J}\mathbb Z\Theta_m. \]
In particular, every $\Theta_m$ has integer Laurent coefficients in
every seed.
\end{proposition}
\begin{proof}
Let $g\in\mathcal U_k\cap\mathcal L_{\mathbb Z,t_0}$.
If $g\in\mathcal L_{\mathbb Z,t}$ and $t'$ is obtained by the
mutation $x_{t,i}'=H_i/x_{t,i}$, Lemma~\ref{lem:recognition-rank-one}
makes each negative slice $g_a$ divisible by $H_i^{-a}$ over $k$.
Faithful flatness descends this divisibility to $\mathbb Q$.
After clearing Laurent monomials, Gauss's lemma and primitivity of
$H_i$ give divisibility over $\mathbb Z$. Thus
$g\in\mathcal L_{\mathbb Z,t'}$. Induction along finite mutation
paths gives $g\in\mathcal U_{\mathbb Z}$; the reverse inclusion follows from
the definitions.

Consequently every $\Theta_m$ belongs to $\mathcal U_{\mathbb Z}$.
For $g\in\mathcal U_{\mathbb Z}$, take its finite expansion
$g=\sum_m c_m\Theta_m$ over $k$. Its least Laurent coefficient is
$c_m$ for the unique summand with least exponent, so $c_m\in\mathbb Z$.
Subtracting $c_m\Theta_m$ and inducting on the number of nonzero
summands proves integral spanning. Independence follows from the
$k$-basis property.
\end{proof}

\section{Valuation tests, reverse traces and mirror geometry}\label{app:sections}

We give finite tests for equality of reciprocal and smallest Newton
polytopes, then derive equations for reverse traces and study how the
mirror changes when a mutable coordinate is split. The last two
subsections treat obstructions to smoothness and interpret realization
in the Cox ring of a projective model.

\subsection{Finite valuation tests for reciprocal polytopes}

Recall that a label $a$ is realized when the reciprocal polytope
$P_a$ equals the smallest Newton polytope $Q_a$ of regular
functions with least exponent $a$.
Initial monomials give cones on which realization can be tested at
the generators. A second test compares the forward and source minima
at finitely many weights by taking the volumes of two polytopes.

\subsubsection{Initial monomials on a cone of Gauss weights}

Let $B\subset k[L]$ be the marked forward algebra of
Theorem~\ref{thm:reciprocal-polytope}, and write $M=L^*$ for the
source character lattice. For a lattice-separating
weight $w$ let $e_w(g)$ be the least Laurent exponent of
$0\ne g\in B$ and $\Gamma_w(B)$ the monoid of these exponents.

\begin{proposition}\label{prop:reciprocal-khovanskii-cone}
Let $\delta\subset L^*_{\mathbb R}$ be a rational polyhedral cone.
Suppose $g_1,\ldots,g_m\in B\setminus\{0\}$ have least exponents
$e_i=e_w(g_i)$ generating $\Gamma_w(B)$, and each $e_i$
minimizes every weight of $\delta$ on $\Newt(g_i)$.
Then $e_w(\Theta_p)$ minimizes every such weight on
$\Newt(\Theta_p)$ for every label $p$, uniquely at weights in
the ambient interior of $\delta$. The family $(P_a)_a$ is Minkowski
additive on $\delta$, and all its integral labels are realized
exactly when any given finite list of integral cone generators is
realized.
\end{proposition}
\begin{proof}
Valuative independence at $w$ gives distinct $e_w(\Theta_p)$.
Choose a product of the $g_i$ with least exponent $e_w(\Theta_p)$.
Its theta expansion contains $\Theta_p$, by the distinct least
exponents. Valuative independence for all weights places the polytope of
this summand inside that of the product. On the latter polytope its
least exponent minimizes throughout $\delta$ by Minkowski addition.
In the ambient interior a different exponent cannot tie, since their
nonzero difference pairs nonnegatively throughout the cone. Thus
every forward support function is linear there and its values on
$\delta$ give Minkowski additivity of $P$. The generator test
follows from Proposition~\ref{prop:single-sum-realization-cone}.
\end{proof}

\begin{corollary}\label{cor:reciprocal-khovanskii-cover}
If cones as in Proposition~\ref{prop:reciprocal-khovanskii-cone}
have dense union in $L^*_{\mathbb R}$, realization at their
generators is equivalent to realization at every integral label.
For a finite cover this is a finite test. Any normalized choice of
the resulting functions forms a basis of $k[U]$ that is valuatively
independent for all marked Gauss valuations.
\end{corollary}
\begin{proof}
The proposition and homogeneity give $P_a=Q_a$ on the rational
points of the union. At each integral direction their support minima
are continuous in $a$, so they agree everywhere. The basis claim is
Corollary~\ref{cor:reciprocal-from-windows}.
\end{proof}

The $g_i$ form a \emph{Khovanskii basis} when their least
exponents generate $\Gamma_w(B)$ \cite{KM}. This can sometimes
be tested on a polynomial presentation.

\begin{lemma}\label{lem:reciprocal-initial-ideal}
Let $T_i\mapsto g_i$ present $B$ with relation ideal $I$,
write $c_i z^{e_i}$ for the least term of $g_i$, and set
$u_i=\langle w,e_i\rangle$ and
$J=\ker(k[T_1,\ldots,T_m]\to k[L],\ T_i\mapsto c_i z^{e_i})$.
Then $\operatorname{in}_u I\subseteq J$ for minimum-weight forms.
Equality implies that the $e_i$ generate $\Gamma_w(B)$ if
$I$ is homogeneous for a positive integral grading, or if
$\sum_i\mathbb N u_i$ meets each bounded interval in a finite set.
\end{lemma}
\begin{proof}
The least terms of a relation cancel, giving the inclusion. Represent
$0\ne f\in B$ by $F(T)$. Under equality, whenever its least form
maps to zero, express that form by initial forms of relations with
multipliers of the required $u$-weight. Subtract the lifted combination.
This increases the least weight and preserves $f$, so the weight is
bounded above by $\nu_w(f)$. Under the grading hypothesis take the
relations homogeneous and lift separately in the original positive
degrees of $F$, where only finitely many monomials occur. The other
hypothesis gives finiteness of possible weights directly. The process
ends with a nonzero least image, a multiple of $z^{e_w(f)}$ because
$w$ separates exponents. Its exponent is a sum of the $e_i$.
\end{proof}

Elimination computes $J$, and a weight-refined homogeneous
Gr\"obner basis computes the initial ideal in the positive grading.

\subsubsection{A volume criterion from finitely many forward counts}

Fix the marking determined by a negative Newton vertex. Let $E=(p_1,\ldots,p_s)$ be a
finite list in $L$, and choose a sublist $E_0$ whose nonnegative
real span is $L_{\mathbb R}$. Write $d=\operatorname{rank}M$ and
$\Delta_0=\{b\in M_{\mathbb R}:\langle b,p\rangle\ge-1\ (p\in E_0)\}$.
If $E_0$ consists of both signs of every vector in a basis of $L$,
then $\Delta_0=[-1,1]^d$ in the dual coordinates.
For $a\in M_{\mathbb R}$ define
\[ h_a(p)=\min_{z\in\Newt(\Theta_p)}\langle a,z\rangle,\qquad \lambda_a(p)=\min_{z\in T_p}\langle a,z\rangle. \]
For $g=h,\lambda$ define the rational polytopes in $M_{\mathbb R}\times\mathbb R^s$
\[ K_g=\{(a,t): |t_j|\le1,\quad g_a(p)\ge-1\ (p\in E_0),\quad t_j\le g_a(p_j)\ (1\le j\le s)\}. \]
For $l\in[-n,n]^s\cap\mathbb Z^s$ let $A_{n,l}$ be the space of
sections on $U$ supported in $n\Delta_0$ with $\nu_{p_j}^x(f)\ge l_j$.

\begin{theorem}\label{thm:finite-valuation-volume}
In each degree $n\ge1$ one has
\[ \sum_l\dim_k A_{n,l}=\#\bigl(nK_\lambda\cap(M\oplus\mathbb Z^s)\bigr). \]
After division by $n^{d+s}$, the limit is $\operatorname{vol}(K_\lambda)
\le\operatorname{vol}(K_h)$ with lattice covolume one. Equality
holds exactly when $T_p=\Newt(\Theta_p)$ for every $p\in E$.
Equivalently, at each integral $a$ one normalized source function
with least exponent $a$ attains $h_a(p)$ simultaneously on $E$.
If $E$ detects every strict inclusion $P_a\subsetneq Q_a$, this
single volume equality is equivalent to realization at every label. Without
a finite detecting list, realization is equivalent to equality for
all finite lists containing both signs of every vector in a fixed basis of $L$.
\end{theorem}
\begin{proof}
A normalized section on the smallest polytope $Q_a$ belongs to
$A_{n,l}$ precisely when $(a,l)\in nK_\lambda$, by
Proposition~\ref{prop:source-polar-polytope}. These sections form a
basis of this space by Theorem~\ref{thm:minimum-source-newton-polytope},
giving the exact count. The inequality $\lambda_a(p)\le h_a(p)$ follows
from $P_a\subseteq Q_a$ by duality and continuity. Also
$h_a(p)\le\langle a,p\rangle$ since $z^p$ occurs in $\Theta_p$.
Thus $K_\lambda\subseteq K_h\subseteq\Delta_0\times[-1,1]^s$, and both
polytopes have interior since all inequalities are strict at
$(0,-\tfrac12\mathbf1)$. Lattice-point asymptotics give the claimed
limit, and equality of the volumes means $K_\lambda=K_h$.

For any integral $a$, the point $(a,(h_a(p_j))_j)$ belongs to
$nK_h$ for $n$ sufficiently large. If $K_\lambda=K_h$, its membership
in $nK_\lambda$ forces $\lambda_a(p)=h_a(p)$ for every $p\in E$.
Conversely, equality at every integral $a$ gives equality of the
support functions by homogeneity and continuity, hence $K_\lambda=K_h$.
For each $p\in E$, equality of these functions of $a$ is equivalent
to $T_p=\Newt(\Theta_p)$.
A section on $Q_a$ attains $\lambda_a(p)$ everywhere, and every section
with least exponent $a$ has value at most $\lambda_a(p)$. This proves
the simultaneous-attainment assertion. A detecting list gives the
global test by Corollary~\ref{cor:uniform-divisor-realization}. Any
strict inclusion gives distinct support minima at some integral weight,
which we can adjoin to this fixed list.
\end{proof}

\subsection{Degrees of exchange maps and reverse traces}

For a general exchange map, the degree of the generic fiber is a
product of mixed volumes. For a reverse trace, the contact and
evaluation conditions give a finite scheme whose degree is the trace.

\subsubsection{The map given by the ordinary mutation partners}

Write $G_i$ for the ordinary exchange at $x_i$. As in
Corollary~\ref{cor:component-linear-reciprocity}, draw an arrow
$j\to i$ when $x_j$ occurs in $G_i$. The strongly connected
components group together vertices that can reach one another by directed paths.

\begin{proposition}\label{prop:ordinary-partner-degree}
Fix polynomial exchange supports containing the origin, over an
algebraically closed field of characteristic zero. For general
coefficients the rational map
$x\mapsto(G_i(x)/x_i)_i$, over the field of the designated units,
is birational precisely when every exchange is jointly affine-linear
in the variables of its own strongly connected component.
\end{proposition}
\begin{proof}
The equations of a fiber over $t=(t_i)$ are $G_i(x)-t_i x_i=0$.
Order the strongly connected components by dependence. In each
component solve after adjoining the coordinates of its predecessors.
For general parameters each system is dominant and generically
finite. Indeed, the Newton polytopes of the projected equations
contain the segments $[0,e_i]$, so their mixed volume is positive.
By Bernstein's theorem the degree of the system is the normalized
mixed volume of these polytopes \cite{Bernstein}. For a component of
size $d$, this is the coefficient of $s_1\cdots s_d$ in the
volume of the Minkowski sum $\sum_j s_j P_j$ of its projected
Newton polytopes $P_j$, with lattice covolume one. Here each support
is projected onto the coordinates of the current component. The coefficients of its distinct projected monomials
use disjoint sets of general parameters from the current
equation, and the diagonal coefficient is the independent fiber
parameter $-t_i$. The generic count thus applies over the
predecessor field. Solving in this order multiplies the degrees.

For a component of size $d>1$, the sum of any nonempty proper
subcollection of its polytopes, indexed by $I$, spans the
coordinate directions $e_i$ for $i\in I$ and at least one
further direction, by strong connectivity. This is an essential
tuple in the sense of \cite{EG}. If its mixed volume is one,
Esterov and Gusev place lattice translates of all its polytopes
inside a unimodular simplex \cite[Theorem~2]{EG}. We use their
edge-tree description \cite[Proposition~3]{EG}. Each translated pair
$0,e_i$ therefore gives an oriented edge of that simplex.
These $d$ edges form a tree on the $d+1$ vertices, since an
undirected circuit would give a linear dependence among their
directions. Each exponent of the $i$-th equation is the
difference from the tail of edge $i$ to a vertex. Its coordinates
are nonnegative. Hence the path from that tail to the vertex
uses only forward edges, since its signed edge expansion is its
expansion in the basis $e_j$. If the equation involves
coordinate $j$, the tail of $j$ is consequently reachable from
the tail of $i$. Strong connectivity forces all tails to
coincide. The tree is a star and its allowed differences are
$0,e_1,\ldots,e_d$. This proves affine-linearity. A singleton
component gives one linear equation because $G_i$ is independent
of $x_i$. Conversely the componentwise linear systems have a
unique generic solution by Cramer's rule, as in
Theorem~\ref{thm:principal-minor-reciprocity}.
\end{proof}

\subsubsection{Reverse traces in a projective embedding}

A projective compactification gives equations for reverse traces.
More generally, the same construction computes Keel--Yu traces on any
counting source $U$ with a marked torus. For two or three contact orders
$P_1,\ldots,P_s$, choose a smooth projective simple normal crossings
compactification $\overline U$ on which each nonzero $P_i$ has divisorial
center $D_{j(i)}$ and is $m_i$ times the corresponding primitive
valuation. Write $D_1,\ldots,D_t$ for all the boundary components.
We may choose $\overline U$ so its boundary supports a very ample divisor
$H=\sum_j h_jD_j$ with every $h_j>0$. Indeed start with the
compactification of \cite[Lemma~19.2]{KY}, extract the contact rays by
boundary subdivisions as in \cite[Definition~1.1]{KY}, and subtract
sufficiently small positive exceptional multiples from the pullback of
the ample divisor at each blowup. A sufficiently divisible positive
multiple is integral and very ample.

\begin{lemma}\label{lem:projective-trace-equations}
Embed $\overline U\subset\mathbb P^N$ by $H$, and set
$d=\sum_{i:P_i\ne0}h_{j(i)}m_i$. The augmented trace of these contacts
is the degree over the generic field of marked points and evaluation
of the following finite scheme. Take basepoint-free tuples of $N+1$
binary forms of degree $d$, up to common scalar. Require their map
to have image in $\overline U$ and the prescribed value at the evaluation mark.
At the $i$th nonzero contact require its value to lie in the open
stratum of $D_{j(i)}$ and its contact order to be $m_i$. An input
of order zero is an interior mark. Equations for $\overline U$ and its boundary
thus determine the augmented trace from the projective model.
\end{lemma}
\begin{proof}
Fix the generic smooth ordered marked projective line and the generic
interior evaluation, as in \cite[Definition~1.1]{KY}. For a counted map
$f$ the pullback of each boundary component is exactly its prescribed
sum of contact marks. Hence $\deg f^*H=d$, and the embedding expresses
$f$ uniquely by the stated tuple up to scalar. Conversely the evaluation
condition prevents the image from lying in any boundary component.
At the prescribed marks the required contacts account for all $d$
zeros of $f^*H$, counted with multiplicity. After division by the
prescribed factors the section has relative degree zero on the fixed
projective line, so it comes from the base. Its evaluation at the
interior mark is a unit. Since every $h_j$ is positive, this also
excludes further boundary intersection in families.

The tuples form the usual open parameter scheme for maps from the fixed
line to projective space. Substitution in the homogeneous equations for
$\overline U$ imposes the image condition coefficient by coefficient. In charts
trivializing a boundary divisor, the first $m_i$ Taylor coefficients
of its local equation after substitution vanish, and the coefficient of
order $m_i$ is nonzero. These jet conditions agree on overlaps.
They give precisely the fixed-domain incidence scheme of the trace.
Its generic evaluation fiber is finite \'etale by
\cite[Proposition~3.13]{KY}. Taking its degree sums the contributions
of all numerical curve classes, as required. Over an effective field
the degree can be computed by elimination in the indicated open
parameter scheme.
\end{proof}

To distinguish numerical curve classes, choose additional very ample
divisors $H_j$ spanning $N^1(\overline U)_{\mathbb Q}$, the space
of rational divisor classes modulo numerical equivalence, and integers
$C_j$ with $C_jH-H_j$ ample. Every degree $d_j=\deg f^*H_j$
lies between $0$ and $C_jd$. For each tuple in these finite ranges,
adjoin binary forms of degree $d_j$ for the additional projective
embeddings, imposing the equations of the graph of their restrictions
to $\overline U$. A map has a unique graph lift. The degrees therefore give
disjoint parts of the original fiber and determine the numerical class
in each part. When $\overline U$ refines a specified compactification, group
these degrees by pushforward of the curve classes, as in
\cite[Proposition~17.3]{KY}.

\subsubsection{Reverse traces from Laurent generators}

\begin{proposition}\label{prop:laurent-generator-traces}
Let $Y$ be smooth affine over an algebraically closed field $k$ of characteristic zero,
with a marked open torus whose logarithmic volume extends to a nowhere-vanishing form on $Y$. Suppose a lattice basis
$w_1,\ldots,w_n$ of torus characters is regular on $Y$.
A finite set of Laurent generators $g_1,\ldots,g_h$ of $k[Y]$
gives polynomial equations and open conditions computing each
augmented trace of two or three contacts. In particular these data
are available for every smooth forward mirror in
Theorem~\ref{thm:automatic-reverse-chart}.
\end{proposition}
\begin{proof}
For each distinct primitive ray $p$ spanned by a nonzero contact, choose a generator
$g$ with $\nu_p(g)<0$, possible by Proposition~\ref{prop:coercive}.
In the toric chart of this ray write $g=s^{-l}h$ near its divisor,
where $s$ is a primitive normal coordinate and $l>0$.
Use the open subset $h\ne0$ of the divisor orbit.
Glue the torus and this subset in each ray chart to $Y$ along
the common torus. The resulting smooth scheme $W$ is separated.
Indeed the charts for distinct rays form part of a separated toric
variety. A valuation centered on one of these open orbits, with generic
point in the common torus, gives $s$ positive value and $h$
value zero. It cannot also have a center on $Y$, where $g$
is regular. This is the valuative criterion for the remaining
chart intersections. The torus form trivializes $K_W(D)$, with
$D$ the union of these open subsets.

Fix finite contact marks $a_i$ and evaluation $b$ at infinity
as in Theorem~\ref{thm:incidence}. Write $d_{ij}=P_i(w_j)$ and
$d_j=-\sum_i d_{ij}$. A negative $d_j$ gives zero. Otherwise set
\[ W_j(t)=w_j(b)A_j(t)\prod_i(t-a_i)^{d_{ij}}, \qquad A_j\text{ monic of degree }d_j,\quad A_j(a_i)\ne0. \]
Require every $g_k(W_1(t),\ldots,W_n(t))$ to be regular on the
complement of the marks. To impose this, clear the marked factors and
the monic denominator $\prod_j A_j^{m_j}$, where $m_j$ clears
negative powers of $w_j$ in $g_k$, and set the remainder modulo
that denominator equal to zero. Shared or repeated roots are allowed.
At each nonzero contact require the sum of the least-weight terms
of its chosen $g$ to be nonzero at the leading units of the $W_j$
at that mark. These are polynomial
open conditions after clearing the marked denominators.

The equations give maps into $Y$ off the marks because the $g_k$
generate its coordinate ring. At a nonzero mark the specified character
orders give a unique extension into the toric ray chart, of the required
multiplicity, and the open condition places it in the chosen open subset.
A zero contact here lies in the torus at the mark, which is a dense
open of the allowed interior. Conversely every fixed-domain map to
this pair of the indicated type has the stated unique coordinate expressions.
To identify the generic degree with the trace, we use the
comparison in Theorem~\ref{thm:incidence}. A projective model
extracting these valuations and $W$ agree on $Y$ and near the
generic point of each added divisor. On each component dominating evaluation,
the pullback of the logarithmic tangent bundle has degree zero
and is generated at the evaluation mark, hence is trivial. Evaluation
at each contact is therefore dominant on its divisor, and the excluded
proper subsets of those divisors, or of $Y$ at a zero contact,
do not contribute generically. The generic maps on the projective
model obey the same avoidance by \cite[Lemmas~3.6, 3.9 and 3.12]{KY}.
Thus the generic reduced degree on this open locus equals the trace.
For effectively given coefficients, elimination computes this degree.

On a smooth forward mirror the characters of a full-dimensional
straight cone are regular. To find a lattice basis in this cone, extend a primitive interior vector to a basis and add sufficiently
large multiples of it to the other basis vectors.
\end{proof}

\subsection{Splitting a mutable coordinate}\label{subsec:split-coordinate}

Suppose a seed as in Theorem~\ref{thm:corner} has a mutable coordinate
$x$. Replace it by $h\ge2$ mutable coordinates $y_1,\ldots,y_h$
and write $Y=\prod_i y_i$. At each $y_i$ take the old ordinary
exchange at $x$, and in every other ordinary exchange replace $x$
by $Y$. Keep the designated units. The new ordinary exchanges are
irreducible. Indeed on the full torus $y_1,\ldots,y_h$ may be
replaced by $Y,y_1,\ldots,y_{h-1}$ by a unimodular change of
coordinates. Irreducibility survives Laurent polynomial extension, and
none of the new ordinary exchanges has a mutable variable as a factor.
Write $\pi:L'\to L$ for the map of source cocharacter lattices that
adds the $h$ entries in these directions and leaves the others fixed.

\begin{lemma}\label{lem:split-source}
Let $A$ be the full upper algebra of the original seed, with its
designated units invertible. The new full upper algebra has the presentation
$A'=A[y_1,\ldots,y_h]/(Y-x)$.
For $A$ of finite type over $k$, $\Spec A'$ is smooth
exactly when $\Spec A$ and its divisor $x=0$ are smooth.
\end{lemma}
\begin{proof}
The LP starfish theorem \cite[Theorem~1.3]{Branching} describes each
upper algebra by its first star. Grade the split Laurent ring by $\mathbb Z^h/\mathbb Z\mathbf1$
using the exponents of the $y_i$. In each grade choose the unique
representative $\gamma$ with $\min_i\gamma_i=0$. Each homogeneous
Laurent polynomial is uniquely $y^\gamma f(Y,\ldots)$. At an
unsplit direction the ordinary exchange involves the $y_i$ only
through $Y$. The adjacent Laurent condition preserves the grade
and is exactly the old one on $f$, since LP normalizers are units in the
unchanged Laurent coordinates. In the term of $Y$-degree $n$ in
$f$, the split directions require the old exchange to the powers
$[{-n-\gamma_i}]_+$. Their simultaneous condition is divisibility
by the largest power, namely $[{-n}]_+$. Thus the new upper algebra
is the direct sum of $y^\gamma A$. Reducing monomials by $Y=x$
gives exactly this module and multiplication in the asserted
presentation.

The open where $h-1$ split coordinates are invertible is
$\Spec A$ times a torus, so smoothness of $\Spec A'$ implies
smoothness of $\Spec A$. If $\Spec A$ is smooth, the defining
hypersurface is smooth wherever at most one $y_i$ vanishes.
At two or more zero coordinates its differential is $-dx$. Every point of $x=0$ lifts to such a point. This proves the
criterion. For simultaneous splittings of distinct coordinates, the same
Jacobian calculation requires linear independence of the
$dx_i$ for every set of vanishing split divisors at their common
points, equivalently smooth intersections of the expected codimensions.
\end{proof}

\begin{theorem}\label{thm:split-coordinate}
Suppose $x=0$ is smooth on the original counting source. Let $g$
be a regular function whose Newton polytope has a vertex $v$ negative
in every mutable coordinate. Use the chart of
Proposition~\ref{prop:negative-vertex-chart} determined by $v$.
On the split source use the corresponding vertex of the pullback of
$g$ under $x=Y$. For a label $m'\in L'$, write
$A'_{m'}$ and $A_{\pi m'}$ for the new and old polynomial
factors. Denote their mutable variables by $w_i$ for the split
coordinates and by $z$ for the remaining coordinates. Then
\[ A'_{m'}(z,w_1,\ldots,w_h)=A_{\pi m'}(z,w_1+\cdots+w_h). \]
For the marking supplied by a polynomial corner, the sum of the degrees
in the split directions of each monomial is at most
$D_x(\pi m')=-\tau_{G_x}(\pi m')$, where $G_x$ is the polynomial
corner representative of the exchange at $x$.
\end{theorem}
\begin{proof}
We compare the incidence equations before and after splitting.
For a list of new input contacts $P'_a$, let $\delta_i$ be
their defects in the $y_i$ directions and let $d=\sum_i\delta_i$.
In the old problem use the contacts $\pi P'_a$ and the evaluation
coordinate $b_x=\prod_i b_{y_i}$. The other parameters agree. When
all defects are nonnegative and the designated units balance, we claim
that the new trace is
\begin{equation}\label{eq:split-trace}
\mathsf T'(P'_1,\ldots,P'_l) =\binom{d}{\delta_1,\ldots,\delta_h} \mathsf T(\pi P'_1,\ldots,\pi P'_l).
\end{equation}
We may compare the reduced generic fibers over a common algebraically
closed extension of their parameter fields. In the old equations the
monic polynomial for $x$ is the product of the new monic polynomials
for the $y_i$. The rational function for $x$ is likewise the
product of those for the $y_i$. All conditions at unmarked roots
then agree. In particular the common exchange is divisible by each
of the pairwise coprime new polynomials exactly when it is divisible
by their product. At a root assigned to $y_i$ its other factors are
units, so nonvanishing of the adjacent coordinate is equivalent in
the two descriptions.

A positive mutable-axis contact upstairs projects to the corresponding
positive axis downstairs. A projected contact can also lie on such an
axis when the original contact does not. The extra condition there holds
at every counted map for generic parameters. Indeed, at a nonzero
upstairs contact the map of open toric divisor orbits to the projected
orbit is dominant because its pullback on character lattices is the injection induced by $\pi^*$ on characters annihilating the projected contact. The
nonzero exchange on the projected orbit therefore pulls back to a
nonzero Laurent polynomial on the upstairs orbit. Generic contact
evaluation avoids its zero locus, by the same transversality used
in Theorem~\ref{thm:incidence}. A contact projecting to zero gives
an interior marked point, distinct from all unmarked roots. Conversely
every partition of the $d$ distinct old roots into ordered parts of
sizes $\delta_i$ satisfies the new marked conditions. It gives unique
monic factors and hence one new solution with the specified evaluation.
This proves the trace identity, including $d=0$.

The corresponding vertex of the pullback has exponent $v_x<0$ at
each split coordinate and agrees with $v$ in the other coordinates.
If $C$ is the normal cone of $v$, the new normal cone is $\pi^{-1}C$.
The correction bound at $m'$ is the old bound at $\pi m'$, with
the split shifts added together. To compute a given coefficient choose $b'$
sufficiently deep in $\pi^{-1}C$ and take $b=\pi b'$, so that
all bounded outputs of the translated products and the test indices
and their predecessors lie in the respective straight cones.
In the triangular trace equations an old coefficient $c_e$
whose correction in direction $x$ is $j$ is replaced, under the
proposed substitution, by
$\binom{j}{e_1,\ldots,e_h}c_e$ for $\sum_i e_i=j$.
For a correction $\boldsymbol d=(d_i)$ in the split
directions and $D=\sum_i d_i$, the trace multiplying this coefficient
is the trace of the corresponding pair before splitting, multiplied by
$\binom{D-j}{d_1-e_1,\ldots,d_h-e_h}$, interpreted as zero if any
difference is negative. Summing at fixed $j$ yields the common factor
\[ \sum_{\sum_i e_i=j} \binom{j}{e_1,\ldots,e_h} \binom{D-j}{d_1-e_1,\ldots,d_h-e_h} =\binom{D}{d_1,\ldots,d_h}. \]
This is the coefficient identity for two powers of $\sum_i w_i$.
By~\eqref{eq:split-trace}, the trace of the triple acquires the same factor. The original identity also holds outside its correction bound, by
pairing the translated product with the negative test label and
extending its coefficient array by zero. Thus the proposed polynomial satisfies
every equation in the new triangular system. Its unit diagonal proves
the result.
\end{proof}

\begin{corollary}\label{cor:split-reciprocal-windows}
Split any collection of mutable coordinates, allowing one factor for an
unchanged coordinate. Assume both varieties are smooth and use the
marking of Theorem~\ref{thm:split-coordinate}. On characters,
$\pi^*$ repeats the old exponent on all factors of a coordinate.
For a character $a'$ of the split source, define $a_i$ as the
minimum of its entries on the factors of the $i$-th coordinate,
and keep its frozen entries. With $b=a'-\pi^*a$ one has
\[ P'_{a'}=b+\pi^*P_a,\qquad Q'_{a'}=b+\pi^*Q_a. \]
Multiplication by the initial monomial $y^b$ after pullback gives
an isomorphism
\[ \{f\in A:\Supp(f)\subseteq P_a\} \ \xrightarrow{\ \sim\ }\ \{f'\in A':\Supp(f')\subseteq P'_{a'}\},\qquad f\longmapsto y^b f(Y). \]
It preserves the coefficients at $a,a'$. In particular, simultaneous
reciprocal realization, including the existence of a basis satisfying
the minimum formula, is equivalent for the two sources.
\end{corollary}
\begin{proof}
Apply Theorem~\ref{thm:split-coordinate} one coordinate at a time.
Each intermediate variety is smooth. Indeed, in the final variety invert
all but one factor of every coordinate still to be split. This open
is the product of the intermediate variety with a torus by Lemma~\ref{lem:split-source}.
The weight of the sum of the characters corresponding to one
original coordinate is the minimum of their weights. The expansion has nonnegative
integral coefficients, so substitution in the old polynomial factor
gives, for every split label $p'$,
\[ \nu^z_{a'}(\Theta'_{p'}) =\langle b,p'\rangle+\nu^z_a(\Theta_{\pi p'}). \]
Taking minimum support functions gives the identity for $P'_{a'}$.
For $Q'_{a'}$, choose a separating weight in the lifted cone and
its image for the old marking. The component of any split section
in the grade of its least exponent $a'$ is $y^b f(Y)$, with
$f\in A$ of least exponent $a$. Its polytope contains
$b+\pi^*Q_a$, and is contained in the polytope of the full section.
Conversely a normalized section on $Q_a$ pulls back after multiplication
by $y^b$ to one on $b+\pi^*Q_a$. This proves the second identity.
Every exponent
in $b+\pi^*P_a$ has the same grade modulo $\pi^*M$, and the
representative of this degree with minimum zero on the factors of each original coordinate is $b$. The graded description
in Lemma~\ref{lem:split-source} gives the claimed isomorphism of spaces of sections.
Every old label occurs by taking $a'=\pi^*a$. The basis assertion
follows from Corollary~\ref{cor:reciprocal-from-windows}.
\end{proof}

For example, partition the mutable indices into sets $I_j$.
Assume that the ordinary exchanges agree within each $I_j$ and
that exchanges outside $I_j$ depend on its variables only through
their product $Y_j=\prod_{i\in I_j}y_i$. Replacing these products by
single variables gives an LP seed. Its exchange representatives
exclude their own variables and have no variable factors. A
nontrivial factorization would pull back to a factorization in the
original seed. If the resulting seed satisfies
Corollary~\ref{cor:component-linear-reciprocity}, the normalized
section on the smooth split source is $y^bf_a(Y)$, with $f_a$
the principal-minor product of Theorem~\ref{thm:principal-minor-reciprocity}.
In particular an exchange can contain products of arbitrarily many
variables from its own strongly connected component.

\subsection{Obstructions to smoothness of the forward mirror}

One obstruction comes from a rigid toric singularity in a degeneration.
Another comes from a local equation for the mirror when several exchanges
depend linearly on a common coordinate.

\subsubsection{Rigid toric singularities}

\begin{corollary}\label{cor:rigid-independent-exchanges}
In the setting of Proposition~\ref{prop:principal-open-smoothability},
choose a nonempty set $I$ of mutable directions whose ordinary
exchanges depend only on the complementary coordinates. In the Laurent
ring of those coordinates, group the exchanges that differ by a unit.
For each class $I_l$ choose one representative with Newton polytope
$P_l$. Let $e_l$ be the standard basis vector for this class in a
separate lattice with one coordinate per class. Form their Cayley
polytope $Q=\operatorname{conv}\bigcup_l(\{e_l\}\times P_l)$.
If the augmented mirror of $U$ is smooth, every face $G$ of $Q$
with primitive edges and only triangular two-dimensional faces is a
unimodular simplex in its saturated affine lattice, the intersection
of its affine span with the ambient lattice. Here an edge is primitive
if its only lattice points are its two vertices. The face may be
$Q$ itself. In particular this test applies to every face of each
individual exchange Newton polytope in every mutable rank.
\end{corollary}
\begin{proof}
Invert the mutable coordinates outside $I$ and denote the resulting
principal open by $V$. Localization commutes with the finite
first-star intersection. At an inverted direction, every initial
Laurent polynomial belongs to the localized adjacent chart by
substituting the exchange for the product of the two partners.
The remaining charts impose exactly the first-star conditions for
$I$ over the complementary Laurent ring.

Translate a vertex of each $P_l$ to zero using one common generic
exposing weight. The translations act on $Q$ by a lattice shear.
Choose a weight basis in this common chamber by extending an interior
primitive vector to a basis and adding sufficiently large multiples
of it to the other basis vectors. In this polynomial corner set
$\tau_l=1+\sum_{i\in I_l}Z_i$. Corollary~\ref{cor:independent-face}
gives every forward function on $V$ in the form
$Z^p z^m\prod_l\tau_l^{-\min_{v\in P_l}\langle v,m\rangle}$.
Let $A$ be the toric algebra with monomials $\tau^a z^m$
satisfying $a_l+\langle v,m\rangle\ge0$ for all $l,v$, initially
with independent $\tau_l$. This expression identifies
$B_k(V)$ with
\[ A[Z_i^{\pm1}:i\in I]\big/ (\tau_l-1-\textstyle\sum_{i\in I_l}Z_i)_l . \]
Indeed the minimal allowed $a_l$ occur in the formula and the
$Z_i^{\pm1}$ also occur. Elimination of one $Z_i$ per class
identifies the quotient with a localization of a polynomial
extension of $A$, since the $\tau_l$, the $z$-coordinates and the remaining $Z_i$
are algebraically independent in the forward chart.

Write $\sigma=\mathbb R_{\ge0}Q$, so $A=k[\sigma^\vee\cap M']$ for the character lattice $M'$ of the independent $\tau_l,z$. All primitive rays of $\sigma$ have height one for the sum of the class coordinates.
For a face $G$, localize at a character exposing
$\gamma=\mathbb R_{\ge0}G$ in $\sigma$. This gives the toric
affine open for $\gamma$. After splitting off a torus it is the
cone for $G$ in the saturated lattice spanned by $\gamma$.
Its closed orbit meets the image of the indicated open defining
$B_k(V)$. In fact, $\tau_l$ is either zero or a nonconstant
character on this orbit, so one may choose the point with all
$\tau_l\ne1$, then choose nonzero $Z_i$ in each class with
sum $\tau_l-1$.

Primitive edges give regular two-dimensional cones. Indeed the
ray through one endpoint has height one, and the difference of the
endpoints is primitive in the height-zero lattice. The cone for
$G$ is therefore Gorenstein and smooth in codimension two.
Altmann's theorem for height sections with triangular two-faces
\cite[Corollary~2.8(1)]{Altmann} gives $T^1=0$.
The toric algebra is defined over $\mathbb Q$, and the finite
polynomial presentation computing $T^1$ commutes with field
extension, so this assertion holds over $k$.
Polynomial extensions and localizations preserve this vanishing.
If $G$ is not a unimodular simplex, the cone is singular on its
closed orbit by the toric smoothness criterion. The localized
$B_k(V)$ is then both singular and infinitesimally rigid.
Lift the exposing character to the Rees algebra of
Proposition~\ref{prop:principal-open-smoothability} and localize there.
This still gives a finite type flat family, with the indicated
special fiber and smooth locus over $t\ne0$ if $B_k(U)$ is
smooth. The Fitting-ideal argument of
Proposition~\ref{prop:principal-open-smoothability} applied to the
localized family rules out smoothness of $B_k(U)$. Taking $I=\{i\}$
gives the single-exchange assertion.
\end{proof}

\subsubsection{Singularities along a linear exchange direction}

Fix a coordinate $x_i$, mutable or a designated unit, in a family
of smooth full upper LP sources satisfying Theorem~\ref{thm:incidence}
with a common polynomial corner. Let $N_i$ be the set of mutable $j$ for which $x_i$
occurs in the ordinary exchange $G_j$. Let $y$ denote all the coordinates outside $\{i\}\cup N_i$, including the remaining designated units. Suppose that for $j\in N_i$
the exchange has the form
\[ G_j=B_j(y)+a_j A_j(y)x_i+ \sum_{j\to k}c_{jk} C_{jk}(y)x_k , \qquad k\in N_i\setminus\{j\}. \]
In this formula $j\to k$ means that $x_k$ occurs in $G_j$.
The arrows point from an exchange to the coordinates in $N_i$ it
involves. The coefficients $B_j,A_j$ are nonzero.
The scalars $a_j,c_{jk}$ vary independently over the field of
the coefficient functions and, if $i$ is mutable, of $G_i$. Assume the standing
source hypotheses on a nonempty open subset of this parameter
torus. The exchanges with indices outside $N_i$ can otherwise be arbitrary.

\begin{proposition}\label{prop:generic-linear-partner}
Write $Z_h=\Theta_{e_h}$ for the coordinate characters in the
polynomial-corner chart of the augmented mirror. On a nonempty
Zariski open subset of the parameters this chart gives
\[ T_i=\Theta_{-e_i}=Z_i^{-1}\prod_{j\in N_i}(1+Z_j). \]
\end{proposition}

\begin{proof}
Theorem~\ref{thm:corner} writes $Z_i T_i$ as a multiaffine
polynomial in $Z_j$, $j\in N_i$, with constant term one.
Fix $S\subseteq N_i$ and write $\mathbf1_S=\sum_{j\in S}e_j$.
Use the graph induced on $S$. For $W\subseteq S$, say that $W$
is \emph{reached from all of $S$} if every vertex has a directed path to
$W$, allowing paths of length zero. We compute two traces.

First consider the triple $(e_i,-e_i,-\mathbf1_S)$, at marks
$0,1,h$ and evaluation at infinity. For $S\ne\varnothing$ the
unmarked root in direction $j\in S$ is $u_j$. The coordinates
are $a t/(t-1)$ in direction $i$,
$b_j(t-u_j)/(t-h)$ for $j\in S$, and constant otherwise.
Consequently the remainder equations are
\[ B_j(u_j-1)(u_j-h)+A_j u_j(u_j-h) +(u_j-1)\sum_{k\in S,\ j\to k} C_{jk}(u_j-u_k)=0 . \]
In these equations $A_j,C_{jk}$ now denote the coefficients after
including their scalar parameters and the evaluation factors of
their moving coordinates. The nonzero constant $B_j$ absorbs all
coordinates outside $\{i\}\cup S$, which are constant along the
parameter line. Over the field of these data and the evaluation
parameters, the indicated $A_j,C_{jk}$ are independent.

These quadrics have $2^{|S|}$ distinct affine solutions generically.
After setting all the active $C_{jk}$ to zero, each affine equation
has the two simple roots $h,B_j/(A_j+B_j)$, and their ordinary
homogenizations have no common projective point at infinity. Properness
and the Jacobian criterion give the same count on a neighborhood of
this coefficient point.

For a solution write $V=\{j:u_j\ne h\}$. Generically there is no
edge from $S\setminus V$ to $V$. To see this, fix the root
coordinates. For each $j\in V$, the $j$-th equation gives an independent
condition on its own parameters whenever it is consistent. For $u_j\ne0,h$
the coefficient of $A_j$ is nonzero. At $u_j=0$ the constant
is $hB_j\ne0$. Each vertex outside $V$ with an edge into $V$ contributes another independent linear condition. Only $|V|$ root
coordinates are free, giving the assertion by dimension. The same
count with one fewer free coordinate excludes a zero, a one, or a
collision among the roots in $V$. When there is no such edge, the
equations outside $V$ vanish identically and the remaining equations
are those for $V$, treating $b_k(t-h)/(t-h)=b_k$ as constant
for $k\notin V$.

We verify the open conditions at these remaining roots. At a fixed
$u_j\ne0,1,h$, solving its equation for $A_j$ leaves the
active roots and the $C_{jk}$ free. A repeated zero of
$F_j(X(t))$ at $t=u_j$ imposes a further proper condition. Indeed for
$C_{jk}=0$ and
$A_j=-B_j(u_j-1)/u_j$, the function
$B_j+A_j t/(t-1)$ has a simple zero there. When $i$ is mutable,
the ordinary exchange at the positive mark $e_i$ is nonzero generically. The constant coordinates are generic and the active coordinates at zero are
$b_j u_j/h$, with the roots free in this parametrization. The
coefficients of $G_i$ are independent of the parameters in these equations.
For a designated unit the two orders in this direction balance and
no exchange condition is imposed there.
Thus the number of solutions on a fixed nonempty
$V$ is the triple trace $T(V)$ for that induced set. Set
$T(\varnothing)=1$, which also counts the all-$h$ solution.
We have proved
\[ 2^{|S|}=\sum_{V\subset S:\ \text{no edge from }S\setminus V \text{ to }V}T(V). \]
Every subset $W\subset S$ determines the set $V$ of vertices
having a directed path to $W$. This $V$ has the condition in
the sum, and within $V$, $W$ is reached from every vertex.
Conversely these two properties characterize $V$. Induction in
the last identity shows that $T(S)$ is the number of $W\subset S$
reached from all of $S$.

For the pair $(\mathbf1_W,-\mathbf1_S)$, the active unmarked roots
are indexed by $L=S\setminus W$. With marks zero and one, write
$w_j=u_j-1$. Clearing the single denominator yields the linear system
\[ d_j w_j-\sum_{k\in L,\ j\to k}C_{jk} w_k =-\sum_{k\in W,\ j\to k}C_{jk}\qquad(j\in L), \]
where $d_j$ is the value of $F_j$ at infinity. Its diagonal
entries vary independently of the off-diagonal $C_{jk}$, using
the independent scalar multiplying the $x_i$-term at its nonzero
evaluation. The matrix is generically invertible.
If a vertex cannot reach $W$, the equations on all such vertices
form a homogeneous subsystem, so their roots equal the mark one.
If every vertex can reach $W$, every $w_j\ne0$ generically.
Indeed specialize the $C_{jk}$ to positive rational numbers
and take all $d_j=d$ sufficiently large. The convergent geometric
series for the inverse matrix has a positive contribution for each
path to $W$, giving $w_j<0$. For nonzero $w_j$ the equations
solve uniquely for the $d_j$ in terms of the roots and off-diagonal
entries. Thus a coincidence among roots or with the other mark
is nongeneric. Each cleared exchange is linear in $t$ with
nonzero slope $d_j$, giving a simple zero. At a positive coordinate
mark, which can occur only for a singleton $W$, its exchange equation
is absent from the equations and is generically nonzero there.
The pair trace is consequently one if $W$ is reached from all
of $S$ and zero otherwise.

The monomials in $Z_iT_i$ are straight theta functions.
Taking the trace against $\theta_{-\mathbf1_S}$, the terms outside
$S$ vanish by the degree condition of Theorem~\ref{thm:incidence}. The computations show that
assigning coefficient one to each subset of $S$ gives exactly the
triple trace. The diagonal pair trace, at $W=S$, is one, so induction
on $|S|$ proves the formula for $T_i$.

\end{proof}

\begin{corollary}\label{cor:generic-linear-exchange-singularity}
Under the hypotheses of Proposition~\ref{prop:generic-linear-partner}, choose parameters in its open subset.
Replace $x_i$ by $G_i/x_i$ when $i$ is mutable and by
$x_i^{-1}$ when it is a unit, fixing all other coordinates.
Suppose that this source torus admits the Laurent realization of
Theorem~\ref{thm:corner} in its coordinate labels.
Then
\[ B_k(U)[Z_h^{-1}:h\ne i] =k[Z_h^{\pm1}:h\ne i][Z_i,T_i]. \]
The sole relation is $Z_i T_i=\prod_{j\in N_i}(1+Z_j)$.
In particular $|N_i|\ge2$ forces a singularity on the forward mirror.
\end{corollary}
\begin{proof}
In the second torus the standard coordinate labels give
the functions $T_i$ and $Z_h$ for $h\ne i$ in the original marking, by the tropical
change of coordinates and the nonzero constant term of $G_i$ in
the mutable case, or simply by inversion in the unit case.
The Laurent realization there places the localized mirror in
both Laurent rings with $Z_i$ or $T_i$ invertible. It contains
all the indicated generators. Their rank-one intersection gives the
presentation. If $j,k\in N_i$ are distinct, the points with
$Z_i=T_i=0,\ Z_j=Z_k=-1$ satisfy the equation and all its first
derivatives vanish. This proves singularity.
\end{proof}

\subsection{A projective model with finitely generated Cox ring}

Let $A$ be a full upper LP algebra over
$k[u_1^{\pm1},\ldots,u_c^{\pm1},s_1,\ldots,s_l]$,
where $k$ is a field. Write $M$ for the character lattice of the
initial torus after inverting the polynomial frozen variables,
$L=M^*$, and $A^\circ=A[(\prod_j s_j)^{-1}]$.
For a variety $X$ with free divisor class group, its \emph{Cox ring}
is the direct sum of the spaces $H^0(X,\mathcal O_X(D))$ over the
divisor classes, with representatives chosen additively so that
multiplication of sections respects the grading.

\begin{lemma}\label{lem:gauss-multisection-ring}
Suppose $A^\circ$ has a regular function with a Newton vertex
negative in every mutable coordinate. For any finite list of integral
Gauss weights $p_1,\ldots,p_s\in L$, let $t_1,\ldots,t_s$ be
formal variables and $t^d=\prod_j t_j^{d_j}$. Consider the multisection ring
\[ \bigoplus_{d\in\mathbb Z^s} \{h\in A:\nu_{p_j}(h)\ge-d_j\text{ for all }j\}\,t^d . \]
It is finitely generated over $k$. Its generators can be obtained by Hilbert bases and finite first-star
linear systems once the regular function and its vertex are specified.
\end{lemma}
\begin{proof}
Choose a full-dimensional lattice polytope $\Delta$ containing
zero in its interior. Adjoin a degree $n\ge0$ and restrict the
component of degree $(n,d)$ to functions supported in $n\Delta$.
Choose a separating weight $\omega$ in the interior of the
minimizing normal cone of the negative vertex.
Let $H$ be the cone of polynomial-frozen inequalities from
Corollary~\ref{cor:polynomial-polyhedral-newton-sections}, with
$H=M_{\mathbb R}$ when $l=0$. Sections of $A$ supported in
$n\Delta$ are exactly the sections of $A^\circ$ supported in
$n(\Delta\cap H)$. By Theorem~\ref{thm:polyhedral-newton-sections}
the triples $(n,d,a)$ attained with least exponent $a$ are
precisely the lattice points of the rational polyhedral cone given
by $a\in nD_\omega(\Delta\cap H)$, computed in $A^\circ$, and
\[ \langle a,v\rangle+d_j\ge0 \qquad(1\le j\le s,\ v\in\operatorname{Vert}(T_{p_j})). \]
Indeed the section with smallest polytope $Q_a$ attains all these
bounds simultaneously. At degree zero in $n$, the allowed
exponents consist only of zero. Lift the Hilbert generators to
monic sections in their degrees. Subtracting products with the
same least exponent generates every component, since each
subtraction strictly increases the least exponent within the
finite set $n\Delta\cap M$. Forgetting the auxiliary grading by
setting its formal variable to one gives a surjection onto the multisection ring,
since every finite Laurent support fits in some $n\Delta$.
\end{proof}

\begin{theorem}\label{thm:exchange-cox-model}
Let $A$ be the full upper algebra of any LP seed over
$k[u_1^{\pm1},\ldots,u_c^{\pm1}]$. There is a projective normal
locally factorial variety $X$ with reduced boundary $B$ and
$K_X+B\sim0$. Its open $X\setminus B$ is the locus in $\Spec A$
where at most one initial mutable coordinate vanishes, and
$\Gamma(X\setminus B,\mathcal O)=A$. The group $\Cl(X)=\Pic(X)$
is free of finite rank, generated by the strict toric boundary.
Over an algebraically closed field of characteristic zero the pair is
log canonical. If $A$ has a regular function with a Newton vertex
negative in every mutable coordinate, the Cox ring is finitely
generated. In this case $X$ is a \emph{Mori dream space} in the
sense of Hu and Keel \cite[Definition~1.10]{HuKeel} over an
algebraically closed field of characteristic zero.
\end{theorem}
\begin{proof}
We adapt the toric blowup construction of Gross, Hacking and Keel
\cite[\S3, Theorem~3.9]{GHK} to ordinary LP exchanges.
Choose a smooth projective toric variety $Y=Y_\Sigma$ whose fan
contains the mutable rays $e_i$. A barycentric subdivision, obtained
by smooth projective toric blowups, makes the divisors $D_{e_i}$
disjoint. Let $Z_i\subset D_{e_i}$ be the closure of $F_i=0$ in
the open orbit. This hypersurface is integral, since $F_i$ remains
an irreducible nonunit there. The divisor $Z_i$ is Cartier
on the smooth $D_{e_i}$, hence regularly embedded in $Y$ with
codimension two. Form $\pi:X=\operatorname{Bl}_{\coprod_iZ_i}Y\to Y$.
Each exceptional divisor $E_i$ is an integral projective-line bundle
over $Z_i$. The standard blowup charts are hypersurfaces in smooth
affine-line bundles over open subsets of $Y$, so $X$ is Cohen--Macaulay.
It is regular away from the $E_i$ and at their generic points,
hence normal by Serre's criterion.

Write $B_\rho$ for the strict toric boundary divisors and
$B=\sum_\rho B_\rho$. At the generic centers,
$\pi^*D_{e_i}=B_{e_i}+E_i$; every other toric divisor pulls back
to its strict transform. Thus all $B_\rho,E_i$ are Cartier.
The canonical divisor formula and the toric logarithmic volume give
$K_X+B=\pi^*(K_Y+\sum_\rho D_\rho)=0$.
Over an algebraically closed field of characteristic zero this
crepant equality also gives log canonicity.

Every toric stratum of codimension at least two lies in a boundary
divisor other than the disjoint $D_{e_i}$, so its preimage lies in $B$.
On the toric chart consisting of the initial torus and the open
orbit of $D_{e_i}$, blow up $(x_i,F_i)$ and delete the strict
transform of $x_i=0$. This leaves
$\operatorname{Spec}C_i[x_i,y_i]/(x_i y_i-F_i)$, where
$C_i=k[u^{\pm1},x_j^{\pm1}:j\ne i]$.
Together with the initial torus these charts cover $X\setminus B$.
Proposition~\ref{prop:upper-factorial} identifies them with
$D(\prod_{j\ne i}x_j)\subset\operatorname{Spec}A$, compatibly
on their common torus. Their coordinate rings intersect in $A$
by the LP starfish theorem. In mutable rank zero the open is
the initial torus.

Since the torus is factorial with units $k^\times x^M$, its divisor
sequence identifies $\Cl(X)$ with the cokernel of
\[ M\longrightarrow\mathbb Z^{\Sigma(1)}\oplus\mathbb Z^r,\qquad m\longmapsto\bigl((\langle m,\rho\rangle)_\rho,(\langle m,e_i\rangle)_i\bigr). \]
Here the rays are primitive. Subtracting the $e_i$-th boundary
coordinate from the $i$-th exceptional coordinate identifies this
cokernel with $\Cl(Y)\oplus\mathbb Z^r$, free of rank
$|\Sigma(1)|-\operatorname{rank}M+r$.
The principal relations for $x_i$ express $[E_i]$ in the strict
boundary classes. These Cartier divisors generate $\Cl(X)$,
so $X$ is locally factorial.

The blowup is an isomorphism at the generic point of each toric
boundary divisor, so $\ord_{B_\rho}(h)=\nu_\rho(h)$.
Global sections of $\mathcal O_X(\sum_\rho d_\rho B_\rho)$ are
precisely the $h\in A$ with $\nu_\rho(h)\ge-d_\rho$.
Under the negative-vertex hypothesis,
Lemma~\ref{lem:gauss-multisection-ring} gives finite generation
of the boundary multisection ring. Choose an additive lift of the
free group $\Cl(X)$ to boundary divisors. For a basis of the kernel
of the boundary class map, choose rational functions with those
principal divisors and extend multiplicatively. Multiplication by
these functions identifies each boundary degree with its chosen
lift and gives a graded surjection onto the Cox ring.
Over an algebraically closed field of characteristic zero the free
finitely generated Picard group has trivial connected component, so
$\Pic(X)_{\mathbb Q}=N^1(X)_{\mathbb Q}$.
Hu and Keel's criterion gives the Mori dream space assertion
\cite[Proposition~2.9]{HuKeel}.
\end{proof}

Over an algebraically closed field of characteristic zero, crepancy
identifies the valuations of log discrepancy zero with those of the
toric pair. These are the monomial valuations along its strata, hence
the marked Gauss valuations \cite[Lemma~2.2]{BL}. The pair thus has
maximal boundary, meaning that its dual complex has dimension
$\dim X-1$. Its Kawamata log terminal locus, where all log
discrepancies are positive, and $X\setminus B$ agree in
codimension one, so normality gives the same ring of regular functions.
Blum and Liu give a valuatively independent basis for these weights
\cite[Theorem~1.2(1)]{BL}. Theorem~\ref{thm:minimum-source-newton-polytope}
also identifies the smallest Newton polytopes and constructs a basis
over the original field by finite first-star equations.

\begin{corollary}\label{cor:cox-reciprocal-face}
In the setting of Theorem~\ref{thm:reciprocal-polytope}, fix a
reciprocal polytope $P_a$ and a primitive integral weight
$q\in\operatorname{int}C$. Choose the fan of
Theorem~\ref{thm:exchange-cox-model} refining the normal fan of
$P_a$ and containing the ray of $q$. Define
$D_a=-\sum_{\rho\in\Sigma(1)} \min_{b\in P_a}\langle b,\rho\rangle\,B_\rho$.
Choose finite homogeneous generators of $\operatorname{Cox}(X)$.
The label $a$ is realized on the source if and only if $[D_a]$
belongs to the rational cone generated by the degrees of those
generators whose sections do not vanish identically along $B_q$.
For a fixed fan and ray the same finite cone tests all labels whose
reciprocal support functions are linear on every cone of that fan.
\end{corollary}
\begin{proof}
Insert the specified rays by projective toric refinements, then take
a smooth refinement and a barycentric subdivision. The ray inequalities
cut out $P_a$, so $H^0(X,\mathcal O_X(nD_a))=V(nP_a)$ for $n\ge1$.
A section $h$ has order $\nu_q(h)-n\langle a,q\rangle$ along $B_q$.
Since $q$ exposes $a$, order zero is equivalent to a nonzero
coefficient at $x^{na}$. Then $na$ is the least exponent at every
separating weight in the straight cone, so
Theorem~\ref{thm:minimum-source-newton-polytope} gives
$nQ_a=Q_{na}\subseteq nP_a$. Hence $P_a=Q_a$.
Conversely, a section realizing $a$ has order zero along $B_q$.

Express such a section as a homogeneous polynomial in the Cox
generators. Some monomial is nonzero at the generic point of $B_q$,
so its factors all have order zero there. Conversely, clear
denominators in a rational conic expression for $[D_a]$ and take
the corresponding product of generators. It is a section of a line
bundle of class $n[D_a]$, for some $n>0$, nonzero at the generic
point of $B_q$. Identifying this line bundle with
$\mathcal O_X(nD_a)$ gives the required section of order zero.
\end{proof}

\subsection*{Acknowledgements}

We thank Helge Ruddat for pointing out the connections with logarithmic
and tropical curve counts.

\section*{Declaration of generative AI and AI-assisted technologies}

The author used ChatGPT and Codex (OpenAI) to explore arguments, check proofs
and references, and edit the manuscript. The author takes responsibility for its content.

\enlargethispage{2\baselineskip}
\bibliographystyle{amsplain}

\end{document}